\documentclass[12pt]{article}
\usepackage{amsmath}
\usepackage{amssymb}
\usepackage{epic}
\numberwithin{equation}{section}
\usepackage{epsfig}

\title{Random and Integrable Models in Mathematics and
Physics}

\author{
 Pierre van
Moerbeke\thanks{ Series of five Lectures at the
``Program on Random Matrices, Random Processes and
Integrable Systems" at the ``Centre de Recherche
Math\'ematique", Montr\'eal, Canada (June 20-July 8,
2005). Department of Mathematics, Universit\'e
Catholique de Louvain, 1348 Louvain-la-Neuve, Belgium
and Brandeis University, Waltham, Mass 02454, USA.
E-mail: Pierre.Vanmoerbeke@uclouvain.be and
vanmoerbeke@brandeis.edu. The support of a National
Science Foundation grant \# DMS-07-06287, a European
Science Foundation grant (MISGAM), a Marie Curie Grant
(ENIGMA) and ``P\^ole d'attraction inter-universitaire"
grants is gratefully acknowledged.}}

\date{}

\catcode`\@=11
\let\c@equation=\relax
\newcounter{equation}[subsection]

\catcode`\@=12

\newcommand{\MAT}[1]{\left(\begin{array}{*#1c}}
\newcommand{\mat}{\end{array}\right)}
\newcommand{\qed}{\leavevmode\unskip\nobreak\penalty200\hskip2pt\null
\nobreak\hfill\rule{1.1ex}{1.1ex}
\medbreak }

\newcommand{\rg}{\rightarrow}
\newcommand{\lrg}{\longrightarrow}

\newcommand{\DR}{{\cal D}}
\newcommand{\ER}{{\cal E}}
\newcommand{\BE}{{\mathbb E}}
\newcommand{\BP}{{\mathbb P}}

\newcommand{\FR}{{\cal F}}

\newcommand{\HR}{{\cal H}}

\newcommand{\LR}{{\cal L}}

\newcommand{\BC}{{\mathbb C}}

\newcommand{\BN}{{\mathbb N}}

\newcommand{\BY}{{\mathbb Y}}
\newcommand{\BZ}{{\mathbb Z}}
\newcommand{\Sg}{\Sigma}
\newcommand{\iy}{\infty}
\newcommand{\pl}{\partial}
\newcommand{\al}{\alpha}

\newcommand{\gs}{{\bf s}}
\newcommand{\no}{\nonumber}

\newcommand{\ba}{{\backslash}}

\newenvironment{proof}{\medskip\noindent{\it Proof:\/} }{\qed}
\newenvironment
        {remark}{\medskip\noindent\underline{\it Remark:\/} }{\medbreak}
\newenvironment
        {example}{\medskip\noindent\underline{\it Example:\/} }{\medbreak}

\newcommand{\om}{\omega}

\newcommand{\la}{\langle}
\newcommand{\ra}{\rangle}
\newcommand{\ga}{\gamma}

\newcommand{\dt}{\delta}
\newcommand{\Dt}{\Delta}
 \newcommand{\vr}{\varepsilon}
\newcommand{\sg}{\sigma}
\newcommand{\BR}{{\mathbb R}}
\newcommand{\lb}{\lambda}
\newcommand{\Lb}{\Lambda}

\newcommand{\BJ}{{\mathbb J}}

\newcommand{\diag}{\operatorname{diag}}
\newcommand{\Res}{\operatorname{Res}}

\def\be#1\ee{\begin{equation}#1\end{equation}}
\def\bea#1\eea{\begin{eqnarray}#1\end{eqnarray}}
\def\bean#1\eean{\begin{eqnarray*}#1\end{eqnarray*}}

\newcommand{\Tr}{\operatorname{\rm Tr}}

\newtheorem{definition}{Definition}[section]
\newtheorem{theorem}[definition]{Theorem}
\newtheorem{lemma}[definition]{Lemma}
\newtheorem{corollary}[definition]{Corollary}
\newtheorem{proposition}[definition]{Proposition}

\catcode `!=11

\newdimen\squaresize
\newdimen\thickness
\newdimen\Thickness
\newdimen\ll! \newdimen \uu! \newdimen\dd! \newdimen \rr! \newdimen
\temp!

\def\sq!#1#2#3#4#5{%
\ll!=#1 \uu!=#2 \dd!=#3 \rr!=#4
\setbox0=\hbox{%
 \temp!=\squaresize\advance\temp! by .5\uu!
 \rlap{\kern -.5\ll!
 \vbox{\hrule height \temp! width#1 depth .5\dd!}}%
 \temp!=\squaresize\advance\temp! by -.5\uu!
 \rlap{\raise\temp!
 \vbox{\hrule height #2 width \squaresize}}%
 \rlap{\raise -.5\dd!
 \vbox{\hrule height #3 width \squaresize}}%
 \temp!=\squaresize\advance\temp! by .5\uu!
 \rlap{\kern \squaresize \kern-.5\rr!
 \vbox{\hrule height \temp! width#4 depth .5\dd!}}%
 \rlap{\kern .5\squaresize\raise .5\squaresize
 \vbox to 0pt{\vss\hbox to 0pt{\hss $#5$\hss}\vss}}%
}
 \ht0=0pt \dp0=0pt \box0
}

\def\vsq!#1#2#3#4#5\endvsq!{\vbox to \squaresize{\hrule
width\squaresize height 0pt%
\vss\sq!{#1}{#2}{#3}{#4}{#5}}}

\newdimen \LL! \newdimen \UU! \newdimen \DD! \newdimen \RR!

\def\vvsq!{\futurelet\next\vvvsq!}
\def\vvvsq!{\relax
  \ifx     \next l\LL!=\Thickness \let\continue=\skipnexttoken!
  \else\ifx\next u\UU!=\Thickness \let\continue=\skipnexttoken!
  \else\ifx\next d\DD!=\Thickness \let\continue=\skipnexttoken!
  \else\ifx\next r\RR!=\Thickness \let\continue=\skipnexttoken!
  \else\def\continue{\vsq!\LL!\UU!\DD!\RR!}%
  \fi\fi\fi\fi
  \continue}

\def\skipnexttoken!#1{\vvsq!}

\def\place#1#2#3{\vbox to 0pt{\vss
\rlap{\kern#1\squaresize
  \raise#2\squaresize\hbox{$#3$}}
\vss}}

\def\Young#1{\LL!=\thickness \UU!=\thickness \DD! = \thickness \RR! =
\thickness \vbox{\smallskip\offinterlineskip
\halign{&\vvsq! ##
\endvsq!\cr #1}}}

\def\blank{\omit\hskip\squaresize}
\catcode `!=12

\begin{document}
\maketitle

\tableofcontents

\vspace*{2cm}

During the last 15 years or so, and since the pioneering
work of E. Wigner, F. Dyson and M. L. Mehta, random
matrix theory, combinatorial and percolation questions
have merged into a very lively area of research,
producing an outburst of ideas, techniques and
connections; in particular, this area contains a number
of strikingly beautiful gems. The purpose of these five
Montr\'eal lectures is to present some of these gems in
an elementary way, to develop some of the basic tools
and to show the interplay between these topics. These
lectures were written to be elementary, informal and
reasonably self-contained and are aimed at researchers
wishing to learn this vast and beautiful subject. My
purpose was to explain these topics at an early stage,
rather than give the most general formulation.
Throughout, my attitude has been to give what is
strictly necessary to understand the subject. I have
tried to provide the reader with plenty of references,
although I may and probably will have forgotten some of
them; if so, my apologies!

 As we now know, random matrix theory has
reached maturity and occupies a prominent place in
mathematics, being at the crossroads of many subjects:
number theory (zeroes of the Riemann zeta functions),
integrable systems, asymptotics of orthogonal
polynomials, infinite-dimensional diffusions,
communication technology, financial mathematics, just to
name a few. Almost 200 years ago A. Quetelet tried to
establish universality of the normal distribution
(mostly by empirical means). Here we are, trying to
prove universality of the many beautiful statistical
distributions which come up in random matrix theory and
which slowly will find their way in everyday's life.

This set of five lectures were given during the first
week of a random matrix 2005-summer school at the
``Centre de Recherches Math\'ematiques" in Montr\'eal;
about half of them are devoted to combinatorial models,
whereas the remaining ones deal with related random
matrix subjects. They have grown from another set of ten
lectures I gave at Leeds (2002 London Mathematical
Society Annual Lectures), and semester courses or
lecture series at Brandeis University, at the University
of California (Miller Institute, Berkeley), at the
Universiteit Leuven (Francqui chair, KULeuven) and at
the Universit\'e de Louvain (UCLouvain).

I would like to thank many friends, colleagues and
graduate students in the audience(s), who have
contributed to these lectures, through their comments,
remarks, questions, etc... , especially Mark Adler, Ira
Gessel, Alberto Gr\"unbaum, Luc Haine, Arno Kuijlaars,
Walter Van Assche, Pol Vanhaecke, and also Jonathan
Del\'epine, Didier Vanderstichelen, Tom Claeys, Maurice
Duits, Maarten Vanlessen, Aminul Huq, Dong Wang and many
others....

Last, but not least, I would like to thank John Harnad
for creating such a stimulating (and friendly)
environment during the 2005-event on ``random matrices"
at Montr\'eal. Finally, I would label it a success if
this set of lectures motivated a few young people to
enter this exciting subject.

\newpage

\section{Permutations, words, generalized permutations and percolation}

%

\subsection{Longest increasing subsequences in permutations,
words and generalized permutations} (i) {\bf
Permutations} $\pi:=\pi_n$ of $1,\ldots,n$ are given by

$$
S_n\ni \pi_n =\left(\begin{array}{ccc}
1&\ldots &n\\
 j_1&\ldots&j_n
\end{array}\right),~~~~1\leq j_1,\ldots,j_n \leq n~\mbox{all distinct
integers}
$$
 with $\pi_n(k)=j_k$. Then
 $$
 \#S_n=n!.
 $$
An {\em increasing subsequence} of $\pi_n \in S_n$ is a
sequence $1\leq i_1<...< i_k\leq n$, such that $\pi_n
(i_1)<...<\pi_n(i_k)$. Define \be L(\pi_n) =  \mbox{
length of a longest (strictly) increasing subsequence of
$\pi_n$ }.
 \ee
 Notice that there may be many longest (strictly) increasing
 subsequences!

\noindent {\bf Question:}~Given uniform probability on
$S_n$, compute
  $$
   P^{n} (L(\pi_n )\leq k,~\pi_n \in S_n) =\frac{\#\{
    L(\pi_n )\leq k,~\pi_n \in S_n\}}{n!}=~?
    $$
    $\hspace{9cm}\mbox{(Ulam's problem 1961)}$

\begin{example} for $\pi_7=
 \left(\begin{array}{cccccccc}
1&2&3&4&5&6&7\\
\underline{3}&1&\underline{4}&2&\underline{6}
  &\underline{7}&5
  \end{array}
  \right)$, we have
  $L(\pi_7)=4$. A longest increasing sequence is
  underlined; it is not necessary unique.
  \end{example}

\vspace{.6cm}

(ii) {\bf Words} $\pi:=\pi_n^q$ of length $n$ from an
alphabet $1,\ldots, q$ are given by integers
$$
S^q_n\ni \pi_n^q =\left(\begin{array}{ccc}
1 &\ldots&  n\\
 j_1 & \ldots & j_n
\end{array}\right),~~~~1\leq j_1,\ldots,j_n \leq q~
$$
 with $\pi_n^q(k)=j_k$. Then
 $$
 \#S^q_n=q^n.
 $$
 An {\em increasing subsequence} of $\pi_n^q \in S^q_n$ is given by a sequence
$1\leq i_1<\ldots< i_k\leq n$, such that $\pi_n^q (i_1)
 \leq \ldots \leq
\pi_n^q(i_k)$. As before, define \be L(\pi^q_n) = \mbox{
length of the longest weakly increasing subsequence of
$\pi^q_n$ }.
 \ee

 \noindent{\bf Question:}~Given uniform probability on
 $S_n^q$, compute
 $$P^q_n
  \left(L(\pi_n^q )\leq k,~\pi_n^q  \in
   S^q_n
   \right)=\frac{\#\left\{L(\pi_n^q )\leq k,~\pi_n^q  \in
   S^q_n\right\}}{q^n}=~?$$

\begin{example}
 for $
\pi=\left(\begin{array}{lllll}
1&2&3&4&5\\
2&\underline{1}&\underline{1}&3&\underline{2}\end{array}\right)\in
S_5^3 $, we have $L(\pi)=3$. A longest increasing
sequence is
  underlined.
\end{example}

  \vspace{.6cm}

(iii) {\bf Generalized permutations} $\pi:=\pi_n^{p,q}$
are defined by an array of integers
$$\mbox{GP}_n^{p,q}\ni
\pi_n^{p,q} =\left(\begin{array}{ccc}
i_1&\ldots&i_n\\
 j_1&\ldots&j_n
\end{array}\right),
$$
%
subjected to
$$
1\leq i_1\leq i_2\leq \ldots\leq i_n\leq p
~~~\mbox{ and }~~~1\leq j_1,\ldots, j_n\leq q$$
$$\mbox{with}~~
i_k =i_{k+1}\mbox{~implying~}j_k\leq j_{k+1}.
$$
Then
$$
\#\mbox{GP}_n^{p,q}=\left(\begin{array}{c}
  pq+n-1\\ n
\end{array}\right).$$
 An {\em increasing subsequence} of a generalized permutation $\pi$ is defined as
$$
 \left(  \begin{array}{c}
i_{r_1}\\
j_{r_1}
\end{array} ,\ldots, \begin{array}{c}
i_{r_m}\\
j_{r_m}
\end{array} \right)\subset\pi
$$
with $r_1\leq \ldots \leq r_m$ and $j_{r_1}\leq
j_{r_2}\leq \ldots\leq j_{r_m}$.
Define
$$
L(\pi):=\mbox{~length of the longest weakly increasing
subsequence of $\pi$}.
$$

\begin{example}
For $\left(\begin{array}{llllllllll}
1&1&1&2&2&2&3&3&4&4\\
 {2}& {3}& {3}&\underline{1}&
 \underline{2}&
\underline{2}&1&\underline{2}&1&\underline{3}
\end{array}\right)\in \mbox{GP}_{10}^{4,3}$, we have $L(\pi)=5$.

\end{example}

For more information on these matters, see Stanley
\cite{Stanley1}.

\newpage

\subsection{Young diagrams and Schur
polynomials}\label{subsection1.3}

 Standard references to this
subject are MacDonald\cite{MacDonald},
Sagan\cite{Sagan}, Stanley\cite{Stanley1}, Stanton and
White \cite{Stanton}. To set the notation, we remind the
reader of a few basic facts.

\begin{itemize}
  \item A {\em partition }$\lb$ of $n$ (noted $\lb \vdash n$)
  or a {\em Young diagram} $\lb$ of weight $n$ is represented
  by a sequence of integers $\lambda_1\geq
  \lambda_2 \geq ...\geq \lambda_{\ell}\geq 0$,
  such that $n=|\lambda|:=\lambda_1+...+\lambda_{\ell}$;
  $n=|\lambda|$ is called the weight.
 A {\em dual Young diagram} $ \lambda^{\top}=(
  (\lambda^{\top})_1 \geq (\lambda^{\top})_2 \geq...)$ is the
  diagram obtained by flipping the diagram $\lambda$
  about its diagonal; set
  \be\lambda_i^{\top}:=
  (\lb^{\top})_i
   =
    \mbox{length of $i$th column of $\lb$}.
     \label{TransposePartition}\ee
    Clearly $|\lambda|=|\lambda^{\top}|$.
  For future use, introduce the following notation:
  \bea
   \BY&:=&\{ \mbox{all partitions}~\lb\}  \nonumber  \\
   \BY_n&:=&\{ \mbox{all partitions}~\lb\vdash n\}
   \nonumber  \\
   \BY ^p&:=&\{ \mbox{all partitions},~\mbox{with}~
      \lb^{\top}_1\leq p\}\no\\
      \BY_n^p&:=&\{ \mbox{all partitions}~\lb\vdash n,~\mbox{with}~
      \lb^{\top}_1\leq p\}
      \label{Y-definitions}\eea

  \item A {\em semi-standard Young tableau} of shape $\lambda$
  is an array of
   integers $a_{ij}>0$ placed in box $(i,j)$ in the Young diagram $\lambda$, which are
   weakly increasing from left to right {\em and} strictly
   increasing from top to bottom.

  \item A {\em standard Young tableau} of shape $\lambda
    \vdash n$ is an array of
   integers $1,...,n=|\lb|$ placed in the Young diagram, which are
   strictly
   increasing from left to right {\em and} from top to bottom. For
   $\lb \vdash n$, define
   $$
   f^{\lb}:=\# \left\{
   \begin{array}{l}
    \mbox{standard tableaux of shape $\lb$}\\
     \mbox{filled with integers $1,\ldots,|\lb|$}
     \end{array}
            \right\}.
  $$

 \item The {\em Schur function $\tilde{\bf s}_{\lambda}$}
  associated with a Young diagram $\lambda \vdash n$
  is a symmetric function in the variables $x_1,x_2,...$,
  (finite or infinite), defined as
%
  \be
  \tilde\gs_{\lambda}(x_1,x_2,...)=
  \sum_{
  {\mbox{\tiny semi-standard}}\atop
     {\mbox{\tiny tableaux $P$ of shape $\lb$}   }
   }  \prod_i x_i ^{\# \{\mbox{times $i$
  appears in P}\}}
  \label{Schur polynomial}\ee
  It equals a polynomial $\gs_{\lambda}(t_1,t_2,\ldots)$
  (which will be denoted without the tilde) in the symmetric variables $kt_k =\sum_{i\geq 1}
 x_i^k$,
 \be
  \tilde\gs_{\lambda}(x_1,x_2,...)=\gs_{\lambda}(t_1,t_2,\ldots)
  = \det \left(  {\bf s}_{\lb_i-i+j}(t)\right) _{1\leq i,j\leq m}
  ,\ee
for any $m\geq n$. In this formula $s_i(t)=0$ for $i<0$
and $s_i(t)$ for $i\geq 0$
 is defined as
 $$e^{\sum^{\iy}_{1}t_iz^i}:=\sum_{i\geq 0} {\bf
   s}_i(t_1,t_2, \ldots)z^i.$$
   Note for $\lb\vdash n$,
   $$
   \tilde{\bf s}_{\lb}(x_1,x_2,\ldots)= f^{\lb}~ x_1\ldots x_n+\ldots.
   =f^{\lb}\frac{t_1^{|\lb|}}{|\lb|!}+\ldots$$

 \item Given two partitions $\lb \supseteq\mu$,
 (i.e.   $\lb_i\geq \mu_i$), the diagram $\lambda \ba \mu$
 denotes the diagram obtained by removing $\mu$ form $\lb$.
  The {\em skew-Schur polynomial ${\bf s}_{\lambda\ba \mu}$} associated with a Young
  diagram $\lambda \ba \mu\vdash n$ is a symmetric function in the
  variables $x_1,x_2,...$, (finite or infinite), 
  defined by
  \bean
  \tilde\gs_{\lambda \backslash \mu}(x_1,x_2,...)&=&
  \sum_{
  {{\mbox{\tiny semi-standard}}\atop
     {\mbox{\tiny skew-tableaux $P$  }   }}\atop
     {\mbox{\tiny of shape $\lb\backslash \mu$}   }
   }  \prod_i x_i ^{\# \{\mbox{times $i$
  appears in P}\}}\\&&\\
  &=& \det \left(  {\bf s}_{\lb_i-\mu_j-i+j}(t)\right) _{1\leq i,j\leq n}
  .
  \eean

 Similarly, for $\lb\ba \mu \vdash n$,
   $$
   \tilde{\bf s}_{\lb\ba \mu}(x_1,x_2,\ldots)= f^{\lb\ba \mu}~ x_1\ldots x_n+\ldots.
   $$

\item The {\em hook length} of the $i,j$th box is
defined
 by\footnote{$h^{\lb}_{ij}:=\lb_i +  \lb_j^{\top}-i-j+1$
 is the {\em hook length} of the $i,j$th box in the
 Young diagram; i.e., the number of boxes covered by the hook formed
 by drawing a horizontal line emanating from the center
 of the box to the right and a vertical line emanating
 from the center of the box to the bottom of the diagram.
  }
$h^{\lb}_{ij}:=\lb_i + \lb^{\top}_j-i-j+1$. Also define
  \bea h^{\lb}&:=&\prod_{(i,j)\in
\lb}h^{\lb}_{ij} \nonumber\\ &=& \frac{\prod_1^m
(m+\lb_i-i)!}{\Dt_m(m+\lb_1-1,\ldots,m+\lb_m-m)},~~\mbox{for}~
m \geq \lb_1^{\top}. \label{hook}\eea


\item   The {\em number of standard Young
   tableaux} of a given shape $\lambda=
   (\lambda_1\geq...\geq\lambda_m)$ is given
   by\footnote{using the Vandermonde determinant
   $\Dt_m(z_1,\ldots,z_m)=\prod_{1\leq i<j\leq m} (z_i-z_j)$
  }
  \bea f^{\lb}&=&\# \left\{\begin{array}{c}
  \mbox{standard ~tableaux  of  shape~
$\lb$}\\
\mbox{filled with the integers $1, \ldots, |\lb|$}
\end{array}
 \right\}\nonumber\\
&=&\mbox{coefficient of $x_1 \ldots x_n$ in
$\tilde\gs_{\lb}(x)$} \nonumber
\\&=&
\frac{|\lb|!}{u^{|\lb|}}~\tilde\gs_{\lb}(x)\Bigl|_{\sum_{i\geq
1} x_i^k=\dt_{k1} u}
 =\frac{|\lb|!}{u^{|\lb|}}\gs_{\lb}(u,0,0,\ldots)
\nonumber\\
&=&\frac{|\lb
|!}{h^{\lb}}  \nonumber\\
&=& |\lb |
!\det\left(\frac{1}{(\lb_{i}-i+j)!}\right)_{1\leq
i,j\leq m} 
\nonumber\\
  &=&  |\lb|! ~
  \frac{\Delta_{m}(m+\lb_1-1,\ldots,m
 +\lb_{m}-m)}
{\displaystyle{\prod_1^{m}} (m+\lb_i-i)!},
\quad\mbox{for any $m\geq\lb^{\top}_{1}$}.\nonumber \\
 \label{number standard tableaux} \eea
In particular, for any $m\geq\lb^{\top}_{1}$ and
arbitrary $u\in \BR$,

\be
  \gs_{\lb}(u,0,0,\ldots)
=u^{|\lb|}\frac{f^{\lb}}{{|\lb|}!}
 = u^{|\lb|}
  \frac{\Dt_m (m+\lb_1-1,\ldots, m+\lb_m-m)}
  {\prod_1^m (m+\lb_i-i)!} .
   \label{number standard tableaux1} \ee

\item   The {\em number of semi-standard Young
   tableaux} of a given shape $\lambda\vdash n $,
   filled with numbers $1$ to $q$ for $q \geq 1$:
  \bean
 \lefteqn{ \# \left\{\begin{array}{l}\mbox{semi-standard
tableaux of shape $\lb$} \\ \mbox{filled with numbers
from $1$ to $q$}
\end{array}
\right\}}\nonumber\\
 &=&
\tilde\gs_{\lb}(\overbrace{1,\ldots,1}^q,0,0,\ldots)
  \no\\\no\\
&=&\gs_{\lb}(q, \frac{q}{2},\frac{q}{3},\ldots)\nonumber\\
 &=&
\prod_{(i,j)\in\lb} \frac{j-i+q}{h^{\lb}_{i,j}}
\eean\bea &=& \left\{ \begin{array}{l}
  \displaystyle{
\frac{\Delta_q(q+\lb_1-1,\ldots,q+\lb_q-q)}
{\displaystyle{\prod^{q-1}_{i=1}}i!},~~~\mbox{when $ q
\geq \lb_1^{\top}$,}}\\
\\
 \displaystyle{   0},~~~\mbox{when $ q
< \lb_1^{\top}$,}
 \end{array} \right.\nonumber\\
 \nonumber\\
 \label{number semi-standard tableaux}\eea
 using the fact that
\be \prod_{(i,j)\in\lb} (j-i+q)=\frac {\prod_{i=1}^q
(q+\lb_i-i)!}{\prod_1^{q-1}i!}. \ee

\item Pieri's formula: given an integer $r\geq 0$ and
the Schur polynomial $\gs_{\mu} $, the following holds
\be
  \gs_{\lb}\gs_{r}=\sum_{\mu \backslash
\lb=\mbox{\tiny horizontal strip}\atop {\mbox{\tiny
and}~|\mu \backslash \lb|=r}}\gs_{\mu}.
 \label{Pieri}
  \ee
   Note $\mu
 \backslash \lb$ is an horizontal $r$-strip, when they
 interlace, $\mu_1\geq \lb_1\geq \mu_2\geq\ldots$, and
 $|\mu \backslash \lb|=r$.



\end{itemize}


\subsection{Robinson-Schensted-Knuth correspondence for
 generalized permutations}

Define the set of $p\times q$ integer matrices (see
\cite{Stanley1})
$$
\mbox{Mat}^{p,q}_n
 :=
 \left\{ W= (w_{ij})_{{1\leq i\leq p}\atop {1\leq j\leq q}}
,~~w_{ij}\in \BZ_{\geq 0}~~\mbox{and}~~\sum_{i,j}
w_{ij}=n
 \right\}
  $$

\begin{theorem} \label{Theo: RSK for GP}There is a 1-1 correspondence between
the following three sets:
$$
GP_n^{p,q}\Longleftrightarrow\left\{\begin{array}{l}
\mbox{two semi-standard Young tableaux $(P,Q)$,}\\
\mbox{of same, but arbitrary shape $\lb\vdash n$, filled}\\
 \mbox{resp. with integers $(1,\ldots,p)$ and $(1,\ldots,q)$}\\
%
\end{array}\right\}
 \Longleftrightarrow
    \mbox{Mat}^{p,q}_n
    .
$$
$$
\pi
  \hspace{0.5cm}
  \longleftrightarrow
 \hspace{1.9cm}
 (P,Q)
 \hspace{1.9cm}
 \longleftrightarrow \hspace{1.5cm}
 \hspace{0.5cm}
  W(\pi)=(w_{ij})_{\tiny\begin{array}{l}
                     1\leq i\leq p\\
                     1\leq j\leq q
                     \end{array}}
$$
 where
 $$w_{ij}=\#\left\{\mbox{times that~} \left(\begin{array}{c}
i\\ j
 \end{array}\right)\in\pi\right\}.
$$
 Therefore, we have\footnote{Use the notation ${\tilde \gs}_{\lb}(1^p)
 =\tilde\gs_{\lb}(\overbrace{1,\ldots,1}^p,0,0,\ldots)$.}
 \be
   \left(\begin{array}{c}
  pq+n-1\\ n
\end{array}\right)= \# \mbox{GP}_n^{p,q}
   =\sum_{\lb\vdash n}{\tilde \gs}_{\lb}(1^p)  {\tilde\gs}_{\lb}(1^q)
=\# \mbox{Mat}^{p,q}_n .\label{numberGP}\ee
%
Also, we have equality between the length of the longest
weakly increasing subsequence of the generalized
permutation $\pi$, the length of the first row of the
associated Young diagram and the weight of the optimal
path:
\be L(\pi)=\lb_1=L(W):=\max_{ \mbox{\tiny all such}\atop
                          \mbox{\tiny paths} }
\left\{\sum w_{ij},\quad\begin{array}{l}
\mbox{over right/down}\\
\mbox{paths starting from}\\
\mbox{entry $(1,1)$ to $(p,q)$}
\end{array}\right\}
\label{longest increasing GP} \ee
 \end{theorem}


{\medskip\noindent{\it Sketch of Proof:\/} } Given a
generalized permutation
$$
\pi =\left(\begin{array}{ccc}
i_1&\ldots&i_n\\
 \\
 j_1&\ldots&j_n
\end{array}\right)
,$$
the correspondence constructs two semi-standard Young
tableaux $P,Q$ having the same shape $\lambda$. This
construction is inductive. Namely, having obtained two
equally shaped Young diagrams $P_k,Q_k$ from
 $$
 \left(\begin{array}{ccc}
i_1&\ldots&i_k\\
 \\
 j_1&\ldots&j_k
\end{array}\right), ~~~1\leq k\leq n
 $$
 with the numbers $(j_1,...,j_k)$ in
the boxes of $P_k$ and the numbers $(i_1,...,i_k)$ in
the boxes of $Q_k$, one forms a new diagram $P_{k+1}$,
by creating a {\em new box in the first row of $P$,
containing the next number $j_{k+1}$ }, according to the
following rule:
\begin{description}
  \item[(i)] if $j_{k+1} \geq $ all numbers appearing
  in the first row of $P_k$, then one creates a new box
  with $j_{k+1}$ in that box to the right of the first
  row,
  \item[(ii)] if not, place $j_{k+1}$ in the box
  (of the first row) containing
  the smallest integer $> j_{k+1}$. The integer,
  which was in that box, then gets
  pushed down to the second row of $P_k$ according to
  the rule (i) or (ii), and the process starts afresh
  at the second row.
\end{description}
The diagram Q is a bookkeeping device; namely, add a box
(with the number $i_{k+1}$ in it) to $Q_k$ exactly at
the place, where the new box has been added to $P_k$.
This produces a new diagram $Q_{k+1}$ of same shape as
$P_{k+1}$.
The inverse of this map is constructed by reversing the
steps above.

Formula (\ref{numberGP}) follows from (\ref{number
semi-standard tableaux}).\qed


\begin{example} For $n=10$, $p=4$, $q=3$, \bean
\mbox{GP}_n^{p,q}\ni\pi&=&\left(\begin{array}{llllllllll}
1&1&1&2&2&2&3&3&4&4\\
 {2}& {3}& {3}&\underline{1}&
 \underline{2}&
\underline{2}&1&\underline{2}&1&\underline{3}
\end{array}\right)\qquad
,~~\mbox{with}~~L(\pi)=5\\
& & \\
& &
 \Updownarrow\hspace{25mm}
 \\
 & & \\
(P,Q)&=&\left(\overbrace{\squaresize .4cm \thickness
.01cm \Thickness .07cm \Young{
 1&1&1&2&3\cr
 2&2&2\cr 3&3\cr}}^5,
 \squaresize .4cm \thickness .01cm \Thickness .07cm \Young{
 1&1&1&3&4\cr 2&2&2 \cr
 3&4\cr
}\right)\qquad L(\pi)=\lb_1=5\\
& &
 \Updownarrow\\&&\\
 W&=& \left(\begin{array}{ccccc}
{\!\!\scriptstyle
0}{\!\!\!\!\bigcirc\!\!}&\rg&{\!\!\scriptstyle
1}{\!\!\!\!\bigcirc\!\!}& &2\\
 & &\downarrow \\
1& &{\!\!\scriptstyle 2}{\!\!\!\!\bigcirc\!\!}& &0\\
 & &\downarrow \\
1& &{\!\!\scriptstyle
1}{\!\!\!\!\bigcirc\!\!}&\rg&{\!\!\scriptstyle
0}{\!\!\!\!\bigcirc\!\!}\\
 & & & & \downarrow \\
1& &0& &{\!\!\scriptstyle 1}{\!\!\!\!\bigcirc\!\!}
\end{array}\right),
 \mbox{~~with} ~L(\pi)=\sum_{(i,j)\in\{\mbox{path}\}}w_{ij}=
5\\
\eean
\end{example}
%
%
 The RSK algorithm proceeds as follows:

 \newpage



{\footnotesize $$
\begin{array}{ccccccccccccccccccccccccccc}
\mbox{adding}& &
\left(\begin{array}{c}2\\1\end{array}\right)&&
\left(\begin{array}{c}2\\2\end{array}\right)&&
\left(\begin{array}{c}2\\2\end{array}\right)&&\left(\begin{array}{c}3\\1
\end{array}\right)
   \\
P&\begin{array}{l}2\,\,3\,\,3\end{array}& &\begin{array}{l}1\,\, 3\,\, 3\\
2\end{array}&
&\begin{array}{l}1\,\,2\,\,3\\2\,\,3\end{array}&
&\begin{array}{l}1\,\,2\,\,2\\
2\,\,3\,\,3\end{array}& &
\begin{array}{l}1\,\,1\,\,2\\2\,\,2\,\,3\\3\end{array}\\
  & &\Rightarrow& &\Rightarrow& &\Rightarrow& &\Rightarrow\\
\\
Q&\begin{array}{l}1\,\,1\,\,1\end{array}& &\begin{array}{l}1\,\,1\,\,1\\
2\end{array}&
&\begin{array}{l}1\,\,1\,\,1\\2\,\,2\end{array}&
&\begin{array}{l}1\,\,1\,\,1\\
2\,\,2\,\,2\end{array}& &
\begin{array}{l}1\,\,1\,\,1\\2\,\,2\,\,2\\3\end{array}\\
\end{array}
$$}

\vspace*{2cm}

{\footnotesize$$
\begin{array}{ccccccc}
\left(\begin{array}{c}3\\2\end{array}\right)&
&\left(\begin{array}{c}4\\1\end{array}\right)&
&\left(\begin{array}{c}4\\3\end{array}\right)\\
\\
   &
    \begin{array}{l}1\,\,1\,\,2\,\,2\\
1\,\,2\,\,3\\
3\end{array}
& & \begin{array}{l}1\,\,1\,\,1\,\,2\\
1\,\,2\,\,2\\
3\,\,3\end{array} & & \begin{array}{l}1\,\,1\,\,1\,\,2\,\,3\\
1\,\,2\,\,2\\
3\,\,3\end{array} \\
\Rightarrow& &\Rightarrow& &\Rightarrow\\
& \begin{array}{l}1\,\,1\,\,1\,\,3\\
2\,\,2\,\,2\\
3\end{array} & & \begin{array}{l}1\,\,1\,\,1\,\,3\\
2\,\,2\,\,2\\
3\,\,4\end{array} & & \begin{array}{l}1\,\,1\,\,1\,\,3\,\,4\\
2\,\,2\,\,2\\
3\,\,4  \end{array}
  \end{array}
~~~~~~~~~=\left(\begin{array}{l}{P}\\{Q}\end{array}\right)
$$}

\noindent yielding the set $(P,Q)$ of semi-standard
Young tableaux above.


\subsection{The Cauchy identity}

\begin{theorem} \label{Theorem 1.2}Using the customary change of variables
$\sum_{k\geq 1}x^i_k=it_i,~~\sum_{\ell\geq
1}y^i_{\ell}=is_i$, we have
$$
 \prod_{i,j\geq 1}\frac{1}{1-x_iy_j} =\sum_{\lb}\tilde
\gs_{\lb}(x)\tilde \gs_{\lb}(y) =\sum_{\lb} \gs_{\lb}(t)
\gs_{\lb}(s)= e^{\sum_{i\geq 1}i t_is_i}
 .$$

\end{theorem}

\proof On the one hand, to every $\pi\in
GP=\bigcup_{n,p,q} GP_n^{p,q}$, we associate a monomial,
as follows
\be
 \pi\lrg\prod_{\left(\begin{array}{l}i\\j\end{array}\right)\in\pi}x_iy_j
\hspace{2cm}\mbox{(with multiplicities)}
 \ee
 Therefore, taking into account the multiplicity of $
 \left(\begin{array}{l}i\\j\end{array}\right)\in \pi$,

\be
  \sum_{\pi\in
GP}\prod_{\left(\begin{array}{l}i\\j\end{array}\right)
 \in\pi}x_iy_j= \prod_{i,j \geq 1}(1+x_iy_j+x^2_i
y^2_j+x_i^3y_j^3+\ldots)=\prod_{i,j\geq
1}\frac{1}{1-x_iy_j}. \label{product}\ee
One must think of the product on the right hand side of
(\ref{product}) in a definite order, as follows:
 \bean
  \prod_{i,j\geq
1}\frac{1}{1-x_iy_j}&=&
  ~~~\frac{1}{1-x_1y_1}~
  \frac{1}{1-x_1y_2}~
  \frac{1}{1-x_1y_3}\ldots \\
 &&
 \times \frac{1}{1-x_2y_1}~
  \frac{1}{1-x_2y_2}~
  \frac{1}{1-x_2y_3}\ldots \\
  &&\\&& \times\ldots,
 \eean
 and similarly for the expanded version. Expanding out all the products,
 \be
  \prod_{ j \geq 1}(1+x_1y_j+x^2_1
y^2_j+x_1^3y_j^3+\ldots)\prod_{j \geq 1}(1+x_2y_j+x^2_2
y^2_j+x_2^3y_j^3+\ldots)\ldots,\label{1.5.3} \ee leads
to a sum of monomials, each of which can be interpreted
 as a generalized permutation, upon respecting the prescribed order. Vice-versa each
 generalized permutation can be found among those
 monomials.
   As an example
illustrating identity (\ref{product}), the monomial
$x_1y_2x_1^2y_3^2~x_3^2y_2^2~x_3y_3~x_4y_1~x_4y_2$,
appearing in the expanded version of (\ref{1.5.3}), maps
into the generalized permutation
 $$\left(
\begin{array}{cccccccc}
 1&1&1&3&3&3&4&4\\
 2&3&3&2&2&3&1&2
\end{array}\right),
 $$
 and vice-versa. On the other hand, to every couple of
semi-standard Young tableaux $(P,Q)$, we associate
%
%
$$
(P,Q)\lrg\prod_ix_i^{\#\mbox{\tiny{times $i$ appears in
$Q$}}} \prod_jy_j^{\#\mbox{\tiny{times $j$ appears in
$P$}}} .
 $$
Therefore, the Robinson-Schensted-Knuth construction,
mapping the generalized permutation $\pi$ into two
semi-standard Young tableaux $(P,Q)$ of same shape
$\lb$, implies
 $$
\prod_{\left(\begin{array}{l}i\\j\end{array}\right)\in\pi}x_iy_j=
\prod_ix_i^{\#\mbox{\tiny{times $i$ appears in $Q$}}}
\prod_jy_j^{\#\mbox{\tiny{times $j$ appears in $P$}}}.
$$
Then, summing up over all $\pi \in GP$,
 using the fact that RSK is a bijection, and using the
 definition of Schur polynomials,
 one computes
\bean \lefteqn{
 \sum_{\pi\in
GP}\prod_{\left(\begin{array}{l}i\\j\end{array}\right)
 \in\pi}x_iy_j=}\\
 &=& \sum_{{\mbox{\tiny all~}(P,Q)\mbox{\tiny
with}}\atop \mbox{\tiny shape~}P=\mbox{\tiny shape~}{Q}}
\prod_ix_i^{\#\mbox{\tiny{times $i$ appears in $Q$}}}
\prod_jy_j^{\#\mbox{\tiny{times $j$ appears in $P$}}}
\\
& & \\
&=& \sum_{\lb}
 \sum_{{\mbox{\tiny all~}(P,Q)\mbox{\tiny
with}}\atop \mbox{\tiny shape~}P=\mbox{\tiny
shape~}{Q=\lb}} \prod_ix_i^{\#\mbox{\tiny{times $i$
appears in $Q$}}} \prod_jy_j^{\#\mbox{\tiny{times $j$
appears in $P$}}}
 \\
&& \\
  &=&\sum_{\lb}\left(
  \sum_{  \mbox{\tiny all $Q$ with}\atop \mbox{\tiny{shape
$Q=\lb$}}     } \prod x_i^{ \#\mbox{\tiny{times $i$
appears in $Q$}} } \right) \left( \sum_{
  \mbox{\tiny all $P$ with}\atop \mbox{\tiny{shape
$P=\lb$}}     } \prod y_j^{ \#\mbox{\tiny{times $j$
appears in $P$}} }
\right)\\
& &\\
&=&\sum_{\lb}\tilde\gs_{\lb}(x)\tilde\gs_{\lb}(y), \eean
using the definition (\ref{Schur polynomial}) of the
  Schur polynomial. %
The proof is finished by observing that
\bean \sum_{i\geq 1}it_is_i&=&\sum_{i\geq
1}i\frac{\displaystyle{\sum_{k\geq 1}}x_k^i}{i}
\frac{\displaystyle{\sum_{\ell\geq 1}}y_{\ell}^i}{i}\\
 & & \\
&=&\sum_{k,\ell\geq 1}\sum_{i\geq 1}\frac{(x_ky_{\ell})^i}{i}\\
 & & \\
&=&\log\prod_{k,\ell\geq 1}(1-x_ky_{\ell})^{-1}, \eean
ending the proof of Theorem \ref{Theorem 1.2}.\qed


\subsection{Uniform Probability on Permutations, Plancherel Measure
and Random Walks}

\subsubsection{Plancherel measure}\label{sect 1.5.1}

In this section, one needs
$$
\mbox{Mat}^{n,n}_n(0,1)
 :=
 \left\{ W= (w_{ij})_{{1\leq i,j\leq n}}
,
~\left|\begin{array}{l}\mbox{with
exactly one $1$}\\
\mbox{in each row and column}\\
\mbox{and otherwise all zeros}\end{array}\right.
 \right\}
  $$
See \cite{Stanley1,BOO,Borodin1,Bo-Ok} and references
within.
\begin{proposition}For {\em permutations}, we
have a 1-1 correspondence between
$$
S_n  \Longleftrightarrow\left\{\begin{array}{c}
\mbox{two standard Young tableaux $(P,Q)$,}\\
\mbox{of same shape $\lb$ and size $n$,}\\
\mbox{each filled with numbers $(1,\ldots,n)$}\\
\end{array}\right\}
 \Longleftrightarrow
    \mbox{Mat}_n^{n,n}( 0,1 ).
$$

\be \pi_n
 \hspace{1.5cm}\longleftrightarrow
 \hspace{1.5cm}
 (P,Q)
 \hspace{1.5cm}
 \longleftrightarrow \hspace{1.5cm}
  W(\pi)=(w_{ij})_{i,j\geq 1},
\label{RSK-correspondence(permutations)}\ee
%
Uniform probability $P^{n }$ on $S_n $ induces a
probability $\tilde P^{n }$ (Plancherel measure) on
Young diagrams $\BY^{ }_n$, given by ($m:=\lb_1^{\top}$)
 \bea
 \tilde P^{n}\left( \lb \right)=
\frac{1}{n!} \#\left\{
\begin{array}{l}
  \mbox{permutations\!\!}\\
  \mbox{leading to }\\
   \mbox{shape $\lb$}
  \end{array}
  \right\}
 &=&
  \frac{(f^{\lb})^2}{n!} \no\\
  &=&n! ~{\bf s}_{\lb}(1,0,\ldots)^2
   \no\\ \no\\
  &=&
   n!  \frac{\Dt_m
(m+\lb_1-1,\ldots, m+\lb_m- m)^2}{\displaystyle
\left(\prod_1^m   (m+\lb_i-i)!\right)^2  }
 \no\\
 \no\\&&\label{permutationPlancherel}
  \eea
%
%
%
%
%
and so
$$\#S_n=\sum_{\lb\vdash n}\left(f^{\lb}\right)^2=n!  \\
.$$
Finally, the length of the longest increasing
subsequence in permutation $\pi_n$, the length of the
first row of the partition $\lb$ and the weight of the
optimal path $L(W)$ in the percolation matrix $ W(\pi)$
are all equal: $$ L(\pi_n)=\lb_1=L(W):=\max_{
\mbox{\tiny all such}\atop
                          \mbox{\tiny paths} }
\left\{\sum w_{ij},\quad\begin{array}{l}
\mbox{over right/down}\\
\mbox{paths starting from}\\
\mbox{entry $(1,1)$ to $(n,n)$}
\end{array}\right\}.$$
Hence
$$P^n(L(\pi)\leq \ell)=
 \sum_{\begin{array}{c}
  \lb \in \BY_n \\
   \lb_1 \leq \ell
   \end{array}
  }\frac{(f^{\lb})^2}{n!}=
n!\sum_{\begin{array}{c}
  \lb \in \BY_n \\
   \lb_1 \leq \ell
   \end{array}
  } ~{\bf s}_{\lb}(1,0,\ldots)^2
. %
$$

\end{proposition}

\proof A \underline{\em permutation} is a generalized
permutation, but with integers $i_1,\ldots,i_n$ and
$j_1,\ldots,j_n$ all distinct and thus both tableaux $P$
and $Q$ are standard.

Consider now the uniform probability $P^{n }$ on {\em
permutations} in $S_n$; from the RSK correspondence
%
%
%
%
%
we have the one-to-one correspondence, given a fixed
partition $\lb$,
$$
 \left\{
\begin{array}{l}
  \mbox{permutations}\\
  \mbox{leading to}\\
  \mbox{the shape $\lb$}
  \end{array}
  \right\}
   \Longleftrightarrow
 \left\{
\begin{array}{l}
  \mbox{standard tableaux}\\
   \mbox{of shape $\lb$,}\\
  \mbox{filled with integers}\\
  \mbox{$1,\ldots,n$}
  \end{array}
  \right\}
   \times
  \left\{
\begin{array}{l}
  \mbox{standard tableaux}\\
  \mbox{of shape $\lb$} \\
  \mbox{filled with integers}\\
  \mbox{ $1,\ldots,n$}
  \end{array}
  \right\}
$$
 and thus, using formulae (\ref{number standard tableaux})
  and (\ref{number semi-standard tableaux}) and noticing that
 $\tilde\gs_{\lb}(1^q)=0$ for $\lb_1^{\top}>q$,
 $$
 \tilde P^{n }\left( \lb \right)=
\frac{1}{ n!} \#\left\{
\begin{array}{l}
  \mbox{permutations leading}\\
  \mbox{to the shape $\lb$}
  \end{array}
  \right\}
 =\frac{(f^{\lb})^2}{ n!},~~~~~~~
 \lb\in\BY_n,
$$
with
$$
\sum _{\lb\in\BY_n } \tilde P^{n}(\lb)=1.
$$
Formula (\ref{permutationPlancherel}) follows
immediately from the explicit values (\ref{number
standard tableaux}) and (\ref{number semi-standard
tableaux}) for $f^{\lb}$. \qed

\begin{example} For permutation $
\pi_5 =\left(\begin{array}{cccccc}
1&2&3&4&5\\
 5&\underline{1}&4&\underline{3}&2
\end{array}\right)\in S_5$, the RSK algorithm gives

{\footnotesize $$
\begin{array}{cccccccccccccccccccccccc}
P\Longrightarrow&&&& 5 &&&& 1 &&&& 1&4 &&&& 1&3 &&&& 1&2 \\
                &&&&   &&&& 5 &&&& 5&  &&&& 4&  &&&& 3&  \\
                &&&&   &&&&   &&&&  &  &&&& 5&  &&&& 4&  \\
                &&&&   &&&&   &&&&  &  &&&&  &  &&&& 5&  \\
                &&&&   &&&&   &&&&  &  &&&&  &  &&&&  &  \\
Q\Longrightarrow&&&& 1 &&&& 1 &&&& 1&3 &&&& 1&3 &&&& 1&3 \\
                &&&&   &&&& 2 &&&& 2&  &&&& 2&  &&&& 2&  \\
                &&&&   &&&&   &&&&  &  &&&& 4&  &&&& 4&  \\
                &&&&   &&&&   &&&&  &  &&&&  &  &&&& 5&
\end{array}
$$}

Hence

$ \displaystyle{\pi_5}$

$\Downarrow$

\vspace*{-2cm}
 $$
  \displaystyle{
(P,Q)=\left(\!\!\!\begin{array}{c}
\overbrace{\,\,\,\,}^2\\
\MAT{1} \squaresize .4cm \thickness .01cm \Thickness
.06cm \Young{
 1&d 2 \cr
 3\cr  4\cr d 5\cr}
\mat ,\!\! \\
\mbox{\tiny standard}
\end{array}
\begin{array}{c}
\\
\MAT{1}\squaresize .4cm \thickness .01cm\Thickness .07cm
\Young{
 1& 3\cr
2\cr 4\cr 5\cr}
\mat\\
\mbox{\tiny standard}
\end{array}
\!\!\!\right)}
 \Longrightarrow\left(\begin{array}{ccccccccc}
{\!\!\scriptstyle 0}{\!\!\!\!\bigcirc\!\!}& &0& &0& &0& &1\\
\downarrow& & & & & & & & \\
{\!\!\scriptstyle 1}{\!\!\!\!\bigcirc\!\!}& &0& &0& &0& &0\\
\downarrow& & & & & & & & \\
{\!\!\scriptstyle
0}{\!\!\!\!\bigcirc\!\!}&\rg&{\!\!\scriptstyle
0}{\!\!\!\!\bigcirc\!\!}&\rg&{\!\!\scriptstyle
0}{\!\!\!\!\bigcirc\!\!}&\rg&{\!\!\scriptstyle
1}{\!\!\!\!\bigcirc\!\!}&\rg&{\!\!\scriptstyle
0}{\!\!\!\!\bigcirc\!\!}\\
 & & & & & & & &\downarrow \\
0& &0& &1& &0& &{\!\!\scriptstyle 0}{\!\!\!\!\bigcirc\!\!}\\
 & & & & & & & &\downarrow \\
0& &1& &0& &0& & {\!\!\scriptstyle
0}{\!\!\!\!\bigcirc\!\!}
\end{array}\right)
 $$

\end{example}

\bigbreak

\remark The Robinson-Schensted-Knuth correspondence has
the following properties
\begin{itemize}
  \item $\pi\mapsto (P,Q)$, then $\pi^{-1}\mapsto (Q,P)$
  \item length (longest increasing subsequence of $\pi$) $=\#$ (columns
  in $P$)
  \item length (longest decreasing subsequence of $\pi$) $=\#$ (rows
  in $P$)
  \item $\pi^2=I$, then $\pi\mapsto (P,P)$
  \item $\pi^2=I$, with $k$ fixed points, then $P$ has exactly $k$
  columns of odd length.
  \be \ee

  \end{itemize}


\subsubsection{Solitaire game}

With Aldous and Diaconis \cite{AD}, consider a deck of
cards $1,\ldots,n$ thoroughly shuffled and put those
cards one at a time into piles, as follows:
\begin{enumerate}
\item[(1)]  a low card may be placed on a higher card,
or can be put into a new pile to the right of the
existing pile. \item[(2)] only the top card of the pile
is seen. If the card which turns up is higher than the
card showing on the table, then start with that card a
new pile to the right of the others.
\end{enumerate}

 \noindent \underline{\em Question}: What is the optimal strategy
which minimizes the number of piles?

\noindent\underline{\em Answer}: put the next card
always on the leftmost possible pile!

  \example  Consider a deck of 7 cards, appearing in the order
3, 1, 4, 2, 6, 7, 5. The optimal strategy is as follows:

$$
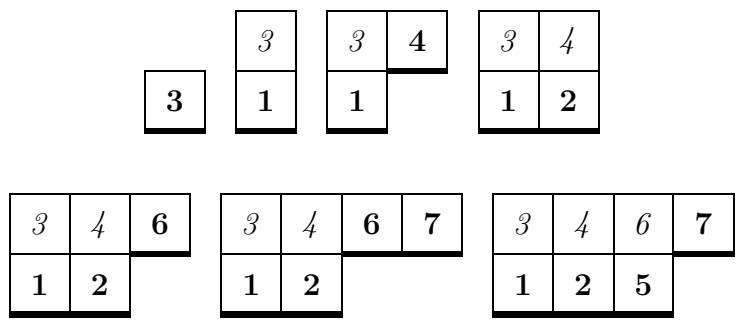
\squaresize .8cm \thickness .01cm \Thickness .07cm
\Young{  d   {\bf  3}\cr}
\quad
\squaresize .8cm \thickness .01cm \Thickness .07cm
\Young{ {\it 3} \cr d {\bf 1}\cr}
\quad
  \squaresize
.8cm \thickness .01cm \Thickness .07cm \Young{
  {\it 3}
  &d {\bf 4}\cr
d {\bf 1}\cr}
\quad
 \squaresize .8cm \thickness .01cm \Thickness
.07cm \Young{
 {\it 3}
 &{\it 4}
  \cr
 d {\bf 1}&d {\bf 2}\cr}
 $$

$$
\squaresize .8cm \thickness .01cm \Thickness .07cm
\Young{
 {\it 3}
  &{\it 4}
  &d {\bf 6}\cr
 d {\bf 1}&d {\bf 2}\cr}\quad
\squaresize .8cm \thickness .01cm \Thickness .07cm
\Young{
 {\it 3}
  &{\it 4}
  &d {\bf 6}&d {\bf 7}\cr
 d {\bf 1}&d {\bf 2}\cr}\quad
\squaresize .8cm \thickness .01cm \Thickness .07cm
\Young{
 {\it 3}
 &{\it 4}
  &{\it 6}
   &d {\bf 7}\cr
 d {\bf 1}&d {\bf 2}&d {\bf 5}\cr}
 $$
This optimal strategy leads to 4 piles!
For a deck of 52 cards, you will find in the average
between 10-13 piles and having 9 piles or less occurs
approximately $5\%$ of the times. It turns out that, for
a given permutation,
$$\mbox{number of piles}=\mbox{length of longest
increasing sequence}.$$








\subsubsection{Anticipating large $n$ asymptotics}

Anticipating the results in section \ref{subsect 4.2}
and \ref{subsect 8.5}, given the permutations of
$(1,\ldots,n)$, given a percolation matrix of size $n$
and given a card game of size $n$, the numbers fluctuate
about $2\sqrt{n}$ like
\bean
 L(\pi_n)
  &=&{\mbox{length of the longest increasing
subsequence}} \\
 &=&\mbox{weight of the optimal path in the
$(0,1)$-percolation matrix}
 \\
 &=& \mbox{number of piles in
the solitaire game}\\
 &\simeq& 2\sqrt{n} + n^{1/6} ~{\cal F}
 ,\eean
  where $ {\cal F}$ is a probability distribution,
  the Tracy-Widom distribution, with
   $$E(\FR)=-1.77109 ~~\mbox{and}~~ \sigma(\FR)=0.9018.$$
In particular
 {for $n=52$ cards,
 \bean E( L(\pi_{52})) &\simeq & 2 \sqrt{52}+(52)^{1/6}
 (-1.77109)= 11.0005
 .\eean}
 The Tracy-Widom distribution will be discussed
 extensively in section \ref{subsect 8.5}.

  \subsubsection{A transition probability and Plancherel measure}

\begin{proposition} (\cite{Vershik and Kerov1,BOO})
$P_n$ on $\BY_n$ can be constructed from $P_{n-1}$ on
$\BY_{n-1}$, by means of a transition probability, as
follows

$$
P_n(\mu)=\sum_{\lb\in\BY_{n-1}}P_{n-1}(\lb)
 p(\lb,\mu),\quad\mu\in
\BY_n
$$
where

$$
p(\lb,\mu):=\left\{\begin{array}{ll}
\frac{f^{\mu}}{f^{\lb}}\frac{1}{\vert\mu\vert}&\mbox{if
$\lb\in\BY_{n-1}$
and $\mu\in\BY_n$ are}\\
&\mbox{such that $\mu=\lb+ \square $}\\
0&\mbox{otherwise}
\end{array}
\right.
$$
is a transition probability, i.e.

$$
\sum _{\mu\in\BY_n\atop{\mu=\lb+\square}}
p(\lb,\mu)=1,\quad\mbox{for fixed $\lb$.}
$$
\end{proposition}

\proof Indeed, for fixed $\mu$, one computes

\bean \sum_{\lb\in\BY_{n-1}}P_{n-1}(\lb)
 p(\lb,\mu)
   &=&
\sum_{\lb\in\BY_{n-1} \atop \mbox{\tiny such
that~}\mu=\lb+\square\in\BY_n}
\frac{(f^{\lb})^2}{(n-1)!}
 \left(\frac{f^{\mu}}{f^{\lb}}
 \frac{1}{\vert\mu\vert}\right)
 \\
&=&\frac{f^{\mu}}{n!}\sum_{\lb\in\BY_{n-1} \atop
\mbox{\tiny such
that~}\mu=\lb+\square\in\BY_n}f^{\lb}\\
&=&\frac{(f^{\mu})^2}{n!}
 = P_n(\mu). \eean Indeed, given a
standard tableau of shape $\lb$, filled with
$1,\ldots,n-1$, adjoining a box to $\lb$ such as to form
a Young diagram $\mu$ and putting $n$ in that box yield
a new standard tableau (of shape $\mu$).

That $p(\lb,\mu)$ is a transition probability follows
from Pieri's formula (\ref{Pieri}), applied to $r=1$,
upon putting $t_i=\dt_{1i}$:
$$
\sum_{\mu\in\BY_n~\mbox{\tiny with}\atop{\mu=\lb+
\square}}p(\lb,\mu)=\sum_{\lb\in \BY_{n-1}~\mbox{\tiny
with}\atop{\lb +
\square\in\BY_n}}\frac{f^{\lb+\square}}{\vert\lb
+\square\vert} \frac{1}{f^{\lb}}
 =\sum_{\lb\in \BY_{n-1}~\mbox{\tiny
with}\atop{\lb +
\square\in\BY_n}}\frac{f^{\lb+\square}}{\vert\lb
+\square\vert !} \frac{|\lb|!}{f^{\lb}}=1.
$$\qed


\begin{corollary}\label{cor:Pieri} The following probability
$$
P_n(\mu_1\leq x_1,\ldots,\mu_k\leq x_k)
$$
decreases, when $n$ increases.
\end{corollary}

\proof Indeed
\bean
\lefteqn{P_n(\mu_1\leq x_1,\ldots,\mu_k\leq x_k)}\\
&=&\sum_{\mu\in\BY_n\atop{\mbox{\tiny all~}\mu_i\leq
x_i}}\sum_{\lb\in\BY_{n-1}}P_{n-1}(\lb)p(\lb,\mu)\\
& & \\
&=&\renewcommand{\arraystretch}{0.5}
\begin{array}[t]{c}
\sum\\
{\scriptstyle \lb\in\BY_{n-1},\mu\in\BY_n}\\
{\scriptstyle \mbox{\tiny such that}}\\
{\scriptstyle \mu=\lb+\square}\\
{\scriptstyle \mbox{\tiny all~} \mu_i\leq x_i}
\end{array}
\renewcommand{\arraystretch}{1}P_{n-1}(\lb)p(\lb,\mu)\\
&\leq&\renewcommand{\arraystretch}{0.5}
\begin{array}[t]{c}
\sum\\
{\scriptstyle \lb\in\BY_{n-1}}\\
{\scriptstyle \mbox{\tiny such that}}\\
{\scriptstyle \mu=\lb+\square}\\
{\scriptstyle \mbox{\tiny all~} \mu_i\leq x_i}
\end{array}
\renewcommand{\arraystretch}{1}P_{n-1}(\lb)=P_{n-1}(\lb_1\leq
x_1,\ldots,\lb_k\leq x_k),\eean
 proving Corollary \ref{cor:Pieri}.\qed

\newpage

  \subsubsection{Random walks}

Consider $m$ random walkers in $\BZ$, starting from
distinct points $x:=(x_1<\ldots<x_m)$, such that, at
each moment, only one walker moves either one step to
the left or one step to the right. Notice that $m$
walkers in $\BZ$, obeying this rule, is tantamount to a
random walk in $\BZ^m$, where at each point the only
moves are
 $$
 \pm e_1,\ldots, \pm e_m,
 $$
with all possible moves equally likely. That is to say
the walk has at each point $2m$ possibilities and
 thus at time $T$ there are $ (2m)^{T} $ different paths. Denote by $P_x$ the probability for such a walk,
 where $x$ refers to the initial condition. {\em Requiring these walks
  not to intersect turns out to be closely
 related to the problem of longest increasing subsequences in random permutations}, as is
 shown in the Proposition below. For skew-partitions,
 see section \ref{subsection1.3}. For references, see
 \cite{Stanley1,AvM-random walks}.

\begin{proposition}

\be 
 P_x\left( \begin{array}{l}
 \mbox{that $m$ walkers in $\BZ$,}\\
 \mbox{reach $y_1<\ldots <y_m$}\\
 \mbox{in $T$ steps, without}\\
 \mbox{ever intersecting}
 \end{array}
   \right)
   =
   \frac{1}{(2m)^{T}}
   \left(  {T}
    \atop{T_L~~ T_R }
    \right)
   \sum_{
  {{\lb~\mbox{\tiny with}~\lb \supset \mu,\nu}\atop
 {   { {|\lb\ba \mu|=T_L}\atop{|\lb\ba \nu|=T_R}}
 \atop {\lb_1^{\top} \leq m}
   }    }
   }
     f^{\lb\ba \mu}  f^{\lb\ba \nu}
\label{1.5.6}\ee
%
 where $\mu$, $\nu$ are fixed partitions defined by the points
 $x_i$ and $y_i$,
 \bean \mu_k &=& k-1-x_k,~~~
        \nu_k = k-1-y_k  \\
 T_L&=&\frac{1}{2} (T+\sum_1^m(x_i-y_i))=
   \frac{1}{2} (T-|\mu|+|\nu|)\\
      T_R&=&\frac{1}{2} (T-\sum_1^m(x_i-y_i))=
   \frac{1}{2} (T+|\mu|-|\nu|)\\
   T&=&T_L+T_R,~~ \sum_1^m (x_i-y_i)=T_L-T_R.
 \eean



\noindent In particular, close packing of the walkers at
times $0$ and $T$ implies \bea
 \lefteqn{P_{1,\ldots,m}\left( \begin{array}{c}
 \mbox{that $m$ walkers in $\BZ$ reach}\\
 \mbox{$1,\ldots,m$ in $2n$ steps,}\\
 \mbox{without ever intersecting}
 \end{array}
   \right)}\no\\
   &&\hspace{2cm}=\frac{1}{(2m)^{2n}} \left(  {2n}\atop{n}\right)
   \sum_{{\lb\vdash n}\atop{\lb_1^{\top} \leq m}}
   \left( f^{\lb}\right)^2
   \no\\ \no\\
 & &\hspace{2cm}=\frac{(2n)!}{ n! } ~~\frac{\#\{
    \pi_n \in S_n  ~~\mbox{such that}~~
     L(\pi_n )\leq m \}}{(2m)^{2n}}
 \no\\ \label{1.5.7}\eea

\end{proposition}

\noindent$\hspace*{5.6cm} {1~2~~3~4~~5~6 }$
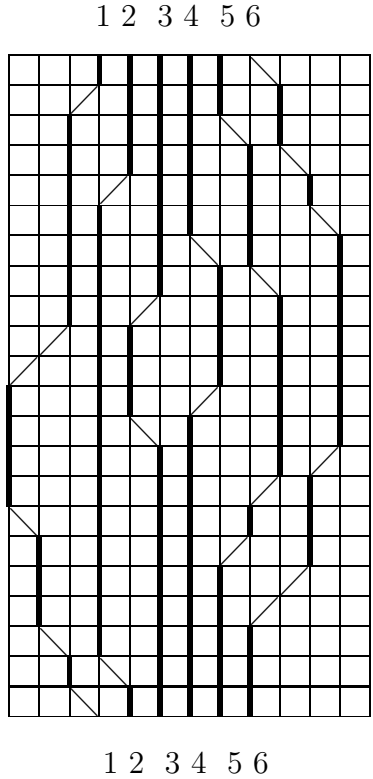
$$
 \squaresize 0.4cm \thickness .01cm\Thickness
.06cm \Young{
 & &r &r &r &r &r & &\diagdown & & &\cr
 & &\diagup &r &r&r &r & &r & & &\cr
 & r& & r&r & r& & \diagdown& r& & &\cr
 & r& & r&r & r& & r& & \diagdown& &\cr
 & r&  &\diagup& r&r & &r &  &r & &\cr
 & r& r& & r& r& & r& & &\diagdown &\cr
 & r&r & &r & &\diagdown &r & & &r &\cr
 & r& r& &r & &r & &\diagdown & & r& \cr
 & r& r& &\diagup & &r & &r & &r &\cr
 &\diagup &r &r & & &r & &r & & r &\cr
  \diagup& & r&r &  & &r & &r & &r &\cr
l & & r& r& & &\diagup & &r & &r &\cr l & & r&
&\diagdown&r & & & r& & r& \cr
 l& & r& & r& r& & & r& &\diagup &\cr
 l& & r&  &r &  r&  & & \diagup& r& &\cr
 \diagdown& & r& &r &r & & r& &r & &\cr
 r& & r& &r &r & &\diagup  & & r& &\cr
 r& & r& &r &r & r& & &\diagup & &\cr
 r& & r& &r &r & r& &\diagup  & & &\cr
 & \diagdown&r &  & r& r& r& r& & & &\cr
 &  & l&\diagdown & r& r& r& r& & & &\cr
 & &\diagdown&r & r& r& r& r& & & &  \cr
 }
$$
$\hspace*{5.7cm} {1~2~~3~4~~5~6 }$

{\it Six non-intersecting walkers leaving from and
returning to $1,\ldots,6$.}

\vspace{1cm}

\proof 
{\em Step 1}: Associate to a given walk
 a sequence of $T_L$ $L$'s and $T_R$ $R$'s:
 \be
 L~~R~~R~~{R~~L}~~R~~L~~L~~R~\ldots ~R~,
 \label{LR's}\ee
thus recording the nature of the move, left or right, at
the first instant, at the second instant, etc...

If the $k^{th}$ walker is to go from $x_k$ to $y_k$,
then
$$
y_k-x_k=\#\left\{ \begin{array}{c}
  \mbox{right moves }\\
  \mbox{for $k^{th}$ walker}
  \end{array} \right\}
  -
  \#\left\{ \begin{array}{c}
  \mbox{left moves }\\
  \mbox{for $k^{th}$ walker}
   \end{array} \right\}
  $$
  and so, if
  $$
  T_L:=\#\left\{ \begin{array}{c}
  \mbox{left moves }\\
  \mbox{for all $m$ walkers}
  \end{array} \right\}
~~~ \mbox{and} ~~~
  T_R:=\#\left\{ \begin{array}{c}
  \mbox{right moves }\\
  \mbox{for all $m$ walkers}
  \end{array} \right\},
 $$
 we have, since at each instant exactly one walker moves,
  $$
  T_R+T_L= T  ~~\mbox{and}~~
  T_R-T_L= \sum_1^m (y_k-x_k),
  $$
  from which
  $$
  T_{\left\{L\atop R\right\}}=\frac{1}{2} \left( T\pm \sum_1^m (x_k-y_k)\right)
  .$$
  Next, we show there is a {\em
canonical} way to map a walk, corresponding to
(\ref{LR's}) into one with left moves only at instants
$1,\ldots,T_L$ and then right moves at instants
$T_L+1,\ldots,T_L+T_R=T$, thus corresponding to a
sequence
 \be
 \overbrace{L~~L~~L~\ldots ~L}^{T_L}~~\overbrace{R~~R~~R~\ldots
 ~R}^{T_R}~.
 \ee
This map, originally found by Forrester
 \cite{Forrester}
 takes on here a slightly different (but canonical) form.
 Indeed, in a typical sequence, as (\ref{LR's}),
 \be
 L~~R~~R~~\underline{R~~L}~~R~~L~~L~~R~\ldots ~R~,
 \ee
 consider the first sequence $R~L$ (underlined) you encounter,
 in reading from left to right. It corresponds to one of the following three
 configurations (in the left column),
%
%
\bean
 \begin{array}{l}
L\\
R \end{array}
 {\displaystyle {\diagdown}\atop{~~~
\big|}}~~{\displaystyle {\big|}\atop{ \big|}}
       ~~{\displaystyle {\big|}\atop{ \big|}}
       ~~{\displaystyle {~~~ \big|}\atop{ \diagup}}
~~~~~~~~~~
  &\Rightarrow&
~~~~~~~~~~
 \begin{array}{l}
R\\
L \end{array}
 {\displaystyle { \big|}\atop{~~ \diagdown}}
 ~~{\displaystyle {\big|}\atop{ \big|}}
       ~~{\displaystyle {\big|}\atop{ \big|}}
       ~
{\displaystyle {~~\diagup}\atop{ \big|}}\\
&&\\&&\\
 \begin{array}{l}
L\\
R \end{array}
  {\displaystyle {~~~ \big|}\atop{ \diagup}}
       ~~{\displaystyle {\big|}\atop{ \big|}}
       ~~{\displaystyle {\big|}\atop{ \big|}}
       ~~{\displaystyle{\diagdown}\atop{~~~ \big|}}
~~~~~~~~~~
  &\Rightarrow &
~~~~~~~~~~
 \begin{array}{l}
R\\
L \end{array}
 {\displaystyle {~~\diagup}\atop{ \big|}}
 ~~{\displaystyle {\big|}\atop{ \big|}}
       ~~{\displaystyle {\big|}\atop{ \big|}}
       ~{\displaystyle { \big|}\atop{~~ \diagdown}}
 \\
 &&\\
 &&
 \\
  \begin{array}{l}
L\\
R \end{array}
  {\displaystyle {~~~ \big|}\atop{~~~ \big| }}
      ~~~ . ~ ~~{\displaystyle {\big|}\atop{ \big|}}
       ~~{\displaystyle {\big|}\atop{ \big|}}
       ~~{\displaystyle{\diagdown}\atop{\diagup}}
~~~~~~~~~~
  &\Rightarrow&
~~~~~~~~~~
 \begin{array}{l}
R\\
L \end{array}
 {\displaystyle {~~\big|}\atop{ ~~\big|}}
 ~~{\displaystyle {\diagup }\atop{ \diagdown }}
       ~~{\displaystyle {\big|}\atop{ \big|}}
       {\displaystyle { ~~\big|}\atop{~~ \big|}}
       ~~
 \eean
which then can be transformed into a new
 configuration $L~R$, with same beginning and end,
  thus yielding a new sequence;
 in the third case
 the reflection occurs the first place it can.
 So, by the moves above, the original configuration
  (\ref{LR's}) can be transformed in a new one.
  In the new sequence, pick again the first sequence $RL$,
  reading from left to right, and use again one of the
  moves. So, this leads again to a new sequence, etc...
 \bean
 &L~~R~~R~~\underline{R~~L}~~R~~L~~L~~R~\ldots ~R&\\
 &L~~R~~\underline{R~~L}~~R~~R~~L~~L~~R~\ldots ~R& \\
 &L~~\underline{R~~L}~~R~~R~~R~~L~~L~~R~\ldots ~R&\\
 &L~~L~~R~~R~~R~~\underline{R~~L}~~L~~R~\ldots ~R&\\
 &L~~L~~R~~R~~\underline{R~~L}~~L~~L~~R~\ldots ~R&\\
 &\vdots&
 \eean
%
%
%
 \be
 \underbrace{L~~L~~L~\ldots ~L}_{T_L}~~
  \underbrace{R~~R~~R~\ldots ~R}_{T_R}~;
 \ee
  Since this procedure is invertible, it gives a {\em one-to-one} map between all the left-right walks
  corresponding to a given sequence, with $T_L$ $L$'s and $T_R$
  $R$'s
  \be
 L~~R~~R~~{R~~L}~~R~~L~~L~~R~\ldots ~R~;
 \ee
and all the walks corresponding to \be
 \overbrace{L~~L~~L~\ldots ~L}^{T_L}
 ~~\overbrace{R~~R~~R~\ldots ~R}^{T_R}~.
 \label{LRwalk}\ee
Thus, a walk corresponding to (\ref{LRwalk}) will map
into $\left( {T}
    \atop{T_L~~ T_R }
    \right)$ different walks, corresponding to the $\left( {T}
    \atop{T_L~~ T_R }
    \right)$ number of permutations of
    $T_L$ $L$'s and $T_R$ $R$'s.

    \bigbreak

 {\em Step 2:}
 To two standard tableaux $P$, $Q$ of shape
$\lb=(\lb_1\geq\ldots\geq\lb_m>0)$ we associate a random
walk, going with (\ref{LRwalk}), in the following way.
%
Consider the situation where $m$ walkers start at
$0,1,\ldots,m-1$.

\medbreak

 $$
\raisebox{15mm}{\begin{tabular}{l}
\raisebox{1mm}{1$^{\rm st}$ walker}\\
\raisebox{3mm}{2$^{\rm nd}$ walker}\\
\raisebox{3mm}{3$^{\rm rd}$ walker}
 \end{tabular}}\overbrace{\squaresize .7cm \thickness .01cm
\Thickness .07cm \Young{ c_{11}&c_{12}&c_{13}&d c_{14}&
d c_{15}&d
c_{16} \cr c_{21}&c_{22}&c_{23}
 \cr c_{31}&d c_{32}& d c_{33}\cr
 d c_{41}\cr}}^{\mbox{instants of left move}}
 $$

\noindent $\bullet$ The 1$^{\rm st}$ walker starts at 0
and moves to the left, only at instants
$$\mbox{$c_{1i}=$ content of box$(1,i)\in P$}$$
and thus has made, in the end, $\lb_1$ steps to the
left.

$$\vdots$$

%
\noindent $\bullet$ The $k^{\rm th}$ walker $(1\leq
k\leq m)$ starts at $k-1$ and moves to the left, only at
instants
$$\mbox{$c_{ki}=$ content of box$(k,i)\in P$.}$$
and thus has made, in the end, $\lb_2$ steps to the
left.
$$\vdots$$
\noindent $\bullet$ Finally, walker $m=\lb_1^{\top}$
walks according to the contents of the last row.

\bigbreak

Since the tableau is standard, filled with the numbers
$1,\ldots,n$ the walkers never meet and at each moment
exactly one walker moves, until instant
$n=\vert\lb\vert$, during which they have moved from
position
$$0~<~1~<~\ldots~ <~  k-1~< ~\ldots ~ <~m-1$$
to position
$$
-\lb_1+0<-\lb_2+1<\ldots <-\lb_k+k-1<\ldots < -\lb_m
+m-1
$$
That is to say the final position is given by unfolding
the right hand boundary of $\lb$, the horizontal (fat)
segments refer to gaps between the final positions and
vertical segments refer to contiguous final positions.

In the same fashion, one associates a similar walk to
the other tableau $Q$, with the walkers also moving
left. These walkers will have reached the same position
as in the first case, since the final position only
depends on the shape of $P$ or $Q$. Therefore, reversing
the time for the second set of walks, one puts the two
walks together, thus yielding $m$ non-intersecting
walkers moving the first half of the time to the left
and then the second half of the time to the right, as in
the example below. Therefore the number of ways that $m$
walkers start from and return to $0,\ldots , m-1$,
without ever meeting each other, by first moving to the
left and then to the right, is exactly $\sum_{{\lb\vdash
n}\atop{\lb_1^{\top} \leq m}} \left( f^{\lb}\right)^2$.

  $$
  \displaystyle{
(P,Q)=\left(\!\!\!\begin{array}{c}
\overbrace{\,\,\,\,}^2\\
\MAT{1} \squaresize 0.5cm \thickness .01cm \Thickness
.06cm \Young{
 1&d 2 \cr
 3\cr  4\cr d 5\cr}
\mat ,\\
\mbox{\tiny standard}
\end{array}
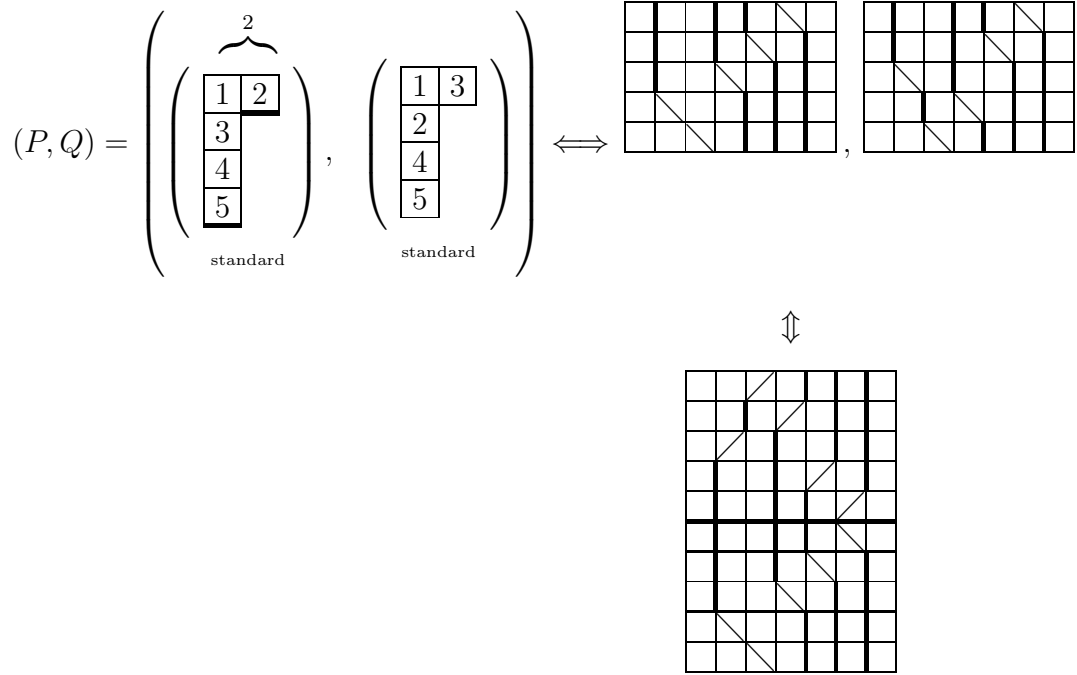
\begin{array}{c}
\\
\MAT{1}\squaresize 0.5cm \thickness .01cm\Thickness
.07cm \Young{
 1& 3\cr
2\cr 4\cr 5\cr}
\mat\\
\mbox{\tiny standard}
\end{array}
\!\!\!\right)}
  \Longleftrightarrow~
   \squaresize 0.4cm \thickness .01cm\Thickness
.05cm \Young{
   r& &r&r& &\diagdown& \cr
   r& &r& &\diagdown&r& \cr
   r& & &\diagdown&r&r& \cr
    &\diagdown& &r&r&r& \cr
     & &\diagdown&r&r&r& \cr
  }
  ~,~
\squaresize 0.4cm \thickness .01cm\Thickness .05cm
\Young{
 r& &r&r& &\diagdown& \cr
 r& &r& &\diagdown&r& \cr
 &\diagdown&r& &r&r&   \cr
 & r& &\diagdown&r&r& \cr
 & &\diagdown&r&r&r& \cr
  }
$$
$$\hspace*{7cm}\Updownarrow$$
 $$  \hspace*{7cm}
 \squaresize 0.4cm \thickness .01cm\Thickness
.05cm \Young{
  & &\diagup&r&r&r& \cr
  & r& &\diagup&r&r& \cr
   &\diagup&r& &r&r&\cr
   r& &r& &\diagup&r& \cr
   rd& d&rd&rd&d &d\diagup & d\cr
   r& &r&r& &\diagdown& \cr
   r& &r& &\diagdown&r& \cr
   r& & &\diagdown&r&r& \cr
    &\diagdown& &r&r&r& \cr
     & &\diagdown&r&r&r& \cr
  }$$



 More generally, an analogous argument shows that the number of ways that
 walkers
 leave from $x_1<\ldots<x_m$ and end up at
 $y_1<\ldots<y_m$ at time $T$, without intersection and by first moving to the
left and then to the right, is given by
 \be
   \renewcommand{\arraystretch}{0.5}
\begin{array}[t]{c}
{\displaystyle \sum} \\
{\scriptstyle \lb \vdash \frac{1}{2} (T+|\mu|+|\nu|)      }\\
{\scriptstyle \lb^{\top}_1\leq m}
\end{array}
\renewcommand{\arraystretch}{1}
f^{\lb \ba \mu}f^{\lb \ba \nu}
 \label{last}\ee
%
On the other hand, there are $\left( {T}
    \atop{T_L~~ T_R }
    \right)$
    sequences of $T_L$ $L$'s and $T_R$ $R$'s,
which combined with (\ref{last}) yields formula
(\ref{1.5.6}).


In the close packing situation, one has $\mu_k=\nu_k=0$
for all $k$, and so $\mu=\nu=\emptyset $ and
$T_L=T_R=T/2$. With these data, (\ref{1.5.7}) is an
immediate consequence of (\ref{1.5.6}).   \qed


\subsection{Probability measure on words}

Remember from section 1.1 words $\pi:=\pi_n^q$ of length
$n$ from an alphabet $1,\ldots, q$ and from section
\ref{subsection1.3}, the set $\BY_n^q$ of partitions
$\lb\vdash n,~\mbox{with}~
      \lb^{\top}_1\leq q\}$. Also define the set of $n\times
q$ matrices,
$$
  \tilde{\mbox{Mat}}^{n,q}_n(0,1)
 :=
 \left\{ W= (w_{ij})_{{1\leq i,j\leq n}}
,
~\left|\begin{array}{l}\mbox{with
exactly one $1$}\\
\mbox{in each row}\\
\mbox{and otherwise all zeros}\end{array}\right.
 \right\}
  $$
  For references, see \cite{Stanley1,TW-words}.

\begin{proposition}
In particular, for {\em words}, we have the 1-1
correspondence
$$
S_n^{q}\Longleftrightarrow\left\{\begin{array}{c}
\mbox{semi-standard and standard Young tableaux}\\
\mbox{$(P,Q)$ of same shape and of size $n$, filled}\\
 \mbox{resp., with integers $(1,\ldots,q)$ and $(1,\ldots,n),$}\\
%
\end{array}\right\}
 \Longleftrightarrow
    \tilde {\mbox{Mat}}_n^{n,q}( 0,1 ).
$$
\be \pi
 \hspace{1.5cm}\longleftrightarrow
 \hspace{1.5cm}
 (P,Q)
 \hspace{1.5cm}
 \longleftrightarrow \hspace{1.5cm}
  W(\pi)=(w_{ij})_{i,j\geq 1}.
\label{RSK-correspondence(words)}\ee
%
Uniform probability $P^{n,q}$ on $S_n^q$ induces a
probability $\tilde P^{n,q}$ on Young diagrams $\lb \in
\BY^{(q)}_n$, given by
 \bea
 \tilde P^{n,q}\left( \lb \right)=
\frac{1}{q^n} \#\left\{
\begin{array}{l}
  \mbox{words in $S_n^q$}\\
  \mbox{leading
  }\\
  \mbox{to shape $\lb$}
  \end{array}
  \right\}
 &=&\frac{f^{\lb}\tilde s_{\lb}(1^q)}{q^n}
 \no\\
 &=& \frac{n!}{q^n}  {\bf s}_{\lb}(1,0,\ldots)
  ~\gs_{\lb}(q,\frac{q}{2},\frac{q}{3},\ldots) \no\\
 &=&
  \frac{n!}{q^n\displaystyle{\prod_1^{q-1}}i!}
\frac{\Dt_
q(q+\!\lb_1\!-\!1,\ldots,q\!+\!\lb_q\!-\!q)^2}
 {\displaystyle{\prod_1^q}(q+\lb_i-i)!}.
 \no\\&&\label{wordPlancherel}\eea
%
%
Also,
\be
 \# S^q_n =\sum_{\lb\vdash n} f^{\lb}~ \tilde s_{\lb}(1^q)
=  q^n.
 \ee
%
Finally, given the correspondence
(\ref{RSK-correspondence(words)}), the length $L(\pi)$
of the longest weakly increasing subsequence of the word
$\pi$ equals
$$ L(\pi)=\lb_1=L(W):=\max_{ \mbox{\tiny
all such}\atop
                          \mbox{\tiny paths} }
\left\{\sum w_{ij},\quad\begin{array}{l}
\mbox{over right/down}\\
\mbox{paths starting from}\\
\mbox{entry $(1,1)$ to $(n,q)$}
\end{array}\right\},$$
and thus
$$
P^{n,q}(L(\pi)\leq \ell)=
\sum_{\begin{array}{l}\lb_1\leq \ell \\
                    \lb \in \BY_n^{(q)}\end{array}}
 \frac{n!}{q^n}  {\bf s}_{\lb}(1,0,\ldots)
  ~\gs_{\lb}(q,\frac{q}{2},\frac{q}{3},\ldots) .$$

\end{proposition}

\proof A word is a special instance of generalized
permutation, where the numbers $i_1,\ldots,i_n$ are all
distinct. Therefore the RSK correspondence holds as
before, except that $Q$ becomes a standard tableau;
thus, a word maps to a pair of arbitrary Young tableaux
$(P,Q)$, with $P$ semi-standard and $Q$ standard and
converse. Also the correspondence with integer matrices
is the same, with the extra requirement that the matrix
contains $0$ and $1$'s, with each row containing exactly
one $1$.
%

Consider now the uniform probability $P^{n,q}$ on {\em
words} in $S^q_n$; from the RSK correspondence
%
%
%
%
%
we have the one-to-one correspondence, given a fixed
partition $\lb$,
$$
 \left\{
\begin{array}{l}
  \mbox{words in $S_n^q$}\\
  \mbox{leading
  }\\
  \mbox{to shape $\lb$}
  \end{array}
  \right\}
   \Longleftrightarrow
 \left\{
\begin{array}{l}
  \mbox{semi-standard}\\
   \mbox{tableaux of shape}\\
  \mbox{$\lb$, filled with}\\
  \mbox{integers $1,\ldots,q$}
  \end{array}
  \right\}
   \times
  \left\{
\begin{array}{l}
  \mbox{standard tableaux}\\
  \mbox{of shape $\lb$} \\
  \mbox{filled with integers}\\
  \mbox{ $1,\ldots,n$}
  \end{array}
  \right\}
$$
 and thus, using formulae (\ref{number standard tableaux})
  and (\ref{number semi-standard tableaux}) and noticing that
 $\tilde\gs_{\lb}(1^q)=0$ for $\lb_1^{\top}>q$,
 $$
 \tilde P^{n,q}\left( \lb \right)=
\frac{1}{q^n} \#\left\{
\begin{array}{l}
  \mbox{words leading}\\
  \mbox{to the shape $\lb$}
  \end{array}
  \right\}
 =\frac{\tilde\gs_{\lb}(1^q)f^{\lb}}{q^n},~~~~~~~
 \lb\in\BY^{q}_n,
$$
with
$$
\sum _{\lb\in\BY_n^{(p)}} \tilde P^{n,q}(\lb)=1.
$$
Formula (\ref{wordPlancherel}) follows immediately from
an explicit evaluation (\ref{number standard tableaux})
and (\ref{number semi-standard tableaux}) for $f^{\lb}$
and $\gs_{\lb}(1^q)$.\qed

\begin{example}
For word
$$
\pi=\left(\begin{array}{lllll}
1&2&3&4&5\\
2&\underline{1}&\underline{1}&3&\underline{2}\end{array}\right)\in
S_5^3
$$
the RSK algorithm gives

{\footnotesize $$
\begin{array}{cccccccccccccccccccccc}
2&&&&1&&&&1&1&&&&1&1&3&&&&1&1&2\\
 &&&&2&&&&2& &&&&2& & &&&&2&3& \\
 \\
1&&&&1&&&&1&3&&&&1&3&4&&&&1&3&4\\
 &&&&2&&&&2& &&&&2& & &&&&2&5&
\end{array}
$$}
Hence
%
$$
\pi\Longrightarrow\begin{array}{c}
\stackrel{3}{\longleftrightarrow}\\
\left(\squaresize .4cm \thickness .01cm \Thickness .07cm
\Young{
 1&1&2\cr
 2&3\cr
 },\right.\\
\mbox{\tiny semi-}  \\
 \mbox{\tiny standard}
 \end{array}
 \begin{array}{c}
 \\
\left.\squaresize .4cm \thickness .01cm \Thickness .07cm
\Young{
 1&3&4\cr
 2&5\cr
 }\right)\\
 \\
\mbox{\tiny standard}
 \end{array}
\Longrightarrow\left(
 \begin{array}{ccccc}
{\!\!\scriptstyle 0}{\!\!\!\!\bigcirc\!\!}& &{\!\!\scriptstyle 1}& &0\\
\downarrow  & & & & \\
{\!\!\scriptstyle 1}{\!\!\!\!\bigcirc\!\!}& &0& &0\\
\downarrow  & & & & \\
{\!\!\scriptstyle 1}{\!\!\!\!\bigcirc\!\!}&\rg&
{\!\!\scriptstyle 0}{\!\!\!\!\bigcirc\!\!}&\rg&
{\!\!\scriptstyle 0}{\!\!\!\!\bigcirc\!\!}\\
 & & & &\downarrow \\
0& &0& &{\!\!\scriptstyle 1}{\!\!\!\!\bigcirc\!\!}\\
 & & & & \downarrow \\
0& &1& &{\!\!\scriptstyle 0}{\!\!\!\!\bigcirc\!\!}
\end{array}\right)
,$$
and $L(\pi)=\lb_1=L(W)=3$.
 \end{example}


 \subsection{Generalized Permutations, Percolation and growth
 models}

 The purpose of this section is to show that uniform probability, generalized
 permutations, percolation, queuing and growth models
 are intimately related. Their probabilities are
 ultimately given by the same formulae.

\subsubsection{Probability on Young diagrams induced
by uniform probability on generalized permutations}
 \begin{proposition} \label{pro:1.8}(Johansson
\cite{Johansson0})
 Uniform probability $P^{p,q}_n$ on $GP_n^{p,q}$ induces a
probability $\tilde P^{p,q}_n$ on Young diagrams
$\BY_n^{\min(p,q)}$, given by
 \bean
\tilde P_n^{p,q}(\lb)&=&  \frac{1}{\# \mbox{GP}_n^{p,q}}
\#\left\{
\begin{array}{l}
  \mbox{generalized permutations}\\
  \mbox{leading to shape $\lb$}
  \end{array}
  \right\}
  \\
  &=& \frac{\tilde{\bf s}_{\lb}(1^p)\tilde{\bf
s}_{\lb}(1^q)}
            {  \# \mbox{GP}_n^{p,q}   }       \\
 &=& \frac{1}{\# \mbox{GP}_n^{p,q} }\prod^{q-1}_{j=0}
 \frac{1}{j!(p-q+j)!}\\
   && \hspace{1cm}\Dt_q(q+\lb_1-1,...,q+\lb _q-q)^2
  \prod^q_{i=1}\frac{(p+\lb_i-i)!}
 {(q+\lb_i-i)!},~\mbox{for} ~q\leq p,
\eean
 with\footnote{There is no loss of generality in
 assuming $q\leq p$.}
 \be
 \# \mbox{GP}_n^{p,q}
 =
 \sum_{\lb\in \BY_n^{\min(p,q)}} \tilde\gs_{\lb}(1^p)~
  \tilde \gs_{\lb}(1^q)
=\left(\begin{array}{c}
  pq+n-1\\ n
\end{array}\right).\ee

\end{proposition}

\proof By the RSK correspondence, we have the one-to-one
correspondence,
$$
 \left\{
\begin{array}{l}
  \mbox{$\pi \in \mbox{GP}_n^{p,q}$ leading}\\
  \mbox{to the shape $\lb$}
  \end{array}
  \right\}\longleftrightarrow
 \left\{
\begin{array}{l}
  \mbox{semi-standard}\\
   \mbox{tableaux of shape}\\
  \mbox{$\lb$, filled with}\\
  \mbox{integers $1,\ldots,q$}
  \end{array}
  \right\}
  \times
  \left\{
\begin{array}{l}
  \mbox{semi-standard}\\
   \mbox{tableaux of shape}\\
  \mbox{$\lb$, filled with}\\
  \mbox{integers $1,\ldots,p$}
  \end{array}
  \right\}
$$
 and thus, for $\lb\in
\BY_n^{\min (q,p)}$, we have\footnote{Remember the
notation from section 1.2: $\lb\in \BY_n^{\min (q,p)}$
means the partition $\lb\vdash n$ satisfies
      $\lb^{\top}_1\leq p,q$. }
 \bean
  &&
\tilde P^{p,q}_n(\lb)=\frac{1}{\# \mbox{GP}^{p,q}_n}
\#\left\{
\begin{array}{l}
  \mbox{$\pi \in \mbox{GP}_n^{p,q}$ leading}\\
  \mbox{to the shape $\lb$}
  \end{array}
  \right\}= \frac{1}{\# \mbox{GP}^{p,q}_n}
\tilde\gs_{\lb}(1^q)\tilde\gs_{\lb}(1^p);
 \eean
%
     when $q<\lb_1^{\top} $, we have automatically
      $\tilde\gs_{\lb}(1^q)=0$,
and thus
$$
\sum _{\lb\in\BY_n^{\min (q,p)}} \tilde P^{p,q}_n(\lb)
=1.
$$
Notice that\footnote{Note
$ \prod_{1\leq i<j\leq q}(j-i)=\prod^{q-1}_{j=0}j! $ .}
for $m\geq \lb_1^{\top}$,
 \bea
 \tilde\gs_{\lb}(1^m)=
 \frac{\Dt_m(m+\lb_1-1,...,m+\lb_m-m)}
  {\displaystyle{\prod_1^{m-1}}i!}
  =
  \prod_{1\leq i<j\leq m}\frac{\lb_i-\lb_j+j-i}{j-i}\no\\
  \eea

 Without loss of generality, assume $p\geq q$;
 then $\lb_i=0$ if $q<i\leq p$. Setting
$h_j=\lb_j+q-j$ for $j=1,...,q$, we have $h_1>...
>h_q=\lb_q\geq 0$. We now compute, using (\ref{number semi-standard
tableaux}),
\bean \lefteqn{\tilde\gs_{\lb}(1^p)\tilde\gs_{\lb}(1^q)
}
 \\&=& \prod_{1\leq i<j\leq
q}\left(\frac{\lb_i-\lb_j+j-i}{j-i}\right)^2
\prod_{i=1}^q\prod_{j=q+1}^p\frac{\lb_i+j-i}{j-i}\\
&=&\prod_{1\leq i<j\leq q}\frac{(h_i-h_j)^2}{(j-i)^2}
\prod_{i=1}^q\prod_{j=q+1}^p\frac{h_i+j-q}{j-i}\\
&=&\frac{1}{\displaystyle\prod_{i=1}^q\prod^p_{j=q+1}(j-i){\prod_{1\leq
i<j\leq q}}(j-i)^2 }
 \prod_{1\leq
i<j\leq q}(h_i-h_j)^2\prod^q_{i=1}
 \frac{(h_i+p-q)!}{h_i!}
 \\
&=&
 \prod_{j=0}^{q-1}\frac{1}{j!(p-q+j)!}
 \prod_{1\leq i<j\leq
q}(h_i-h_j)^2\prod^q_{i=1}\frac{(h_i+p-q)!}{h_i!}
 ,\eean
using
\bean
\lefteqn{\prod_{1\leq i<j\leq q}(j-i)\prod^{q}_{i=1}\prod^p_{j=q+1}(j-i)}\\
& & \\
&=&\prod^{q-1}_{0}j!\prod^{q}_{i=1}(q+1-i)(q+2-i)(q+3-i)...(p-i)\\
& & \\
&=&q!\frac{(q+1)!}{1!}\frac{(q+2)!}{2!}...\frac{(p-1)!}{(p-q-1)!}1!
...
(q-1)!\\
& & \\
&=&\prod^{q-1}_{j=0}(p-q+j)!
 .\eean
 This ends the proof of Proposition \ref{pro:1.8}.\qed
%
%
%
%


\subsubsection
{Percolation model with geometrically distributed
entries}\label{section 1.7.2}

Consider the ensemble
 $$
 \mbox{Mat}^{(p,q)}=\{p\times q\mbox{~matrices $M$
with
  entries $M_{ij}=0,1,2,\ldots$}\}
 $$
 with {\em independent and geometrically
distributed entries}, for fixed $0<\xi<1$,
 $$
P(M_{ij}=k)=(1-\xi)\xi^k,~~~~~k=0,1,2,\ldots.
$$

\vspace*{1cm}

 \hspace*{3.5cm}$(1,1)$ \hspace*{3.6cm}$(1,q) $
$$
  \squaresize .4cm \thickness .02cm \Thickness
.05cm \Young{
     l&&&&&&&&&\cr
     l &&&&&&&&&\cr
     ur&ld  & & &&&&&&\cr
      &&l&& &&&&&\cr
      &&l&& &&&&&\cr
      &&u&u&u&ur&&&&\cr
      &&&&&&u&ur&d&\cr
      & &&&&&&&r& d \cr
  }
$$
\hspace*{4.0cm}$(p,1)$ \hspace*{3.6cm}$(p,q) $

   \vspace*{1cm}

\begin{theorem}(Johansson \cite{Johansson0}) \label{theorem: percolation}
 Then
 $$L(M):=\max_{ \mbox{\tiny
all such}\atop
                          \mbox{\tiny paths} }
\left\{\sum M_{ij},\quad\begin{array}{l}
\mbox{over right/down}\\
\mbox{paths starting from}\\
\mbox{entry $(1,1)$ to $(p,q)$}
\end{array}\right\} $$
 %
%
has the following distribution, 
assuming $q\leq p$,
\bean 
 P(L(M)\leq\ell)&=&\sum_{\lb\in
\BY^{\min(q,p)}\atop{\lb_1\leq\ell}}
(1-\xi)^{pq}\xi^{\vert\lb\vert}\tilde\gs_{\lb}
 (1^q)\tilde\gs_{\lb}(1^p)
 \\
&=&Z^{-1}_{p,q} \sum_{h\in \BN^q \atop{\max(h_i)\leq
\ell
+q-1}}\Dt_q(h_1,\ldots,h_q)^2\prod^q_{i=1}\frac{(h_i+p-q)!}{h_i!}\xi^{h_i}
\eean
 where
\be
Z_{p,q}=\xi^{\frac{q(q-1)}{2}}(1-\xi)^{-pq}q!\prod^{q-1}_{j=0}j!(p-q+j)!
\label{constant1}\ee
\end{theorem}

%

\proof Then the probability that $M$ be a given matrix
$A=(a_{ij})_{1\leq i\leq p\atop{1\leq j\leq q}}$ equals

\bean
 P\left(M=(a_{ij})_{1\leq i\leq
p\atop{1\leq j\leq q}}\right)&=&
  \prod_{1\leq i,j\leq p\atop{1\leq
j\leq q}}P(M_{ij}=a_{ij}),\mbox{~using independence,}\\
&=& \prod_{1\leq i\leq p\atop{1\leq j\leq q}}(1-\xi)\xi^{a_{ij}}\\
&=&
 (1-\xi)^{pq}\xi^{\sum_{1\leq i\leq
p\atop{1\leq j\leq q}}a_{ij}}\\&&\\
 &=&
  (1-\xi)^{pq }\xi^{\vert A\vert}.
  \eean
This probability only depends on the total weight
$|A|=\sum_{i,j} a_{ij}$. Hence the
matrices\footnote{$M\in \mbox{Mat}_{n}^{pq}\subset
\mbox{Mat}^{pq}$, means that $\sum M_{ij}=n$.} $M \in
\mbox{Mat}_{n}^{pq}$ have all equal probability and, in
particular, due to the fact that, according to Theorem
\ref{Theo: RSK for GP}, the matrices in
$\mbox{Mat}_{n}^{pq}$ are in one-to-one correspondence
with generalized permutations of size $n$, with
alphabets $1,\ldots,p$ and $1,\ldots,q$, one has
  $$
  P(|M|=n)
=\renewcommand{\arraystretch}{0.5}
\begin{array}[t]{c}
\sum\\
{\scriptstyle \mbox{\tiny all~}a_{ij}}\\
{\scriptstyle \mbox{\tiny with} }\\
{\scriptstyle \Sg~a_{ij}=n}
\end{array}
\renewcommand{\arraystretch}{1}
P\left(M=(a_{ij})_{1\leq i\leq p\atop{1\leq j\leq
q}}\right)=
 (\#\mbox{GP}^{p,q}_n)(1-\xi)^{qp}\xi^n.
$$
%
We now compute
  \bean
\lefteqn{P(L(M)\leq\ell ~\Big\vert ~\vert M\vert =n)}\\
&=&\frac{\#\{M\in \mbox{Mat}^{p,q}_n,~L(M)\leq \ell  \}
}
 {\# \mbox{Mat}^{p,q}_n}
 \\
 &=&\frac{\#\{\pi \in \mbox{GP}^{p,q}_n,~L(\pi)\leq \ell  \} }
 {\# \mbox{GP}^{p,q}_n}
  \\
  &=&P_n^{p,q}(\lb_1\leq\ell)
 \\
&=&\frac{1}{\#\mbox{GP}^{p,q}_n}\sum_{\lb_1\leq\ell\atop{\vert\lb\vert
=n}}\tilde\gs_{\lb} (1^q)\tilde\gs_{\lb}(1^p). \eean
Hence,
\bean P\left(L(M)\leq\ell
\right)&=&\sum^{\iy}_{n=0}P\left(L(M)\leq\ell
~\Big\vert~ \vert
M\vert =n\right)~P(\vert M\vert =n)\\
\\
&=&\sum^{\iy}_{n=0}\sum_{\lb_1\leq\ell\atop{\vert\lb\vert
=n}}
\frac{1}{\#\mbox{GP}^{p,q}_n}\tilde\gs_{\lb}(1^q)\tilde\gs_{\lb}(1^p)
 (\#\mbox{GP}^{p,q}_n) (1-\xi)^{pq}\xi^n\\
\\
&=&\sum^{\iy}_{n=0}\sum_{\lb_1\leq\ell\atop{\vert\lb\vert
=n}}
\tilde\gs_{\lb}(1^q)\tilde\gs_{\lb}(1^p)(1-\xi)^{pq}\xi^{|\lb|}\\
\\
&=&\sum_{\lb\in\BY\atop{\lb_1\leq\ell }}
\tilde\gs_{\lb}(1^q)\tilde\gs_{\lb}(1^p)(1-\xi)^{pq}\xi^{|\lb|}.
\eean
Now, using the expression 
for $\tilde\gs_{\lb}(1^p)\tilde\gs_{\lb}(1^q)$ in
Proposition \ref{pro:1.8}, one computes, upon setting
$h_i=q+\lb_i-i$, and noticing that $\ell\geq \lb_1\geq
  \ldots\lb_q\geq 0$ implies $\ell+q-1\geq
h_1> \ldots> h_q\geq 0$,
\bean
\lefteqn{P(L(M)\leq\ell)}\\
\\
&=&\sum_{\lb\in\BY^q\atop{\lb_1\leq\ell }}
(1-\xi)^{pq}\xi^{|\lb|}\tilde\gs_{\lb}(1^p)\tilde\gs_{\lb}(1^q)
\eean\bean
&=&\!\!\!\sum_{ {\ell+q-1\geq h_1>\ldots
>h_q\geq 0}}\frac{\xi^{\displaystyle{\sum_1^q
h_i-\frac{q(q-1)}{2}}}}{(1-\xi)^{-pq}\displaystyle{\prod
^{q-1}_{j=0}}j!(p-q+j)!}\Dt_q(h_1,\ldots,h_q)^2\prod^q_{i=1}
\frac{(h_i+p-q)!}{h_i!}
\\
&=&Z^{-1}_{p,q}\sum_{h\in\mathbb
N^q\atop{\max(h_i)\leq\ell +q-1}}\Dt_q(h_1,\ldots,h_q)^2
\prod^q_{i=1}\frac{(h_i+p-q)!}{h_i!}\xi^{h_i}, \eean
since the expression in the sum is symmetric in
$h_1,\ldots,h_q$. The normalization 
 $
Z_{p,q}$ is as announced, ending the proof of Theorem
\ref{theorem: percolation}.\qed


\subsubsection{Percolation model with exponentially distributed entries}

\begin{theorem}\cite{Johansson0}\label{theorem: percolation-exp}
Consider the ensemble
 $$
 \mbox{Mat}^{p,q}=\{p\times q\mbox{~matrices $M$
with
 $\BR^+$-entries}\}
 $$
 with {independent and exponentially
distributed entries},
 $$
P(M_{ij}\leq t)=1-e^{-t},~~~~~t\geq 0~~\mbox{and}~~\in
\BR.
$$
%
 Then
 \bean
 L(M)&=&\max_{  \mbox{\tiny all such}\atop
                               \mbox{\tiny paths}
                               }
\left\{\sum M_{ij},\qquad\begin{array}{l}
\mbox{over right/lower}\\
\mbox{paths starting from}\\
\mbox{entry $(1,1)$ to $(p,q)$}
\end{array}\right\}
\eean
%
has the following distribution, 
(assuming $q\leq p$, without loss of generality),
\bean 
 P(L(M)\leq t)&=&\frac{\int_{(0,t)^q}\Dt_q(x_1,\ldots,x_q)^2
 \prod^q_{i=1}x_i^{p-q}e^{-x_i}dx_i}{
  \int_{(0,\iy)^q}\Dt_q(x_1,\ldots,x_q)^2
 \prod^q_{i=1}x_i^{p-q}e^{-x_i}dx_i}.
\\ \\
 &=&\frac{1}{Z_n}\int_{M\in {\cal H}_q,~ \mbox{\tiny Spectrum}(M)\leq t} (\det M)^{p-q}e^{-\Tr M}dM,
 \eean

\end{theorem}

\remark It is remarkable that this percolation problem
coincides with the probability that the spectrum of an
Hermitian matrix does not exceed $t$, where the matrix
is taken from a (positive definite) random Hermitian
ensemble with the Laguerre distribution, as appears in
the second formula; this ensemble will be discussed much
later in section \ref{subsubsect 8.3.1}.

\bigbreak

\proof For fixed $0<\xi <1$, let $X_{\xi}$ have a
geometric distribution

$$
P(X_{\xi}=k)=(1-\xi)\xi^k\qquad 0<\xi <1,\quad
k=0,1,2,\ldots;
$$
then in distribution
$$
(1-\xi)X_{\xi}\lrg Y,\quad \mbox{for~}\xi\rg 1
$$
where $Y$ is an exponential distributed random variable.
Indeed, setting $\vr :=1-\xi$,
\bean
P\left((1-\xi)X_{\xi}\leq t\right)&=&P(\vr X_{1-\vr}\leq t)\\
\\
&=&\sum_{0\leq k\leq t/\vr}\vr(1-\vr)^k\\
\\
&=&\sum_{0\leq k\leq
t/\vr}\vr(1-\vr)^{\frac{1}{\vr}k\vr}
  ~~~~\mbox{(Riemann sum)}\\
\\
&\rg&\int^t_0ds~e^{-s}=P(Y\leq t) \eean
 Then, setting $\xi-1-\vr$,
$t=\ell \vr$, $\vr h_i=x_i$ in formula (\ref{constant1})
of Theorem \ref{theorem: percolation} and letting
$\vr\rg 0$ , one computes
\bean
\lefteqn{\lim_{\vr\rg 0}P_{\vr}(L(M)\leq t)}\\
\\
&=&\lim_{\vr\rg 0}\frac{Z^{-1}}{\vr^{q(q-1)+q(p-q)+q}}
\\
&&~~~\sum_{\max(h_i\vr)\leq\ell\vr
+(q-1)\vr}\Dt_q(h_1\vr,\ldots,h_q\vr)^2
\prod^q_{i=1}\left(\prod_{k=1}^{p-q}\vr(h_i+k)\right)(1-\vr)^{\frac{1}{\vr}
h_i\vr}\vr^q\\
\\
&=&\frac{1}{Z'}\int_{[0,t]^q}\Dt_q(x_1,\ldots,x_q)^2
\prod^q_{i=1}x_i^{p-q}e^{-x_i}dx_i
 \\
 &=&\frac{1}{Z_n}\int_{M\in {\cal H}_q,~
 \mbox{\tiny Spectrum}(M)\leq t} (\det M)^{p-q}e^{-\Tr
 M}dM.
 \eean
  This last identity will be shown later in
  Proposition \ref{Proposition 7.2}, ending the proof of Theorem \ref{theorem: percolation-exp}.\qed


\subsubsection{Queuing problem}

Consider servers $1,2,\ldots,q$ waiting on customers
$1,\ldots,p$, with a ``{\it first--in first--out}"
service rule; see e.g., \cite{Bary} and references
within. At first, the system is empty and then $p$
customers are placed in the first queue. When a customer
$k$ is done with the $i^{\rm th}$ server, he moves to
the queue waiting to be served by the $i+1^{\rm st}$
server. Let
$$
V(k,\ell)=\mbox{~service time of customer $k$ by server
$\ell$}
$$
be all geometrically independently distributed
$$
P(V(k,\ell)=t)=(1-\xi)\xi^t,\quad t=0,1,2,\ldots~.
$$
\vspace*{0.51cm}

$$
\begin{array}{ccccccc}
\mbox{servers:}&\fbox{$1$}~&~\fbox{$2$}~&~\fbox{$3$}~&
   ~\fbox{$4$}~ &~\ldots~&~\fbox{$q$}
  \\ \\
  \mbox{customers:}&\dashbox{3}(18,18){$5$}
    &\dashbox{3}(18,18){$2$}
     & &\dashbox{3}(18,18){$1$}   &
   \\
  & \dashbox{3}(18,18){$6$}&\dashbox{3}(18,18){$3$}& &\\
  &\dashbox{3}(18,18){$7$}&\dashbox{3}(18,18){$4$}\\
 & \dashbox{3}(18,18){$8$}&\\
  &\vdots&\\
 & \dashbox{3}(18,18){$p$}&\\
   \end{array}
$$
~~~~~~~~~~~~~~~~~~{\it Servers $1,2,\ldots,q$ waiting on
customers
   $1,\ldots,p$}

\vspace*{1cm}

\begin{theorem}\label{Theo: queuing}

The distribution of
%
$$
D(p,q)=\left\{\begin{array}{l}
\mbox{departure time for the last }\\
\mbox{customer $\dashbox{3}(18,18){$p$}$ at the last
server $\fbox{$q$}$}
\end{array}\right\}
$$
 is given by (assuming $q\leq p$)
 \bea  
 P(D(p,q)\leq\ell)&=&\sum_{\lb\in
\BY^{\min(q,p)}\atop{\lb_1\leq\ell}}
(1-\xi)^{pq}\xi^{\vert\lb\vert}\gs_{\lb}
 (1^q)\gs_{\lb}(1^p)
 \no\\
&=&Z^{-1}_{q,p} \sum_{h\in \BN^q \atop{\max(h_i)\leq
\ell
+q-1}}\Dt_q(h_1,\ldots,h_q)^2\prod^q_{i=1}\frac{(h_i+p-q)!}{h_i!}\xi^{h_i}
\no\\
\label{distribution queuing}\eea

\end{theorem}

\proof We show the problem is equivalent to the
percolation problem discussed in Theorem \ref{theorem:
percolation}. Indeed:

\underline{\em Step 1}: {\em Setting $D(p,0)=D(0,q)=0$
for all $p,q $, we have for all $p,q\geq 1$, \be D(p,q
)=\max\bigl(D(p-1,q),D(p,q-1)\bigr)+V(p,q
).\label{waiting time} \ee}
Indeed, if $D(p-1,q)\leq D(p,q-1)$, then customer $p-1$
has left server $q$ by the time customer $p$ reaches
$q$, so that customer $p$ will not have to queue up and
will be served immediately. Therefore
$$
D(p,q)=D(p,q-1)+V(p,q).
$$
Now assume
$$
D(p-1,q)\geq D(p,q-1);
$$
then when customer $p$ reaches server $q$, then customer
$p-1$ is still being served at queue $q$. Therefore the
departure time of customer $p$ at queue $q$

$$
D(p,q )=D(p-1,q)+V(p,q).
$$
In particular
$$
D(p,1)=D(p-1,1)+V(p,1)
$$
and
$$
D(1,2)=D(1,1)+V(1,2),
$$
establishing (\ref{waiting time}).

\medbreak

\underline{\em Step 2}: we now prove {\em $$ D(p,q
)=\max_{ \mbox{\tiny all such}\atop
                               \mbox{\tiny paths}
                               } \left\{
\sum V(i,j)\Bigr|\begin{array}{ll}
&\mbox{over right/down paths}\\
&\mbox{from entry $(1,1)$ to $(p,q )$}
\end{array}
\right\},
$$
where the paths are taken in the random matrix
$$
V=(V(i,j))_{1\leq i\leq p\atop{1\leq j\leq q}} .$$} By a
straightforward calculation,
 \bean
D(p,q )&=&\max(D(p-1,q)+V(p,q ),D(p,q-1)+V(p,q ))\\
&=&\max\left(\max_{\mbox{\tiny all
paths}\atop{(0,0)\rg(p-1,q)}}\sum_{\mbox{\small path}}
V(i,j)+V(p,q),\right.\\
&&\hspace*{1cm}\left.\max_{\mbox{\tiny all
paths}\atop{(0,0)\rg(p,q-1)}}\sum_{\mbox{\small path}}
 V(i,j)+V(p,q )\right)\\
&=&\max_{\mbox{\tiny all paths}\atop{(0,0)\rg(p,q
)}}\left(\sum_{\mbox{\small path}}V(i,j)\right), \eean
ending the proof of Theorem \ref{Theo: queuing}.\qed


 \subsubsection{Discrete polynuclear growth
models}

 Consider geometric i.i.d. random variables $\om(x,t)$, with $x \in \BZ, t\in \BZ_{+}$,
 $$P(\om {(x,t)}=k)=(1-\xi)\xi^k,~~~k\in \BZ_+
,$$ except
 $$\om(x,t)=0\mbox{ ~if $t-x$ is even or $|x|>t$}.
$$
Define inductively a growth process, with a height curve
$h(x,t)$, with $x \in \BZ, t\in \BZ_{+}$, given by
$$h(x,0)=0, 
$$
 $$
h(x,t+1)=\max \Bigl(h(x-1,t), h(x ,t), h(x+1,t)\Bigr)
+\om (x,t+1)
 .$$
 For this model, see \cite{Johansson2}, \cite{Spohn} and
 \cite{Krug}.

\begin{theorem}\label{Theo: polynuclear}
The height curve at even sites $2x$ at times $2t-1$ is
given by
$$ 
  h(2x,2t - 1)= \max_{  \mbox{\tiny all
such}\atop
                               \mbox{\tiny paths}
                               }  \left\{\sum V({i,j}),
\begin{array}{l}
\mbox{over right/down\!\!}\\
\mbox{paths starting from\!\!}\\
\mbox{entry $(1,1)$ to $(t+x,t-x)$}
\end{array}\right\}
$$
where
$$V(i,j):=\om(i-j,i+j-1).$$
Thus $h(2x,2t - 1)$ has again the same distribution as
in (\ref{constant1}).
\end{theorem}


\proof It is based on the fact that, setting
 $$
  G(q,p) :=  h(q-p,q+p-1),
  $$
 one computes
 \bean
 G(q,p)& =&\max (G(q-1,p), G(q,p-1)) +V(q,p)\\
 &&\\&=& \max \left\{\sum V_{ij},~\begin{array}{l}
\mbox{over all right/down}\\
\mbox{paths starting from}\\
\mbox{entry $(1,1)$ to $(q,p)$}
\end{array}\right\}
 .\eean
So
 $$
   h(2x,2t-1)= G(t+x,t-x),
$$establishing Theorem \ref{Theo: polynuclear}.\qed
%

\bigbreak

 The figure below gives an example of such a
growth process and also it shows that $h(2x,2t+1)$ is
given by the maximum of right/down paths starting at the
upper-left corner and going to site $(t+x,t-x)$, where
$x$ is the running variable along the anti-diagonal.

\newpage



\newpage

\vspace*{-1cm}

  $t=1$
$$
  \squaresize .8cm \thickness .02cm \Thickness
.07cm \Young{
   * \cr * \cr
  }
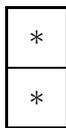
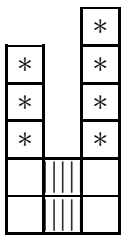
$$
\vspace{-.5cm}
$$ \mbox{\tiny 
 $0$} $$

$t=2$
$$
  \squaresize .5cm \thickness .02cm \Thickness
.07cm
 \Young{
  \blank &\blank &  *  \cr
   * &\blank &  *  \cr
    * &\blank&   *        \cr
 * &  \blank    &  *    \cr
     &||||| & \cr
  &  ||||| &     \cr
  }
$$
\vspace{-0.5cm}
$$ \mbox{\tiny 
 $-1~~~0~~~~1$} $$

$t=3$

$$
  \squaresize .5cm \thickness .02cm \Thickness
.07cm \Young{ * & \blank&\blank & \blank & *\cr
    * &\blank & * & \blank &   *           \cr
  *  & \blank  &  & |\!|\!|\!|\!|\!|\!|\!|\!|\!|\!|  &          \cr
  & |\!|\!|\!|\!|\!|\!|\!|\!|\!|\!|  &   & |\!|\!|\!|\!|\!|\!|\!|\!|\!|\!|  &     \cr
  &  |\!|\!|\!|\!|\!|\!|\!|\!|\!|\!|  & & |\!|\!|\!|\!|\!|\!|\!|\!|\!|\!|   &     \cr
 &  |\!|\!|\!|\!|\!|\!|\!|\!|\!|\!| &         &|\!|\!|\!|\!|\!|\!|\!|\!|\!|\!|  &    \cr
  & |\!|\!|\!|\!|\!|\!|\!|\!|\!|\!|   &|||||  & |\!|\!|\!|\!|\!|\!|\!|\!|\!|\!| & \cr
 & |\!|\!|\!|\!|\!|\!|\!|\!|\!|\!| &  ||||| &|\!|\!|\!|\!|\!|\!|\!|\!|\!|\!|   &  \cr
  }
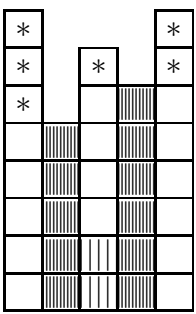
$$
\vspace{-0.5cm}
$$ \mbox{\tiny 
 $-2~~ -1~~~0~~~~1~~~~2$}~~ $$

%


 \hspace*{3.3cm}$0$ \hspace*{5.3cm}$(2t,0)$
$$
  \squaresize .6cm \thickness .02cm \Thickness
.07cm \Young{
     l&&&&&&&&\cr
     l &&&&&&&\cr
     ur&ld  & & &&&\cr
      &&l&& &\cr
      &&l&& \cr
      &&u&u\cr
      &&\cr
      &   \cr
       \cr
  }
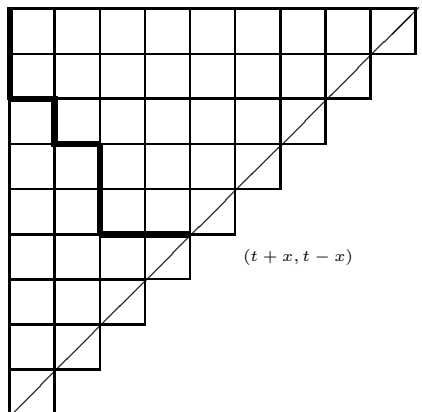
$$

\vspace*{-.9cm}\hspace{3.46cm}
\begin{picture}(12,1.6) 
\put(0,2){\line(1, 1){153}}
\end{picture}

\vspace*{-.1cm} \hspace*{2.6cm}  $(0,2t)$

\vspace*{-3.0cm} \hspace*{6.5cm} $\mbox{\tiny
$(t+x,t-x)$}$



\newpage

 \section{Probability on partitions, Toeplitz and
   Fredholm determinants}



Consider variables $x=(x_1,x_2,\ldots)$ and $y=(y_1,
y_2,\ldots)$, and the corresponding symmetric functions
$$
t=(t_1, t_2, \ldots) \mbox{ and } s=(s_1,s_2,\ldots)
$$
with
$$
kt_k=\sum_{i\geq 1} x_i^k \mbox{ and } ks_k=\sum_{i\geq
1} y_i^k.
$$
Following Borodin-Okounkov-Olshanski (see
\cite{Borodin1,Borodin2,BOO}), given arbitrary $x,y$,
define the (not necessarily positive) {\em probability
measure} on $\lb\in \BY$,
 \be P_{x,y} (\lb) :=\frac{1}{Z}
\tilde\gs_{\lb} (x) \tilde\gs_{\lb} (y)=\frac{1}{Z}
 \gs_{\lb} (t)  \gs_{\lb} (s)
 \label{generalized probability}\ee
with
$$
Z=\prod_{1\leq i ,j}(1-x_iy_j)^{-1}= e^{\sum_1^{\iy} k
t_k s_k}.
$$
Indeed, by Cauchy's identity,
 \bean
\sum_{\lb\in\BY}P_{x,y}(\lb)&=&\sum_{\lb}\gs_{\lb}(x)\gs_{\lb}(y)\prod_{1\leq
i,j}(1-x_iy_j)\\
& & \\
&=&\prod_{i,j\geq 1}(1-x_iy_j)^{-1}\prod_{i,j\geq
1}(1-x_iy_j)=1. \eean
%
%
The main objective of this section is to compute
\be P_{x,y}(\lambda_1 \leq n  ),\label{2.0.2} \ee
 which then will be specialized in the next section
 to specific $x$'s and $y$'s or, what is the
 same, to specific $t$'s and $s$'s.
 This probability (\ref{2.0.2}) has three different
expressions, one in terms of a determinant of a Toeplitz
determinant, another in terms of an integral over the
unitary group and still another in terms of a Fredholm
determinant. The Toeplitz representation enables one to
compute in an effective way the probability
(\ref{2.0.2}), whereas the Fredholm representation is
useful, when taking limits for large permutations.
In the statement below, we need the Fredholm determinant
of a kernel $K(i,j)$, with $i,j\in \BZ$ and restricted
to $[n, \iy)$; it is defined as
\be
  \det\left(I-K(i,j)\Big\vert_{[n,n+1,...]}\right)
   :=
 \sum^{\iy}_{m=0}(-1)^m\sum_{{n\leq x_1<...<
x_m}\atop{x_i \in \BZ}}\det(K(x_i,x_j))_{1\leq i,j\leq
m}. \label{def: Fredholm det}\ee

Now, one has the following statement
\cite{Gessel,Borodin1,Bo-Ok,TW-permutations,AvM2}:

\begin{theorem}\label{Theo:2.1} Given the
``probability measure" (\ref{generalized probability}),
%
the following probability has three different
expressions:
\bea
  P\left(\lb~\mbox{with}~\lb_1\leq n\right)
&=&Z^{-1}\det\left(\oint_{S^1}\frac{dz}{2\pi
iz}z^{k-\ell}e^{-\sum_1^{\iy}
(t_iz^i+s_iz^{-i})}\right)_{1\leq k,\ell\leq n}\no\\
& & \no\\
 &=&
 Z^{-1} \int_{U(n)}
  e^{-Tr~\sum_{i\geq 1}(t_iX^i+s_i\bar X^i)}dX\no
  \\
  &=&\det\left(I-K(k,\ell)\Bigg\vert_{[n,n+1,...]}\right)
,\label{formulaTheo:2.1}\eea
where $K(k,\ell)$ is a kernel \bea
K(k,\ell)&=&\left(\frac{1}{2\pi i}\right)^2
 \oint_{|w|=\rho<1}
 \oint_{|z|=\rho^{-1}>1}
  \frac{dz~
dw}{z^{k+1}w^{-\ell}} \frac{e^{V(z)-V(w)}}{z-w},~~~
\mbox{for
 $k,\ell \in \BZ$}%
 \no\\ \no\\ \no\\
 &=&\frac{1}{k-\ell}\left(\frac{1}{2\pi
i}\right)^2\oint\!\!\oint _{{|w|=\rho<1}\atop
{|z|=\rho^{-1}>1}}\frac{dz~dw}{z^{k+1}
w^{-\ell}}\frac{z\frac{d}{dz}V(z)-w\frac{d}{dw}V(w)}{z-w}
 ~e^{V(z)-V(w)}
\no\\
&& \hspace*{7cm} \mbox{for
 $k,\ell \in \BZ$, with $k\neq \ell$},
 \no\\
 \label{kernel-non-diagonal}%
\eea
 with
$$
V(z):=V_{t,s}(z):=-\sum_{i\geq 1}(t_iz^{-i}- s_iz^i).
$$

\end{theorem}



\subsection{Probability on partitions
 expressed as Toeplitz determinants}

 In this subsection, the first part of Theorem 2.1 will
 be reformulated as Proposition \ref{prop:Toeplitz}
 and also demonstrated:
 (see e.g. Gessel \cite{Gessel},
 Tracy-Widom \cite{TW-permutations},
 Adler-van Moerbeke \cite{AvM2})

\begin{proposition}\label{prop:Toeplitz}

Given the ``probability measure"

$$
P(\lb)=Z^{-1}\gs_{\lb}(t)\gs_{\lb}(s),\qquad
Z=e^{\sum_{i\geq 1} it_is_i},
$$
the following holds
\bean P\left(\lb~\mbox{with}~\lb_1\leq n\right)
&=&Z^{-1}\det\left(\oint_{S^1} \frac{dz}{2\pi iz}
z^{k-\ell}
 e^{-\sum_1^{\iy}(t_i z^i+s_i z^{-i})}\right)_{1\leq k, \ell\leq n}
\eean
and
$$
 P(\lb~\mbox{with}~\lb_1^{\top}\leq n)
=Z^{-1}\det\left(\oint_{S^1} \frac{dz}{2\pi iz}
z^{k-\ell}
 e^{\sum_1^{\iy}(t_i z^i+s_i z^{-i})}\right)_{1\leq k, \ell\leq n}
  $$

\end{proposition}

\proof Consider the semi-infinite Toeplitz matrix
$$
m_{\iy}{(t,s)}=(\mu_{k\ell})_{k,\ell\geq
0},\mbox{~with~}\mu_{k\ell}{(t,s)}=\oint_{S^1}
z^{k-\ell}e^{\sum^{\iy}_{1}
(t_jz^j-s_jz^{-j})}\frac{dz}{2\pi iz}.
$$
Note that
\bean
 \frac{\pl\mu_{k\ell}}{\pl t_m}&=&\oint_{S^1}
z^{k-\ell +m}e^{\sum^{\iy}_{1}
(t_jz^j-s_jz^{-j})}\frac{dz}{2\pi iz}=\mu_{k+m,\ell}\\
& & \\
\frac{\pl\mu_{k\ell}}{\pl s_m}&=&-\oint_{S^1} z^{k-\ell
-m}e^{\sum (t_jz^j-s_jz^{-j})}\frac{dz}{2\pi
iz}=-\mu_{k,\ell +m},
\eean
 with initial condition
$\mu_{k\ell}(0,0)=\dt_{k\ell}$. In matrix notation, this
amounts to the system of differential
equations\footnote{The operator $\Lb$ is the
semi-infinite
 shift matrix, with zeroes everywhere, except for $1$'s just above
 the diagonal, i.e., $(\Lambda v)_n=v_{n+1}$.
  $I_{\iy}$ is the semi-infinite identity matrix.
}
$$
\frac{\pl m_{\iy}}{\pl
t_i}=\Lb^im_{\iy}\mbox{~and~}\frac{\pl m_{\iy}}{\pl
s_i}=-m_{\iy}(\Lb^{\top})^i,\mbox{~with initial
condition~}m_{\iy}(0,0)=I_{\iy}.
$$
The solution to this initial value problem is given by
\newline ({\bf i})
\be m_{\iy}{(t,s)}=(\mu_{k\ell}(t,s))_{k,\ell\geq 0}
 \label{solution1}\ee
  and
\newline ({\bf ii})
 \be
m_{\iy}(t,s)=e^{\sum^{\iy}_1t_i\Lb^i}m_{\iy}(0,0)e^{-\sum^{\iy}_1s_i\Lb^{\top_i}
}, \label{solution2}\ee
  where
\bean e^{\sum_{1}^{\iy}t_i \Lb^i}=\sum_0^{\iy}
\gs_i(t)\Lb^i&=& \left(\begin{array}{ccccc}
1&\gs_1(t)&\gs_2(t)&\gs_3(t)&...\\
0&1&\gs_1(t)&\gs_2(t)& ...\\
0&0&1&\gs_1(t)&...\\
0&0&0&1& \\
\vdots&\vdots&\vdots&\vdots&
\end{array}\right)=(\gs_{j-i}(t))_{1\leq i<\iy\atop{1\leq j<\iy}}.
 \eean
 Then, by the uniqueness of solutions of ode's, the two solutions
 coincide, and in particular the $n\times n$ upper-left blocks
of (\ref{solution1}) and (\ref{solution2}), namely
 \be
 m_n(t,s)=E_n(t)m_{\iy}(0,0)E^{\top}_n(-s),
 \ee
 where
 \bean
E_n(t)
&=&\left(\begin{array}{ccccccc}
1&\gs_1(t)&\gs_2(t)&\gs_3(t)&...&\gs_{n-1}(t)&...\\
0&1&\gs_1(t)&\gs_2(t)& ...&\gs_{n-2}(t)&...\\
\vdots& \\
& & & & &\gs_1(t)&...\\
0& &...& &0&1&...
\end{array}\right)=(\gs_{j-i}(t))_{1\leq i<n\atop{1\leq j<\iy}}
 .\eean
 Therefore the determinants coincide:
\be \det
m_n(t,s)=\det(E_n(t)m_{\iy}(0,0)E^{\top}_n(-s)).
 \label{Toeplitz}\ee

 We shall need to expand the right hand side of (\ref{Toeplitz}) in
 ``Fourier series", which is based on the following
 Lemma:

\begin{lemma}\label{Lemma:Toeplitz}  Given the semi-infinite initial condition
$m_{\iy}(0,0)$, the expression below admits an expansion
in Schur polynomials, 
 \be
\det(E_n(t)m_{\iy}(0,0)E^{\top}_n(-s))
 = \sum_{\lb,~ \nu \atop
  \lb^{\top}_1,~  \nu^{\top}_1 \leq n}
    \det(m_{}^{\lb,\nu}) \gs_{\lb}(t) \gs_{\nu}(-s),
 ~~ \mbox{for $n>0$}, \ee
where the sum is taken over all Young diagrams $\lb$ and
$\nu$, with first columns $\lb^{\top}_1$ and
$\nu^{\top}_1\leq n $ and where
  \be
   m_{}^{\lb,\nu}:=\left( \mu_{\lb_{i}-i+n,\nu_{j}-j+n}
    \right)_{1\leq i,j\leq n}.
 \ee

\end{lemma}

\proof The proof of this Lemma is based on the
{Cauchy-Binet formula},
  which affirms that
given
two matrices $\renewcommand{\arraystretch}{0.5}
\begin{array}[t]{c}
A\\
{\scriptstyle (m,n)}
\end{array}
\renewcommand{\arraystretch}{1}$, $\renewcommand{\arraystretch}{0.5}
\begin{array}[t]{c}
B\\
{\scriptstyle (n,m)}
\end{array}
\renewcommand{\arraystretch}{1}$, for $n$ large $\geq m$

\bea \det(AB)&=&\det\left
(\sum_ia_{\ell i}b_{ik}\right)_{1\leq k,\ell\leq m}\no\\
&=&\sum_{1\leq i_1<...< i_m\leq
n}\det(a_{k,i_{\ell}})_{1\leq k,\ell\leq
m}\det(b_{i_k,\ell})_{1\leq k,\ell\leq m}.
\label{Cauchy-Binet}\eea
 Note that every decreasing sequence $
\iy>k_n>\ldots>k_1\geq 1$ can be mapped into a Young
diagram $\lb_1\geq \lb_2\geq...\geq \lb_n\geq 0$, by
setting $k_j=j+\lb_{n+1-j}$. Relabeling the indices
$i,j$ with $1\leq i,j\leq n$, by setting
$j':=n-j+1,~i':=n-i+1$, we have $1\leq i',j'\leq n$ and
$k_j-i=\lb_{j'}-j'+i'$ and $k_i-1=\lb_{i'}-i'+n$. In
other terms, the sequence of integers $
\iy>k_n>\ldots>k_1\geq 1$ leads to a partition
$\lb_1=k_n-n>\lb_2=k_{n-1}-n+1>\ldots>\lb_n=k_1-1\geq
0$. The same can be done for the sequence $1\leq
\ell_1<...<\ell_n<\iy$, leading to the Young diagram
$\nu$, using the same relabeling.
 Applying the Cauchy-Binet formula twice, expression (\ref{Toeplitz}) leads to:
  \bean
 \lefteqn{
 \det \left( E_n(t)~m_{\iy}(0,0)
  ~E_n^{\top}(-s)\right)}\\
 &=& 
    \sum_{1\leq k_1<...<k_n<\iy}\det
    \left(\gs_{k_j-i}(t)\right)_{1\leq i,j \leq n}
    \det
    \left((m_{\iy}(0,0)E_n^{\top}(-s))_{k_i,\ell}
    \right)_{1\leq i,\ell\leq n}\\
 &=& \hspace{-.3cm}
    \sum_{1\leq k_1<...<k_n<\iy}\det
    \left(\gs_{k_j-i}(t)\right)_{1\leq i,j \leq n}
    \det
    \left( (\mu_{k_i-1,j-1})_{ 1\leq i\leq n\atop
    1\leq j<\iy }
      (\gs_{i-j}(-s))_{ 1\leq i<\iy \atop 1\leq j\leq n
     })\right) \\
 &=& 
   \sum_{1\leq k_1<...<k_n<\iy}\det
    \left(\gs_{k_j-i}(t)\right)_{1\leq i,j \leq n}\\
 &&~~~~~~  \sum_{1\leq \ell_1<...<\ell_n<\iy}\det
    \left( \mu_{k_i-1,\ell_j-1}\right)_{1\leq i,j \leq n}
    \det
    \left( \gs_{\ell_i-j}(-s)\right)_{1\leq i,j \leq n}
    \\
 &=&  
   \sum_{\lb\in \BY \atop  \lb^{\top}_1\leq n}\det
    \left(\gs_{\lb_{j'}-j^{\prime}+i'}(t)
    \right)_{1\leq i',j^{\prime} \leq n}
    \\
 &&~~~~~~  \sum_{\nu\in \BY \atop  \nu^{\top}_1\leq n}\det
    \left( \mu_{\lb_{i'}-i'+n,\nu_{j'}-j'+n}
    \right)_{1\leq i',j' \leq n}
    \det
    \left( \gs_{\nu_{i'}-i^{\prime}+j'}(-s)
    \right)_{1\leq i',j' \leq n}
    \\
 &=&
    \sum_{\lb,~ \nu\in \BY \atop   \lb^{\top}_1,~  \nu^{\top}_1 \leq n}
    \det(m_{}^{\lb,\nu}) \gs_{\lb}(t) \gs_{\nu}(-s),
\eean which establishes Lemma \ref{Lemma:Toeplitz}.\qed

\newpage

Continuing the proof of Proposition \ref{prop:Toeplitz},
apply now Lemma \ref{Lemma:Toeplitz} to
$m_{\iy}(0,0)=I_{\iy}$, leading to
$$
\det m^{\lb,\nu}=\det(\mu_{\lb_i-i+n,\nu_j-j+n})_{1\leq
i,j\leq n}\neq 0
$$
if and only if
$$
 \lb_i-i+n=\nu_i-i+n  ~~\mbox{~for all~}1\leq i\leq n,
$$
 i.e., $\lb=\nu$, in which case
$$
\det m^{\lb,\lb}=1.
$$
Therefore, from (\ref{Toeplitz}), it follows
$$
\sum_{\lb\in \BY\atop \lb^{\top}_1  \leq
n}\gs_{\lb}(t)\gs_{\lb}(-s)=\det\left(
\int_{S^1}\frac{dz}{2\pi
iz}z^{k-\ell}e^{\sum_1^{\iy}(t_iz^i-s_iz^{-i})}\right)_
{1\leq k,\ell\leq n}.
$$
 So, we have, changing
$s\mapsto -s$,
 \bean
P\left(\lb ~\mbox{with}~ 
 ~\lb^{\top}_1\leq
n\right)&=&  Z^{-1}
 \sum_{\lb\in \BY\atop \lb_1  \leq
n}\gs_{\lb}(t)\gs_{\lb}(s) \\
 &=& Z^{-1}\det\left(\oint_{S^1}
\frac{dz}{2\pi iz} z^{k-\ell}
 e^{\sum_1^{\iy}(t_i z^i+s_i z^{-i})}\right)_{1\leq k, \ell\leq n}
, \eean
and, using $\gs_{\lb}(-t)=(-1)^{|\lb|}
\gs_{\lb^{\top}}(t)$, we also have
  \bean
P\left(\lb~\mbox{with} 
 ~\lb^{ }_1\leq
n\right)&=&  Z^{-1}
 \sum_{\lb\in \BY\atop \lb^{}_1  \leq
n}\gs_{\lb}(t)\gs_{\lb}(s) \\
 &=&  Z^{-1}
 \sum_{\lb\in \BY\atop \lb^{\top}_1  \leq
n}\gs_{\lb^{\top}}(t)\gs_{\lb^{\top}}(s) \\
 &=&  Z^{-1}
 \sum_{\lb\in \BY\atop \lb^{\top}_1  \leq
n}\gs_{\lb}(-t)\gs_{\lb}(-s) \\
 &=& Z^{-1}\det\left(\oint_{S^1}
\frac{dz}{2\pi iz} z^{k-\ell}
 e^{-\sum_1^{\iy}(t_i z^i+s_i z^{-i})}\right)_{1\leq k, \ell\leq n}
, \eean where
$$
Z=e^{\sum_1^{\iy} it_i s_i},
$$
ending the proof of Proposition \ref{prop:Toeplitz}.
\qed

\newpage


\subsection{The calculus of infinite wedge spaces}

\newcommand{\fem}{f_{\emptyset}}

The material in this section can be found in Kac
\cite{Kac} and Kac-Raina \cite{KacR} and many specific
results are due to Borodin-Okounkov-Olshanski (see
\cite{Borodin1,Borodin2,BOO}). Given a vector space
$V=\displaystyle{\oplus_{j\in\BZ}}\mathbb C v_j$ with
inner-product $\la v_i,v_j\ra =\dt_{ij}$, the infinite
wedge space $V^{\iy}=\Lambda^{\iy}V$ is defined as

$$
V^{\iy}=\mbox{span}\left\{\begin{array}{ll}
 &s_1>s_2>\ldots\\
v_{s_1}\wedge v_{s_2}\wedge v_{s_3}\wedge\ldots,& \\
& s_k=-k\mbox{~for~}k>\!\!\!>0
\end{array}\right\}
$$
containing the ``vacuum"
$$
\fem=v_{-1}\wedge v_{-2}\wedge\ldots\,.
$$
The vector space $V^{\iy}$ comes equipped with an
inner-product $\la \, ,\,\ra$, making the basis vectors
$v_{s_1}\wedge v_{s_2}\wedge\ldots$ orthonormal.
To each $k\in\BZ$, we associate two operations, a {\em
wedging} with $v_k$ and a {\em contracting operator},
removing a $v_k$,
\bean
& &\psi_k:V^{\iy}\rg V^{\iy}:f\longmapsto\psi_k(f)=v_k\wedge f\\
& & \\
& &\psi^*_k:V^{\iy}\rg
V^{\iy}:v_{s_1}\!\wedge\!\ldots\!\wedge\!
v_{s_i}\!\wedge\!\ldots \longmapsto\sum_i (-1)^{i+1}\la
v_k,v_{s_i}\ra v_{s_1}\wedge  \ldots\wedge\hat
v_{s_i}\wedge\ldots
 \eean
Note that
\bean
\begin{array}{ll}
\psi_k(f)=0, &\mbox{if $v_k$ figures in $f$}\\
\\
\psi^*_k(v_{s_1}\wedge\ldots)=0,&\mbox{if
$k\not\in(s_1,s_2,\ldots).$}
\end{array}
\eean
Define the shift
$$
\Lambda^r:=\sum_{k\in\BZ}\psi_{k+r}\psi^*_k,\qquad
r\in\BZ
$$
acting on $V^{\iy}$ as follows

\bean
  \Lambda^r v_{s_1}\wedge v_{s_2}\wedge\ldots
&=&
 ~v_{s_1+r}\wedge v_{s_2}\wedge v_{s_3}\wedge\ldots
 \\&& + v_{s_1} \wedge
v_{s_2+r}\wedge v_{s_3}\wedge\ldots
 \\&& +v_{s_1}\wedge v_{s_2}\wedge
v_{s_3+r}\wedge\ldots \\&& +\ldots \eean
One checks that
\be
[\Lambda^r,\psi_k]=\psi_{k+r},[\Lambda^r,\psi_k^*]=-\psi^*_{k-r}
\label{commutation relation}\ee
\be [\Lambda^k,\Lambda^{\ell}]=\ell~\dt_{k,-\ell} \ee
 and hence
\be [\sum_{i\geq 1}t_i\Lambda^i,\sum_{j\geq
1}s_i\Lambda^{-i}]=-\sum_{i\geq 1}it_is_i. \ee


\begin{lemma} {(Version of the Cauchy identity)}
\be
 e^{\sum_{i\geq 1}t_i\Lambda^i}e^{-\sum_{j\geq
1}s_i\Lambda^{-j}}
 =
  e^{\sum_{i\geq 1}it_is_i}
   e^{-\sum_{j\geq 1}s_j\Lambda^{-j}}e^{\sum_{i\geq
1}t_i\Lambda^i} \label{Cauchy1}\ee

\end{lemma}

\proof When two operators $A$ and $B$ commute with their
commutator $[A,B]$, then according to Kac \cite{Kac}, p
308,
$$
e^Ae^B=e^Be^Ae^{[A,B]}.
$$
Setting $A=\displaystyle{\sum_{i\geq i}}t_i\Lambda^i$
and $B=-\displaystyle{\sum_{j\geq 1}}s_j\Lambda^{-j}$,
we find \bean e^{\sum_{i\geq
1}t_i\Lambda^i}e^{-\sum_{j\geq
1}s_i\Lambda^{-j}}&=&e^{-\sum_{j\geq
1}s_j\Lambda^{-j}}e^{\sum_{i\geq 1}t_i\Lambda^i}
e^{-\left[\sum_{i\geq 1}t_i\Lambda^i,\sum_{j\geq 1}s_j\Lambda^{-j}\right]}\\
&=&e^{-\sum_{j\geq 1}s_j\Lambda^{-j}}e^{\sum_{i\geq
1}t_i\Lambda^i}e^{\sum_{i\geq
1}it_is_i}\\
&=&e^{\sum_{i\geq 1}it_is_i}~e^{-\sum
s_j\Lambda^{-j}}e^{\sum t_i\Lambda^i}  .\eean\qed


It is useful to consider the generators of $\psi_i$ and
$\psi^*_i$:

\be \psi(z)=\sum_{i
\in\BZ}z^i\psi_i,\qquad\psi^*(w)=\sum_{j\in\BZ}w^{-j}\psi_j^*.
\ee
 From (\ref{commutation relation}), it follows that

\bean
  [ \Lambda^r,\psi(z)]
   &=&\sum_kz^k[\Lambda^r,\psi_k]
    =\sum_k z^k\psi_{k+r}=\frac{1}{z^r}\psi(z),
  \eean
$$
 [ \Lambda^r ,\psi^*(w)]
  = -\frac{1}{w^r}\psi^*(w).
 $$
The two relations above lead to the following, by taking
derivatives in $t_i$ and $s_i$ of the left hand side and
setting all $t_i=s_i=0$:

\bea
e^{\pm\sum_1^{\iy}t_i\Lambda^i}\left\{\begin{array}{c}
\psi(z)\\
\psi^*(w)
\end{array}\right\}e^{\mp\sum_1^{\iy}t_i\Lambda^i}&=& \left\{\begin{array}{c}
e^{\pm\sum_1^{\iy}t_r/z^r}\psi(z)\\
e^{\mp\sum t_r/w^r}\psi^*(w)
\end{array} \right.\nonumber\\
& &\nonumber\\
e^{\pm\sum_1^{\iy}s_i\Lambda^{-i}}\left\{\begin{array}{c}
\psi(z)\\
\psi^*(w)
\end{array}\right\}e^{\mp\sum_1^{\iy}s_i\Lambda^{-i}}&=&
\left\{\begin{array}{c}
e^{\pm\sum_1^{\iy}s_rz^r}\psi(z)\\
e^{\mp\sum_1^{\iy} s_rw^r}\psi^*(w)
\end{array} \right.
\label{conjugation}\eea
To each partition
$\lb=(\lb_1\geq\lb_2\geq\ldots\geq\lb_m\geq 0)$,
associate a vector
$$
f_{\lb}:=v_{\lb_1-1}\wedge v_{\lb_2-2}\wedge\ldots\wedge
v_{\lb_m-m}\wedge v_{-(m+1)}\wedge
v_{-(m+2)}\wedge\ldots\in V^{\iy}.
$$
The following lemma holds:

\begin{lemma}\label{Lemma:2.5}
\be e^{\sum_{i\geq 1} t_i \Lb^i} f_{\lb} = \sum_{\mu \in
\BY\atop{\mu\supset\lb}}
 {\gs}_{\mu\backslash \lb}(t) f_{\mu},
 \quad e^{\sum_{i\geq 1}
t_i\Lb^{-i}} f_{\lb} = \sum_{\mu\in\BY\atop{\mu\subset
\lb}} \gs_{\lb\backslash \mu} (t) f_{\mu}.
 \label{Nakayama}\ee
 In particular,
  \be
  e^{\sum_{i\geq 1} t_i \Lb^i} f_{\emptyset} =
   \sum_{\mu \in
\BY }
 {\gs}_{\mu   }(t) f_{\mu}
 \qquad \mbox{and}\qquad e^{\sum_{i\geq 1}
t_i\Lb^{-i}}  \fem=\fem.
 \label{Nakayama1}\ee

\end{lemma}

\proof First notice that a matrix $A\in GL_{\iy}$ acts
on $V^{\iy}$ as follows

$$
A(v_{s_1}\wedge v_{s_2}\wedge\ldots)=\sum_{s'_1>
s'_2>\ldots}\det\left(A ^{s_1,s_2,\ldots}
_{s'_1,s'_2,\ldots}\right)v_{s'_1}\wedge
v_{s'_2}\wedge\ldots ,
$$
where

$$
A^{s_1,s_2,\ldots}_{s'_1,s'_2,\ldots}=\left\{\begin{array}{l}
\mbox{matrix located at the intersection of the}\\
\mbox{rows $s'_1,s'_2,\ldots$ and columns
$s_1,s_2,\ldots$ of $A$.}
\end{array}
\right\} .$$
 Here the rows and columns of the bi-infinite matrix are labeled by

\newpage

  $$    1~~~~~0~~~      -1    ~-2~-3 $$
$$
 \left(\begin{array}{cccccccccc}
  *&&*&&*&&*&&*\\ \\
  *&&*&&*&&*&&*\\  \\
  *&&*&&\bullet&&*&&*\\  \\
  *&&*&&*&&*&&*\\  \\
  *&&*&&*&&*&&*\\
  \end{array}\right)
$$

\vspace{-4.8cm}

$$\hspace{-6.6cm}\begin{array}{c}
  1\\ \\0\\ \\-1\\ \\-2\\ \\-3
  \end{array}
  $$

\vspace{1cm}

\noindent Hence, for the bi-infinite matrix
$e^{\sum_1^{\iy}t_i\Lb^i}$,
$$
\det \left(\left(e^{\sum_1^{\iy}t_i\Lb^i}\right)
 ^{s_1,s_2,\ldots}_{s'_1,s'_2,\ldots}\right)
 =\det \left(\gs_{s'_i-s_j}\right)$$

Setting $s_i=\lb_i-i$ and defining $\mu_i$ in the
formula below by $ s'_i=\mu_i-i$, one checks
\bean
e^{\sum_1^{\iy}t_i\Lb^i}f_{\lb}&=&e^{\sum_1^{\iy}t_i\Lb^i}(v_{s_1}\wedge
v_{s_2}\wedge\ldots)\\
 \\
&\stackrel{*}{=}&\sum_{s'_1> s'_2>\ldots}\det
\left(\left(e^{\sum_1^{\iy}t_i\Lb^i}\right)
 ^{s_1,s_2,\ldots}_{s'_1,s'_2,\ldots}\right)v_{s'_1}\wedge
v_{s'_2}\wedge\ldots\\
 \\
&=&\sum_{s'_1>s'_2>\ldots}\det\left(\gs_{s'_i-s_j}(t)\right)_{1\leq
i,j<\iy}v_{s'_1}\wedge
v_{s'_2}\wedge\ldots\\
 \\
&=&\sum_{\mu_1-1>\mu_2-2>\ldots}
 \det\left(\gs_{(\mu_i-i)-(\lb_j-j)}(t)\right)_{1\leq
i,j\leq\iy}v_{\mu_1-1}\wedge v_{\mu_2-2}\wedge\ldots\\
 \\
&=&\sum_{\mu\in\BY\atop{\mu\supset\lb}}\gs_{\mu\backslash\lb}(t)f_{\mu}
\eean
 The second identity in (\ref{Nakayama}) is shown in the
same way. Identities (\ref{Nakayama1}) follow
immediately from (\ref{Nakayama}), ending the proof of
Lemma \ref{Lemma:2.5}. \qed

We also need:
\begin{lemma}\label{Lemma:2.6}
 For
  \be
   \Psi^*_{b_k}=\displaystyle{\sum_{i\geq
1}}b_{ki}\psi^*_{-i}~~\mbox{and}~~
\displaystyle{\Psi_{a_k}=\sum_{i\geq 1}a_{ki}\psi_{-i}},
 \label{Psi's}\ee
   the following identity holds:
 \be
\left\langle \Psi_{a_1} \ldots \Psi_{a_m}
\Psi^{\ast}_{b_m}\ldots \Psi^{\ast}_{b_1}
f_{\emptyset},f_{\emptyset}\right\rangle
 =  \det\Bigl(\left\langle
\Psi_{a_k}  \Psi^{\ast}_{b_{\ell}} f_{\emptyset},
f_{\emptyset}\right\rangle\Bigr)_{1\leq k, \ell\leq m}
 \ee

\end{lemma}

\proof First one computes for the $\Psi$'s as in
(\ref{Psi's}),
\bean
 \Psi^*_{b_m}\ldots \Psi^*_{b_1}\fem
&=&\sum_{i_1>\ldots
>i_m\geq
1}(-1)^{\sum^m_1(i_k+1)}\det\left(b_{k,i_{\ell}}\right)_{1\leq
k,\ell\leq m}\\ \\&
&\hspace*{-.1cm}v_{-1}\wedge\ldots\wedge\hat
v_{-i_m}\wedge\ldots\wedge\hat
v_{-i_{m-1}}\wedge\ldots\wedge\hat
v_{-i_2}\wedge\ldots\wedge\hat v_{-i_1}\wedge\ldots
\eean
Then acting with the $\Psi_{a_k}$ as in (\ref{Psi's}),
it suffices to understand how it acts on the wedge
products appearing in the expression above, namely:
\bean
\lefteqn{\Psi_{a_1}\ldots\Psi_{a_m}v_{-1}\wedge\ldots\wedge\hat
v_{-i_m} \wedge\ldots\wedge \hat v_{-i_{m-1}}
\wedge\ldots\wedge \hat v_{-i_1}
\wedge\ldots }~~~~~~~~~~~~~~~~~~~~~~~~~~~~~~\\
 \\
&=&(-1)^{\sum^m_1(i_k+1)}\det\left(a_{k,i_{\ell}}\right)_{1\leq
k,\ell\leq m}\fem \eean
 Thus, combining the two, one finds, using the
 Cauchy-Binet formula (\ref{Cauchy-Binet}) in the last equality,
\bean
\lefteqn{\Psi_{a_1}\ldots\Psi_{a_m}\Psi^*_{b_m} \ldots\Psi^*_{b_1}\fem}\\
&=&\sum_{i_1>\ldots >i_m\geq 1}
\det\left(a_{k,i_{\ell}}\right)_{1\leq k,\ell\leq
m}\det\left(b_{k,i_{\ell}}\right)_{1\leq
k,\ell\leq m}\fem\\
&=& \det\left(\sum_ia_{\ell i }b_{ki}\right)_{1\leq
k,\ell\leq m}\fem   .
 \eean
In particular for $m=1$,
$$
\Psi_{a_{\ell}}\Psi^*_{b_k}\fem =\sum_ia_{\ell
i}b_{ki}~\fem.
$$
Hence
\bean
\la\Psi_{a_1}\ldots\Psi_{a_m}\Psi^*_{b_m}\ldots\Psi^*_{b_1}\fem,\fem\ra
&=&\det\left(\sum_ia_{ki}b_{\ell i}\right)_{1\leq k,\ell\leq m} \\
&=&\det\left(\la\Psi_{a_k}\Psi^*_{b_{\ell}}\fem,\fem\ra\right)_{1\leq
k,\ell\leq m}, \eean
 ending the proof of Lemma \ref{Lemma:2.6}. \qed


\subsection{Probability on partitions
 expressed as Fredholm determinants}

Remember the definition (\ref{def: Fredholm det}) of the
Fredholm determinant of a kernel $K(i,j)$, with $i,j\in
\BZ$ and restricted to $[n, \iy)$. This statement has
appeared in the literature in some form (see
Case-Geronimo \cite{CaseG}) and a more analytic
formulation has appeared in Basor-Widom \cite{Basor}.
The proof given here is an ``integrable" one, due to
Borodin-Okounkov \cite{Bo-Ok}.

\begin{proposition}\label{theo:Fredholm}
 $$
P(\lb  ~\mbox{with}~\lb_1 \leq n) =
\det\left(I-K(k,\ell)\Big\vert_{[n,n+1,\ldots]}\right),
$$
where $K(k,\ell)$ is a kernel
 \bea
K(k,\ell)&=&\left(\frac{1}{2\pi i}\right)^2
 \oint_{|w|=\rho<1}
 \oint_{|z|=\rho^{-1}>1}
  \frac{dz~
dw}{z^{k+1}w^{-\ell}} \frac{e^{V(z)-V(w)}}{z-w}
 \no\\&&\no\\&=&\frac{1}{k-\ell}
 \left(\frac{1}{2\pi
i}\right)^2\oint\!\!\oint _{{|w|=\rho<1}\atop
{|z|=\rho^{-1}>1}}\frac{dz~dw}{z^{k+1}
w^{-\ell}}\frac{z\frac{d}{dz}V(z)-w\frac{d}{dw}V(w)}{z-w}
 e^{V(z)-V(w)},\no\\
  && \hspace*{8cm}\mbox{for $k\neq \ell$,}\no\\
 \label{kernel}
 \eea
 with
  $$
    V(z)=-\sum_{i\geq
1}(t_iz^{-i}-  s_iz^i).
$$
\end{proposition}

The proof of Proposition \ref{theo:Fredholm} hinges on
the following Lemma, due to Borodin-Okounkov
(\cite{Okounkov0,Okounkov1,Bo-Ok}):

\begin{lemma}\label{Theo:det} If $X=\{x_1,\ldots,x_m\}\subset \BZ$ and
 $ S(\lb):=\{ \lb_1-1,\lb_2-2,\lb_3-3,\ldots\}$
$$
P\left(\lb~\big|~ S(\lb)\supset X  \right)=\frac{1}{Z}
\sum_{\lb \mbox{\tiny{such that}}\atop{S(\lb)\supset X}}
s_{\lb}(t) s_{\lb} (s) = \det
\left(K(x_i,x_{j})\right)_{1\leq i,j\leq m}.
$$

\end{lemma}

\medskip\noindent{\it Proof of Proposition \ref{theo:Fredholm}:\/}
Setting $A_k:=\{\lb\mid k\in S(\lb)\}$ with $ S(\lb):=\{
\lb_1-1,\lb_2-2,\lb_3-3,\ldots\}$, one computes (the
$x_i$ below are integers)
\bean
\lefteqn{P(\lb  ~\mbox{with}~\lb_1\leq n)}\\
&=&P(\lb  ~\mbox{with all $\lb_i\leq n$)}\\
&=&P(\lb~\mbox{with}~ S(\lb)\cap\{n,n+1,n+2,...\}=\emptyset)\\
&=&1-P(\lb~\mbox{with}~ S(\lb)\mbox{~contains some $k$ for $k\geq n$})\\
&=&1-P(A_n\cup A_{n+1}\cup ...)\\
&=&1-\sum_{n\leq i}P(A_i)+\sum_{n\leq i<j}P(A_i\cap
A_j)-\sum_{n\leq i<j<k}P(A_i\cap A_j\cap
A_k)+...\\
& & \hspace*{5cm}\mbox{using Poincar\'e's formula}\\
&=&\sum^{\iy}_{m=0}(-1)^m\sum_{n\leq x_1<...< x_m}P(\lb
~\mbox{with}~
\{x_1,...,x_m\}\subset S(\lb))\\
&=&\sum^{\iy}_{m=0}(-1)^m\sum_{n\leq x_1<...<
x_m}\det(K(x_i,x_j))_{1\leq
i,j\leq m}\\
&=&\det\left(I-K(i,j)\Big\vert_{[n,n+1,...]}\right),
\eean
 from which Proposition \ref{theo:Fredholm} follows.
\qed


\medskip\noindent{\it Proof of Lemma \ref{Theo:det}:\/}
Remembering the probability measure introduced in
(\ref{generalized probability}), we have that
$$
P\left(\lb~\big|~ S(\lb)\supset X  \right)=\frac{1}{Z}
\sum_{\lb \mbox{\tiny{such that}}\atop{S(\lb)\supset X}}
s_{\lb}(t) s_{\lb} (s) .$$
Next, from the wedging-contracting operation
 \bean
 \psi_x\psi_x^{\ast}f_{\lb} &=&  f_{\lb},~~\mbox{if $x\in S_{\lb}$}
 \\
 &=&  0,~~~~\mbox{if $x\notin S_{\lb}$},
 \eean
 and using both relations (\ref{Nakayama}),
  one first computes:
\bean \lefteqn{  \left\langle e^{\sum_1^{\iy} s_i
\Lb^{-i}}\prod_{x\in X}\psi_x \psi_x^{\ast}
e^{\sum_1^{\iy} t_i \Lb^i} f_{\emptyset},
f_{\emptyset}\right\rangle}\\
 &=&  \left\langle e^{\sum s_i \Lb^{-i}}\prod_{x\in X}\psi_x
\psi_x^{\ast} \sum_{\lb} \gs_{\lb}(t) f_{\lb},
f_{\emptyset}\right\rangle\\
 &=& \sum_{\lb} {\bf s}_{\lb}(t)  \left\langle e^{\sum s_i
\Lb^{-i}}\prod_{x\in X}\psi_x \psi_x^{\ast} f_{\lb},
f_{\emptyset}\right\rangle\\
 &=& \sum_{\lb\atop{S(\lb) \supset X}} {\bf s}_{\lb}(t)
\left\langle e^{\sum s_i \Lb^{-i}} f_{\lb},
f_{\emptyset}\right\rangle\\
 &=& \sum_{\lb\atop{S(\lb) \supset X}} {\bf s}_{\lb}(t)
\sum_{\mu\subset\lb} \left\langle
{\bf s}_{\lb\backslash\mu}(s) f_{\mu},f_{\emptyset}\right\rangle\\
 &=& \sum_{\lb\atop{S(\lb) \supset X}} {\bf s}_{\lb}(t)
\sum_{\mu\subset\lb} {\bf
s}_{\lb\backslash\mu}(s)\left\langle
 f_{\mu},f_{\emptyset}\right\rangle
  = \sum_{\lb\mbox{\tiny{~such that}}\atop{S(\lb) \supset X}} {\bf
s}_{\lb}(t)    {\bf s}_{\lb}  (s).
 \eean
Using this fact, one further computes (in the
exponentials below the summation $\Sigma$ stands for
$\Sigma_1^{\iy}$)
\bean \lefteqn{\frac{1}{Z}\sum_{\lb\mbox{\tiny{~such
that}}\atop{S(\lb) \supset X}} {\bf s}_{\lb}(t)  {\bf
s}_{\lb}
(s)}\\
&=&\frac{1}{Z} \left\langle e^{\sum  s_i
\Lb^{-i}}\prod_{x\in X}\psi_x \psi_x^{\ast} e^{\sum  t_i
\Lb^i} f_{\emptyset},
f_{\emptyset}\right\rangle\\
&=&
 \frac{1}{Z} \left\langle e^{\sum  s_i
\Lb^{-i}}\psi_{x_m} \ldots \psi_{x_1}
\psi^{\ast}_{x_1}\ldots \psi^{\ast}_{x_m}
 e^{\sum  t_i \Lb^i} f_{\emptyset},
f_{\emptyset}\right\rangle\\
&&\hspace*{4cm}\mbox{using } \psi_{x_i}
\psi^{\ast}_{x_j}= -\psi^{\ast}_{x_j}
\psi_{x_i}\mbox{~for~}i\neq j\\
&=&
 \frac{1}{Z} \left\langle e^{\sum  s_i
\Lb^{-i}}\psi_{x_m} \ldots \psi_{x_1}
\psi^{\ast}_{x_1}\ldots \psi^{\ast}_{x_m}
 e^{\sum t_i \Lb^i} e^{-\sum s_i\Lb^{-i}} f_{\emptyset},
e^{-\sum  t_i\Lb^{-i}}f_{\emptyset}\right\rangle\\
&&\hspace*{9cm}\mbox{using (\ref{Nakayama1})
}
%
 \\&=&
  \frac{1}{Z} \left\langle e^{-\sum  t_i \Lb^{i}} e^{\sum
s_i\Lb^{-i}} \psi_{x_m} \ldots \psi_{x_1}
\psi^{\ast}_{x_1}\ldots \psi^{\ast}_{x_m}
 e^{\sum t_i \Lb^i} e^{-\sum s_i\Lb^{-i}} f_{\emptyset},
f_{\emptyset}\right\rangle\\
 &=&
 \left\langle e^{-\sum  t_i \Lb^{i}} e^{\sum s_i\Lb^{-i}}
\psi_{x_m} \ldots \psi_{x_1} \psi^{\ast}_{x_1}\ldots
\psi^{\ast}_{x_m} e^{-\sum s_i\Lb^{-i}} e^{\sum t_i
\Lb^i} f_{\emptyset},
f_{\emptyset}\right\rangle\\
 && \hspace{5cm}\mbox{using Cauchy's identity (\ref{Cauchy1}) }
\eean\bean \\
&=&\left\langle \Psi_{x_m} \ldots \Psi_{x_1}
\Psi^{\ast}_{x_1}\ldots \Psi^{\ast}_{x_m} f_{\emptyset},
f_{\emptyset}\right\rangle , \eean
 where
\bean \Psi_{k } &=& e^{-\sum  t_i\Lb^i} e^{\sum
s_i\Lb^{-i}} \psi_k e^{-\sum s_i\Lb^{-i}} e^{\sum t_i\Lb^i}\\
\Psi_{k}^{\ast} &=& e^{-\sum  t_i\Lb^i} e^{\sum
s_i\Lb^{-i}} \psi_k^{\ast} e^{-\sum s_i\Lb^{-i}}
 e^{\sum t_i\Lb^i}
 . \eean
Then, using Lemma \ref{Lemma:2.6}, the expression above
equals
\bean  \left\langle \Psi_{x_m} \ldots \Psi_{x_1}
\Psi^{\ast}_{x_1}\ldots \Psi^{\ast}_{x_m}
f_{\emptyset},f_{\emptyset}\right\rangle
 &=&  \det\left(\left\langle
\Psi_{x_k}  \Psi^{\ast}_{x_{\ell}}
f_{\emptyset},f_{\emptyset}\right\rangle\right)_{1\leq k, \ell\leq m}\\
& =&\det\left( K(x_k, x_{\ell})\right)_{1\leq k,
\ell\leq m} \eean where

$$
K(k,\ell)=\left\la \Psi_{k}\Psi^{\ast}_{\ell}
 f_{\emptyset},f_{\emptyset}\right\ra=
\left\la e^{-\Sg t_i\Lb^i}e^{\Sg
s_i\Lb^{-i}}\psi_{k}\psi^*_{\ell}e^{-\Sg s_i\Lb^{-i}}
e^{\Sg t_i\Lb^i}f_{\emptyset },f_{\emptyset }
 \right\ra.
$$
Using
 \bean
 \left\la\psi_i \psi_j^{\ast}f_{\emptyset},
 f_{\emptyset}\right\ra  &=& 1~~~\mbox{if}~ i=j<0  \\
                         &=& 0~~~\mbox{otherwise}~
   \eean
and setting
 \be
 V(z)=-\sum_{i\geq
1}(t_iz^{-i}-s_iz^i)
 ,\ee
 the generating function of the $K(k,\ell)$ takes on the following
 simple form:
\bean \lefteqn{\sum_{k,\ell \in\BZ}z^{k}w^{-\ell
}K(k,\ell)}
 \\
&=&
 \left\la e^{-\sum t_i\Lb^i}e^{\sum s_i\Lb^{-i}}\psi(z)\psi^*(w)e^{-\Sg
s_i\Lb^{-i}}
e^{\sum t_i\Lb^i}f_{\emptyset},f_{\emptyset}\right\ra\\
&=&
 \left\la e^{-\sum t_i\Lb^i}e^{\sum s_i\Lb^{-i}}\psi(z)
 e^{-\sum s_i\Lb^{-i}}e^{\sum t_i\Lb^i}\right.\\
  &&\hspace{3cm}
 .\left.e^{-\sum t_i\Lb^i}e^{\sum s_i\Lb^{-i}}\psi^*(w)e^{-\sum
s_i\Lb^{-i}}
e^{\sum t_i\Lb^i}f_{\emptyset},f_{\emptyset}\right\ra\\
 &=&
 e^{\sum s_iz^i}e^{-\sum t_iz^{-i}}e^{-\sum
 s_iw^i}e^{\sum
t_iw^{-i}}\left\la\psi(z)\psi^*(w)f_{\emptyset},f_{\emptyset}
\right\ra
~~\mbox{using (\ref{conjugation})}\\
&=&e^{V(z)-V(w)}\la\psi(z)\psi^*(w)f_{\emptyset},f_{\emptyset}\ra\\
&=&e^{V(z)-V(w)}\sum_{i,j\in\BZ}z^iw^{-j}\la\psi_i\psi_j^*f_{\emptyset},f_{\emptyset}\ra\\
&=&e^{V(z)-V(w)}\sum_{i\geq 1}\left(\frac{w}{z}\right)^i
\eean\bean&=&e^{V(z)-V(w)}\left(\frac{w}{z}\right)\left(1+\frac{w}{z}+...\right)\\
&=&e^{V(z)-V(w)}\left(\frac{w}{z}\right)\frac{1}{1-\frac{w}{z}}\qquad\mbox{~for~
}|w|<|z|\\
&=&e^{V(z)-V(w)}\frac{w}{z-w},
 \eean
  and so
\bean K(k,\ell)
&=&\left(\frac{1}{2\pi i}\right)^2\oint_{|w|=\rho<1}
 \oint_{|z|=\rho^{-1}>1}\frac{dz~dw}{z^{k+1}
w^{-\ell}}\frac{e^{V(z)-V(w)}}{z-w}, \eean ending the
proof of Proposition \ref{Theo:det}.\qed

  Formula (\ref{kernel-non-diagonal}) in the remark is obtained by setting $z\mapsto tz$,
$w\mapsto tw$ in (\ref{kernel}),
$$
K(k,\ell)=\left(\frac{1}{2\pi
i}\right)^2\frac{1}{t^{k-\ell}}\oint\oint\frac{dz~dw}{z^
{k+1} w^{-\ell+1}}\frac{e^{V(tz)-V(tw)}}{(z/w {-1})}.
$$
 and taking $\displaystyle{\frac{\pl}{\pl t}\Bigg\vert_{t=1}}$ of both
sides, using the $t$-independence of the left hand side,
yielding (\ref{kernel-non-diagonal}) . \qed%
%


\subsection{Probability on partitions
 expressed as U(n) integrals}

This section deals with formula (\ref{formulaTheo:2.1})
in Theorem \ref{Theo:2.1}, stating that
$P(\lb~\bigr|~\lb_1\leq n)$ can be expressed as a
Unitary matrix integral. First we need a Lemma, whose
proof can be found in \cite{Carter}:

\begin{lemma}\label{lemma2.9}
If $f$ is a symmetric function of the eigenvalues
$u_1,\ldots,u_n$ of the elements in U$(n)$, then

$$
\int_{U(n)}f=\frac{1}{n!} \int_{(S^1)^n}|\Dt(u)|^2
 f(u_1,\ldots,u_n)\prod_1^n \frac{du_j}{2\pi i u_j}
  .$$
\end{lemma}
\proof Set $G:=U(n)$ and
$T:=\{\diag(u_1,\ldots,u_n),~~\mbox{with}~u_k=e^{i\theta_k
}\in S^1\}$, and let $\frak{t}$ and $\frak{g}$ denote
the Lie algebras corresponding to $T$ and $G$.  An
element in the quotient $G/T$ can be identified with
$gT$, because $g'\in gT$ implies $g^{-1}g'\in T$, and
thus there is a natural map

$$
T\times (G/T)\rg G:(t,gT)\mapsto gtg^{-1}.
$$
Note that the Jacobian $J$ of this map (with respect to
the invariant measures on $T$, $G/T$ and $G$) only
depends on $t\in {\frak t}$, because of invariance of
the measure under conjugation. For a function $f$ as
above
$$
\int_{G}f=\frac{1}{n!} \int _{T\times
(G/T)}f(gtg^{-1})J(t)dt
d(gT)=\frac{\mbox{vol}(G/T)}{n!}\int_T f(t)J(t) dt .$$

Denote the tangent space to $G/T$ at its base points by
$\frak{t}^{\bot}$, the orthogonal complement of
$\frak{t}$ in $\frak{g}=\frak{t}\oplus\frak{t}^{\bot}$.
Consider infinitesimal changes $t(1+\vr\xi)$ of $t$ with
$\xi\in\frak{t}$; also, an infinitesimal change of $1\in
G$, namely $1\mapsto 1+\vr\eta$, with
$\eta\in\frak{t}^{\bot}$. Then

\bean
gtg^{-1}-t&=&(1+\vr\eta)t(1+\vr\xi )(1-\vr\eta+{\bf O}(\vr^2))-t\\
&=&\vr(t\xi+\eta t-t\eta)+{\bf O}(\vr^2)\\
&=:&\vr\rho +{\bf O}(\vr^2) \eean and so
$$
t^{-1}\rho=\xi+(t^{-1}\eta t-\eta)\in
\frak{t}\oplus\frak{t}^{\bot}={\frak g}.
$$
Thus the Jacobian of this map is given by
 $$J(t)=\det(A(t^{-1})-I),$$
where $A(t)$ denotes the adjoint action by means of the
diagonal matrix $t$, denoted by $t=\diag
(u_1,\ldots,u_n)$, with $|u_i|=1$. Let $E_{jk}$ be the
matrix with $1$ at the $(j,k)$-th entry and $0$
everywhere else. Then
 $$(A(t^{-1})-I)E_{jk}=(u_j^{-1}u_k-1)E_{jk}$$
 and thus the matrix $A(t^{-1})-I$ has $n(n-1)$ eigenvalues
 $(u_j^{-1}u_k-1)$, with $1\leq j,k\leq n$ and $j\neq k$.
 Therefore $$\det(A(t^{-1})-I)= \prod_{j\neq k}
 (u_j^{-1}u_k-1)=\prod_{j<k} |u_j-u_k|^2,$$
using the fact that $u_j^{-1}=\bar u_j$.\qed

\begin{proposition}\label{prop:Unitary} Given the ``probability measure"

$$
P(\lb)=Z^{-1}s_{\lb}(t)s_{\lb}(s),\qquad
Z=e^{\sum_{i\geq 1} it_is_i},
$$
the following holds
\bean P\left(\lb~\Bigl|~\lb_1\leq n\right)
 &=&
 Z^{-1} \int_{U(n)} e^{-Tr~\sum_1^{\iy}(t_iX^i+s_i\bar X^i)}dX. \eean

\end{proposition}

\proof Using Lemma \ref{lemma2.9} and the following
matrix identity,
$$
 \sum_{\sigma \in S_n}\det \left(
a_{i,\sigma(j)}~b_{j,\sigma(j)}\right)_{1\leq i,j\leq n}
=
 \det \left( a_{ik}\right)_{1\leq i,k\leq n}
  \det \left( b_{ik}\right)_{1\leq i,k\leq n},
  $$
one computes
  \bea
\lefteqn { n!\int_{U(n)}e^{-Tr \sum_1^{\iy}
(t_iX^i+s_i\bar X^i)}  dX}\nonumber\\
 &=&
 \int_{(S^1)^{n}}|\Dt_n(z)|^{2}
 \prod_{k=1}^n
\left(e^{-\sum_1^{\iy}(t_i z_k^i+s_iz_k^{-i})}
 \frac{ dz_k}{2\pi i z_k}\right)\nonumber\\
&=&
 \int_{(S^1)^{n}}\Dt_n(z)\Dt_n(\bar z)
 \prod_{k=1}^n
\left(e^{-\sum_1^{\iy}(t_i z_k^i+s_iz_k^{-i})}
 \frac{dz_k}{2\pi i z_k}\right)\nonumber\\
&=&
 \int_{(S^1)^{n}}\sum_{\sigma\in S_n}\det\left(z_{\sigma(m)}^{\ell-1}
 \bar z_{\sigma(m)}^{m-1}  \right)_{1\leq \ell,m\leq n}
 \prod_{k=1}^n
\left(e^{-\sum_1^{\iy}(t_i z_k^i+s_iz_k^{-i})}
 \frac{dz_k}{2\pi i z_k}\right)
 \nonumber\\
&=&
 \sum_{\sigma\in S_n} \det\left(\oint_{S^1}z_{k}^{\ell-1}
 \bar z_{k}^{m-1}
e^{-\sum_1^{\iy}(t_i z_k^i+s_iz_k^{-i})}
 \frac{dz_k}{2\pi i z_k}\right)_{1\leq \ell,m\leq n}
  \nonumber\\
  &=&
 n!\det\left(\oint_{S^1}z^{\ell-m}
e^{-\sum_1^{\iy}(t_i z^i+s_iz^{-i})}
 \frac{dz}{2\pi i z}\right)_{1\leq \ell,m\leq n}
 =n!\det m_n(t,s)
  \nonumber
 \eea
ending the proof of Proposition \ref{prop:Unitary}.\qed

 \newpage

\section{Examples}

\subsection{Plancherel measure and Gessel's theorem}
\label{section3.1}

The point of this section will be to restrict the
probability
 \be P_{x,y} (\lb) :=\frac{1}{Z}
\tilde\gs_{\lb} (x) \tilde\gs_{\lb} (y)=\frac{1}{Z}
 \gs_{\lb} (t)  \gs_{\lb} (s)
 \label{generalized probability1}\ee considered in section
 2 (formula (\ref{generalized probability}))
 to the locus $\LR_1$, defined for real $\xi>0$
(expressed in $(x,y)$ and $(s,t)$ coordinates),
  \bea  \LR_1&=&
   \{(x,y) ~~\mbox{such that } \sum_{i\geq 1} x_i^k=\sum_{i\geq 1} y_i^k = \dt_{ki} \sqrt{\xi}  \}
   \no\\
 &=&
 \{\mbox{all $s_k=t_k=0$, except $t_1=s_1=\sqrt{\xi}$}\}
 .\label{locus1}\eea
 The reader is referred back to subsection \ref{sect 1.5.1}
 for a number of basic formulae.

\begin{theorem} \cite{Gessel,TW-permutations,Bo-Ok,AvM2} For the permutation group $S_k$, the generating function
for distribution of the length of the longest increasing
subsequence

$$
\tilde P^k(L(\pi_k)\leq n)
$$
is given by
\bean e^{-\xi}\sum_{k=0}^{\iy}\frac{\xi^k}{k!}\tilde
P^{k}(L(\pi_k)\leq n)
&=&e^{-\xi}\det\left(\oint_{S^1}\frac{dz}{2\pi
iz}z^{j-\ell}e^{
\sqrt{\xi}(z+z^{-1})}\right)_{1\leq j,\ell\leq n}\\
\\
&=&e^{-\xi}\det\left(J_{j-\ell}(2\sqrt{-\xi})\right)_{1\leq j,\ell\leq n}\\
\\
&=&e^{-\xi}\int_{U(n)}e^{\sqrt{\xi}\Tr(X+X^{-1})}dX\\
\\
&=&\det\left(I-K(j,\ell)\Big\vert_{[n,n+1,\ldots]}\right)
\eean where $J_i(z)$ is the Bessel function and

\bea
K(k,\ell)&=&\frac{\sqrt{\xi}\left(J_k(2\sqrt{\xi})J_{\ell
+1}(2\sqrt{\xi})-J_{k+1}(2\sqrt{\xi})
J_{\ell}(2\sqrt{\xi})\right)}
 { k-\ell },~~~(k\neq \ell)\no\\
 &=&
  \sum_{n=1}^{\iy} J_{k+n}(2\sqrt{\xi})
                   J_{\ell+n}(2\sqrt{\xi}).
 \eea

\end{theorem}

\proof For an arbitrary partition $\lb \in \BY$, and
using formula (\ref{number standard tableaux1}), the
restriction of $P_{x,y}$ to the locus $\LR_1$, as in
(\ref{locus1}), reads as follows:
 \bean
P^{\xi}(\lb):=P_{x,y}(\lb)\Bigr\vert_{\LR_1}
&=&e^{-\sum_{k\geq 1} kt_k s_k}  \gs_{\lb} (t)
 \gs_{\lb} (s)\Bigr\vert_ {{t_i=\sqrt{\xi}
\dt_{i1}}
\atop{s_i=\sqrt{\xi} \dt_{i1} }} \\
&=&e^{-\xi} \xi^{\vert\lb\vert/2}
\frac{f^{\lb}}{\vert\lb\vert!} \xi^{\vert\lb\vert/2}
\frac{f^{\lb}}{\vert\lb\vert!}\\
&=&e^{-\xi}\frac{\xi^{\vert\lb\vert}}{\vert\lb\vert
!}\frac{(f^{\lb})^2}{\vert\lb\vert !}
 \\
 &=&e^{-\xi}\frac{\xi^{\vert\lb\vert}}{\vert\lb\vert
!}\tilde P^{(n)}(\lb),~~~\mbox{for $n=|\lb|$},
 \eean
where $\tilde P^{(n)}(\lb)$ can be recognized as {\em
Plancherel measure} on partitions in $\BY_n$, as defined
in section \ref{sect 1.5.1},
 $$ \tilde P^n(\lb)=\frac{(f^{\lb})^2}{n !}, \quad
\lb\in\BY_n
 .$$
 It is clear that 
%
%
 %
$$
P^{\xi}(\lb)=e^{-\xi}\frac{\xi^{|\lb|}}{|\lb|!}
\left(\frac{(f^{\lb})^2}{|\lb|!}\right),
\qquad\lb\in\BY,
$$
 is a genuine probability ($\geq 0$), called {\em Poissonized Plancherel measure}.

 We now compute
\bean
 P^{\xi}(\lb\mbox{~with~}\lb_1\leq
n)&=&e^{-\xi}\sum_{\lb\in\BY\atop{\lb_1\leq n}}
\frac{\xi^{|\lb|}}{|\lb|!}\frac{(f^{\lb})^2}{|\lb|!}\\
\\
&=&e^{-\xi}\sum_{k=0}^{\iy}\frac{\xi^k}{k!}\sum_{|\lb|=k\atop{\lb_1\leq
n}}
\frac{(f^{\lb})^2}{|\lb|!}\\
\\
&=&e^{-\xi}\sum_{k=0}^{\iy}\frac{\xi^k}{k!}\tilde
P^{(k)}(L(\pi_k)\leq n), \eean
and thus using Theorem \ref{Theo:2.1}\footnote{The
Bessel function $J_n(u)$ is defined by
 $$ e^{u(t-t^{-1})}=\sum^{\iy}_{-\iy}t^nJ_n(2u)
$$
and thus
$$ e^{-\sqrt{\xi}(z+z^{-1})}=e^{\sqrt{-\xi}(( iz)-( iz)^{-1})}=\sum
( iz)^nJ_n(2\sqrt{-\xi}).
$$},
 \bean
{P^{\xi}(\lb\mbox{~with~}\lb_1\leq n)}
  &=&
  Z^{-1}\det\left(\oint_{S^1}\frac{dz}{2\pi
iz}z^{k-\ell}e^{-\sum_1^{\iy}(t_iz^i +s_iz^{-i})
}\right)_{1\leq k,\ell\leq n}\Bigg\vert_{\LR_1}\\
\\
&=&e^{-\xi}\det\left(\oint_{S^1}\frac{dz}{2\pi
iz}z^{k-\ell}e^{-\sqrt{\xi}(z+z^{- 1})
}\right)_{1\leq k,\ell\leq n}\\
\\
&=&e^{-\xi}\int_{U(n)}e^{\sqrt{\xi}\Tr(X+X^{-1})}dX\\
&=&e^{-\xi}\det\left(J_{k-\ell}(2\sqrt{-\xi})\right)_{1\leq
k,\ell\leq n}, \eean where we used the fact that for a
Toeplitz matrix
 $$
\det(a^{k-\ell}c_{k-\ell})_{1\leq k,\ell\leq
n}=\det(c_{k-\ell})_{1\leq k,\ell\leq n},~~~a\neq 0.
 $$
It also equals
$$
P^{\xi}(\lb\mbox{~with~}\lb_1\leq
n)=\det\left(I-K(i,j)\Big\vert_{[n,n+1,...]}\right)
$$
where $K(i,j)$ is given by formula
(\ref{kernel-non-diagonal}), where
$V(z)=\sqrt{\xi}(z-z^{-1})$. Since
$$
\frac{z\frac{d}{dz}V(z)-w\frac{d}{dw}V(w)}{z-w}
 = \sqrt{\xi}\left(1-\frac{1}{wz}\right)
, $$
one checks
\bea
 (k-\ell)K(k,\ell)
  &=&\frac{\sqrt{\xi}}{(2\pi i)^2}
\renewcommand{\arraystretch}{0.5}
\begin{array}[t]{c}
\displaystyle{\oint\oint}\\
{\scriptstyle |z|=c_1}\\
{\scriptstyle |w|=c_2}
\end{array}
\renewcommand{\arraystretch}{1}dz~dw
\left(\frac{e^{\sqrt{\xi}(z-z^{-1})}}{z^{k+1 }}
\frac{e^{\sqrt{\xi}(w^{-1}-w)}}{w^{-\ell}}\right.\no\\
& &\hspace*{2cm}
 -\left.\frac{e^{\sqrt{\xi}(z-z^{-1})}}{z^{k+2}}
\frac{e^{\sqrt{\xi}(w^{-1}-w)}}{w^{-\ell +1}}\right)\no\\
& & \no\\
&=&\sqrt{\xi}\left(J_k(2\sqrt{\xi})J_{\ell+1}(2\sqrt{\xi})
 -J_{k+1}
 (2\sqrt{\xi})J_{\ell}(2\sqrt{\xi})\right)\no\\
& & \no\\
&=&(k-\ell)\sum^{\iy}_{n=1}J_{k+n}(2\sqrt{\xi})J_{\ell+n}(2\sqrt{\xi}).
\label{Bessel kernel}\eea
The last equality follows from the recurrence relation
of Bessel functions
$$
J_{k+1}(2z)=\frac{k}{z}J_k(2z)-J_{k-1}(2z).
 $$
 Indeed, subtracting the two expressions
  \bean
J_{\ell +1}(2z)J_k(2z)&=&
\frac{\ell}{z}J_k(2z)J_{\ell}(2z)- J_{\ell-1}(2z)J_k(2z)
\\
J_{k
+1}(2z)J_{\ell}(2z)&=&\frac{k}{z}J_k(2z)J_{\ell}(2z)-J_{k-1}(2z)J_{\ell}(2z).
,\eean
 one finds
\bean
(\ell -k)J_k(2z)J_{\ell}(2z)&=&z (J_k(2z)J_{\ell+1}(2z)-J_{k+1}(2z)J_{\ell}(2z) )\\
\\
& &-z(J_{k-1}(2z)  J_{\ell}(2z)-J_k(2z)J_{\ell -1}(2z) )
\eean implying
$$
z(J_k(2z)J_{\ell
+1}(2z)-J_{k+1}(2z)J_{\ell}(2z))=(k-\ell)\sum^{\iy}_{n=1}J_{k+n}(2z)
J_{\ell +n}(2z),
$$
thus proving (\ref{Bessel kernel}). \qed



\remark  ~Incidentally, the fact that $P^{\xi}$ is a
probability shows that Plancherel measure itself is a
probability; indeed, for all $\xi$,
$$
e^{-\xi}\sum_{n=0}^{\iy} \frac{\xi^n}{n!} = 1=
 \sum_{\lb\in\BY}P^{\xi}(\lb)=e^{-\xi}
 \sum^{\iy}_{n=0}\frac{\xi^n}{n!}\sum_{\vert\lb
\vert =n} \frac{(f^{\lb})^2}{\vert\lb\vert !}
 =
 e^{-\xi}
 \sum^{\iy}_{n=0}\frac{\xi^n}{n!}
 \sum_{\lb
\in \BY_n}\tilde P^n(\lb);
$$
comparing the extremities leads to $ \sum_{\lb \in
\BY_n}\tilde P^n(\lb)=1$.


\subsection{Probability on random words}

Here also, restrict the probability (\ref{generalized
probability1}) to the locus below, for $\xi>0$ and
$p\in\BZ_{\geq 1}$,
 \bean
  \LR_2 &=&\left\{\sum_{i\geq 1}x_i^k=\dt_{k1}\xi,\quad y_1=\ldots
 =y_p=\beta,\quad\mbox{all other~}y_j=0\right\} \\
 &=&\{t_k=\delta_{k1}\xi,~~ks_k=p\beta^k  \}
   \eean
Recall from section 1.5.4 the probability $\tilde
P^{k,p}$ on partitions, induced from uniform probability
$P^{k,p}$ on words. This was studied by Tracy-Widom
\cite{TW-words} and Borodin-Okounkov \cite{Bo-Ok}; see
also Johansson\cite{Johansson1} and Adler-van Moerbeke
\cite{AvM2}.

\begin{theorem}  For the set of words $S^p_n$, the generating
function for the distribution of the length of the
longest increasing subsequence, is given by (setting
$\beta=1$)
\begin{eqnarray*}
 e^{-p\xi} \sum^{\iy}_{k=0} \frac{(p\xi)^k}{k!} \tilde
P^{k,p} (L(\pi_k)\leq n)
&=&
e^{-p\xi}\det\left(\oint_{S'}\frac{dz}{2\pi
iz}z^{k-\ell}e^{
\xi z^{-1}}(1\!+\! z)^p\right)_{1\leq k,\ell\leq n}\\
&=&e^{-p\xi} \int_{U(n)} e^{\xi Tr\bar M} \det (I+M)^p dM\\
&=&\det \Bigl(I-K(j,k) \Bigr )_{_{(n, n+1,\ldots)}},
\end{eqnarray*}
with $K(j,k)$ a Christoffel-Darboux kernel of Charlier
polynomials:
\bea
\lefteqn{(j-k)K(j,k)}\no\\
&=&\frac{\xi}{(2\pi
i)^2}\renewcommand{\arraystretch}{0.5}
\begin{array}[t]{c}
\displaystyle{\oint\oint}\\
{\scriptstyle |z|=c_1<\frac{1}{|\xi|}}\\
{\scriptstyle |w|=c_2}
\end{array}
\renewcommand{\arraystretch}{1}dz~dw\left(p\frac{(1-\xi
z)^{-p-1}e^{-\xi z^{-1}}}{z^{j+1}}\frac{(1-\xi
w)^{p-1}e^{\xi w^{-1}}}
{w^{-k}}\right.\no\\
& &\hspace*{2cm}-\left.\frac{(1-\xi z)^{-p}e^{-\xi
z^{-1}}}{z^{j+2}}\frac{(1-\xi w)^{p}e^{\xi
w^{-1}}}{w^{-k+1}}
\right)\no\\
&=&(p)_{j+1}e^{-\xi^2}\left(\frac{_1F_1(-p,j+1;\xi^2)}{j!}~
\frac{_1F_1(-p+1,k+2;\xi^2)}{(k+1)!}\right.\no\\
& & \no\\
&
&\hspace*{2cm}-\left.\frac{_1F_1(-q+1,j+2;\xi^2)}{(j+1)!}~\frac{_1F_1(-q,k+1;\xi^
2)}{k!} \right)\label{3.1.4} \eea
where
$$
(a)_j:=a(a+1)...(a+j-1),\qquad (a)_0=1
$$
and $_1F_1(a,c;x)$ is the {\em confluent hypergeometric
function}:
\bean \frac{1}{2\pi
i}\oint_{|z|=c_1<\frac{1}{|\xi|}}(1-\xi z)^{-p}e^{\xi
z^{-1}}
\frac{dz}{z^{m+1}}&=&\frac{(p)_m}{m!}e^{-\xi^2}\, _1F_1(1-p,m+1;\xi^2)\\
 & & \\
\frac{1}{2\pi i}\oint_{|w|=c_2}(1-\xi w)^qe^{\xi w^{-1}}
w^{m-1}dw &=&\frac{1}{m!} \, _1F_1(-q,m+1;\xi^2). \eean
These functions are related to
 Charlier polynomials\label{footnote11}\footnote{\label{footnote11} Charlier polynomials $P(k;\al)
 $, with $k\in \BZ_{\geq 0}$, are
 discrete orthonormal polynomials defined by
 the orthonormality condition
  $$\sum_{k=0}^\iy P_n(k;\al)P_m(k;\al)w_\al(k)=\dt_{nm},~~
  \mbox{for}~~
   w_\al(k)=e^{-\al}\frac{\al^k}{k!},$$
  with generating function
   $$
   \sum_{n=0}^{\iy}\al^{n/2} \frac{1}{\sqrt{n!}}
    P_n(k;\al)w^n=e^{-\al w}(1+w)^k.$$
   }
   .

\end{theorem}

\proof This proof will be very sketchy, as more details
will be given for the percolation case in the next
section.
%
%
%
%
One now computes
\begin{eqnarray*}
 {P^{\xi,p} (\mbox{$\lb$ with $\lb_1\leq n$})}
&=&
   \sum_{{\lb}\atop{\lb_1\leq n}}
   e^{-\sum_{i\geq 1} it_is_i}\gs_{\lb} (t)
    \gs_{\lb}(s)
     \Bigr\vert_{\LR_2}\\
&=&
 e^{-p\xi \beta}\sum_{\lb   \atop{
\lb_1\leq n}}
 {\gs_{\lb}(\xi,0,\ldots)\gs_{\lb}
 (p\beta,\frac{p\beta^2}{2},\frac{p\beta^3}{3},
 \ldots )} \\
 &=&
 e^{-p\xi \beta}\sum_{\lb \in\mathbb Y^p \atop{
\lb_1\leq n}}
 \frac{(p\xi\beta)^{\vert\lb\vert}}{\vert\lb\vert!}
\frac{\gs_{\lb}(1,0,\ldots)\gs_{\lb}
 (p,\frac{p}{2},\frac{p}{3},
 \ldots )}{p^{\vert\lb\vert}}\\
&=&e^{-p\xi\beta} \sum_{k=0}^{\iy}
\frac{(p\xi\beta)^k}{k!} \sum_{ \lb\in \mathbb
Y^p\atop{\vert\lb\vert=k\atop{\lb_1\leq n}}
}\frac{f^{\lb} s_{\lb} (1^p)}{p^{|\lb|}}\\
&=& e^{-p\xi\beta}\sum^{\iy}_{k=0}
\frac{(p\xi\beta)^k}{k!} \tilde P^{k,p} (L(\pi_k)\leq n)
.\end{eqnarray*}
In applying Theorem \ref{Theo:2.1}, one needs to compute
%
$$
e^{-\sum_1^{\iy} (t_i z^i+s_i z^{-i}) }\Big\vert_{\LR_2}
= e^{-\xi z} e^{-p \sum_1^{\iy}\frac{1}{i}
\left( \frac{\beta}{z}\right)^i}
  = e^{-\xi z} \left(1-\frac{\beta}{z}\right)^p
$$
and
$$
\left. e^V\right|_{\LR_2}=e^{-\sum_1^{\iy}(t_i
z^{-i}-s_i z^i) }\Big\vert_{\LR_2} = e^{-\xi z^{-1}}
(1-\beta z)^{-p} .$$
%
Therefore
\begin{eqnarray*}
P^{\xi,p}(\lb_1\leq n) &=&Z^{-1} \left.\det \left(
\oint_{S^1} \frac{dz}{2\pi iz} z^{k-\ell}
e^{-\sum_1^{\iy} (t_iz^i+s_i z^{-i})}
\right)_{1\leq k,\ell\leq n}\right\vert_{\LR_2}\\
&=&e^{-\xi p\beta}\det\left(\oint_{S^1}\frac{dz}{2\pi
iz}z^{k-\ell}e^{-\xi z}(1-\beta z^{-1})^p\right)_{1\leq k,\ell\leq n}\\
&=&e^{-\xi p\beta}\det\left(\oint_{S^1}\frac{dz}{2\pi
iz}z^{k-\ell}e^{\xi z^{-1}}(1+\beta z)^p\right)_{1\leq k,\ell\leq n}\\
&&\hspace{4cm} \mbox{the change of variable
$z\longmapsto
-z^{-1}$}\\
&=&e^{-\xi p\beta} \int_{U(n)} e^{\xi Tr\bar M} \det
(1+\beta M)^p dM.
\end{eqnarray*}
Then, one computes

$$
\frac{z\frac{d}{dz}V(z)-w\frac{d}{dw}V(w)}{z-w}=\frac{\beta
p}{(1-\beta z)(1-\beta w)}-\frac{\xi}{zw}.
$$
and
$$
e^{V(z)-V(w)}=e^{-\xi z^{-1}}(1-\beta z)^{-p}
    e^{ \xi w^{-1}}(1-\beta w)^{ p},
    $$
leading after an appropriate rescaling, using $\beta$ to
formula (\ref{3.1.4}). From footnote \ref{footnote11},
the confluent hypergeometric functions turn out to be
Charlier polynomials. \qed

\subsection{Percolation}\label{section 3.3}

Considering now the locus
$$\LR_3:=\left\{ kt_k=q\xi^{k/2}
   ,~~
   ks_k=p\xi^{k/2}  \right\}
   ,$$
   one is led to the probability appearing in the
   generalized permutations and percolations (sections
   1.5.4 and 1.5.5), namely
   $$
   P(L(M)\leq \ell)
   = \sum_{\tiny \begin{array}{c}
   \lb   \\
  \lb_1\leq \ell \end{array}}
(1-\xi)^{pq}\xi^{\vert\lb\vert}
 \gs_{\lb}(q,\frac{q}{2},\ldots)
 \gs_{\lb}(p,\frac{p}{2},\ldots).
 $$
We now state: (see \cite{Bo-Ok, Johansson0})
\begin{theorem}\label{Theo:3.3}
Assuming $q>p$, we have \bean
 \lefteqn{P(L(M)\leq \ell)}\\
 \\&=&
 \!\!\frac{(1-\xi)^{pq}}{{\ell}!} \int_{(S^1)^{\ell}}|\Dt_{\ell}(z)|^2
 \prod_{j=1}^{\ell}
  (1+\!\sqrt{\xi}z_j)^q  (1+\!\sqrt{\xi}\bar z_j)^p
~ \frac{ dz_j }{2\pi i z_j}
  \\ \\ \\
  &=&
     \displaystyle{ (1\!-\!\xi)^{pq} \!\!\int_{U(\ell)} }
    \det (1\! +\!  \sqrt{\xi}M)^q
                 \det (1\!+\!  \sqrt{\xi}\bar M)^p dM
\\
&=&
  \det\left(I-K(i,j)\right)\Bigg\vert_{[\ell,\ell+1,...]}
\eean
where ($C$ is a constant depending on $p,q,\xi$)
$$
K(i,j)= C \frac{ {M_p(p+i )M_{p-1}(p+j)-M_{p-1}(p+i)M_p(
p+i)}}{i-j} ,$$
 where
 $$M_p(k):=M_p(k;q-p+1,\xi)
  \mbox{~~with~~} k\in \BZ_{\geq 0}.$$ are
  Meixner polynomials\footnote{Meixner polynomials are discrete
  orthogonal polynomials defined on $\BZ_{\geq 0}$ by
  the orthogonality relations
  $$ \sum_{k=0}^{\iy}
  \left({\beta+k-1}\atop {k}\right) \xi^k
   M_p(k;\beta,\xi)M_{p'}(k;\beta,\xi)
   =\frac{\xi^{-p} }{ (1-\xi)^{\beta}}
\left({\beta+p-1}\atop {p}\right)^{-1}\dt_{pp'},~~~\beta
>0 ,$$ and with generating function \be
\left(1-\frac{z}{\xi}\right)^x(1-z)^{-x-\beta }=
 \sum_{p=0}^\iy\frac{(\beta )_p
 z^p}{p!}M_p(x;\beta ,\xi).
 \label{genfunction}\ee}.
\end{theorem}

\proof Using the restriction to locus $\LR_3$, one finds
 $$
 e^{-\sum_1^{\iy}
(t_jz^j+s_jz^{-j})}\Bigr|_
 {\LR_3}=(1-\sqrt{\xi}z)^q
    (1-\sqrt{\xi}z^{-1})^p
    $$
and
    $$
\left. e^V\right|_{\LR_3}=e^{-\sum_1^{\iy}(t_i
z^{-i}-s_i z^i) }\Big\vert_{\LR_3}
=(1-\sqrt{\xi}z^{-1})^q
     (1-\sqrt{\xi}z)^{-p},$$
and thus
$$
V(z)=q\log(1-\sqrt{\xi}~z^{-1})-p\log(1-\sqrt{\xi}~z).
$$
Hence
$$
\frac{z\frac{d}{dz}V(z)-w\frac{d}{dw}V(w)}{z-w}=\frac{p\sqrt{\xi}}{(1-\sqrt{\xi}~z)(1
-\sqrt{\xi}~w)}-
\frac{q\sqrt{\xi}}{zw(1-\sqrt{\xi}~z^{-1})(1-\sqrt{\xi}~w^{-1})}
$$
and so
\bean
 \lefteqn{\frac{1}{z^{k+1}
w^{-\ell}}\frac{z\frac{d}{dz}V(z)-w\frac{d}{dw}V(w)}{z-w}
 ~e^{V(z)-V(w)}}\\
 &=&
 \sqrt{\xi}\left(p\frac{(1-\sqrt{\xi}~z)^{-p-1}(1-\sqrt{\xi}~z^{-1})^q}{z^{k+1}}\frac{(1-\sqrt{
\xi}~w)^{p-1}(1-\sqrt{\xi}~
w^{-1})^{-q}}{w^{-\ell}}\right.\\
\\
& &  \left.
-q\frac{(1-\sqrt{\xi}~z)^{-p}(1-\sqrt{\xi}~z^{-1})^{q-1}}{z^{k+2}}\frac{(1-\sqrt{\xi}
~w)^p(1-\sqrt{\xi}~ w^{-1})^{-q-1}}{w^{-\ell
+1}}\right)\eean
Expanding in Laurent series, one finds
\bean
\lefteqn{(1-\eta z)^{-\beta}(1-\eta z^{-1})^{\beta'}}\\
\\
&=&\left(\sum^{\iy}_{i=0}\frac{(-\beta)(-\beta
-1)\ldots(-\beta -i+1)}{i!}(-\eta)^iz^i\right)\\
\\
& & \qquad\left(\sum^{\iy}_{j=0}\frac{\beta'(\beta'-1)
\ldots(\beta' -j+1)}{j!}(-\eta)^j\frac{1}{z^j}\right)\\
\\
&=&\sum_{m\in\BZ}z^m\sum_{i-j=m\atop{i,j\geq
0}}\frac{\beta(\beta +1)\ldots(\beta
+i-1)}{i!}\frac{(-\beta')(-\beta'+1)\ldots(-\beta'+j-1)}{j!}\eta^{i+j}\\
\\
&=& \sum_{m\in\BZ}z^m\sum_{j=0}^{\iy}\frac{\beta(\beta
+1)\ldots(\beta +j+m-1)}{(j+m)!}\frac{(-\beta')
\ldots(-\beta'+j-1)}{j!}\eta^{2j+m}. \eean
The coefficient of $z^m$ for $m\geq 0$ reads \bean
\lefteqn{ \sum_{j=0}^{\iy}\frac{\beta(\beta
+1)\ldots(\beta +j+m-1)(-\beta')
\ldots(-\beta'+j-1)}{(j+m)! }\frac{\eta^{2j+m}}{j!}
 ~~~~~~~~~~~~~~~~~~~~~~~~~}\\
\\& &~~~~~~~~~~= \frac{(\beta)_m}{m!}\eta^m ~_2F_1(\beta
+m,-\beta';m+1;\eta^2),
~~~~~~~~~~~~~~~~~~~~~~~~~~~\eean 
which is Gauss' Hypergeometric
function\footnote{$~_2F_1(a,b;c;z):=\sum_0^{\iy}
 \frac{(a)_n(b)_n}{(c)_n}\frac{z^n}{n!}$. Notice that when $a=-m<0$ is an
 integer, then $~_2F_1(a,b;c;z)$ is a polynomial of degree $m$.} and
thus \bean \oint_{z=0}(1-\eta z)^{-\beta}(1-\eta
z^{-1})^{\beta'}\frac{dz}{2\pi iz^{k+1}}&=&
\frac{(\beta)_k}{k!}\eta^k~_2F_1(-\beta',\beta
+k;k+1;\eta^2).
 \eean
One does a similar computation for the other piece in
the kernel, where one computes the coefficient of
$z^{-m}$. Using this fact, using a standard linear
transformation formula\footnote{In general one has
$~_2F_1(a,b;c;z)=(1-z)^{-b}
 ~_2F_1(b,c-a;c;\frac{z}{z-1})$.} for the
hypergeometric function, and using a polynomial property
for hypergeometric functions\footnote{When $p$ is a
positive integer ($<q$) and $k$ an integer, then
$~_2F_1(-p,-q;k+1;z)$ is a polynomial of degree $p$ in
$z$, which satisfies $~_2F_1(-p,-q;k+1;z)
 =\frac{\Gamma(k+1)\Gamma(1-q+p)}{(k+p)!\Gamma(1-q)}
(-z)^p ~_2F_1(-p,-k-p;q-p+1;\frac{1}{z})$.} and the
formula for $0<p<q$
$$
\frac{\Gamma(1+p-q)}{\Gamma(1-q)}=(-1)^{p+1}
\frac{(q-1)!}{(q-p-1)!},$$
one finds Meixner polynomials (assume $q> p> 0,~x\geq 0$
and $k\geq 0$) (see Nifikorov-Suslov-Uvarov \cite{NSU}
and Koekoek-Swartouw \cite{Koekoek}):
\bean \lefteqn{  \frac{1}{2\pi
i}\oint_{|\sqrt{\xi}|<|z|=c_1 }
 \frac{(1-\sqrt{\xi}~z)^{-p-1}(1-\sqrt{\xi}~z^{-1})^q}
{z^{k+1}}dz
 }\\
 &=&\frac{(p+1)_k}{k!}
 {\xi^{k/2}}~_2F_1(p+1+k,-q;k+1;\xi)
 \\
 &=&\frac{(p+1)_k}{k!}
 {\xi^{k/2}}{(1-\xi)^q}
~_2F_1\left(-p,-q;k+1;\frac{\xi}{\xi-1}\right)
 \\
 &=&
 \frac{(p+1)_k \Gamma(1+p-q)}{(k+p)!\Gamma(1-q)}
\frac{\xi^{p+k/2}}{(1-\xi)^{p-q}}
 ~_2F_1(-p,-k-p; q-p+1;\frac{\xi -1}{\xi})
\\
&=&\frac{ \Gamma(1+p-q)}{p!\Gamma(1-q)}
\frac{\xi^{p+k/2}}{(1-\xi)^{p-q}}M_p(p+k;q-p+1;\xi)
 \\
 &=&
  \left({q-1}\atop {p}\right)
   \frac{\xi^{p+k/2}}{(1-\xi)^{p-q}}M_p(p+k;q-p+1;\xi)
   ,\eean
where
\bean
 M_p(x;\beta,\xi)
 =~_2F_1(-p,-x;\beta;\frac{\xi-1}{\xi})
 &=& \left(\frac{1-\xi}{\xi}\right)^p
 \frac{1}{(\beta)_p}x^p+\ldots
 ,~~x\in
 \BZ_{\geq 0}
  \\
  &=:&a_px^p+\ldots\mbox{,~~with }~\beta=q-p+1
 \eean
 are Meixner polynomials in $x$, satisfying the following
orthogonality properties:
 $$ \sum_{x=0}^{\iy}
  w(x)
   M_p(x;\beta,\xi)M_m(x;\beta,\xi)
   =\frac{p!}{\xi^p(1-\xi)^\beta (\beta)_p}
   \dt_{pm}=:h_p\dt_{pm}
,$$
with weight
$$
w(x):=\frac{(\beta)_x}{x!}
\xi^x=\frac{(1+q-p)_x}{x!}\xi^x .$$
  Using these facts, one
computes
\bean \lefteqn{(k-\ell)K(k,\ell)=\frac{\sqrt{\xi}}{(2\pi
i)^2}\int\!\!\int_{|\sqrt{\xi}|<|z|=c_1\atop{|w|=c_2 <
|\frac{1}{\sqrt{\xi}}|}}dz~dw}\\
\\
&
&\left(p\frac{(1-\sqrt{\xi}~z)^{-p-1}(1-\sqrt{\xi}~z^{-1})^q}{z^{k+1}}\frac{(1-\sqrt{
\xi}~w)^{p-1}(1-\sqrt{\xi}~
w^{-1})^{-q}}{w^{-\ell}}\right.\\
\\
& &  \left.
-q\frac{(1-\sqrt{\xi}~z)^{-p}(1-\sqrt{\xi}~z^{-1})^{q-1}}{z^{k+2}}\frac{(1-\sqrt{\xi}
~w)^p(1-\sqrt{\xi}~
w^{-1})^{-q-1}}{w^{-\ell +1}}\right)\\
\\
&=&\frac{a_{p-1}}{a_p h_{p-1}}\sqrt{w(p+k)w(p+\ell)}
\left( {M_p(p+k )M_{p-1}(p+\ell)-M_{p-1}(p+k)M_p(
p+\ell)}\right) , \eean
where $$M_p(p+k):=M_p(p+k;q-p+1,\xi).$$
 We also have, using Propositions 2.2 and 2.10 restricted to the
 locus $\LR_3$,
%
%
%
%
\bean\lefteqn{\BP(L(M)\leq \ell)}\\ \\
 &=&
  \sum_{\tiny \begin{array}{c}
   \lb   \\
  \lb_1\leq \ell \end{array}}
(1-\xi)^{pq}\xi^{\vert\lb\vert}
 \gs_{\lb}(q,\frac{q}{2},\ldots)
 \gs_{\lb}(p,\frac{p}{2},\ldots) \\ \\
&=&  e^{-\sum_1^{\iy} kt_ks_k}
 \sum_{\tiny \begin{array}{c}
   \lb   \\
  \lb_1\leq \ell \end{array}}\gs_{\lb}(t_1,t_2,\ldots)
 \gs_{\lb}(s_1,s_2,\ldots)\Bigr|_
 {\begin{array}{l}kt_k=q\xi^{k/2}
   \\
   ks_k=p\xi^{k/2}\end{array}}
   \\ \\
   & {=}&\left.
 (1\!-\!\xi)^{pq}
 \det\!\!\left(\oint_{S^1}
  z^{\al-\al'} e^{-\sum_1^{\iy}
(t_jz^j+s_jz^{-j})}
 \frac{dz}{2\pi iz}
 \right)_{1\leq \al,\al'\leq \ell}
  \right|_
 {\begin{array}{l}kt_k = q\xi^{k/2}
   \\
   ks_k\!=\!p\xi^{k/2}\end{array} 
                           }
                             \\
&  {=}
 &(1-\xi)^{pq}
 \det \!\left(\!\!
  \oint _{S^1} \!\!
   z^{\al-\al'}
   (1-\!\sqrt{\xi}z)^q
   (1-\!\sqrt{\xi}z^{-1})^p 
   \frac{dz}{2\pi iz}\right)
 _{1 \leq \al,\al'\leq \ell} 
 \\&=&
 \!\!\frac{(1-\xi)^{pq}}{{k}!} \int_{(S^1)^{\ell}}|\Dt_{\ell}(z)|^2
 \prod_{j=1}^{k}
  (1+\!\sqrt{\xi}z_j)^q  (1+\!\sqrt{\xi}\bar z_j)^p
~ \frac{ dz_j }{2\pi i z_j}
  \\ \\
  &=&
      (1\!-\!\xi)^{pq}  \int_{U(\ell)}
    \det (1  +  \sqrt{\xi}M)^q
                 \det (1 +   \sqrt{\xi}\bar M)^p dM
 ,\eean
thus ending the proof of Theorem 3.3.\qed

\newpage

\section{Limit theorems}

\subsection{Limit for Plancherel measure}












Stanislaw Ulam \cite{Ulam} raised in 1961 the question:
how do you compute the probability $P^{(n)} (L  (\pi_n)
\leq k)$ that the length $L:=L^{(n)}$ of the longest
increasing sequence in a random permutation is smaller
than $k$. What happens for very large permutations,
i.e., when $n\rightarrow \iy$?
By Monte-Carlo simulations, Ulam conjectured that
\be c:=\lim_{n\rightarrow \infty} \frac{E^{(n)} (L
)}{\sqrt{n}}
 \label{Ulamlimit}\ee
 exists $(E^{(n)}$ denote the
expectation with respect to $P^{(n)}$). A much older
argument of Erd\"os and Szekeres \cite{Erdos and
Szekeres} implied that $E^{(n)}(L ) \geq \frac{1}{2}
\sqrt{n-1}$ and so $c\geq\frac{1}{2}$. Numerical
computation by Baer and Brock \cite{Baer} suggested
$c=2$. Hammersley \cite{Hammersley} showed the existence
of the limit (\ref{Ulamlimit}); in 1977, Logan and Shepp
\cite{Logan and Shepp} proved $c\geq 2$ and, at the same
time, Vershik and Kerov \cite{Vershik and Kerov2} showed
$c=2$. More recently other proofs have appeared by
Aldous and Diaconis \cite{AD-Ulam}, Sepp\"al\"ainen
\cite{Seppalainen} and Johansson \cite{Johansson-Ulam}.
Meanwhile, Gessel \cite{Gessel} found a generating
function for the probability (with respect to $n$) and
connected this problem with Toeplitz determinants.
Monte-Carlo simulations by Odlyzko and Rains \cite{Od}
suggested
$$
c_0\sim \lim_{n\rightarrow \infty} \frac{\mbox{Var}
L^{(n)}}{n^{1/3}} \quad c_1=\lim \frac{E^{(n)}
(L^{(n)})-2\sqrt{n}}{n^{1/6}},
$$
with $c_0 \sim 0.819$ and $c_1\sim -1.758$.

In this section, we explain some of the ideas underlying
this problem. It is convenient to write a partition
$\lb\in\BY$ in $(u,v)$-coordinates, as shown in the
picture below.


\newpage
 $\hspace*{4.2cm}u_1$

 \vspace*{1.9cm}

 $\hspace*{11cm} y$

 $\hspace*{1.6cm}u_2\vspace*{-1.4cm}\hspace*{2cm}
  (x,y)\hspace*{2.7cm} (x,F_{\tilde \lb}(x))=(u_1,\psi(u_1))$


$$
  x~~~\squaresize .2cm \thickness .02cm \Thickness
.08cm  \Young{
      &&&&&&&&&&&&&&&&&&&&&&&&&&&&&&&&&&&&&\cr
      &&&&&&&&&&&&&&&&&&&&&&&&&&&&&\cr
      &&&&&&&&&&&&&&&&&&&&&&\cr
      &&&&&&&&-&-&-&-&-&-&-&-&-&-&-&-&-&-\cr
      &&&&&&&&|&&&&&&&&&&\cr
      &&&&&&&&|&&&&&&&&\cr
      &&&&&&&&|&&&&&&\cr
      &&&&&&&&|&&&&&&\cr
      &&&&&&&&|&&&&&& \cr
      &&&&&&&&|&&&&&&\cr
      &&&&&&&&|&&&&&\cr
      &&&&&&&&|&&&&\cr
      &&&&&&&&|&& \cr
      &&&&&&&&|\cr
      &&&&&&\cr
      &&&&&\cr
      &&&\cr
      &&&\cr
      &&\cr
      &&\cr
      &&\cr
      &&\cr
      &\cr
      &\cr\cr\cr
  }
$$


\vspace*{-.8cm}\hspace{-2.46cm}
\begin{picture}(12,1.6) 
\put(70,70){\line(1, 1){153 }}
 \put(198,198){\line(1, -1){71}}
\put(146,146){\line(1, -1){120}}
 \put(129,129){\line(1,-1){63}}
\end{picture}\hspace*{9cm} {\it v}

\vspace*{-2cm}\hspace*{4cm} $(F_{\tilde\lb}^{-1}(y),y)=
 (u_2,\psi(u_2))$

\vspace*{2cm}

 $$\mbox{Figure 4.1}$$


Remember from (\ref{hook}) the hook length $h^{\lb}$,
defined as
$$
h^{\lb}= \prod_{(i,j)\in\lb}
 h^{\lb}_{ij},~~~\mbox{with}~~h^{\lb}_{ij}=\mbox{~hook
length}=\lb_i+\lb_j^{\top}-i-j-1, $$
 where $\lb^{\top}_j$
is the length of the $j^{\rm th}$ column. Also remember
from (\ref{number standard tableaux}) the formula
$f^{\lb}$ expressed in terms of $h^{\lb}$. Then one has
the following Theorem, due to Vershik and Kerov
\cite{Vershik and Kerov1}; see also Logan and Shepp
\cite{Logan and Shepp}. A sketchy outline of the proof
will be given here.

Consider Plancherel measure (see section
\ref{section3.1}), which using formula (\ref{number
standard tableaux}), can be written in terms of the hook
length,
$$
\tilde
P^{(n)}(\lb)=\frac{(f^{\lb})^2}{n!}=\frac{n!}{(h^{\lb})^2}\mbox{~for~}|\lb
|=n ,$$ and define the function \be
\Omega(u):=\left\{\begin{array}{ll}\frac{2}{\pi}(u\mbox{{\rm
~arcsin}}\frac{u}{2}+\sqrt{4-u^2}),&\mbox{for~}
|u|\leq 2\\
\\
|u|,&\mbox{for~}|u|\geq 2.
\end{array}\right.
\label{Omega-function}\ee

\begin{theorem} {\rm (Vershik--Kerov-Logan-Shepp)}\cite{Vershik and Kerov1,Logan and Shepp}
Upon expressing $\lb$ in $(u,v)$-coordinates, define the
subset of partitions, for any $\vr >0$,
$$
\BY_n(\vr):=\left\{\lb\in\BY\mid\sup_{u}\left|
 \frac{1}{\sqrt{n}}\lb(u\sqrt{n})-\Omega(u)\right|
<\vr\right\}.
$$
Then, for large $n$, Plancherel measure concentrates on
a Young diagram whose boundary has shape $\Omega(u)$;
i.e.,
$$
\lim_{n\rightarrow \iy} \tilde P^{(n)}(\BY_n(\vr))=1.
$$
Moreover, for the length $\lb_1$ of the first row, one
has
 $$
 \lim_{n\rightarrow \iy} \tilde P^{(n)}\left( \left|\frac{\lb_1}{2\sqrt{n}}-1\right|<\vr
   \right)=1,
   $$ and thus for uniform measure on permutations
   (remembering
   $L(\pi_n)=$ the length  of the longest increasing
   sequence in $\pi_n$), one has
   $$
 \lim_{n\rightarrow \iy}  P^{(n)}\left( \left|\frac{L(\pi_n)}{2\sqrt{n}}-1\right|<\vr
   \right)=1.
   $$

\end{theorem}



\medskip\noindent{\it Brief outline of the Proof:\/}
In a first step, the following expression will be
estimated, using Stirling's formula:\footnote{$\log
n!=n\log n+\log\sqrt{2\pi n}-n+\ldots$ for
$n\nearrow\iy$.}
\bean
\lefteqn{-\frac{1}{n}\log P_n(\lb)}\\
\\
&=&-\frac{1}{n}\log\frac{n!}{(h^{\lb})^2}\\
\\
&=&\frac{2}{n}\log\prod_{(i,j)\in\lb}h^{\lb}_{ij}-\frac{1}{n}\log n!\\
\\
&=&\frac{2}{n}\log\prod_{(i,j)\in\lb}h^{\lb}_{ij} -\log
n -\frac{1}{n}\log \sqrt{2\pi n}+1+{\bf O}(1),~\mbox{by
Stirling}
 \\
&=&1+\frac{2}{n}(\log\prod_{(i,j)\in\lb}h^{\lb}_{ij}-n\log
n^{1/2})-\frac{1}{n}\log \sqrt{2\pi n}+{\bf O}(1)
,~\mbox{using $|\lb|=n$,}
\\
\\
&=&1+2 \sum_{(i,j)\in\lb}
 \frac{1}{n}\log\frac{h_{ij}^{\lb}}{\sqrt{n}}+{\bf
O}\left(\frac{\log\sqrt{n}}{n}\right)\mbox{(Riemann
sum)}
 \eean\bean
&\longrightarrow &1+2\int\!\!\!\int_{_{\{(x,y),x,y\geq
0,y\leq F_{\tilde\lb}(x)\}}}dx~dy
\log(F_{\tilde\lb}(x)-y+F^{-1}_{\tilde\lb}(y)-x), \eean
 assuming the partition $\tilde\lb
=\frac{\lb}{\sqrt{n}}$ tends to a continuous curve
$y=F_{\tilde\lb}(x)$ in $(x,y)$ coordinates. Then the
Riemann sum above tends to the expression above for
$n\rg\iy$; note $1/n$ is the area of a square in the
Young diagram, after rescaling by $1/\sqrt{n}$ and thus
turns into $dxdy$ in the limit. Clearly, the expression
$F_{\tilde\lb}(x)-y+F^{-1}_{\tilde\lb}(y)-x$ is the hook
length of the continuous curve, with respect to the
point $(x,y)$. In the $(u,v)$-coordinates, this hook
length is particularly simple:
$$
 F_{\tilde\lb}(x)-y+F^{-1}_{\tilde\lb}(y)-x=
 \sqrt{2}(u_2-u_1),~~\mbox{with }u_1\leq u_2,
 $$
 where $(u_1,\psi(u_1))=(x,F_{\tilde\lb}(x))$ and $(u_2,\psi(u_2))
 =(F_{\tilde\lb}^{-1}(y),y)$ and where $v=\psi(u)$
 denotes the curve $y=F_{\tilde\lb}(x)$ in
 $(u,v)$-coordinates.

 Keeping in mind Figure 4.1, consider two points $(u_1,\psi(u_1))$ and
 $(u_2,\psi(u_2))$ on the curve, with $u_1<u_2$. The
 point $(x,y)$ such that the hook (emanating
  from $(x,y)$ parallel to the $(x,y)$-coordinates) intersects the curve
 at the two points $(u_1,\psi(u_1))$ and
 $(u_2,\psi(u_2))$, is given by a point on the line
 emanating from $(u_1,\psi(u_1))$ in the $(1,-1)$-direction;
 to be precise, the point
 $$(x,y)=(u_1,\psi(u_1))+\frac{1}{2}
 (u_2-\psi(u_2)-u_1+\psi(u_1))(1,-1).$$


 So the surface element $dx dy$ is transformed into the
 surface element $du_1 du_2$ by means of the Jacobian,
 $
 dx dy = \frac{1}{2} (1+\psi'(u_1))(1-\psi'(u_2)
 du_1du_2,$
and thus, further replacing $u_1\mapsto \sqrt{2}u_1$ and
 $u_2\mapsto \sqrt{2}u_2$, one is led to
\bean
 \lefteqn{-\frac{1}{n}\log P_n(\lb)}\\
 &\simeq&  1+2\int\!\!\int_{\{(x,y),x,y\geq 0,y\leq
F(x)\}}dx~dy\log(F(x)+F^+(y)-x-y) \\
\\
&=&1+\int\!\!\!\int_{u_1<u_2}du_1du_2(1+\psi'(u_1))(1-\psi'(u_2))
\log( \sqrt{2}(u_2-u_1))
\\
\\
&=&-\frac{1}{2}\int\!\!\!\int_{\BR^2}\log(\sqrt{2}|u_2-u_1|)
f'(u_1)f'(u_2)du_1du_2+2\int_{|u|>2}f(u)
\mbox{~arccosh}|\frac{u}{2}|du\\
\\
&=&\frac{1}{2}\int\!\!\!\int_{\BR^2}\left(\frac{f(u_1)
-f(u_2)}{u_1-u_2}\right)^2du_1du_2+2\int_{|u|>2}f(u)
\mbox{~arccosh}|\frac{u}{2}|du \eean upon setting
$\psi(u)=\Omega(u)+f(u)$, where the function $\Omega(u)$
is defined in (\ref{Omega-function}). The last identity
is obtained by using Plancherel's formula of Fourier
analysis,
$$ \int_{\BR}g_1(v)g_2(v)dv=\frac{1}{2\pi}
 \int_{\BR}\widehat g_1(v)\widehat g_2(v)dv$$
applied to the two functions $$g_1(v):=\int_{\BR
}du_1\log(\sqrt{2}|v-u_1|)
f'(u_1)=\int_{\BR}\frac{f(u_1)}{\sqrt{2}|v-u_1|}du_1$$
and $g_2(v):=f'(v)$.

This shows the expression $-\frac{1}{n}\log P_n(\lb)$
above is minimal (and $=0$), when $f(u)=0$; i.e. when
the curve $\psi(u)=\Omega(u)$ and otherwise the
integrals above $>0$. So, when the integral equals $\vr
>0$, the expression $-\frac{1}{n}\log
P_n(\lb)\simeq\vr>0$ and thus $P_n(\lb)\simeq e^{-\vr
n}$, which tends to $0$ for $n\nearrow \iy$. Only, when
$\vr=0$, is there a chance that $P_n(\lb)=1$; this
happens only when $\psi=\Omega$.\qed

 \subsection{Limit theorem for longest increasing sequences}
 \label{subsect 4.2}


In section \ref{section3.1}, it was shown that a
generating function for the probability of the length
 $L(\pi_k)$ of
the longest increasing sequence in a random permutation
is given in terms of a Bessel kernel:
  \bean
e^{-\xi}\sum_{k=0}^{\iy}\frac{\xi^k}{k!}
P^{k}(L(\pi_k)\leq n)
&=&\det
\bigl(I-K(j,\ell)\bigr)\Big\vert_{[n,n+1,\ldots]} ,\eean
where
\be
K(k,\ell)=
  \sum_{m=1}^{\iy} J_{k+m}(2\sqrt{\xi})
                   J_{\ell+m}(2\sqrt{\xi})
 .\label{BesselKernel}\ee
 In the statement below $A(x)$ is the classical Airy function
$$
 A(x)=\frac{1}{\pi}\int_0^{\iy}\cos\left(\frac{u^3}{3}+xu\right)du
;$$ this function is well known to satisfy the
 ordinary
differential equation $ A'' (x)=xA(x)$ and to behave
asymptotically as
 $$
A(x):=\frac{e^{-\frac{2}{3}x^{3/2}}}{2\sqrt{\pi}x^{1/4}}
(1+O(x^{3/2})),~~~\mbox{as $x\rg \iy$.}$$


\begin{theorem}\label{BDJ}(J. Baik, P. Deift and K.
Johansson \cite{BDJ}) The distribution of the length
$L(\pi_n)$ of the longest increasing sequence in a
random permutation behaves as
$$  \lim_{n\rightarrow
\iy}
    P^{(n)}\left(L(\pi_n)\leq
2 {n}^{1/2}+xn^{1/6}\right)   =\det \left(I-{\bf A}
\chi_{[x,\iy)}\right),
  $$
  where
 \be
{\bf A}(x,y):=\int_0^{\iy}
 du A(x+u)A(y+u)
 = \frac{A(x)A'(y)
           -A^{\prime}(x)A(y)}{x-y}
           .\label{AiryK}\ee

\end{theorem}

\noindent The proof of this Theorem, presented here, is
due to A. Borodin, A. Okounkov and G. Olshanski
\cite{BOO}; see also \cite{BDJ1}. Before giving this
proof, the following estimates on Bessel
functions are needed:
 \begin{lemma}\cite{BOO}\label{Lemma 4.4}
 {\bf (i)} The following holds for $r\rg\iy$,
$$
\big\vert r^{1/3}J_{2r+xr^{1/3}}(2r)-A(x)\big\vert
=O(r^{-1/3}),
$$
uniformly in $x$, when $x\in$ compact $K\subset\BR$.

{\bf (ii)} For any $\dt >0$, there exists $M>0$ such
that for $x,y>M$ and large enough $r$,
$$
\left\vert \sum^{\iy}_{s=1}  J_{2r+xr^{1/3}+s}(2r)
J_{2r+yr^{1/3}+s}(2r)\right\vert <\dt r^{-1/3}.
$$
\end{lemma}

\begin{lemma} \label{Depoissonization lemma}
(de-Poissonization lemma) (Johansson
\cite{Johansson-Ulam}) Given $1\geq F_0\geq F_1 \geq
F_2\geq ...\geq 0$, and
$$F(\xi):=e^{-\xi}\sum_0^{\iy} \frac{\xi^k}{k!} F_k,$$
there exists $C>0$ and $k_0$, such that
$$F(k+4\sqrt{k\log k})-\frac{C}{k^2}\leq F_k\leq
   F(k-4\sqrt{k\log k})+\frac{C}{k^2},
  ~~\mbox{for all}~~ k>k_0 .$$
\end{lemma}

\medskip\noindent{\it Sketch of proof of Theorem \ref{BDJ}:\/}
 Putting the following scaling in
  the Bessel kernel $K(k, \ell)$, as in (\ref{BesselKernel}),
  one obtains, setting $r:=\sqrt{\xi}$,
 \bea
\lefteqn{ \xi^{1/6}K(2\xi^{1/2}+x\xi^{1/6},
 2\xi^{1/2}+y\xi^{1/6})}\no\\ 
 &=&
 (\sum_{k=1}^{N}+\sum_{k=N+1}^{\iy} ) \no\\
 &&\xi^{-1/6}
  \left[\xi^{
  1/6}J_{2\xi^{1/2}+(x+k\xi^{-1/6})\xi^{1/6}}
  (2\sqrt{\xi})\right]
 \left[ \xi^{ 1/6}J_{2\xi^{1/2}+(y+k\xi^{-1/6})\xi^{1/6}}
 (2\sqrt{\xi})\right]
 \no\\
 &=&
 (\sum_{k=1}^{N}+\sum_{k=N+1}^{\iy} ) \no\\
 &&r^{-1/3}
  \left[r^{
  1/3}J_{2r+(x+kr^{-1/3})r^{1/3}}
  (2r)\right]
 \left[ r^{ 1/3}J_{2r+(y+kr^{-1/3})r^{1/3}}
 (2r)\right] \label{4.0.2}
%
%
\eea
Fix $\dt>0$ and pick $M$ as in Lemma \ref{Lemma 4.4}
(ii). Define
$N:=\left[(M-m+1)r^{1/3}\right]=O(r^{\frac{1}{3}})$,
where $m$ is picked such that $x,y\geq m$ (which is
possible, since $x$ and $y$ belong to a compact set).
Then
%
\bea \lefteqn{\left\vert \sum^{\iy}_{k=N+1}
J_{2r+r^{1/3}(x+kr^{-1/3})}(2r)
 J_{2r+r^{1/3}(y+kr^{-1/3})}(2r)\right\vert}\no\\
 &=&
  \left\vert \sum^{\iy}_{s=1}
J_{2r+r^{1/3}(x+M-m+1)+s}(2r)
 J_{2r+r^{1/3}(y+M-m+1)+s}(2r)\right\vert <\dt
r^{-1/3};\no \\ \label{4a}\eea the latter inequality
holds, since
$$
  x+(M-m+1)=M+(x-m)+1>M.
$$
On the other hand,
$$
\left\vert r^{1/3}   J_{2r+r^{1/3}(x+kr^{-1/3})}(2r)
-A(x+kr^{-1/3}) \right\vert =O(r^{-1/3})
$$
uniformly for $x\in$ compact $K\subset\BR$, and all $k$
such that $1\leq k\leq N=\left[(M-m+1)r^{1/3}\right]$.
Indeed, for such $k$,
$$
m\leq x\leq x+kr^{-1/3}\leq x+M-m+1=M+(x-m)+1
$$
and thus, for such $k$, the $x+kr^{-1/3}$'s belong to a
compact set as well.
Since the number of terms in the sum below is
$N=[(M-m+1)r^{1/3}]$, Lemma \ref{Lemma 4.4} (i) implies
\bea  & &\Bigl\vert \sum_{k=1}^{N}\Bigl(
r^{\frac{1}{3}}J_{2r+xr^{1/3}+k} (2r)
r^{\frac{1}{3}}J_{2r+yr^{1/3}+k}(2r)
 %
\no\\
\no\\
&
&\hspace*{3cm} -
 A(x+kr^{-1/3})A(y+kr^{1/3})\Bigr)\Bigr\vert
=O(1) .\label{4b}\eea
But the Riemann sum tends to an integral \be r^{-1/3}
\sum_{k=1}^{[(M-m+1)r^{1/3}]}
A(x+kr^{-1/3})A(y+kr^{-1/3})\rg\int_0^{M-m+1}A(x+t)A(y+t)dt.
\label{4c}\ee Hence, combining estimates (\ref{4a},
\ref{4b}, \ref{4c}) and multiplying with $r^{1/3}$ leads
to
$$
\left\vert r^{\frac{1}{3}}  \sum_{k=1}^{\iy}
J_{2r+xr^{1/3}+k}
(2r)J_{2r+yr^{1/3}+k}(2r)-\!\int_0^{M-m+1}
A(x\!+\!t)A(y\!+\!t)dt\right\vert \leq\dt +o(1).
$$
for $r\rightarrow \iy$. Finally, letting $\dt\rg 0$ and
$M\rg\iy$ leads to the result.
%
%
Thus the expression (\ref{4.0.2}) tends to the Airy
kernel
$$
{\bf A}(x,y):=\int_0^{\iy}
 du A(x+u)A(y+u)
  %
           .$$
 for $r=\sqrt{\xi}\rightarrow \iy$. %
Hence
  \bean \lefteqn{ \lim_{\xi\rightarrow
\iy}e^{-\xi}
  \sum_{n=0}^{\iy}\frac{\xi^n}{n!}P^{}\left(L(\pi_n)\leq
2 {\xi}^{1/2}+x\xi^{1/6}\right)}\\ 
  &=&\lim_{\xi\rightarrow \iy}\left.
\det\left(I-K(\ell,\ell')\Big\vert_{[k,k+1,\ldots]}
 \right)\right|_{ \begin{array}{l}
 k=2 {\xi}^{1/2}+x\xi^{1/6}
 \end{array} }\\
 &=&1+\sum^{\iy}_{k=1} (-1)^k
  \int_{x\leq z_1\leq
...\leq z_k}\det\Bigl({\bf A}(z_i,z_j)\Bigr)_{1\leq
i,j\leq
k}\prod_1^k 
dz_i
\\
&=&\det \left(I-{\bf A} \chi_{[x,\iy)}\right). \eean
Finally, one uses Johansson's de-Poissonization Lemma
\ref{Depoissonization lemma}. From Corollary
\ref{cor:Pieri}, Plancherel measure $P_n(\lb_1\leq
x_1,\ldots,\lb_k\leq x_k)$ decreases, when $n$
increases, which is required by Lemma
\ref{Depoissonization lemma}. It thus follows that
$$\hspace*{-.7cm} {\displaystyle \lim_{n\rightarrow
\iy}
    P^{}\left(L(\pi_n)\leq
2 {n}^{1/2}+xn^{1/6}\right)}  =\det \left(I-{\bf A}
\chi_{[x,\iy)}\right),
  $$
  ending the proof of Theorem \ref{BDJ}.\qed

\newpage
 \newpage

\subsection{Limit theorem for the geometrically distributed
percolation model, when one side of the matrix tends to
$\iy$ }\label{subsect 4.3}

Consider Johansson's percolation model in section
\ref{section 1.7.2}, but with $p$ and $q$ interchanged.
 Consider the ensemble
 $$
 \mbox{Mat}^{(q,p)}=\{q\times p\mbox{~matrices $M$
with
  entries $M_{ij}=0,1,2,\ldots$}\}
 $$
 with {\em independent and geometrically
distributed entries}, for fixed $0<\xi<1$,
 $$
P(M_{ij}=k)=(1-\xi)\xi^k,~~~~~k=0,1,2,\ldots.
$$

 $$L(M):=\max_{ \mbox{\tiny
all such}\atop
                          \mbox{\tiny paths} }
\left\{\sum M_{ij},\quad\begin{array}{l}
\mbox{over right/down}\\
\mbox{paths starting from}\\
\mbox{entry $(1,1)$ to $(q,p)$}
\end{array}\right\} $$
 %
%
has the following distribution, 
assuming $q\leq p$,
\bean 
 P(L(M)\leq\ell)&=&Z^{-1}_{p,q} \sum_{h\in \BN^q \atop{\max(h_i)\leq
\ell
+q-1}}\Dt_q(h_1,\ldots,h_q)^2\prod^q_{i=1}\frac{(h_i+p-q)!}{h_i!}\xi^{h_i}
\eean
 where
\be
Z_{p,q}=\xi^{\frac{q(q-1)}{2}}(1-\xi)^{-pq}q!\prod^{q-1}_{j=0}j!(p-q+j)!
\label{constant2}\ee

Assuming that the number of columns $p$ of the $q\times
p$ random $M$ matrix above gets very large, as above,
the maximal right/lower path starting from (1,1) to
$(q,p)$ consists, roughly speaking, of many horizontal
stretches and $q$ small downward jumps. The $M_{ij}$
have the geometric distribution, with mean and standard
deviation
\bean E(M_{ij})&=&\sum_{k=0}^{\iy}k~P(M_{ij}=k)
 =(1-\xi)\xi
  \sum_1^{\iy}k\xi^{k-1}
   =\frac{\xi}{1-\xi}\\
%
 \sg^2(M_{ij})&=&E(M_{ij}^2)-(E(M_{ij}))^2
  =
  \sum k^2P(M_{ij}=k)-  \left(\frac{\xi}{1-\xi}\right)^2\\
  & &\hspace{3.8cm}=
   \frac{\xi(\xi+1)}{(1-\xi)^2}
    -  \left(\frac{\xi}{1-\xi}\right)^2
    =\frac{\xi }{(1-\xi)^2}
\eean
%
So, in the average,
$$
L(M)\simeq pE(M_{ij})=\frac{p\xi}{1-\xi},\quad\mbox{for
$p\rg\iy$}
$$
with
$$
\sg^2\left(L(M)-\frac{p\xi}{1-\xi}\right)\simeq
p\sg^2(M_{ij})=\frac{p\xi}{(1-\xi)^2}.
$$
Therefore, it seems natural to consider the variables

$$
x_1=\frac{L(M)-\frac{\xi}{1-\xi}p}{\frac{\sqrt{\xi
p}}{1-\xi}}=
\frac{\lb_1-\frac{\xi}{1-\xi}p}{\frac{\sqrt{\xi
p}}{1-\xi}}
$$
and, since $\lb=(\lb_1,\ldots,\lb_q)$, with $q$ finite,
all $\lb_i$ should be on the same footing. So,
remembering from the proof of Theorem \ref{theorem:
percolation} that $h_i$ in (\ref{constant2}) and $\lb_i$
are related by
 $$
 h_i=q+\lb_i-i,
 $$
  we set
$$
x_i=\frac{\lb_i-\frac{\xi}{1-\xi}p}{\frac{\sqrt{\xi
p}}{1-\xi}}.
$$

\hspace*{ -1.2cm}$(1,1)$ \hspace*{11.7cm}$(1,p)$
$$\hspace*{-1.0cm}
  \squaresize .38cm \thickness .02cm \Thickness
.05cm \Young{
     u&u&ur&&&&&&&&&&&&&&&&&&&&&&&&&&&&&&\cr
      &&&u&u&u&ur&&&&&&&&&&&&&&&&&&&&&&&&&&\cr
      &&&&&&&u&u&u&u&u&u&u&ur&&&&&&&&&&&&&&&&&&\cr
      && && &&&&&&&&&&&u&ur&&&&&&&&&&&&&&&& \cr
      &&&&&&&&&&&&&&&&&u&u&u&u&u&u&u&u&u&u&ur&d&d&d&d& d  \cr
      }
$$
\hspace*{ -.4cm}$(q,1)$ \hspace*{11.8cm}$(q,p)$

\begin{theorem} (Johansson \cite{Johansson1}) \label{Theo:4.5} The following limit holds:
$$
\lim_{p\rightarrow \iy}
  P\left( \frac{L(M)-\frac{\xi}{1-\xi}p}{\frac{\sqrt{\xi p}}{1-\xi}}
  \leq y\right)
  =
   \frac{\displaystyle{\int_{(-\iy,
y)^q}\Dt_q(x)^2\prod^q_1 e^{- x_i^2/2}dx_i}}
 {\displaystyle{\int_{\BR^q}}\Dt_q(x
)^2\prod^q_1e^{-x_i^2/2}dx_i} ,$$ which coincides with
the probability that a $q\times q$ matrix from the
Gaussian Hermitian ensemble (GUE) has its spectrum less
than $y$; see section \ref{sect 8.2}.

\end{theorem}

 \proof The main tool here is Stirling's formula\footnote{
Stirling's formula:
$$
n!=\sqrt{2\pi~ n}~e^{n(\log n-1)}\left(1+{
O}\left(\frac{1}{n}\right)\right),\mbox{~for $n\rg\iy$},
$$}.
Taking into account $
 h_i=q+\lb_i-i,
 $ we substitute
$$
\lb_i=\frac{\xi p}{1-\xi}+\frac{\sqrt{\xi p}}{1-\xi}x_i
$$
in (\ref{constant2}). The different pieces will now be
computed for large $p$ and fixed $q$,
  \bean
   \lefteqn{\prod^q_{i=1}(h_i+p-q)!}\\
&=&{\prod^q_{i=1}(p+\lb_i-i)!}\\
&=&\prod^q_{i=1}\left(\frac{1}{1-\xi}(p+x_i\sqrt{p\xi})-i\right)!\\
\\
&=&\left(\frac{2\pi~p}{1-\xi}\right)^{q/2}\prod^q_{i=1}
\left(1+x_i\sqrt{\frac{\xi}{p}}-\frac{i(1-\xi)}{p}\right)^{1/2}
\left(1+{\bf O}\left(\frac{1}{p}\right)\right)\\
\\
&
&e^{\displaystyle{\sum^q_{i=1}\left(\frac{p}{1-\xi}+x_i\frac{\sqrt{\xi
p}}{1-\xi}-i\right)\left(\log\frac{p}{1-\xi}-1+\log\left(1+x_i
\sqrt{\frac{\xi}{p}}-\frac{i(1-\xi)}{p}\right)\right)}}\\
\\
&=&\left(\frac{2\pi p}{1-\xi}\right)^{q/2}
e^{\displaystyle{\frac{\xi}{2(1-\xi)}\sum^q_1x_i^2-
\frac{qp}{1-\xi}+{\bf O}\left(\frac{1}{\sqrt{p}}\right)}}\\
\\
&&\hspace*{1cm}
 e^{\displaystyle{\left(\frac{qp}{1-\xi}+\frac{\sqrt{\xi p}}
  {1-\xi}\sum_1^qx_i-\frac{q(q+1)}{2}\right)\log
\frac{p}{1-\xi}}}\left(1+{\bf
O}\left(\frac{1}{p}\right)\right), \eean upon expanding
the logarithm in powers of $1/\sqrt{p}$.
 Similarly
%
\bean
\lefteqn{\prod^q_{i=1}h_i}\\
&=&\prod^q_{i=1}(q+\lb_i-i)! \\
&=&\prod^q_{i=1}\left(\frac{1}{1-\xi}(p\xi +x_i\sqrt{p\xi})+q-i\right)!\\
\\
&=&\left(\frac{2\pi p\xi}{1-\xi}\right)^{q/2}
 \prod^q_{i=1}\left(1+\frac{x_i}{\sqrt{p\xi}}
+\frac{(q-i)(1-\xi)}{p\xi}\right)^{1/2}
 \left(1+{\bf O}\left(\frac{1}{p}\right)\right)
 \\
\\
& &.e^{\displaystyle{\sum^q_{i=1} \left(\frac{1}{1-\xi}
(p\xi+x_i\sqrt{p\xi})+q-i\right)
\left(\log\frac{p\xi}{1-\xi}-1\right)}}
  \\&&
 .e^{\displaystyle{
  \sum^q_{i=1}\left(\frac{1}{1-\xi}
 (p\xi+x_i\sqrt{p\xi})+q-i\right)
%
%
\log\left(1+\frac{x_i}{\sqrt{p\xi}}
+\frac{(q-i)(1-\xi)}{p\xi}\right) }}\\
&=& 
 \left(\frac{2\pi p\xi}{1-\xi}\right)^{q/2}
%
  e^{\displaystyle{\frac{1}{2(1-\xi)}
\sum^q_1x_i^2-\frac{qp\xi}{1-\xi }
 }}
  \left(1+ {\bf O}\left(\frac{1}{\sqrt{p}}\right)\right)
  \\
 \\
 & &\hspace*{.5cm}
.e^{\displaystyle{\left(\frac{q\xi
p}{1-\xi}+\frac{\sqrt{p\xi}}{1-\xi}\sum_1^qx_i+\frac{q(q-1)}{2}\right)\log\frac{
p\xi}{1-\xi}}} \eean
Also
 \bean
  \prod^q_1\xi^{h_i}=
\prod^q_1\xi^{q+\lb_i-i}&=&e^{\displaystyle{\sum^q_1
\left(\frac{1}{1-\xi}(p\xi+x_i
\sqrt{p\xi})+q-i\right)\log\xi}}\\
\\
&=&e^{\displaystyle{\left(\frac{q\xi p}{1-\xi}+\frac{
\sqrt{p\xi}}{1-\xi}\sum_1^qx_i+\frac{q(q-1)}{2}\right)\log\xi}}
.\eean
%
%
Besides
  \bean
\prod^{q-1}_{j=0}(p-q+j)!  
&=&(2\pi p)^{q/2}\prod^{q-1}_{j=0}
 \left(1-\frac{q-j}{p}\right)^{1/2}
  \left(1+{\bf O}\left(\frac{1}{p}\right)\right)
  \\
 &&\hspace*{5mm}
 e^{\displaystyle{\sum^{q-1}_{j=0}(p-q+j)\left(\log
p-1+\log\left(1-\frac{q-j}{p}\right)\right)}}
  \\
&=&(2\pi p)^{q/2}
 e^{\displaystyle{-pq+qp\log
p-\frac{q(q+1)}{2}\log p}}
 \left(1+{\bf
O}\left(\frac{1}{p}\right)\right)
\\&&\\
%
%
%
 &=&
   {(2\pi  )^{q/2} p^{qp-q^2/2} e^{-qp}}
    \left(1+{\bf
O}\left(\frac{1}{p}\right)\right)
  \eean
   and
  \bean
\Dt_q(h_1,..., h_q)&=&\prod_{1\leq i<j\leq q}(h_i-h_j)\\
&=&\prod_{1\leq i<j\leq q}\left((x_i-x_j)\frac{\sqrt{\xi
p}}{1-\xi}-i+j\right)\\
\\
&=&\left(\frac{\sqrt{\xi
p}}{1-\xi}\right)^{\frac{q(q-1)}{2}}\left(\prod_{1\leq
i<j\leq q}(x_i-x_j)+{\bf
O}\left(\frac{1}{\sqrt{p}}\right)\right) \eean
%
%
%
Remembering the relation $h_i=q+\lb_i-i$, and using the
estimates before in equality $\stackrel{*}{=}$ below,
one computes
\bean
\lefteqn{Z^{-1}_{p,q}\sum_{h\in\mathbb N^q\atop{\max
h_i\leq\ell +q-1}}\Dt_q(h_1,...,
h_q)^2\prod_1^q\frac{(h_i+p-q)!}{h_i!}\xi^{h_i}}\\
\\
&=& \!\!\!\!\!\!\!\! \sum_{h\in\mathbb N^q\atop{\max
h_i\leq\ell
+q-1}}\frac{\xi^{-\frac{q(q-1)}{2}}}{\prod^q_{j=0}j!}(1-\xi)^{pq}
\Dt_q(h_1,...,h_q)^2\prod^q_{i=1}
\frac{(p+\lb_i-i)!~\xi^{h_i}}{(q+\lb_i-i)!(p-q+i-1)!}
\\
\\
&\stackrel{*}{=}&
 \left(\frac{1}{(2\pi  )^{q/2}\prod_1^q i!}\right)
  \sum_{x_i\leq y}
 \left(\frac{1-\xi}{\sqrt{\xi p}}\right)^q
  \Dt(x)^2e^{\displaystyle{-\frac{1}{2}\sum^q_1x_i^2}}
   \left(1+{\bf
O}\left(\frac{1}{\sqrt{p}}\right)\right)
\eean\bean
&\longrightarrow&
  \frac{\displaystyle{\int_{x_i\leq y}\Dt(x)^2
  \prod^q_1e^{-x_i^2/2}
  dx_i}}{\displaystyle{\int_{\BR^q}}\Dt_q(x
)^2\prod^q_1e^{-x_i^2/2}dx_i},\eean
 using Selberg's formula
$$
\int_{\BR^q}
\Dt_q(x)^2\prod^q_1e^{-x_i^2/2}dx_i=(2\pi)^{q/2}\prod^q_1i!
$$
  and noticing that in
   $$
h_i=q+\lb_i-i=\frac{\xi p}{1-\xi}+\frac{\sqrt{\xi
p}}{1-\xi}x_i+q-i
$$
an increment of one unit in $\lb_i$ implies an increment
of $x_i$ by
  $$
 dx_i= \frac{1-\xi}{\sqrt{p\xi}}.
  $$
Therefore, one has
 $$
 \left(\frac{1-\xi}{\sqrt{p\xi}}
   \right)^{q }
  \simeq  \prod^q_1dx_i
 .$$
Finally, for $p$ large, one has
 $$
  h_i\leq
 \ell+q-1=\frac{\xi p}{1-\xi}
  +\frac{\sqrt{\xi p}}{1-\xi}y+q-1
\Longleftrightarrow x_i\leq y.
 $$
  The connection with the spectrum of Gaussian Hermitian matrices
   will be discussed in section \ref{sect 8.2}. This ends the proof of Theorem \ref{Theo:4.5}.
\qed


 \subsection{Limit theorem for the geometrically distributed
percolation model, when both sides of the matrix tend to
$\iy$ }


The model considered here is the percolation model as in
section \ref{subsect 4.3}, whose probability can, by
section \ref{section 3.3}, be written as a Fredholm
determinant of a Meixner kernel; see Theorem
\ref{Theo:3.3}.

\begin{theorem}(Johansson \cite{Johansson0})\label{Theo:4.7}
Given
 $$
  a=\frac{(1+\sqrt{\xi\gamma})^2}{1-\xi}-1,
  ~~~
 \rho =\left(\frac{\xi}{\gamma}\right)^{1/6}
  \frac{
  (\sqrt{\gamma}+\sqrt{\xi})^{2/3}
  (1+\sqrt{\gamma\xi})^{2/3}}{1-\xi},
 $$
the following holds:
 \bean
 \lim_{\begin{array}{c}
  p,q\rightarrow \iy\\
  q/p=\gamma\geq 1~\mbox{fixed}
  \end{array}}
   P\left(
 \frac{L(M_{q,p})- a p}
   { \rho p^{1/3} }
  \leq y
  \right)
&= &{\cal F}(y),
 \eean
 which is the Tracy-Widom distribution (see section \ref{subsect
 8.5}). In other terms, the random variable $L(M_{q,p})$
 behaves in distribution, like
 $$
  L(M_{q,p})\sim ap+\rho p^{1/3} {\cal F}
  .$$
%


\vspace*{1cm}

 \hspace*{1.5cm}$(1,1)$ \hspace*{7.6cm}$(1,p)\rightarrow \iy$
$$
  \squaresize .8cm \thickness .02cm \Thickness
.07cm \Young{
     l&&&&&&&&&\cr
     l &&&&&&&&&\cr
     ur&ld  & & &&&&&&\cr
      &&l&& &&&&&\cr
      &&l&& &&&&&\cr
      &&u&u&u&ur&&&&\cr
      &&&&&&u&ur&d&\cr
      & &&&&&&&r& d \cr
  }
$$
$$\begin{array}{c}
  (q,1)  \\
   \downarrow  \\
   \iy
   \end{array}
  \hspace*{7cm} \begin{array}{c}
  (q,p)  \\
   ~~~\searrow  \\
   ~~~~~~~~~\iy
   \end{array}$$

   \vspace*{1cm}

   \end{theorem}

\medskip\noindent{\it Sketch of Proof:\/}
 Remember from Theorem \ref{Theo:3.3}, the
probability $P(L(M_{q,p})\leq z)$ is given by a Fredholm
determinant of a Christoffel-Darboux kernel composed of
Meixner polynomials. Proving the statement of Theorem
\ref{Theo:4.7} amounts to proving the limit of this
Meixner kernel with the scaling mentioned in the Theorem
tends to the Airy kernel, i.e.,
$$
           \lim_{p\rightarrow \iy}
           \rho p^{1/3}K_p(a p+\rho p^{1/3}\eta,
                           a
                           p+\rho p^{1/3}\eta')={\bf
                           A}(\eta,\eta'):=
                            \frac{A (\eta)A' (\eta')-A' (\eta
                            )A (\eta')}{\eta-\eta'}
                           $$
where\footnote{with $
           d_p=\frac{p!(p+\beta'-1)!}
           {(1-\xi)^{\beta'}\xi^p(\beta'-1)!}
           $ and $\beta'=q-p+1=p\gamma -p+1$.} \bean K_p(x,y)&=&
 -\frac{\xi}{(1-\xi)d_{p-1}^2}
 \sqrt{\left({x+\beta'-1}\atop {x}\right)\xi^x
           \left({y+\beta'-1}\atop {y}\right)\xi^y}
            \\ \\
            &&~~~~\times ~\frac{m_p(x,\beta',\xi)m_{p-1}(y,\beta',\xi)
    -m_p(y,\beta',\xi)m_{p-1}(x,\beta',\xi)}{x-y}
           \eean
           for the Meixner polynomials, which have the following integral representation, a consequence of the
            the generating function
          ( \ref{genfunction}), (slightly
           rescaled)
            $$
           m_p(x,\beta',\xi)= p!\xi^{-x}\oint_{|z|=r<1}
           \frac{(\xi-z)^x}{z^{p}(1-z)^{x+\beta'}}
           \frac{dz}{2\pi iz}$$
Thus, it will suffice to prove the following limit \bea
\lefteqn
{\lim_{p\rg\iy}p^{1/3}\left(\frac{y^p(y-1)^{x+[\gamma
p]-p+1}}{(y-\xi)^x}\right)
\Biggl|_{y=-\sqrt{\frac{\xi}{\gamma}}}}\no\\
 \no\\
 &&~~~
{\oint_{|z|=r<1}\frac{(z-\xi)^x}{z^p(z-1)^{x+[\gamma
p]-p+1}}\frac{dz}{2\pi iz}\Bigg|_{x=\beta p+\rho
p^{1/3}\eta}}=C A(\eta) \label{limit Meixner},\eea where
$$
\al:=\frac{(\sqrt{\xi}+\sqrt{\gamma})^2}{1-\xi},~~~~
 \beta
=\al-\gamma+1=\frac{(\sqrt{\gamma\xi}+1)^2}{1-\xi},$$
 $$\qquad\rho
=\left(\frac{\xi}{\gamma}\right)^{1/6}\frac{(\sqrt{\gamma}+\sqrt{\xi})^{2/3}
 (1+\sqrt{\gamma\xi})^{2/3}}{1-\xi}
$$
$$
C:=
 \gamma^{-1/3}\xi^{-1/6}(\sqrt{\xi}
+\sqrt{\gamma})^{1/3}(1+\sqrt{\xi\gamma})^{1/3} .$$
%
%
In view of the saddle point method, define the function
$F_p(z)$ such that
$$
 e^{pF_p(z)}:=\frac{(z-\xi)^x}{z^p(z-1)^{x+[\gamma
p]-p+1}}\Bigg|_{x=\beta p+\rho p^{1/3}\eta}
 .
$$
Then one easily sees that
 \be
 F_p(z)=F(z)+\rho p^{-2/3}\eta \log\frac{z-\xi}{z-1}
  +(\gamma -\frac{[\gamma p]}{p})\log(z-1)\label{Fp},\ee
where the $p$-independent function $F(z)$ equals,
 \be
 F(z)=\beta\log(z-\xi)-\al\log(z-1)-\log z.\label{F(z)}
  \ee
For the specific values above of $\al$ and $\beta$, one
checks the function $F(z)$ has a critical point at
$z_c:=-\sqrt{ {\xi}/{\gamma}}$, i.e.,
 $$
 F '\left(z_c\right)
 =F ''\left(z_c\right)=0
,~~~~~F '''\left(z_c\right) =
 \frac{2\gamma^{5/2}}{\xi(\sqrt{\gamma}+\sqrt{\xi})(1+\sqrt{\xi\gamma})}
 $$
 and thus
 \be
 F(z)-F \left(z_c\right)=
  \frac{1}{6} (z-z_c)^3
   F'''\left(z_c\right)+
   O((z-z_c)^4)
   .\label{F Taylor}\ee
   Setting
    \be
    z= z_c(1-p^{-1/3}sC)e^{itC p^{-1/3}}\label{z}
    \ee
     in $F_p(z)$ as in (\ref{Fp}), one first computes this substitution
     in $F(z)$, taking into account
     (\ref{F Taylor}),
     \bea
     F(z)&=& F(z_c)+
 \frac{i(-z_cC)^{3}}{6p} \left(t+is \right)
^{3}
F'''(z_c)+O(p^{-4/3}) \no\\
 &=& F(z_c)+ \frac{i}{3p}(t+is)^3+O(p^{-4/3}),
\label{4.4.6}\eea  upon picking -in the last equality-
the constant $C$ such that
  $$
   (-z_cC)^3F'''\left(z_c\right)=2 .
 $$
 Also, substituting the $z$ of (\ref{z}) in the part
 of $F(z)$ (see (\ref{F(z)})) containing $\eta$,
%

      \be
   \rho p^{-2/3} \eta \log\frac{z-\xi}{z-1}
   =    \frac{\rho \eta}{p^{2/3}}
    \log    {\frac {  z_c-\xi}{ z_c-1}}   -
      {\frac
{  \rho   C z_c\left(  \xi-1 \right)   }{
 \left( z_c-1 \right)  \left( z_c-\xi \right) }}\frac{i\eta (t+is)}{p}
+O(p^{-4/3}).\label{4.4.7}\ee
 Thus, adding the two contributions (\ref{4.4.6})
 and (\ref{4.4.7}), one finds
   \be
   pF_p\left(z_c(1-p^{-1/3}sC)e^{itC p^{-1/3}}\right)
   =
   pF_p(z_c)+\frac{i}{3}\Bigl((t+is)^3+3\eta  (t+is)\Bigr)
   +O\left(p^{-1/3}\right)
   \label{subst}\ee
   One then considers two contributions of the contour
   integral about the circle $|z|=r$ appearing in (\ref{limit Meixner}), a first
   one along the arc $(\pi-\delta_p,\pi+\delta_p)$, for
   $\delta_p$ tending to $0$ with $p\rg \iy$ and a
   second one about the complement of
   $(\pi-\delta_p,\pi+\delta_p)$. The latter tends to
   $0$, whereas the former tends to the Airy function
   (keeping $s$ fixed, and in particular $=1$)
  $$  CA(\eta)=\frac{C}{2\pi} \int_{\iy}^{\iy} e^{i(\frac{(t+is)^3}{3}
    +\eta(t+is))} dt,
    $$
    upon noticing that $dz/z=iCp^{-1/3}dt$ under the change
    of variable $z\mapsto t$, given in (\ref{z}), establishing limit (\ref{limit Meixner}) and finally
    the limit of the Meixner kernel and its Fredholm
    determinant. Further details of this proof can be found
    in Johansson \cite{Johansson0}. The Fredholm determinant of the Airy kernel
    is precisely the Tracy-Widom distribution,
    as will be shown in section \ref{subsect 8.5}.
    \qed

\subsection{Limit theorem for the exponentially distributed
percolation model, when both sides of the matrix tend to
$\iy$ }

Referring to the exponentially distributed percolation
model, discussed in section \ref{theorem:
percolation-exp}, we now state

\begin{theorem}\label{Theo:4.6}(Johansson
\cite{Johansson0})
  Given $$
  a = {(1+\sqrt{ \gamma})^2}  ,
  ~~~
 \rho =
  \frac {
  (1+\sqrt{\gamma})^{4/3}}{{\gamma} ^{1/6}}
 ,$$
  the following limit holds:
 \bean \lefteqn{\lim_{\begin{array}{c}
  p,q\rightarrow \iy\\
  q/p=\gamma\geq 1~\mbox{fixed}
  \end{array}}
   P\left(
 \frac{L(M_{ q,p })-
  ( p^{1/2}+q^{1/2})  }
   {
   (p^{1/2}+q^{1/2})(p^{-1/2}+q^{-1/2})^{1/3}}
  \leq y
  \right)}\\
  &=& \lim_{\begin{array}{c}
  p,q\rightarrow \iy\\
  q/p=\gamma\geq 1~\mbox{fixed}
  \end{array}}
   P\left(
 \frac{L(M_{q,p})- a p}
   { \rho p^{1/3} }
  \leq y
  \right)
={\cal F}(u), \eean
which is again the Tracy-Widom distribution.
 Here again $
  L(M_{q,p})$ behaves, 
  after some rescaling and in distribution, like the Tracy-Widom distribution, for large $p$ and $q$ such that $q/p=\gamma\geq 1$:
 $$
  L(M_{q,p}) \sim ap+\rho p^{1/3} {\cal F}
$$

 \end{theorem}

 \proof In Theorem \ref{theorem: percolation-exp},
 one has shown that $P(L(M)\leq t)$ equals the ratio of
 two integrals; this ratio will be shown in section \ref{section 7} on
 random matrices
 (see Propositions \ref{Proposition 7.9}
 and \ref{Proposition 7.4}) to equal
 a Fredholm determinant of a kernel corresponding
  to Laguerre polynomials:
 $$P(L(M)\leq t)=\frac{\int_{[0,t]^q}\Dt_p(x)^2
\prod^p_{i=1}x_i^{q-p}e^{-x_i}dx_i}{\int_{[0,\iy]^q}\Dt_p(x)^2
\prod^p_{i=1}x_i^{q-p}e^{-x_i}dx_i} =\det
(I-K_p^{(\al)}(x,y)\chi_{[t,\iy]})$$
 where
\bean K_p^{(\al)}(x,y)&=&
 \sqrt{\frac{h_{p}}{h_{p-1}}}(xy)^{\al/2}e^{-\frac{1}{2}(x+y)}
 \frac{{\cal L}_p^{(\al)}(x)
{\cal L}_{p-1}^{(\al)}(y)-{\cal L}_p^{(\al)} {\cal
L}_{p-1}^{(\al)}(x)}{x-y};
   \eean
in the formula above, the
$\LR_p^{(\al)}(x)=\frac{1}{\sqrt{h_p}}x^n+\ldots= (-1)^p
\left(\frac{p!}{(p+\al)!} \right)^{1/2}L_p^{\al}(x)$ are
the normalized Laguerre
polynomials\footnote{$L_n^{(\al)}(y)=\sum_{m=0}^n (-1)^m
\left( {n+\al}\atop  n-m
\right)\frac{x^m}{m!}=\frac{e^x}{2\pi
i}\int_C\frac{e^{-xz}z^{n+\al}}{(z-1)^{n+1}}dz$, where
$C$ is a circle about $z=1$.}:
$$
\int^{\iy}_0\LR_n^{(\al)}(x)\LR_m^{(\al)}(x)x^{\al}e^{-x}dx=\delta_{nm}.
$$
Therefore, using the precise values of $a$ and $\rho$
above,
 \bean P\left(
 \frac{L(M_{q,p})- a p}
   { \rho p^{1/3} }
  \leq y
  \right)&=&\det\left(I-{\cal
K}(\xi,\eta)\chi_{[y,\iy)}\right) \eean
with
 $$ {\cal K}(\xi,\eta)=bp^{1/3}K_p^{(\gamma
-1)p}\left(ap+\rho p^{1/3}\xi,ap+\rho p^{1/3}\eta\right)
$$ The result follows from an asymptotic formula for
Laguerre polynomials
$$
\lim_{p\rightarrow \iy} K_p^{(\gamma -1)p}\left(ap+\rho
p^{1/3}\xi,ap+\rho p^{1/3}\eta\right)={\bf A}(\xi,\eta)
,$$ with ${\bf A}(x,y)$ the Airy kernel, as in
(\ref{AiryK}),namely
$$
{\bf A}(x,y):=\frac{A (x)A' (y)-A' (x)A (y)}{x-y}.
$$
The Fredholm determinant of the Airy kernel is the
Tracy-Widom distribution. This ends the proof of Theorem
\ref{Theo:4.6}.\qed



\section{Orthogonal polynomials for a
time dependent weight and the KP equation}

\subsection{Orthogonal polynomials} \label{sect:
Orthogonal polynomials}\label{section 5.1}


The inner-product with regard to the weight $\rho(z) $
over $ \BR$, assuming $\rho(z) $ decays fast enough at
the boundary of its support\footnote{In this section,
the support of the weight $ \rho (z)$ can be the whole
of $\BR$ or any other interval.}
 \be
\la f,g\ra =\int_{\BR}f(z)g(z)\rho (z)dz,
 \label{5.1}\ee
  leads to a moment matrix
   \be
   m_{n} =(\mu_{ij} )_{0\leq
i,j<n}=(\la z^i,z^j\ra )_{0\leq i,j\leq n-1}.
 \label{5.2}\ee
Since the $\mu_{ij}$ depends on $i+j$ only, this is a
{\em H\"ankel matrix}, and thus symmetric. This is
tantamount to the relation $$\Lambda
m_{\iy}=m_{\iy}\Lambda^{\top},$$
  where $\Lb$ denotes
 the semi-infinite shift matrix
  $$\left(\begin{array}{ccccccc}
   0&1&0&0&0&0&\ldots  \\
   0&0&1&0&0&0&\ldots  \\
   0&0&0&1&0&0&\ldots  \\
   \vdots
\end{array}
\right).
 $$
 Define
 $$
 \tau_n:= \det m_n
 .$$
 Consider the
factorization of $m_{\iy}$ into a lower- times an
upper-triangular matrix\footnote{This factorization is
possible as long as $ \tau_n:=\det m_n \neq 0 $ for all
$n\geq 1$.}
 \be
    m_{\iy} =S ^{-1}S^{\top
-1},\label{factorization}
 \ee
  with
 $$   S =\mbox{\,lower triangular with
 non-zero diagonal elements.}
 $$
For any $z\in \BC$, define the semi-infinite column
 \be
 \chi(z):=(1,z,z^2,\ldots)^{\top},
  \ee
%
%
%
and functions  $p_n(z)$ and $q_n(z)$,
 \bea
 p_n(z)&:=&
\left(S\chi(z)\right)_n
   \no \\
 q_n( z)&:=&\left(S^{\top -1} \chi(z^{-1})\right)_n
 .\eea

\vspace*{.5cm}

 \noindent {\bf (1)} {\em  The $p_n(z)$ are
polynomials of degree $n$, orthonormal with regard to
$\rho(z)$, and $
 q_n(z)$ is the Stieltjes transform of $p_n(z)$,}
 $$
 q_n(z)=z  \int_{\BR} \frac{ p_n( u)\rho (u)}{z-u}du
 .$$
 Indeed,
  $$ \Bigl(\la
p_k ,p_{\ell} \ra\Bigr)_{0\leq k,\ell
<\iy}=\int_{\BR}S\chi(z)(S\chi(z))^{\top}\rho (z)dz=S
m_{\iy}S^{\top}=I. $$
 Note that $S\chi(z)(S\chi(z))^{\top}$
  is a
semi-infinite matrix obtained by multiplying the
semi-infinite column $S\chi(z)$ and row
$(S\chi(z))^{\top}$.
%
%
The definition of $q_n$, together with the decomposition
 $S^{\top -1}=Sm_{\iy}$, leads to
\begin{eqnarray*}
q_n(z)
 &=&
  \left( S^{\top
-1}\chi(z^{-1}\right)_n \\
&=&\sum_{j\geq 0} \left( S m_{\iy} \right)_{nj} z^{-j}
 \\
  &=&\sum_{j\geq 0} z^{-j} \sum_{0\leq \ell\leq n}
 S_{n\ell}~ \mu_{\ell j} \\
 &=&  \sum_{j\geq 0 } z^{-j}
\sum_{0\leq \ell\leq n} S_{n\ell} \int_{\BR}
u^{\ell+j}\rho (u)du
\\ &=&  \int_{\BR} \sum_{0\leq \ell\leq n} S_{n\ell}  u^{\ell}
\sum_{j\geq 0}\left(\frac{u}{z} \right)^j\rho(u)du\\
&=&z  \int_{\BR} \frac{ p_n( u)\rho (u)}{z-u}du.\\
\end{eqnarray*}

\vspace*{.5cm}

\noindent{\bf (2)} {\em The orthonormal polynomials
$p_n$ have the following representation}
 \be p_n( z)=
 \frac{1}{\sqrt{\tau_n \tau_{n+1} }}
\det\left(
\begin{array}{lll|l}
  & & &1\\
  &m_n& &z\\
  & & &\vdots\\
 \hline
 \mu_{n,0}&\ldots&\mu_{n,n-1}&z^n
\end{array}
\right).
\label{2}\ee {\em As a consequence, the monic orthogonal
polynomials $\tilde p_n(z)$ are related to $p_n(z)$ as
follows:} \be
 p_n( z)=\sqrt{\frac{\tau_n}{\tau_{n+1}} }
  \tilde p_n(z).
  \label{2a}\ee
Defining $p_n'(z)$ to be the polynomial on the right
hand side of (\ref{2}), it suffices to show that for
$k<n$,
 $$
 \la p'_n,p'_k \ra=0\mbox{~~and~~}
  \la p'_n,p'_n \ra=1,
 $$
thus leading to $p_n=p'_n$. Indeed,
 $$
 \la   p'_n,  z^k \ra=
       \frac {1}{\sqrt{\tau_n\tau_{n+1}}}
  \det\left(
\begin{array}{lll|l}
  & & &\la 1, z^k\ra\\
  &m_n& &\la z, z^k\ra\\
  & & &\vdots\\
 \hline
 \mu_{n,0}&\ldots&\mu_{n,n-1}&\la z^n,z^k\ra
\end{array}
\right)=\left\{ \begin{array}{l}
   0,~~~\mbox{for $k<n$,}\\ \\
   \sqrt{\frac{\tau_{n+1}}{\tau_n}}
    ~~\mbox{for $k=n$,}\end{array}\right.
$$
 and thus for $k=n$,
 $$
 \la   p'_n,  p'_n \ra
  =\sqrt{\frac {\tau_n}{  \tau_{n+1} }}
  \la   p'_n,z^n+\ldots\ra=
  \sqrt{\frac {\tau_n}{  \tau_{n+1} }}
  \la   p'_n,z^n \ra=1 ,
  $$
 from which (\ref{2}) follows. Formula (\ref{2a}) is a
 straightforward consequence.

\vspace*{.5cm}

\noindent{\bf (3)} {\em The monic orthogonal polynomials
$ \tilde p_n$ and their Stieltjes transform have the
following representation}
 \bea
 \tilde p_n( z)&=&
   {z^n}
   \frac
  {\det\left(\mu_{ij}-\frac{1}{z}\mu_{i,j+1}\right)_{0\leq
i,j\leq n-1}
  }
  {\det\left(\mu_{ij} \right)_{0\leq
i,j\leq n-1}
  }\no\\
 %
 %
  \int_{\BR}
  \frac{\tilde p_n(u)\rho(u)}{z-u} du
   &=&
 {z^{-n-1}}
   \frac
  {\det\left(\mu_{ij}+\frac{1}{z}\mu_{i,j+1}
   +\frac{1}{z^2}\mu_{i,j+2}+\ldots\right)_{0\leq
i,j\leq n }
  }
  {\det\left(\mu_{ij} \right)_{0\leq
i,j\leq n-1}
  }\no\\
 \label{2b} \eea

  \proof Setting
   $$
 \vec\mu_j=(\mu_{0j},...,\mu_{n-1,j})\in\BR^n,
   $$
one computes
\bean
\lefteqn{z^n\det\left(\mu_{ij}-z^{-1}\mu_{i,j+1}\right)_{0\leq i,j\leq n-1}}\\
\\
&=&\det\left(z\mu_{i,j}-\mu_{i,j+1}\right)_{0\leq i,j\leq n-1}\\
\\
&=&\det(z\vec\mu_0^{\top}-\vec\mu_1^{\top},z\vec\mu_1^{\top}-\vec\mu_2^{\top},
...,z\vec\mu^{\top}_{n-1}
-\vec\mu_n^{\top})\\
\\
&=&\det\left(\sum_0^{n-1}\frac{z\vec\mu_j^{\top}-
\vec\mu_{j+1}^{\top}}{z^j},
\sum_0^{n-2}\frac{z\vec\mu_{j+1}^{\top}-
\vec\mu_{j+2}^{\top}}{z^j},...,
z\vec\mu_{n-1}^{\top}-\vec\mu_n^{\top}\right)\\
& &\hspace*{4cm}\mbox{by column operations}\\
&=&\det\left(z\vec\mu_0^{\top}-\frac{\vec\mu_n^{\top}}{z^{n-1}},
z\vec\mu_1^{\top}-\frac{\vec\mu_n^{\top}}{z^{n-2}},...,
z\vec\mu_{n-1}^{\top}-\vec\mu_n^{\top}\right)
\eean\bean &=&\frac{1}{z^n}\det\left(\begin{array}{c|c}
z\vec\mu_0-\frac{\vec\mu_n}{z^{n-1}}&0\\
\\
z\vec\mu_1-\frac{\vec\mu_n}{z^{n-2}}&0\\
\vdots&\vdots\\
z\vec\mu_{n-1}-\vec\mu_n &0\\
\\
\vec\mu_n&z^n
\end{array}\right)\\
& &\hspace*{3cm}\mbox{enlarging the matrix by one row and column}\\
&=&\frac{1}{z^n}\det\left(\begin{array}{c|c}
z\vec\mu_0 &z\\
\\
z\vec\mu_1 &z^2\\
\vdots&\vdots\\
z\vec\mu_{n-1} &z^n\\
\\
\vec\mu_n&z^n
\end{array}\right)\\
& &\hspace*{3cm}\mbox{by adding a multiple of the last row to rows $1$ to $n$}\\
&=&\tau_n \tilde p_n(z). \eean


\noindent Setting this time
$$
\vec\mu_j:=(\mu_{0j},...,\mu_{nj})\in\BR^{n+1} ,
 $$ one
computes
\bean
\lefteqn{\det\left(\mu_{ij}+\frac{\mu_{i,j+1}}{z}+\frac{\mu_{i,j+2}}{z^2}+\ldots\right)_{0\leq
i,j\leq n}}\\
\\
&=&\det\left(\sum_0^{\iy}\frac{\vec\mu_j^{\top}}{z^j},
\sum_0^{\iy}\frac{\vec\mu_{j+1}^{\top}}{z^j},...,\sum_0^{\iy}
\frac{\vec\mu_{j+n}^{\top}}{z^j}\right)\\
\\
&=&\det\left(\vec\mu_0^{\top},\vec\mu_1^{\top},\ldots,\vec\mu_{n-1}^{\top},
\sum_0^{\iy}\frac{\vec\mu^{\top}_{j+n}}{z^j}\right)\\
\\
&=&z^n\det\left(\vec\mu^{\top}_0,\vec\mu^{\top}_1,\ldots,\vec\mu^{\top}_{n-1},
\sum_0^{\iy}\frac{\vec\mu^{\top}_j}{z^j}\right)
\eean\bean&=&z^n \det\left(\begin{array}{c|c}
 &\displaystyle{\int_{\BR}\sum^{\iy}_{j=0}\left(\frac{u}{z}\right)^j}\rho(u)du\\
 &\displaystyle{\int_{\BR}\sum^{\iy}_{j=0}\left(\frac{u}{z}\right)^j}u\rho(u)du\\
\vec\mu^{\top}_0,\vec\mu^{\top}_1,\ldots,\vec\mu^{\top}_{n-1}& \\
 &\vdots\\
&\displaystyle{\int_{\BR}\sum^{\iy}_{j=0}\left(\frac{u}{z}\right)^j}u^n\rho(u)du
\end{array}\right)\\
&=&z^n\int_{\BR}\sum_{j=0}^{\iy}\left(\frac{u}{z}\right)^j\rho(u)du\det
\left(\begin{array}{c|c}
 &1\\
 &u\\
\vec\mu_0^{\top},\ldots,\vec\mu^{\top}_{n-1} & \\
 & \vdots\\
  &u^n
\end{array}\right)\\
\\
&=&z^n\tau_n\int_{\BR}\sum_{j=0}^{\iy}\left(\frac{u}{z}\right)^j\tilde
p_n(u)\rho(u)du=z^{n+1}\tau_n \int_{\BR}\frac{\tilde
p_n(u)\rho_{\mu}(u)}{z-u}du \eean

\remark Representation (\ref{2b}) for orthogonal
polynomials $p_n$ can also be deduced from Heine's
representation. However, representation (\ref{2b}) is
much simpler.


  \noindent{\bf (4)} {\em The vectors $p$ and $q$ are
  eigenvectors of the tridiagonal symmetric matrix}
\be L:=S \Lb S ^{-1}.\label{L} \ee
%
Conjugating the shift matrix $\Lambda$ by $S$ yields a
matrix \bean L &=&S \Lb
S ^{-1}\\ &=&S \Lb S^{-1}S^{\top -1}S^{\top}\\
&=&S\Lb m_{\iy}S^{\top},\mbox{\,\,using (\ref{factorization})},\\
&=&Sm_{\iy}\Lb^{\top}S^{\top},\mbox{\,\,using $\Lb
m_{\iy}=m_{\iy}\Lb^{\top}$},\\ &=&S(S^{-1}S^{\top
-1})\Lb^{\top} S^{\top},\mbox{\,\,using again
(\ref{factorization})},\\ &=&(S\Lb S^{-1})^{\top}=L
^{\top}, \eean which is symmetric and thus tridiagonal.
Remembering $\chi(z)=(1,z,z^2,...)^{\top}$, and the
shift $(\Lambda v)_n=v_{n+1}$, we have
 $$ \Lb \chi(z)=z\chi(z)~~\mbox{and}~~
  \Lb^{\top}
  \chi(z^{-1})=z\chi(z^{-1})-ze_1,~~\mbox{with}~e_1=(1,0,0,...)^{\top}.
  $$
  Therefore, $p(z)=S\chi(z)$ and $q(z)=S^{\top
  -1}\chi(z^{-1})$ are eigenvectors, in the sense
  \bean
  Lp&=&S\Lb S^{-1}S\chi(z)=zS\chi(z)=zp  \\
  L^{\top}q&=& S^{\top -1}\Lb^{\top} S^{\top}
   S^{\top -1} \chi(z^{-1})=zS^{\top -1}\chi(z^{-1})
   -zS^{\top -1}e_1=zq-zS^{\top -1}e_1 .
   \eean
Then, using $L=L^{\top}$, one is led to
$$((L-zI)p)_n=0,~~\mbox{for}~ n\geq 0~~~\mbox{and}~~~
 ((L-zI)q)_n=0,~~\mbox{for}~~ n\geq 1.$$

\vspace*{.5cm}

\noindent{\bf(5)} {\em The off-diagonal elements of the
symmetric tridiagonal matrix $L$ are given by }\be
L_{n-1,n}=\sqrt{\frac{h_n}{h_{n-1}}}.\label{5a} \ee
Since $\la \tilde p_n,\tilde p_n\ra =h_n$, one has $
p_n(y)
=\frac{1}{\sqrt{h_n}}\tilde p_n(y). $ From the three
step relation $Lp(y)=yp(y)$, it follows that  \bean
\left(\frac{1}{\sqrt{h_{n-1}}}y^{n}+\ldots\right)
&=&yp_{n-1}(y)=L_{n-1,n}p_n(y)+(\mbox{terms of degree
$\leq n-1$})\\
&=&L_{n-1,n}\left(\frac{1}{\sqrt{h_{n}}}y^{n}+\ldots\right)
,\eean leading to statement (\ref{5a}).

\subsection{Time dependent Orthogonal polynomials and the KP equation}

Introduce now into the weight $\rho(z)$ a dependence on
parameters $t=(t_1,t_2,\ldots)$, as follows
\be
 \rho_t(z):= \rho(z)e^{\sum_{1}^{\iy}t_iz^i}.
  \label{rhot}\ee
  Consider the moment matrix $m_n(t)$, as in
  (\ref{5.2}), but now dependent on $t$, and
  the factorization of $m_{\iy}$ into lower- times
upper-triangular $t$-dependent matrices, as in
(\ref{factorization})
 \bea
  m_{\iy}(t) =(\mu_{ij}(t) )_{0\leq
i,j<\iy}=S ^{-1}(t)S^{\top -1}(t).
 \label{t-factorization}\eea
 The Toda lattice mentioned in the Theorem below will
 require the following Lie algebra splitting
 \be{\frak g\frak l}(n)={\frak s}  \oplus  \frak b,
 \label{Lie-splitting}\ee
into skew-symmetric matrices and
 (lower) Borel matrices.

 Also, one needs in this section the Hirota symbol:
  given a polynomial
 $p(t_1,t_2,...)$ of a finite or infinite
 number of variables and functions $f(t_1,t_2,\ldots)$
  and $g(t_1,t_2,\ldots)$,
 also depending on a finite or infinite number of variables
 $t_i$, define the
symbol \be
 p\left(\frac{\pl}{\pl t_1},\frac{\pl}{\pl
t_2},\ldots\right)f\circ g:= p\left(\frac{\pl}{\pl
y_1},\frac{\pl}{\pl y_2},...\right)f(t+y)g(t-y)
\Bigl|_{y=0}.\label{HirotaSymbol} \ee
The reader is reminded of the elementary Schur
polynomials $e^{\sum^{\iy}_{1}t_iz^i}:=\sum_{i\geq 0}
{\bf s}_i(t)z^i$ and for later use, set for
$\ell=0,1,2,\ldots$, \be{\bf s}_{\ell}(\tilde \pl):={\bf
s}_{\ell}(\frac{\pl}{\pl t_1},\frac{1}{2}\frac{\pl}{\pl
t_2},\ldots).\label{Schur-Hirota}\ee
One also needs {\em Taylor's formula} for a ${\cal
C}^{\iy}$-function $f$:
\be
 f(z+y)=e^{ y \frac{\pl}{\pl z }}f(z)
 , \label{Taylor}\ee
 which is seen by expanding the exponential.
%
%
The following Lemma will also be used later in the proof
of the bilinear relations:
\begin{lemma} \cite{AvM-3,AvM0}\label{Lemma: residue}
If $\oint_{\iy}$ denotes the integral along a small
circle about $\iy$, the following identity holds (formal
identity in terms of power series):\be \int_{\BR}
f(u)g(u)du =\frac{1}{2\pi i}\oint_{\iy} dz
f(z)\int_{\BR}\frac{g(u)}{z-u}du, \ee for holomorphic
$f(z)=\sum_{i\geq 0} a_iz^i$ and $g(z)$, the latter
assumed to have all its moments.\end{lemma}

\proof For holomorphic functions $f$ in $\BC$, \bea
\frac{1}{2\pi i}\oint_{\iy} dz
f(z)\int_{\BR}\frac{g(u)}{z-u}du &=&\Res_{z=\iy}
(\sum_{i\geq 0} a_i z^i)
 \Bigl(\frac{1}{z}\sum_{j \geq 0} z^{-j}\int_{\BR} g(u) u^j du\Bigr)
\nonumber\\
&=& \sum_{i\geq 0} a_i\int_{\BR} g(u)u^i du
 \nonumber\\
 &=& \int_{\BR} g(u)\sum_{i\geq 0} a_i u^i du
 \nonumber\\
 &=& \la f,g\ra,
 \eea
ending the proof of Lemma \ref{Lemma: residue}.\qed

 \bigbreak

 The next Theorem shows that the determinant of the {\em time}-dependent
 moment
 matrices, satisfies the {\em KP hierarchy}, a non-linear hierarchy,
 whereas in the
 next section, it will be shown that these same determinants
 satisfy {\em Virasoro equations}. These two features
 will play an important role in random matrix theory.
 Notice that this result is very robust: it can be
 generalized
 from orthogonal polynomials to multiple orthogonal
 polynomials, from the KP hierarchy to multi-component
 KP hierarchies; see \cite{AvM-Mops}.


\begin{theorem} \cite{AvM-1,AvM-3} \label{Theo: 5.2}The determinants
of the moment matrices, also representable as a multiple
integral\footnote{$\Dt_n(z)=\prod_{1\leq i<j\leq
n}(z_i-z_j)$},
 \be
  \tau_n(t): =\det m_n(t)
=\frac{1}{n!}\int_{E^n}\Delta^{2}_n(z)
\prod^n_{k=1}\rho (z_k) e^{\sum_{i=1}^{\iy}t_iz_k^i}dz_k 
\label{tau-function}\ee
  satisfy
\begin{description}

\item [(i)]  \underline{Eigenvectors of $L$}: The
tridiagonal matrix $L(t)$ admits two independent
eigenvectors: \bean (L(t)p(t;z))_n&=&zp_n(t;z),~n\geq
0\\
 (L(t)q(t;z))_n&=&zq_n(t;z)
 ,~~n\geq 1
. \eean

\begin{itemize}
  \item 
 $p_n(t;z)$
are $n$th degree polynomials  in $z$, depending on $t
\in \BC^{\iy}$, orthonormal with respect to $ \rho_t(z)$
(defined in (\ref{rhot})),
 and enjoying the following representations:
 (define $\chi(z):=(1,z,z^2,...)^{\top}$)
 \be  
p_n(t;z) := (S(t)\chi(z))_n
 =z^nh_{n}^{-1/2}\frac{\tau_n(t-[z^{-1}])}{\tau_{n}
(t)}, ~~~~h_{n}:=\frac{\tau_{n+1}(t)}{\tau_{n}(t)}
\ee 
 \item
$
 q_n(t,z), ~{n\geq 0}$, are Stieltjes transforms of the
 polynomials $p_n(t;z)$
and have the following $\tau$-function representations:
 \bea
q_n(t;z):=z\int_{\BR^n}\frac{p_n(t;u)}{z-u}\rho_t(u)du
&=&\left(S^{\top -1}(t)\chi(z^{-1})\right)_n
\nonumber\\
&=&z^{-n}h_n^{-1/2}
\frac{\tau_{n+1}(t+[z^{-1}])}{\tau_n(t)}.\nonumber\\
\eea
 \end{itemize}


\item [(ii)] The \underline{standard Toda lattice},
i.e., the symmetric tridiagonal matrix
 \be
  L(t):=S(t)\Lb S(t)^{-1}=
\left(\begin{tabular}{lllll} $\frac{\pl}{\pl t_1}\log
\frac{\tau_1}{\tau_0}$ &
 $ \left(\frac{\tau_{0}\tau_{2}}{\tau_{1}^2}\right)^{1/2}$
 & ~~ $0$      &    \\
$\left(\frac{\tau_{0}\tau_{2}}{\tau_{1}^2}\right)^{1/2}$&
$\frac{\pl}{\pl t_1}\log \frac{\tau_2}{\tau_1}$ &
$\left(\frac{\tau_{1}\tau_{3}}{\tau_{2}^2}\right)^{1/2}$&
\\~~~$0$&$\left(\frac{\tau_{1}\tau_{3}}{\tau_{2}^2}\right)^{1/2}$
&
 $\frac{\pl}{\pl t_1}\log \frac{\tau_3}{\tau_2}$&
   \\
 & &  &      $\ddots$\\
\end{tabular}
\right)
 \label{L-matrix}\ee
satisfies the commuting equations\footnote{in terms of
the Lie algebra splitting (\ref{Lie-splitting}).}
 \be\frac{\pl L}{\pl
t_k}=\left[\frac{1}{2}(L^k)_{\frak s},L\right]
 =-\left[\frac{1}{2}(L^k)_{\frak b},L\right].
 \ee

\item[(iii)] The functions $\tau_n(t)$ satisfy the
following {\em bilinear} identity, for $n\geq m+1$, and
all $t,t' \in \BC^{\iy}$, where one integrates along a
small circle about $\iy$, \be
\oint_{z=\iy}\tau_n(t-[z^{-1}])\tau_{m+1}(t'+[z^{-1}])e^{\sum_1^{\iy}
(t_i-t'_i)z^i}z^{n-m-1}dz=0. \label{bilinear
identity}\ee

 \item[(iv)] The \underline{KP-hierarchy}
 \footnote{Remember the Hirota symbol (\ref{HirotaSymbol})
 and the Schur polynomial notation (\ref{Schur-Hirota}).
 } 
 for $k=0,1,2,\ldots$ and for all $n=1,2\ldots$,
  $$
\left({\bf s} _{k+4}\bigl(\frac{\pl}{\pl
t_1},\frac{1}{2}\frac{\pl}{\pl
t_2},\frac{1}{3}\frac{\pl}{\pl
t_3},\ldots\bigr)-\frac{1}{2}\frac{\pl^2}{\pl t_1\pl
t_{k+3}}\right)\tau_n \circ\tau_n=0, $$
 of
which the first equation reads: \be \hspace{-1cm}
\left(\left(\frac{\pl}{\pl t_1}
\right)^4+3\left(\frac{\pl}{\pl
t_2}\right)^2-4\frac{\pl^2}{\pl t_1 \pl
t_3}\right)\log\tau_n+6\left(\frac{\pl^2}{\pl
t^2_1}\log\tau_n \right)^2=0. \ee


\end{description}

\end{theorem}

\remark In order to connect with classical integrable
theory, notice that, when $\tau$ satisfies the equation
above, the function
$$q(t_1,t_2,\ldots):=2\frac{\pl^2 \log \tau_n}{\pl t_1^2}$$
 satisfies the classical Kadomtsev-Petviashvili equation (KP
 equation):
\be
 3\frac{\pl^2 q}{\pl t_2^2}-
\frac{\pl}{\pl t_1}\left(
  4\frac{\pl q}{\pl t_3}-
   \frac{\pl^3 q}{\pl t_1^3}
   -6q \frac{\pl q}{\pl t_1}\right)=0.
   \ee
If $q$ happens to be independent of $t_2$, then $q$
satisfies the Korteweg-de Vries equation
 \be
  4\frac{\pl q}{\pl t_3}=
   \frac{\pl^3 q}{\pl t_1^3}
   +6q \frac{\pl q}{\pl t_1}.
  \ee

\vspace*{.6cm}

\proof Identity (\ref{tau-function}) follows from the
general fact that the product of two matrices can be
expressed as a symmetric sum of determinants\footnote
{Indeed, $
 \sum_{\sigma \in S_n}\det \left(
a_{i,\sigma(j)}~b_{j,\sigma(j)}\right)_{1\leq i,j\leq n}
=
 \det \left( a_{ik}\right)_{1\leq i,k\leq n}
  \det \left( b_{ik}\right)_{1\leq i,k\leq n}
  .$},
in particular the square of a Vandermonde can be
expressed as a sum of determinants:
  $$
\Delta^2(u_1,\ldots,u_n)=\sum_{\sg\in
S_n}\det\left(u_{\sg(k)}^{\ell +k-2}\right)_{1\leq
k,\ell\leq n}. $$
  Indeed,
 \bean n! \tau_n(t)&=& n! \det m_n(t)\\
 &=& \sum_{\sigma \in S_n}\det
 \left( \int_E z_{\sigma(k)}^{\ell +k-2} \rho_t(z_{\sigma(k)})
 dz_{\sigma(k)}  \right)_{1\leq
k,\ell\leq n}\\
 &=&   \sum_{\sigma \in S_n} \int_{E^n}  \det
 \left(z_{\sigma(k)}^{\ell +k-2}   \right)_{1\leq
k,\ell\leq n}\rho_t(z_{\sigma(k)})
 dz_{\sigma(k)}\\
 &=& \int_{E^n}\Delta^{2}_n(z)
\prod^n_{k=1}\rho_t(z_k)dz_k.
  \eean

\noindent{\bf (i)} At first, note
\bean \mu_{ij}(t\mp[z^{-1}])&=&\int_{\BR}
u^{i+j}e^{\sum_1^{\iy}\left(t_i \mp
\frac{z^{-i}}{i}\right)u^i}\rho(u)du\\
&=&
 \int_{\BR}
u^{i+j}\left(1-\frac{u}{z}\right)^{\pm
1}\rho(u)e^{\sum_0^{\iy} t_iu^i}du
 \\  \\
  &=&
  \left\{
   \begin{array}{l}
   \mu_{i,j}(t)-\frac{1}{z}\mu_{i,j+1}(t) \\  \\
   \mu_{i,j}(t)+\frac{1}{z}\mu_{i,j+1}(t)
    +\frac{1}{z^2}\mu_{i,j+2}(t)+\ldots,
   \end{array}\right.
 \eean which by formula (\ref{2b}) of section \ref{sect:
Orthogonal polynomials} leads at once to the following
representation for the monic orthogonal polynomials
$\tilde p_n(t;z)$ and their Stieltjes transforms,
\bea
 \tilde p_n(t;z)  &=&
 z^n\frac{\tau_n(t-[z^{-1}])}
         {\tau_n(t)}
 \no \\
   z\int_{\BR} \frac{ \tilde p_n( t;u)\rho_t (u)}{z-u}du
   &=&
 z^{-n }\frac{\tau_{n+1}(t+[z^{-1}])}
         {\tau_n(t)}
.\label{ortho+stieltjes}\eea
{\bf (ii)} The matrix $L:=S\Lb S^{-1}$ satisfies the
standard Toda lattice. One computes
 $$
 \frac{\pl \mu_{ij}}{\pl t_k}  =
  \mu_{i+k,j}~~\mbox{implying}~~~
  \frac{\pl m_{\iy}}{\pl t_k}=\Lb^km_{\iy}
.$$
Then, using the factorization (\ref{t-factorization})
and the definition (\ref{L}) of $L=S\Lb S^{-1}$ of
section \ref{sect: Orthogonal polynomials}, one computes
  \bean
0&=&S\left(\Lb^km_{\iy}-\frac{\pl m_{\iy}}{\pl
t_k}\right)S^{\top}\\ &=&S\Lb^kS^{-1}-S\frac{\pl}{\pl
t_k}(S^{-1}S^{\top -1})S^{\top}\\ &=&L^k+\frac{\pl
S}{\pl t_k} S^{-1}+S^{\top -1}\frac{\pl S^{\top}}{\pl
t_k}. \eean
 Upon taking the $()_-$ and $()_0$ parts of this equation
  ($A_-$ means the lower-triangular part of the
matrix $A$, including the diagonal and $A_0$ the
diagonal part) leads to
$$
  (L^k)_-+\frac{\pl S}{\pl
t_k} S^{-1}+\left(S^{\top -1}\frac{\pl S^{\top}}{\pl
t_k}\right)_0=0~\mbox{and}~\left(\frac{\pl S}{\pl t_k}
S^{-1}\right)_0=-\frac{1}{2}(L^k)_0.
 $$
  Upon observing that for any
symmetric matrix the following holds,
  $$
 \left(
\begin{array}{ll}
 a&c\\
  c&b
\end{array}
\right)_{\frak b}=
 \left(
\begin{array}{ll}
 a&0\\
  2c&b
\end{array}
\right)
 =2
 \left(
\begin{array}{ll}
 a&c\\
  c&b
\end{array}
\right)_- - \left(
\begin{array}{ll}
 a&c\\
  c&b
\end{array}
\right)_0,$$
 it follows that
 the matrices $L(t)$, $S(t)$ 
and the vector $p(t;z)=(p_n(t;z))_{n\geq 0}=S(t)\chi(z)$
satisfy the (commuting) differential equations and the
eigenvalue problem
 \be \frac{\pl S}{\pl
 t_k}=-\frac{1}{2}(L^k)_{\frak
 b}S,~~~~L(t)p(t;z)=zp(t;z),
 \ee
and thus the tridiagonal matrix $L$ satisfies the
 standard Toda lattice equations
 \bean
 \frac{\pl L}{\pl
t_k}=
 \frac{\pl }{\pl
t_k} S\Lb S^{-1}&=&
 \frac{\pl S}{\pl
t_k}S^{-1}S\Lb S^{-1}-S\Lb S^{-1}\frac{\pl S}{\pl t_k}
S^{-1}\\
&=&-\left[\frac{1}{2}(L^k)_{\frak b},L\right]
,\eean
 with $p(t;z)$ satisfying
 $$
 \frac{\pl p}{\pl t_k}
 =
  \frac{\pl S}{\pl t_k}\chi(z)
  =-\frac{1}{2}(L^k)_{\frak
 b}S\chi(z) =-\frac{1}{2}(L^k)_{\frak b} p.
$$
The two leading terms of $p_n(t;z)$ look as follows,
upon using (\ref{2}) and (\ref{ortho+stieltjes}) :
 \bea
p_n(t;z)&=&\sqrt{\frac{\tau_n}{\tau_{n+1}}}\tilde p_n(t;z)\\
\no\
&=&z^n\frac{\tau_n(t-[z^{-1}])}{\sqrt{\tau_n\tau_{n+1}}}\\
\no\\
&=&\sqrt{\frac{\tau_n}{\tau_{n+1}}}z^n\left(1-z^{-1}
 \frac{\pl\tau_n/\pl t_1}{\tau_n}+
\ldots\right). \label{6a}\eea
Thus, $z^n$ admits the following representation in terms
of the orthonormal polynomials $p_i$:
\bea  z^n&=& \sqrt{\frac{\tau_{n+1}}{\tau_{n}}}\left(
p_n+\frac{\pl\tau_n/\pl t_1}{\sqrt{\tau_{n}\tau_{n+1}}}
z^{n-1}+O(z^{n-2})\right)\no\\
 &=&\sqrt{\frac{\tau_{n+1}}{\tau_{n}}}p_n+
 \frac{\pl\tau_n/\pl t_1}{\sqrt{\tau_{n-1}\tau_{n}}}
p_{n-1}+O(z^{n-2})
 .\label{6b}\eea
Then, using (\ref{6a}) in $zp_n$ and then using the
representation (\ref{6b}) for $z^n$ and $z^{n+1}$, one
checks that the diagonal entries $b_n$ and non-diagonal
entries $a_n$ of $L$ are given by
 \bean
 b_n=\la zp_n,p_n\ra
 &=&\sqrt{\frac{\tau_n}{\tau_{n+1}}}\left(\la
z^{n+1},p_n\ra -\la z^n,p_n\ra\frac{\pl\tau_n/\pl
t_1}{\tau_n}\right)
 \\ \\
  &=&\frac{\pl\tau_{n+1}/\pl t_1}{\tau_n}-\frac{\pl\tau_n/\pl t_1}{\tau_n}\\
\\
&=&\frac{\pl}{\pl t_1}\log\frac{\tau_{n+1}}{\tau_n}
\eean and \bean
 a_n =
  \la zp_n,p_{n+1}\ra
 =\sqrt{\frac{\tau_n}{\tau_{n+1}}}
 \la z^{n+1}+\ldots,p_{n+1}\ra
 =\sqrt{\frac{\tau_n\tau_{n+2}}{\tau^2_{n+1}}}, \eean
establishing (\ref{L-matrix}).

\medbreak

\noindent{\bf (iii)}  {\em The bilinear identity:} The
functions $\tau_n(t)$ satisfy the following identity,
for $n\geq m+1, ~ t,t' \in \BC^{\iy}$, where one
integrates along a small circle about $\iy$, \be
\oint_{z=\iy}\tau_n(t-[z^{-1}])\tau_{m+1}(t'+[z^{-1}])e^{\sum_1^{\iy}
(t_i-t'_i)z^i}z^{n-m-1}dz=0. \label{bilinear
identity1}\ee
%
%
Indeed, using the $\tau$-function representation for the
monic orthogonal polynomials and their Stieltjes
transform (\ref{ortho+stieltjes}), one checks:
 \bean
  \lefteqn{
 \frac{1}{\tau_n(t)\tau_m(t')}
   \oint_{z=\iy}\tau_n(t-[z^{-1}])\tau_{m+1}(t'+[z^{-1}])
    e^{\sum_1^{\iy}
(t_i-t'_i)z^i}z^{n-m-1}dz }\\
&=&
 \oint_{z=\iy}z^n\frac{\tau_n(t-[z^{-1}])}{\tau_n(t)}
  z^{-m}\frac{\tau_{m+1}(t'+[z^{-1}])}{\tau_m(t')}
   e^{\sum_1^{\iy}
(t_i-t'_i)z^i} \frac{dz}{z}
\\
&=& \oint_{z=\iy}~
 dz ~e^{\sum_1^{\iy}(t_i-t'_i)z^i} \tilde p_n(t;z) \int_{\BR}\frac{\tilde
p_m(t';u)}{z-u}e^{\sum_1^{\iy}t'_iu^i}\rho(u)du
 \\
 &=&
2\pi i \int_{\BR} e^{\sum_1^{\iy}(t_i-t'_i)z^i}\tilde
p_n(t;z) \tilde p_m(t';z)e^{\sum_1^{\iy}
t'_iz^i}\rho(z)dz, \mbox{\,using
Lemma \ref{Lemma: residue},}\nonumber\\
&=&2\pi i\int_{\BR}\tilde p_n(t;z)\tilde
p_m(t';z)e^{\sum_1^{\iy} t_iz^i}\rho(z)dz
=0,\mbox{\,\,when $m\leq n-1$,}
 \eean
by orthogonality, establishing (\ref{bilinear
identity1}).

\noindent{\bf (iv)}  {\em The KP hierarchy:} Setting
$n=m+1$ in (\ref{bilinear identity}), shifting $t\mapsto
t-y,t'\mapsto t+y$, evaluating the residue, Taylor
expanding in $y_k$ (see (\ref{Taylor})) and using the
notation
$$\tilde \pl =\left(\frac{\pl}{\pl
t_1},\frac{1}{2}\frac{\pl}{\pl
t_2},\frac{1}{3}\frac{\pl}{\pl t_3},\ldots\right)
 ,$$
one computes the following residue about $z=\iy$,
%
%
\bean 0&=&
 \frac{1}{2\pi i} \oint_{z=\iy} dz~ e^{-\sum_1^{\iy} 2y_i
z^i}\tau_n(t-y-[z^{-1}])
 \tau_n(t+y+[z^{-1}])\\
 &=&
 \frac{1}{2\pi i} \oint dz~ e^{-\sum_1^{\iy} 2y_i
z^i}
 e^{\sum_1^\iy  \frac{z^{-i}}{i}\frac{\pl}{\pl u_i}}
 e^{\sum_1^\iy y_k \frac{\pl}{\pl u_k}}
 \tau_n(t-u)\tau_n(t+u)\Bigr|_{u=0}\\
 &=&
 \frac{1}{2\pi i} \oint dz~ e^{-\sum_1^{\iy} 2y_i
z^i}
 e^{\sum_1^\iy  \frac{z^{-i}}{i}\frac{\pl}{\pl t_i}}
 e^{\sum_1^\iy y_k \frac{\pl}{\pl t_k}}
 \tau_n(t)\circ\tau_n(t)
%
\eean\bean &=&\frac{1}{2\pi i}\oint dz
~\left(\sum_0^{\iy} z^{i}\gs_i(-2y)\right)
\left(\sum_0^{\iy} z^{-j}\gs_j(\tilde \pl)\right)
 e^{\sum_1^{\iy} y_k\frac{\pl}{\pl t_k}}\tau_n\circ\tau_n
 \\
%
&=& e^{\sum_1^{\iy} y_k\frac{\pl}{\pl t_k}}
 \sum_0^{\iy} {\bf s}_i(-2y){\bf s}_{i+1}(\tilde \pl)\tau_n\circ\tau_n\\
&=& \left(1+\sum_1^{\iy}y_j \frac{\pl}{\pl t_j}
 +O(y^2)\right)\left(\frac{\pl}{\pl t_1}+
 \sum_1^{\iy}{\bf s}_{i+1}(\tilde \pl)(-2y_i +O(y^2))
 \right) \tau_n\circ\tau_n\\
 &=&\left( \frac{\pl}{\pl t_1}+
 \sum_{1}^{\iy} y_{k}\left(\frac{\pl}{\pl t_{k}}
 \frac{\pl}{\pl t_1}
   -2{\bf s}_{k+1}(\tilde \pl)
 \right)\right)\tau_n
 \circ\tau_n+O(y^2),
\eean for arbitrary $y_k$, implying
 $$
 \frac{\pl}{\pl t_1}\tau\circ \tau=0
 ~~\mbox{and}~~ \left(\frac{\pl^2}{\pl t_{k}\pl t_1} -2{\bf s}_{k+1}(\tilde \pl)
 \right) \tau_n
 \circ\tau_n=0~~\mbox{for~} k=1,2,\ldots.
 $$
 Taking into account the fact that trivially
$(\pl / \pl t_1)\tau \circ \tau=0$ and that the equation
above is trivial for $k=1$ and $k=2$, one is led to the
KP hierarchy:
%
$$
\left({\bf s} _{k+4}\bigl(\frac{\pl}{\pl
t_1},\frac{1}{2}\frac{\pl}{\pl
t_2},\frac{1}{3}\frac{\pl}{\pl
t_3},\ldots\bigr)-\frac{1}{2}\frac{\pl^2}{\pl t_1\pl
t_{k+3}}\right)\tau_n \circ\tau_n=0, \mbox{~for} ~
 k=0,1,2,\ldots.$$
 In particular, for $k=0$, one computes
 $$\gs_4(t):=\frac{t_1^{4}}{4!}
  +\frac{1}{2}
  { t_2} {  t_1}%
 ^{2}+  {  t_3} {  t_1}+\frac{1}{2}  {t_2}^{2}
 + t_4,
 $$
leading to
the first equation in the hierarchy \be 
\left(\left(\frac{\pl}{\pl t_1}
\right)^4+3\left(\frac{\pl}{\pl
t_2}\right)^2-4\frac{\pl^2}{\pl t_1 \pl
t_3}\right)\log\tau_n+6\left(\frac{\pl^2}{\pl
t^2_1}\log\tau_n \right)^2=0. \ee This ends the proof of
Theorem \ref{Theo: 5.2}.\qed

\remark As mentioned earlier this method is very robust
and can be generalized to other integrals, besides
(\ref{tau-function}), upon using multiple orthogonal
polynomials. Such integrals with appropriate
 multiple time-deformations lead to $\tau$-functions for
multi-component KP hierachies; see Adler-van
Moerbeke-Vanhaecke \cite{AvM-Mops}.


 \newpage

 \newpage

\section{Virasoro constraints}

\subsection{Virasoro constraints for $\beta$-integrals}

Consider
 weights $\rho(z)dz=e^{-V(z)}dz$
 with rational logarithmic derivative and $E$ a disjoint
 union of intervals:
 $$
  -\frac{\rho'}{\rho}=V'(z)=\frac{g}{f}=
 \frac{\sum_0^{\iy} b_iz^i}{\sum_0^{\iy} a_iz^i}~~~\mbox{
    and     }~~~ E=\bigcup_1^{r}~[c_{2i-1},c_{2i}] \subset F
    \subseteq \BR,
 $$
  where $F=[A,B]$ is an interval such that
   $$
    \lim_{z\rightarrow
A,B}f(z)\rho(z)z^k=0\mbox{\,\,for all\,\,}k\geq 0.
 $$
 Consider an integral $I_n(t,c;\beta)$,
 generalizing (\ref{tau-function}), where $t:=(t_1,t_2,...)$ and $c=(c_1,c_2,...,c_{2r})$;
 namely with a Vandermonde\footnote{$\Dt_n(z)=\prod_{1\leq
i<j\leq n}(z_i-z_j)$.}
 to the power $2\beta>0$ instead of a square, and omitting the
 $n!$ appearing in (\ref{tau-function}),
%
 \be
I_n(t,c;\beta):=\int_{E^n}|\Dt_n(z)|^{2\beta}\prod_{k=1}^n
\left(e^{\sum_1^{\iy}t_i z_k^i}\rho(z_k)dz_k\right)
~\mbox{for} ~~n>0.\label{Vir:int}\ee Then the following
Theorem holds:


 \begin{theorem}\label{thm: Virasoro constraints} {\em (Adler-van Moerbeke \cite{AvM2})}
 ~~The multiple integrals $I_n(t):= \\ I_n(t,c;\beta)$
%
 with $I_0=1$, satisfy the following
 Virasoro constraints\footnote{When $E$ equals the whole range $F$, then the
  the first term, containing the partials
  with respect to the $c_i$'s, are absent in the
  formulae (6.1.3).}
  for all $k\geq -1$:
\be  
 \left(-\sum_1^{2r} c^{k+1}_if(c_i)\frac{\pl}{\pl c_i}
 + \sum_{i\geq
0}\left( a_i~  \BJ_{k+i,n}^{(2)}(t,n)-b_i ~
\BJ_{k+i+1,n}^{(1)}(t,n)\right) \right)I_n(t) = 0,
 \label{Vir:constraints}\ee
 where $\BJ^{(2)}_{k,n}(t,n)$ and $\BJ^{(1)}_{k,n}(t,n)$ are
 combined differential and multiplication (linear) operators.
 For all $n\in \BZ$, the operators
 $
  \BJ_{k,n}^{(2)}(t,n) $
and $
  \BJ_{k,n}^{(1)}(t,n)$ form a Virasoro and a
Heisenberg algebra respectively, interacting as follows
\bea \left[ \BJ_{k,n}^{(2)},~
 \BJ_{\ell,n}^{(2)} \right] &=&(k-\ell)~
 \BJ_{k+\ell,n}^{(2)} +c\left(
\frac{k^3-k}{12} \right)\dt_{k,-\ell}\nonumber\\
\left[ \BJ_{k,n}^{(2)},~ \BJ_{\ell,n}^{(1)} \right]
&=&-\ell ~ ~
 \BJ_{k+\ell,n}^{(1)}+c'k(k+1)\delta_{k,-\ell}.
\nonumber\\ \left[ \BJ_{k,n}^{(1)},~
 \BJ_{\ell,n}^{(1)}  \right]
&=&\frac{k}{2\beta}\delta_{k,-\ell}, \eea with ``central
charge"
 \be
c=1-6\left( \beta ^{1/2} - \beta ^{-1/2} \right)^2 ~~~
\mbox{and}~~~c'=\frac{1}{2}\left(\frac{1}{\beta}-1
\right). \ee

\end{theorem}

\medskip\noindent\underline{\it Remark 1:\/}
The operators $  \BJ_{k,n}^{(2)}=
\BJ_{k,n}^{(2)}(t,n)$'s are defined as follows: (the
normal ordering symbol ``$: ~:$" means: always pull
differentiation to the right, ignoring commutation
rules)
 \be
  \BJ_{k,n}^{(2)}= {\beta} \sum_{i+j=k}
:  \BJ_{i,n}^{(1)} ~  \BJ_{j,n}^{(1)}: +\left(1- {\beta}
\right)\left((k+1)~  \BJ_{k,n}^{(1)}
-k\BJ_{k,n}^{(0)}\right) ,\ee
 in terms of the $  \BJ_{k,n}^{(1)}=
\BJ_{k,n}^{(1)}(t,n)$'s. Componentwise, we have $$ ~
 \BJ_{k,n}^{(1)}(t,n)= ~
 J_k^{(1)}+nJ_k^{(0)}
 ~\mbox{and}~ ~  \BJ_{k,n}^{(0)}=nJ_k^{(0)}= n\dt_{0k}
$$
and
 \bean   \BJ_{k,n}^{(2)}(t,n)  =
{\beta}  J_k^{(2)} + \Bigl(2n\beta
+(k+1)(1-{\beta})\Bigr) ~~
 J_k^{(1)} +
n\Bigl( n  {\beta} +1-\beta\Bigr) J_k^{(0)},
 \eean
  where
  \bea
  J_k^{(0)}&=&\dt_{k0}\no\\
  J_k^{(1)}&=&\frac{\pl}{\pl
t_k}+\frac{1}{2\beta}(-k)t_{-k}
 \label{Vir:Vir}\\
   J^{(2)}_{k}&=&\sum_{i+j=k}\frac{\pl^2}{\pl
 t_{i}\pl t_{j}}+\frac{1}{\beta}\sum_{-i+j=k}it_{i}\frac{\pl}{\pl
 t_{j}}+\frac{1}{4\beta^2}\sum_{-i-j=k}it_{i}jt_{j}.
\nonumber \eea The integer $n$ appears explicitly in
 $  \BJ_{k,n}^{(2)}(t,n)$ to indicate
  the explicit $n$-dependence of the
 $n$th component, besides $t$.

%



\medskip\noindent\underline{\it Remark 2:\/}
 In the case $\beta=1$, the Virasoro generators
  (6.1.6) take on a particularly
  elegant form, namely for ${n \geq 0}$,
  \bean \BJ_{k,n}^{(2)}(t)&=&\sum_{i+j=k}
   :~ \BJ_{i,n}^{(1)}(t) ~
\BJ_{j,n}^{(1)}(t):~=~  J_k^{(2)}(t) + 2n J_k^{(1)}(t)
+ n^2\dt_{0k}\\
 \BJ_{k,n}^{(1)}(t) &=& ~  J_k^{(1)}(t)+
 n\dt_{0k} ,
  \eean
 with\footnote{The expression $J_k^{(1)}=0$ for $k=0$. }
\be J_k^{(1)}=\frac{\pl}{\pl
t_k}+\frac{1}{2}(-k)t_{-k}~,~
J^{(2)}_{k}=\sum_{i+j=k}\frac{\pl^2}{\pl
 t_{i}\pl t_{j}}+ \sum_{-i+j=k}it_{i}\frac{\pl}{\pl
 t_{j}}+\frac{1}{4}\sum_{-i-j=k}it_{i}jt_{j}.
 \ee

\bigbreak

One now establishes the following Lemma:

\begin{lemma}\label{lemma:Vir}Setting
$$
 dI_n(x):=|\Dt_n(x)|^{2\beta}\prod_{k=1}^n
\left(e^{\sum_1^{\iy}t_i x_k^i}\rho(x_k)dx_k\right),$$
the following variational formula holds:
 \be\left.\frac{d}{d\vr}dI_n (x_i\mapsto
x_i+\vr f(x_i)x_i^{k+1} )\right|_{\vr=0}
=\sum^{\iy}_{\ell=0}
 \left(a_{\ell} \BJ^{(2)}_{k+\ell,n}
-b _{\ell} \BJ^{(1)}_{k+\ell +1,n}\right)dI_n.
\label{Vir:diff}\ee
\end{lemma}

\proof Upon setting
 \be {\cal E}(x,t):=\prod_1^n e^{\sum_1^{\iy}t_i
x_k^i}\rho(x_k)
,\ee
 the following two relations hold:
 \bea
  \left(\frac{1}{2}\sum_{{i+j=k}\atop{i,j>0}}\frac{\pl^2}{\pl
t_i\pl t_j}-\frac{n}{2} \delta_{k,0}\right)\ER
&=&\left(\sum_{{1\leq a<\beta\leq n
}\atop{{i,j>0}\atop{i+j=k}}}x^i_{\al}x^j_{\beta}+
\frac{k-1}{2}\sum_{1\leq\al\leq n}x_{\al}^k\right)\ER,
\nonumber\\ \left(\frac{\pl}{\pl
t_k}+n\delta_{k,0}\right)\ER&=&\left( \sum_{1\leq\al\leq
n}x_{\al}^k\right)\ER,~~\mbox{all $k\geq 0$}.
\label{(6.1.11)}\eea
So, the point now is to compute the $\vr$-derivative \be
\frac{d}{d\vr}\left(|\Delta_n(x)|^{2\beta}
 e^{\sum^n_{k=1}(-V(x_k)+\sum^{\iy}
_{i=1}t_ix^i_k)} dx_1...dx_n\right)_{x_i\mapsto x_i+\vr
f(x_i)x_i^{k+1}}\Biggl|_{\vr=0}, \ee which consists of
three contributions:

\medbreak

\noindent\underline{Contribution 1}:

\medbreak

\noindent$\displaystyle{\frac{1}{2\beta|\Delta(x)|^{2\beta}}\frac{\pl}{\pl\vr}\left|\Delta(x+\vr
f(x)x^{k+1})\right|^{2\beta}\Biggl|_{\vr=0}}$ \bea
&=&\sum_{1\leq\al <\ga\leq
n}\frac{\pl}{\pl\vr}\log\left(|x_{\al}-x_{\ga}
+\vr(f(x_{\al})x_{\al}^{k+1}
-f(x_{\ga})x_{\ga}^{k+1})|\right)\Biggl|_{\vr=0}\nonumber\\
&=& \sum_{1\leq\al <\ga\leq
n}\frac{f(x_{\al})x_{\al}^{k+1}-f(x_{\ga})x_{\ga}^{k+1}}
{x_{\al}-x_{\ga}}\nonumber\\
&=& \sum_{\ell=0}^{\iy}a_{\ell} \sum_{1\leq\al <\ga\leq
n}
\frac{x_{\al}^{k+\ell+1}-x_{\ga}^{k+\ell+1}}{x_{\al}-x_{\ga}}
 \nonumber\\
&=& \sum_{\ell=0}^{\iy}a_{\ell}
 \left(\sum_{{i+j=\ell+k}\atop{{i,j>0}\atop {1\leq\al
<\ga\leq n}}}x^i_{\al}x^j_{\ga}+(n-1)\sum_{1\leq\al\leq
n}x_{\al}^{\ell+k}-\frac{n(n-1)}{2}\delta_{\ell+k,0}\right)\nonumber\\
&=&  \ER^{-1}|
 \sum_{\ell=0}^{\iy}a_{\ell}\Biggl(\frac{1}{2}
\sum_{{i+j=k+\ell}\atop{i,j>0}}\frac{\pl^2}{\pl t_i\pl
t_j}-\frac{n}{2}\delta_{k+\ell,0}\nonumber\\ &
&+\left(n\!-\!\frac{k\!+\!\ell\!+\!1}{2}\right)
\left(\frac{\pl}{\pl
t_{k+\ell}}+n\delta_{k+\ell,0}\right)-\frac{n(n\!-\!1)}{2}
\delta_{k+\ell,0}\Biggr)\ER, ~~\mbox{using (\ref{(6.1.11)}),}\nonumber\\
 &=&
\ER^{-1} \sum_{\ell=0}^{\iy}a_{\ell}
  \nonumber\\&&
 \Biggl(\frac{1}{2}
\sum_{{i+j=k+\ell}\atop{i,j>0}}\frac{\pl^2}{\pl t_i\pl
t_j}+\left(n-\frac{k+\ell+1}{2}\right) \frac{\pl}{\pl
t_{k+\ell}} +\frac{n(n-1)}{2}
\delta_{k+\ell,0}\Biggr)\ER.\nonumber\\ && \eea

 \noindent\underline{Contribution 2}: Using
 $f(x)=\sum_0^{\iy}
 a_ix^i$,

\medbreak

\noindent$\displaystyle{\frac{\pl}{\pl\vr}\prod^n_1d(x_{\al}+\vr
f(x_{\al})x_{\al}^{k+1})\Biggl|_{\vr=0}}$
 \bea
&=&\sum^n_1\left(f'(x_{\al})x_{\al}^{k+1}
+(k+1)f(x_{\al})x_{\al}^k\right)\prod^n_1dx_i\nonumber\\
&=&\sum^{\iy}_{\ell=0}(\ell+k+1)a_{\ell}
 \sum^n_{\al=1}x_{\al}^{k+\ell}
\prod^n_1dx_i\nonumber\\
&=&\ER^{-1}\sum^{\iy}_{\ell=0}(\ell+k+1)a_{\ell}
 \left(\frac{\pl}{\pl
t_{k+\ell}}+n\delta_{k+\ell,0}\right)\ER \prod^n_1
~dx_i.
\eea

\medbreak

\noindent\underline{Contribution 3}: again using
 $f(x)=\sum_0^{\iy}
 a_ix^i$,
 \bea
\lefteqn{\frac{\pl}{\pl\vr}\prod^n_{\al=1}\exp \left(-V
\left(x_{\al}+\vr
f(x_{\al})x^{k+1}\right)+\sum^{\iy}_{i=1}t_i
\sum^n_{\al=1}\left( x_{\al}+\vr f(x_{\al})
x_{\al}^{k+1}\right)^i\right)\Biggl|_{\vr=0}}\nonumber\\
&=&\left(-\sum^n_{\al=1}V'(x_{\al})f(x_{\al})x_{\al}^{k+1}+\sum^{\iy}_{i=1}
it_i\sum^n_{\al=1}f(x_{\al})x_{\al}^{i+k}\right)
\ER\nonumber\\
&=&\!\left(-\sum^{\iy}_{\ell=0}b_{\ell}\sum^n_{\al=1}
x_{\al}^{k+\ell+1}+\sum_{{\ell\geq 0}\atop{i\geq
1}}a_{\ell}it_i\sum_{\al=1}^nx_{\al}^{i+k+\ell}\right)
\ER\nonumber\\
&=&\!\Biggl(\!-\sum^{\iy}_{\ell=0}b_{\ell}\left(\frac{\pl}{\pl
t_{k+\ell+1}}\!+n\delta_{k+\ell+1,0}\right)
 \!+\!\sum^{\iy}_{\ell=0}a_{\ell}
\sum^{\iy}_{i=1}it_i\left(\frac{\pl}{\pl
t_{i+k+\ell}}\!+n\delta_{i+k+\ell,0}\right)\Biggr) \ER.
\nonumber\\
\eea
As mentioned, for knowing (\ref{Vir:diff}), we must add
up the three contributions 1, 2 and 3, resulting in:
\bea
 \lefteqn{\left.\frac{\pl}{\pl\vr}dI_n
 (x_i\mapsto
x_i+\vr f(x_i)x_i^{k+1} )
 \right|_{\vr=0}}\\
&=&\left(\sum^{\iy}_{\ell=0}a_{\ell}\left( {\beta}
J^{(2)}_{k+\ell} +(2n\beta +(\ell+k+1)(1- {\beta} ))
J^{(1)}_{k+\ell}+n(( n-1) {\beta}
+1)\delta_{k+\ell,0}\right)\right.\nonumber\\ &
&\quad\quad\quad\quad
\left.-\sum^{\iy}_{\ell=0}b_{\ell}\left(
J^{(1)}_{k+\ell+1}+n\delta_{k+\ell+1,0}\right)\right)
dI_n(x).\nonumber \eea
 Finally, the use of (\ref{Vir:Vir}) ends the
proof of Lemma \ref{lemma:Vir}.\qed

\medbreak

\underline{\sl Proof of Theorem \ref{thm: Virasoro
constraints}}: The change of integration variable $
x_i\mapsto x_i+\vr f(x_i)x_i^{k+1}
  $
   in the
integral (\ref{Vir:int}) leaves the integral invariant,
but it induces a change of limits of integration, given
by the inverse of the map above; namely the $c_i$'s in
$E={\bigcup^r_1}[c_{2i-1},c_{2i}]$, get mapped as
follows $$ c_i\mapsto c_i-\vr f(c_i)c_i^{k+1}+O(\vr^2).
$$ Therefore, setting
$$E^{\vr}=\displaystyle{\bigcup^r_1}[c_{2i-1}-\vr
f(c_{2i-1}) c_{2i-1}^{k+1}+O(\vr^2),c_{2i}-\vr
f(c_{2i})c_{2i}^{k+1} +O(\vr^2)],$$ we find, using Lemma
\ref{lemma:Vir} and the fundamental theorem of calculus,
\begin{eqnarray*}
0&=&\frac{\pl}{\pl\vr}\int_{(E^{\vr})^{2n}}|\Delta_{2n}(x+\vr
f(x)x^{k+1})|\prod^{2n}_{i=1}e^{-V(x_i+\vr
f(x_i)x_i^{k+1} ,t)}d(x_i+\vr f(x_i)x_i^{k+1})\\
&=&\left(-\sum^{2r}_{i=1}c_i^{k+1}f(c_i)\frac{\pl}{\pl
c_i} +\sum^{\iy}_{\ell=0}\left( a_{\ell} ~
\BJ^{(2)}_{k+\ell,2n}-b_{\ell}~ \BJ^{(1)}_{k+\ell+1,2n}
\right)\right)  I_{n}(t,c,\beta),
\end{eqnarray*}
ending the proof of Theorem \ref{thm: Virasoro
constraints}.\qed


\subsection{Examples}

These examples are taken from \cite{ASV,AvM-1,AvM1,vM};
for the Laguerre ensemble, see also \cite{Haine1} and
for the Jacobi ensemble, see \cite{Haine2}.

\subsubsection*{Example 1 (GUE)}
 Here we pick $$\rho(z)=e^{-V(z)}= e^{-z^2},
~~V'=g/f=2z,$$ $$ a_0=1,b_0=0,b_1=2,~\mbox{and all
other}~a_i,b_i=0 .$$
  Define the differential operators
  \be
   {\cal B}_{k}:=
\sum_1^{2r}c_i^{k+1}\frac{\pl}{\pl c_i}
 \label{B-operator}
 \ee
in terms of the end points of the set $
 E=\bigcup_1^{2r}~[c_{2i-1},c_{2i}] \subset \BR.
 $
  From Theorem \ref{thm: Virasoro
constraints}, the integrals \be I_n= \int_{E^n} \Dt_n(z)
^{2}\prod_{k=1}^n
e^{-z_k^2+\sum_{i=1}^{\iy}t_iz_k^i}dz_k
 \ee
satisfy the Virasoro constraints
 \be
  -{\cal B}_{k}I_n=
 \left(- \BJ_{k,n}^{(2)}+2
 \BJ_{k+2,n}^{(1)}\right)I_n,~~k=-1,0,1,... .
  \ee
 The first three constraints have the following
form, upon setting $F=\log I_n$ (this will turn out to
be more convenient in the applications),
 \bean
 -{\cal B}_{-1} F&=&
 \left(
2\frac{\pl}{\pl t_1}-\sum_{i\geq 2} it_i
\frac{\pl}{\pl t_{i-1}}  \right)F - nt_1 \\
 - {\cal B}_0 F  &=&
 \left(
2\frac{\pl}{\pl t_2}- \sum_{i\geq 1} it_i \frac{\pl}{\pl
t_i} \right)F - n^2
\\
 -{\cal B}_1 F &=&
 \left(
2\frac{\pl}{\pl t_3} -2n
\frac{\pl}{\pl t_1}- \sum_{i\geq 1} it_i \frac{\pl}{\pl
t_{i+1}} \right)F \eean
 For later use, take linear combinations
 such that each expression contains the pure
 differentiation term $\pl F/\pl t_i$,
 yielding
  \bea
   -\frac{1}{2}{\cal B}_{-1} F\!&=:&\!{\cal D}_{-1}F=\left(\frac{\pl}{\pl t_1}-\frac{1}{2}
\sum_{i\geq 2}it_i\frac{\pl}{\pl
t_{i-1}}\right)F-\frac{nt_1}{2}\nonumber\\
-\frac{1}{2}{\cal B}_0 F\!&=:&\!{\cal
D}_{0}F=\left(\frac{\pl}{\pl t_2}-\frac{1}{2}
\sum_{i\geq 1}it_i\frac{\pl}{\pl
t_i}\right)F-\frac{n^2}{2}
\nonumber\\ -\frac{1}{2}\left({\cal B}_1\!+\!
 n{\cal B}_{-1}\right)F
 \!&=:&\!{\cal D}_{1}F=\left(\frac{\pl}{\pl t_3}-\frac{1}{2} \sum_{i\geq
1}it_i\frac{\pl}{\pl
t_{i+1}}-\frac{n}{2}
\sum_{i\geq 2}it_i\frac{\pl}{\pl t_{i-1}}\right)F
-\!\frac{n^2t_1}{2}
 .\nonumber\\ \label{GUE1}
\eea

\subsubsection*{Example 2  (Laguerre ensemble)}
 Here, the
weight is
 $$e^{-V}=z^ae^{-z},~~V'=\frac{g}{f}=\frac{z-a}{z}, $$
 $$ a_0=0,~a_1=1,~ b_0=-a,~b_1=1,~\mbox{and all
other}~a_i,b_i=0 .$$
Here define the (slightly different) differential
operators
  \be
   {\cal B}_{k}:=
\sum_1^{2r}c_i^{k+2}\frac{\pl}{\pl c_i}
 \label{B-operatorLaguerre}
 .\ee
 Thus from Theorem \ref{thm: Virasoro constraints}, the
integrals
 \be
I_n= \int_{E^n} \Dt_n(z) ^{2}\prod_{k=1}^n
z_k^{a}e^{-z_k+\sum_{i=1}^{\iy}t_iz_k^i}dz_k
 \ee
satisfy the Virasoro constraints, for $k\geq -1$,
 \be
- {\cal B}_k I_n=
 \left(-  \BJ_{k+1,n}^{(2)} -a
 \BJ_{k+1,n}^{(1)}+ \BJ_{k+2,n}^{(1)}
 \right) I_n.\ee
 Written out, the first three have the form, again upon setting
  $F= \log
I_n$,
 \bean
- {\cal B}_{-1} F&=& \left(\frac{\pl}{\pl
t_1}-\sum_{i\geq 1} it_i \frac{\pl}{\pl t_{i}}
 \right) F -n(n+a)\\
- {\cal B}_0 F   &=&   \left(\frac{\pl}{\pl t_2} -
(2n+a)\frac{\pl}{\pl t_1}-\sum_{i\geq 1} it_i
 \frac{\pl}{\pl t_{i+1}} \right) F \\
  -{\cal B}_1 F  &=& \left(\frac{\pl}{\pl
t_3}-(2n+a)\frac{\pl}{\pl t_2}-\sum_{i\geq 1} it_i
\frac{\pl}{\pl t_{i+2}}-  \frac{\pl^2}{\pl t_1^2} \right) F 
 -
 \left(\frac{\pl F}{\pl t_1} \right)^2;
 \eean
Replacing the operators ${\cal B}_i$ by linear
combinations 
  \bea {\cal D}_1&=&-{\cal
 B}_{-1}\nonumber\\
  {\cal D}_2&=&-{\cal B}_0-(2n+a){\cal B}_{-1}\nonumber\\
  {\cal D}_3&=&-{\cal B}_1-
  (2n+a){\cal B}_0 -(2n+a)^2 {\cal B}_{-1}
   \label{6.2.8}\eea
  yields expressions, each containing a pure derivative
$\pl F/\pl t_i$

  \bea
{\cal D}_1F&=&\frac{\pl F}{\pl t_1}-\sum_{i\geq 1}
it_i\frac{\pl F}{\pl
t_i}-n(n+a)\nonumber\\
 {\cal D}_2F&=&\frac{\pl F}{\pl
t_2}+ \sum_{i\geq 1} it_i \left( -(2n+a)\frac{\pl}{\pl
t_i}-\frac{\pl}{\pl t_{i+1}}\right)F
 - {n} (n+a)(2n+a)
\nonumber\\
 {\cal D}_3F&=&
  \frac{\pl F}{\pl t_3}
  -\sum_{i\geq 1}it_i
\Biggl((2n+a)^2\frac{\pl}{\pl t_i} +(2n+a)\frac{\pl}{\pl
t_{i+1}}+\frac{\pl}{\pl t_{i+2}}\Biggr)F   \nonumber\\&&
 - \left(\frac{\pl^2F}{\pl t^2_1}+
\left(\frac{\pl F}{\pl t_1}\right)^2\right)- {n}
(n+a)(2n+a)^2.
 \label{6.2.9}\eea
Notice the non-linearity in this expression is due to
the fact that one uses $F=\log I_n$ rather than $I_n$.


\subsubsection*{Example 3  (Jacobi ensemble)}
 The
weight is given by
 $$
\rho_{ab}(z):=e^{-V}=(1-z)^{a}(1+z)^{b} ,~~
V'=\frac{g}{f}=\frac{a-b+(a+b)z}{1-z^2}
 $$
 $$
 a_0=1,a_1=0,a_2=-1,b_0=a-b,b_1=a+b
 ,~\mbox{and all
 other}~a_i,b_i=0 .$$
 Here define
 $${\cal B}_k :=
\sum_1^{2r}c_i^{k+1}(1-c_i^2)\frac{\pl}{\pl c_i}.$$
 The integrals
 \be
 \int_{E^n} \Dt_n(z) ^{2}\prod_{k=1}^n
(1-z_k)^{a}(1+z_k)^{b}e^{\sum_{i=1}^{\iy}t_iz_k^i}dz_k
 \ee
satisfy the Virasoro constraints $(k\geq-1)$:
 \bea  - {\cal B}_k I_n
=
 \left( \BJ_{k+2,n}^{(2)}- \BJ_{k,n}^{(2)}+ b_0
\BJ_{k+1,n}^{(1)} +b_1 \BJ_{k+2,n}^{(1)}\right)I_n. \eea
 Introducing the following notation
 $$\sigma=2n+b_1$$
 the first four having the following form, upon setting $F=\log
I_n$,
 \bea
 - {\cal B}_{-1} F&=& \left(\sigma
 \frac{\pl}{\pl
t_1}+\sum_{i\geq 1} it_i \frac{\pl}{\pl t_{i+1}}-
 \sum_{i\geq 2} it_i \frac{\pl}{\pl t_{i-1}}
 \right) F 
 +n(b_0-t_1) \nonumber \\ &&  \nonumber\\
 -{\cal B}_0 F   &=&   \left(
 \sigma  \frac{\pl}{\pl t_2}
+  b_0\frac{\pl}{\pl t_1}
 +\sum_{i\geq 1}it_i
 (\frac{\pl}{\pl t_{i+2}}-\frac{\pl}{\pl t_{i}})
   + \frac{\pl^2 }{\pl t_1^2}
   \right) F \nonumber\\
 &&+\left(\frac{\pl F}{\pl
 t_1}\right)^2-\frac{n}{2}(\sigma -b_1) \nonumber \\  &&
\nonumber\\
 - {\cal B}_1 F  &=& \left(\sigma
 \frac{\pl}{\pl t_3}+b_0 \frac{\pl}{\pl t_2}
 -(\sigma -b_1)\frac{\pl}{\pl t_1}
 +\sum_{i\geq 1} it_i (\frac{\pl}{\pl
t_{i+3}}
 -\frac{\pl}{\pl t_{i+1}})
 \right.   \nonumber \\&& \left.
  +2 \frac{\pl^2}{\pl t_{1}\pl t_2}\right)F
  +2 \frac{\pl F}{\pl t_{1}} \frac{\pl F}{\pl t_{2}}
 \nonumber\\
  - {\cal B}_2 F  &=& \left(\sigma
 \frac{\pl}{\pl t_4}+b_0 \frac{\pl}{\pl t_3}
 -(\sigma -b_1)\frac{\pl}{\pl t_2} +\sum_{i\geq 1} it_i (\frac{\pl}{\pl
t_{i+4}}
 -\frac{\pl}{\pl t_{i+2}})
 \right.   \nonumber \\
&&\hspace{-0.5cm} \left. +   \frac{\pl^2}{\pl t_{2}^2}-
  \frac{\pl^2}{\pl t_{1}^2}+
  2\frac{\pl^2}{\pl t_{1}\pl t_3} \right)F
   + \left(\Bigl( \frac{\pl F}{\pl t_{2}} \Bigr)^2-
  \Bigl( \frac{\pl F}{\pl t_{1}} \Bigr)^2
  +2\frac{\pl F}{\pl t_{1}}\frac{\pl F}{\pl t_{3}}\right)
.\nonumber\\ \eea


\section{Random matrices}\label{section 7}

This whole section \ref{section 7} is a very standard
chapter of random matrix theory. Most of the results can
be found e.g., in Mehta \cite{Mehta}, Deift \cite{Deift}
and others.

\subsection{Haar measure on the space $\HR_n$ of Hermitian matrices
}

Consider the most na\"ive measure (Haar measure)
\be
 dM:=\displaystyle{\prod_1^n dM_{ii}\prod_{1\leq
i<j\leq n}} d\Re e M_{ij}~d\Im m M_{ij} \label{Haar
measure} \ee
 on the space of Hermitian matrices
$$
 \HR_n:=
  \{n\times n ~\mbox{matrices such that} ~M^{\top}=\bar
  M\}
  .$$
The parameters $M_{ij}$ in (\ref{Haar measure}) are
precisely the free ones in $M$. This measure turns out
to be invariant under conjugation by unitary matrices:
(see Mehta \cite{Mehta} and more recently Deift
\cite{Deift})

\begin{proposition}\label{Proposition 7.1} For a fixed $U\in SU(n)$, the map

$$
\HR_n\rg\HR_n:M\longmapsto M'=U^{-1}MU
$$
has the property

$$
dM=dM',\quad\mbox{i.e.~} \Big\vert\det\left(\frac{\pl
M'}{\pl M}\right)\Big\vert =1.
$$
\end{proposition}

\proof Setting $M'=U^{-1}MU$, we have

$$
\Tr M^2=\Tr M^{'2}
$$
and so

$$
\sum_{i,j}M_{ij}M_{ji}=\sum_{i,j}M'_{ij}M'_{ji}.
$$
Working out this identity, one finds
\bea \lefteqn{\sum_1^nM^2_{ii}+2\sum_{1\leq i<j\leq
n}((\Re
e~M_{ij})^2+(\Im m~M_{ij})^2)}\nonumber\\
&=&\sum^n_1M^{'2}_{ii}+2\sum_{1\leq i<j\leq n} ((\Re
e~M'_{ij})^2+(\Im m~M'_{ij})^2). \label{7.1.1}\eea
Setting
$$
\vec M:=(M_{11},...,M_{nn},\Re e~M_{12},...,\Re
e~M_{n-1,n},\Im m~M_{12},..., \Im m~M_{n-i,n}),
$$
identity (\ref{7.1.1}) can be written, in terms of the
usual inner product $\la \, ,\, \ra$ in
$\BR^{\frac{n(n+1)}{2}}$,

$$
\la\vec M,D\vec M\ra =\la\vec{M}',D\vec{M}'\ra
$$
for the $n^2\times n^2$ diagonal matrix with $n$ 1's and
$n(n-1)$ 2's,
$$
D=\left(\begin{array}{ccccccc}
1& \\
 &\ddots& & & &{\bf O}& \\
 & & 1\\
& & &2\\
 &{\bf O}& & & &\ddots&  \\
 & & & & & &2
\end{array}\right).
$$
Let $V$ be the matrix transforming the vectors $M$ into
$M'$
$$
\vec{M}'=V\vec{M}
$$
and so, (\ref{7.1.1}) reads
$$
\la \vec{M'},D\vec{M'}\ra =\la V\vec{M},DV\vec{M}\ra
=\la \vec{M},V^{\top}DV\vec{M}\ra ,
$$
from which it follows that $D=V^{\top}DV$ and so

$$
0\neq\det D=\det(V^{\top}DV)=(\det V)^2\det D,
$$
implying
$$
|\det V|=1.
$$
But for a linear transformation $\vec{M}'=V\vec{M}$ the
Jacobian of the map is $V$ itself and so
$$
\Big\vert\det\frac{\pl \vec{M}'}{\pl M}\Big\vert =|\det
V|=1,
$$
ending the proof of Proposition \ref{Proposition 7.1}.
\qed

\begin{proposition}\label{Proposition 7.2}
The diagonalization $M=UzU^{-1}$ leads to ``polar" or
``spectral" coordinates $M\mapsto  (U,z)$, where
$z=\mbox{diag} (z_1,\ldots,z_n),\quad z_i\in\BR$. In
these new coordinates
\be dM=\Dt^2(z) dz_1\ldots dz_n~dU. \label{polar coord}
\ee
\end{proposition}
\proof Every matrix $M\in\HR_n$ can be diagonalized,
\bean M&=&UzU^{-1},
 \eean
with $U= e^Aze^{-A}\in SU(n)$ and
$$
A=\sum_{1\leq k\leq\ell\leq
n}(a_{k\ell}(e_{k\ell}-e_{\ell
k})+i~b_{k\ell}(e_{k\ell}+e_{\ell k}))\in
su(n),\mbox{~with~}a_{k,\ell},b_{k,\ell}\in\BR.
$$
Then, using the definition of the measure $dM$, as in
(\ref{Haar measure}), and using the fact that $[A,z]$ is
a Hermitian, one computes
\bean
dM\Big\vert_{M=z}=d(e^Aze^{-A})\Big\vert_{A=0}&=&d(z+[A,z]+{\bf
O}(A^2))\Big\vert_{A=0}\\
\\
&=&\prod^n_1dz_i\prod_{1\leq j<k\leq n}d~\Re
e[A,z]_{jk}~d~\Im
m[A,z]_{jk}\Big\vert_{A=0}\\
\\
&=&\prod^n_1dz_i\prod_{1\leq j<k\leq
n}(z_j-z_k)^2\prod_{1\leq j<k\leq
n}da_{jk}~db_{jk}\Big\vert_{A=0}. \eean
Since\footnote{$e_{k\ell}$ denotes the matrix having a
$1$ at the entry $(k,\ell)$ and $0$ everywhere else.}
\bean [A,z]&=&\sum_{1\leq k,\ell\leq
n}(a_{k\ell}[e_{k\ell}-e_{\ell
k},z]+i~b_{k\ell}[e_{k\ell}+e_{\ell k},z])\\
\\
&=&\sum_{1\leq k,\ell\leq n}(a_{k\ell}(z_{\ell}-z_k)(
e_{k\ell}+e_{\ell
k})+i~b_{k\ell}(z_{\ell}-z_k)(e_{k\ell}-e_{\ell k}))
,\eean
 with
$$
\Re e[A,z]_{k,\ell}=a_{k\ell}(z_{\ell}-z_k),\qquad\Im
m[A,z]_{k\ell}=b_{k\ell}(z_{\ell}-z_k),
$$
establishing (\ref{polar coord}) near $M=z$. By Lemma
\ref{Proposition 7.1},
 $
 dM=d(U'^{-1}MU')$ for any unitary matrix $U'$,
 implying the result (\ref{polar coord})
 holds everywhere, establishing Lemma \ref{Proposition 7.2}. \qed

\remark The set of Hermitian matrices is the tangent
space to a symmetric space $G/K=SL(n,\BC)/SU(n)$. The
argument in Theorem \ref{Proposition 7.2} can be
generalized to many other symmetric spaces, as worked
out in van Moerbeke \cite{vM}, page 324-329.
\newpage

\subsection{Random Hermitian ensemble}

Consider the probability distribution on the space
${\cal H}_n$ of Hermitian matrices, in terms of Haar
measure $dM$,
given by
$$
P(M\in dM)=\frac{1}{Z_n}e^{-\Tr V(M)}dM.
$$
Let $z_1\leq z_2\leq\ldots\leq z_n$ be the real
eigenvalues of $M$. Then \bean P(z_1\in
dz_1,\ldots,z_n\in dz_n)&=&P(z_1,\ldots ,z_n)dz_1\ldots
dz_n\\
&:=&\frac{1}{Z_n}\Dt^2(z)\prod ^n_{i=1}e^{-V(z_i)}dz_i
\eean with
$$
Z_n=\int_{z_1\leq\ldots\leq
z_n}P(z_1,\ldots,z_n)\prod_1^ndz_i.
$$

\begin{lemma}\label{Lemma 7.3}
\be  P(M\in{\cal H}_n,~\mbox{spectrum }(M)\subseteq E)
=\frac{\displaystyle{\int_{E^n}\Dt^2(z)\prod^n_1e^{-V(z_i)}}dz_i}{
\displaystyle{\int_{\BR^n}\Dt^2(z)\prod_1^ne^{-V(z_i)}}dz_i}
.\label{ratio integrals}\ee
\end{lemma}
\proof Indeed
\bean
\lefteqn{P(M\in{\cal H}_n,\mbox{spectrum }M\subseteq E)}\\
\\
 &=&\frac{1}{Z_n}\int_{z_{1}<\ldots
<z_{n}}\Dt_n^2\left(z_{1},\ldots,z_{n}\right)\prod^n_{i=1}
\chi_{_{E}}(z_{i})e^{-V(z_{i})}dz_{i}
 \\
 \\
  &=&\frac{1}{Z_nn!}\sum_{\pi\in
S^n}\int_{z_{\pi(1)}<\ldots
<z_{\pi(n)}}\Dt_n^2\left(z_{\pi(1)},\ldots,z_{\pi(n)}\right)\prod^n_{i=1}
\chi_{_{E}}(z_{\pi(i)})e^{-V(z_{\pi(i)})}dz_{\pi(i)}\\
\\
&=&\frac{1}{Z_nn!}\sum_{\pi\in
S^n}\int_{z_{\pi(1)}<\ldots
<z_{\pi(n)}}\Dt_n^2(z_1,\ldots,z_n)\prod^n_{i=1}
\chi_{_{E}}(z_i)e^{-V(z_i)}dz_i
\eean\bean \hspace*{-6cm}&=&\frac{\displaystyle{\int_{
E^n}\Dt_n^2(z)\prod^n_1e^{-V(z_i)}}dz_i}{
\displaystyle{\int_{\BR^n}\Dt_n^2(z)\prod_1^ne^{-V(z_i)}}dz_i}
,\eean
 showing Lemma \ref{Lemma 7.3}.\qed


Let $p_0(z),p_1(z),p_2(z),\ldots $ be orthonormal
polynomials with regard to the weight $\rho(z)$ defined
on $\BR$, as discussed in section \ref{section 5.1},
$$
 \int_{\BR}p_i(z)p_j(z)\rho(z)dz=\dt_{ij},
 $$
and $\tilde p_n(z)$ be the monic orthogonal polynomials,
$$
 \int_{\BR}\tilde p_i(z)\tilde p_j(z)\rho(z)dz=h_i\dt_{ij}.
 $$
 Then,
  \be
    p_n(z)=\frac{1}{\sqrt{h_n}}\left(z^n+\ldots\right)
   =\frac{1}{\sqrt{h_n}}\tilde p_n(z)\label{monic}
 .\ee

\begin{proposition} \label{Proposition 7.4}Setting
$$
Z_n=\int_{\BR^n}\Dt^2_n(z)\prod_1^n\rho(z_k)dz_k,
$$
we have the identity
 \be
Z_n^{-1}\Dt_n^2(z)\prod^n_1\rho(z_k)=\frac{1}{n!}\det(K_n(z_k,z_{\ell}))_{1\leq
k,\ell\leq n} ,\label{identity}\ee where the symmetric
kernel $K_n$ is given by
\bea
K_n(y,z)&=&\sqrt{\rho(y)\rho(z)}\sum^{n-1}_{j=0}p_j(y)p_j(z)
 \quad\mbox{(Christoffel-Darboux)}\no\\
\no\\
&=&\sqrt{\frac{h_n}{h_{n-1}}} \sqrt{\rho(y)\rho(z)}
 \frac{p_n(y)p_{n-1}(z)-p_{n-1}(y)
 p_n(z)}{y-z}.\label{Darboux kernel}
 \eea
The kernel $K_n(y,z)$ has the following ``reproducing"
property
$$
\int_{\BR}K_n(y,u)K_n(u,z)du=K_n(y,z)~~
\mbox{and}~~\int_{\BR}K_n(z,z)dz=n .$$
\end{proposition}

\proof Notice that the Vandermonde $\Dt_n(z)$ can also
be expressed as \\$\det(\tilde p_{i-1}(z_j))_{1\leq
i,j\leq n}$ by row operations, where $\tilde p_{i}(z)$
can be chosen to be any monic polynomial of degree $i$,
and in particular the monic orthogonal polynomial of
degree $i$. Thus, one computes for the normalization
$Z_n$,
\bean
Z_n&=&\int_{\BR^n}\Dt^2(z)\prod_1^n\rho(z_i)dz_i\\
\\
&=&\int_{\BR^n}\det(\tilde p_{i-1}(z_j))_{1\leq i,j\leq
n}\det(\tilde p_{k-1}(z_{\ell}))_{1\leq
k,\ell\leq n}\prod_{i=1}^n\rho(z_i)dz_i\\
\\
&=&\sum_{\pi,\pi'\in
S_n}(-1)^{\pi+\pi'}\prod^n_{k=1}\int_{\BR}\tilde
p_{\pi(k)-1}(z_k)
\tilde p_{\pi'(k)-1}(z_k)\rho(z_k)dz_k\\
\\
&=&\sum_{\pi\in S_n}\prod^n_{k=1} \int_{\BR}\tilde
p^2_{\pi(k)-1}(z_k)\rho(z_k)dz_k~,~~
\mbox{using the orthogonality of the $\tilde p_i$'s,}\\
\\
&=&n!\prod_0^{n-1}\int_{\BR}\tilde p^2_k(z)\rho(z)dz
=n!\prod_0^{n-1}h_k .\eean
Then using the expression obtained for $Z_n$ and $(\det
A)^2=\det(AA^{\top})$ in the third equality, one further
computes
\bean
\lefteqn{Z_n^{-1}\Dt_n^2(z)\prod_1^n\rho(z_k)}\\
\\
&=&\frac{1}{n!\prod_0^{n-1}h_k}\det(\tilde
p_{i-1}(z_j))_{1\leq i,j\leq n}\det
(\tilde p_{k-1}(z_{\ell}))_{1\leq k,\ell\leq n}\prod^n_{k=1}\rho(z_k)\\
\\
&=&\frac{1}{n!}\det\left(\frac{\tilde
p_{i-1}(z_j)}{\sqrt{h_{i-1}}}
\sqrt{\rho(z_j)}\right)_{1\leq i,j\leq
n}\det\left(\frac{\tilde
p_{k-1}(z_{\ell})}{\sqrt{h_{k-1}}}
\sqrt{\rho(z_{\ell})}\right)_{1\leq k,\ell\leq n}\\
\\
&=&\frac{1}{n!}\det\left(\sum^n_{j=1}\frac{\tilde
p_{j-1}(z_k)}{\sqrt{h_{j-1}}} \frac{\tilde
p_{j-1}(z_{\ell})}{\sqrt{h_{j-1}}}\sqrt{\rho(z_k)\rho(z_{\ell})}\right)_{1\leq
k,\ell\leq n}\\
\\
&=&\frac{1}{n!}\det(K_n(z_k,z_{\ell}))_{1\leq k,\ell\leq
n}. \eean
One finally needs the classical Christoffel-Darboux
identity: setting $p_j(y)=0$ for $j<0$, one checks
\bean
(y-z)\sum_0^{n-1}p_j(y)p_j(z)&=&\sum_{j=0}^{n-1}\left(a_{j,j-1}p_{j-1}(y)+a_{jj}p_j(y)+
a_{j,j+1}p_{j+1}(y)\right)p_j(z)
\\
\\
&
&-\sum_{j=0}^{n-1}p_j(y)\left(a_{j,j-1}p_{j-1}(z)+a_{jj}p_j(z)+
a_{j,j+1}p_{j+1}(z)\right)
\\
\\
&=&a_{n-1,n}\left(p_n(y)p_{n-1}(z)-p_{n-1}(y)p_n(z)\right)
, \eean
 and one uses (\ref{monic}). The reproducing property follows
immediately from
 the Christoffel-Darboux representation of the kernel in
 terms of orthogonal polynomials. This proves Proposition \ref{Proposition 7.4}.\qed



\subsection{Reproducing kernels}

     \begin{lemma} \label{Lemma 7.5} Let $K(x,y)$ be a symmetric kernel
     satisfying the reproducing property 
$$
\int_{\BR}K(x,y)K(y,z)dy=K(x,z).
$$
Then
\bea \lefteqn{\int\det\left(K(z_i,z_j)\right)_{1\leq
i,j\leq n}dz_n}
 \no\\
\no\\
&=&\left(\int_\BR
K(z,z)dz-n+1\right)\det(K(z_i,z_j))_{1\leq i,j\leq n-1}
\label{Reproducing},\eea where $dy$ stands for any
measure on $\BR$.\end{lemma}

\proof The proof of (\ref{Reproducing}) is due to J.
Verbaarschot \cite{Verbaarschot}, which proceeds in two
steps:

 \vspace{.1cm}

 \underline{Step 1} : Let

$$
M_n=\left(\begin{array}{cc}
M_{n-1}&m\\
\\
\bar m^{\top}&\gamma
\end{array}\right)
$$
be a $n\times n$ Hermitian matrix, with $M_{n-1}$ a
$n-1\times n-1$ Hermitian matrix. Then $m$ is a
$(n-1)\times 1$ column and $\gamma\in\BR$. Then
$$
\det M_n=\gamma \det M_{n-1}-\bar m^{\top}\tilde
M_{n-1}m.
$$
Indeed given the column $u\in\BC^{n-1}$, one checks

\bean \lefteqn{\left(\begin{array}{cc}
I&0\\
\bar u^{\top}&1\end{array}\right)\left(\begin{array}{cc}
M_{n-1}&m\\
\bar
m^{\top}&\gamma\end{array}\right)\left(\begin{array}{cc}
I&u\\
0&1\end{array}\right)}\\
\\
&=&\left(\begin{array}{ccc}
M_{n-1}& &M_{n-1}u+m\\
\bar u^{\top}M_{n-1}+\bar m^{\top}& &\bar
u^{\top}M_{n-1} u+\bar m^{\top}u+\bar u^{\top}m+\gamma
\end{array}\right)\\
\\
&=&\left(\begin{array}{ccc}
M_{n-1}& &0\\
0& &\gamma -\bar m^{\top}M_{n-1}^{-1}m
\end{array}\right), ~\mbox{upon setting~}u=-M^{-1}_{n-1}m.
\eean
 The latter assumes that $M_{n-1}$ is invertible.
Furthermore, using $M^{-1}_{n-1}=(\det
M_{n-1})^{-1}\tilde M_{n-1}$,

\bean
\det M_n&=&(\gamma-\bar m^{\top}M_{n-1}^{-1}m)\det M_{n-1}\\
&=&(\gamma -\bar m^{\top}\tilde M_{n-1}m(\det M_{n-1})^{-1})\det M_{n-1}\\
&=&\gamma\det M_{n-1}-\bar m^{\top}\tilde M_{n-1}m.
\eean

\underline{Step 2} : Proof of identity
(\ref{Reproducing}). Define the matrix
$$
M_k:=(K(z_i,z_j))_{1\leq i,j\leq k}\quad\quad\gamma
:=K(z_k,z_k)
$$
and
$$ m=\left(\begin{array}{c}
K(z_1,z_k)\\
\vdots\\
K(z_{k-1},z_k)
\end{array}
\right);
$$
$M_k$ is a symmetric matrix, since $K$ is a symmetric
kernel. One finds, upon integration over $\BR$,
\bean
\lefteqn{\int_{\BR}\det(M_n)dz_n}\\
\\
&=&\det(M_{n-1})\int_{\BR}K(z,z)dz
 -\int_{\BR}dz_n\sum^{n-1}_{i,j=1}K(z_n,z_i)(\tilde M_{n-1})_{ij}K(z_j,z_n)\\
\\
&=&\det(M_{n-1})\int_{\BR}K(z,z)dz
-\sum^{n-1}_{i,j=1}(\tilde M_{n-1})_{ij}\int_{\BR}K(z_n,z_i)K(z_j,z_n)dz_n\\
\\
&=&\det(M_{n-1})\int_{\BR}K(z,z)dz
-\sum_{i,j=1}^{n-1}(\tilde M_{n-1})_{i,j}(M_{n-1})_{j,i}\\
\\
&=&\det(M_{n-1})\int_{\BR}K(z,z)dz
 -\det (M_{n-1})\sum^{n-1}_{i,j=1}
 (M^{-1}_{n-1})_{ij}(M_{n-1})_{ji}\\
\\
&=&\det(M_{n-1})\left(\int_{\BR}K(z,z)dz-(n-1)\right)
,\eean establishing Lemma \ref{Lemma 7.5}. \qed

\begin{lemma}\label{Lemma 7.6} If
\bean
 (i)&&    \int_{\BR}
K(x,y)K(y,z)dy=K(x,z) \\ \\
(ii)&&   \int_{\BR} K(z,z)dz=n ,
  \eean
then \be \int\!\!\ldots\!\!\int_{\BR^{n-m}}\det
(K(z_i,z_j))_{1\leq i,j\leq n}dz_{m+1}\ldots dz_n=(n-m)!
\det(K(z_i,z_j))_{1\leq i,j\leq m}
 \label{kernel integration}\ee
\end{lemma}

\proof The proof proceeds by induction on $m$. On the
one hand, assuming (\ref{kernel integration}) to be
true, integrate with regard to $z_m$ and use identity
(\ref{Reproducing}):
\bean \lefteqn{\int\ldots\int_{\BR^{n-m+1}}
 \det (K(z_i,z_j))_{1\leq i,j\leq n}dz_mdz_{m+1}\ldots dz_n}
\eean\bean
&=&(n-m)!\int_\BR dz_m\det(K (z_i,z_j))_{1\leq i,j\leq
m}
\\ \\
 &=&(n-m)!(\int_\BR K(z,z)dz-m+1)
  \det(K(z_i,z_j))_{1\leq i,j\leq m-1}
 \\ \\
&=&(n-m)!(n-m+1)\det(K (z_i,z_j))_{1\leq i,j\leq
m-1}\\
&=& (n-m+1)! \det(K(z_i,z_j))_{1\leq i,j\leq m-1}. \eean
On the other hand for $m+1=n$, the statement follows at
once from Lemma \ref{Lemma 7.5}. This ends the proof of
Lemma \ref{Lemma 7.6}.\hfill$\Box$


\subsection{Correlations and Fredholm determinants}

 For this section, see M. L. Mehta \cite{Mehta}, P. Deift \cite{Deift}
 , Tracy-Widom \cite{TW-Schev} and others. Returning now to the probability distribution on the space ${\cal H}_n$ of Hermitian matrices  (setting
$\rho(z):=e^{-V(z)}$)
$$
P(M\in  dM)=\frac{1}{Z_n}e^{-\Tr V(M)}dM, $$
remember from Lemma \ref{Lemma 7.3} and Proposition
\ref{Proposition 7.4},
 $$
 P(M\in{\cal H}_n,\mbox{spectrum }M\subseteq E)
=\int_{E^n}
 P_n(z)dz_1\ldots dz_n
  $$
 with
 \bea
 P_n(z)=\frac{\displaystyle{\Dt_n^2(z)\prod^n_1\rho(z_i)}dz_i}{
\displaystyle{\int_{\BR^n}\Dt_n^2(z)\prod_1^n\rho(z_i)}dz_i}
 &=& Z_n^{-1}\Dt_n^2(z)\prod^n_1\rho(z_k)\no\\
  &=&\frac{1}{n!}
\det(K_n(z_k,z_{\ell}))_{1\leq k,\ell\leq n}
,~~~~~\label{7.4.1}
 \eea
 with the kernel $K_n(y,z)$
defined in Proposition \ref{Proposition 7.4},
 \be
K_n(x,y)=\sqrt{\rho(x)\rho(y)}\sum^{n-1}_{j=0}p_j(x)p_j(y)
.\label{7.4.2}
 \ee
Let $\BE$ be the expectation associated with the
probability $P$ above. Then one has the following
``classical" Proposition for any subset $E\subset \BR$
  (for which a precise statement
 and proof was given by P. Deift \cite{Deift}):
\begin{proposition}
\label{Proposition 7.7} The 1- and 2-point correlations
have the following meaning\footnote{If $x_1,x_2\in E$,
then it counts for 2 in the second formula.}:
\bean
 \int_E K_n(z,z)dz&=&\BE \, (\#\mbox{~of eigenvalues in $E$)}\\
\\
 \int_{E\times E}\det(K_n(z_i,z_j))_{1\leq i,j\leq
2}dz_1dz_2&=&\BE \, ( \#\mbox{~of pairs of eigenvalues
 in $E$)}
,\eean and thus \be K_n(z,z)=\frac{1}{dz}\BE \,
(\#\mbox{~of eigenvalues in $dz$)}\label{equil measure}
\ee
\end{proposition}

\proof Using (\ref{kernel integration}) for $m=1$ and
(\ref{identity}), one computes
\bean
\lefteqn{\int_E K_n(z,z)dz}\\
\\
&=&\int_{\BR}\chi_{_{E}}(z_1)K_n(z_1,z_1)dz_1\\
\\
&=&\frac{1}{(n-1)!}\int_{\BR}dz_1\chi_{_{E}}(z_1)\int\ldots
\int_{\BR^{n-1}}\det(K_n(z_i,z_j))_{1\leq
i,j\leq n} dz_2\ldots dz_n\\
\\
&=&n\int_{\BR}dz_1\chi_{_{E}}(z_1)\int\ldots
\int_{\BR^{n-1}}\frac{1}{Z_n}\Dt_n^2(z)\prod_1^n\rho(z_k)
dz_2\ldots dz_n\\
\\
&=&\frac{n}{Z_n}\int_{\BR^n}\chi_{_{E}}(z_1)\Dt_n^2(z)
\prod_1^n\rho(z_k)dz_k\\
\\
&=&\frac{1}{Z_n}\int_{\BR^n}
\left(\sum^n_{i=1}\chi_{_{E}}(z_i)\right)
\Dt_n^2(z)\prod_1^n\rho(z_k)dz_k\\
\\
&=&\frac{1}{Z_n}\int_{\BR^n}\#\{i\mbox{~such
that~}z_i\in E\}\Dt_n^2(z)\prod_1^n\rho(z_k)dz_k \\
\\
&=& \BE \, (\#\mbox{~of eigenvalues in $E$)}. \eean A
similar argument holds for the second identity of
Proposition \ref{Proposition 7.7}.\qed


Consider disjoint intervals $E_1,\ldots,E_m$ and
integers $1\leq n_1,\ldots,n_m\leq n$ and set
$n_{m+1}:=n-\sum_1^m n_i$. Then for the $n\times n$
Hermitian ensemble with $P_n$ as in (\ref{7.4.1}), one
has:
\bea \lefteqn{P(\mbox{exactly $n_i$ eigenvalues}
 \in E_i, 1\leq i\leq m)}\no\\
\no\\
&=&\left(\begin{array}{c}n\\
  n_1,...,n_m,n_{m+1}\end{array}\right)
  \int_{\BR^n}\prod^{n_1}_{i=1}\chi_{_{E_1}}(x_i)
 \prod^{n_1+n_2}_{i=n_1+1}\chi_{_{E_2}}(x_i)...
  \no\\
  \no\\
& &\prod^{n_1+...+n_m}_{i=n_1+...+n_{m-1}+1}
 \chi_{_{E_m}}(x_i)\prod^n_{i=\sum_1^mn_k+1}
 \chi_{_{(\cup^m
_{i=1}E_i)^c}}(x_i)P_n(x)dx_1\ldots dx_n.\no\\
 \label{7.4.3}\eea
This follows from the symmetry of $P_n(x)$ under the
permutation group; the multinomial coefficient takes
into account the number of ways the event occurs.

\begin{lemma} \label{Lemma 7.8}The following identity holds
\bean  \lefteqn {\hspace*{-1cm}\int_{\BR^n}\prod^n_{k=1}
 \left(1+\sum^m_{i=1}\lb_i\chi_{_{E_i}}(x_k)\right)
  P_n(x)dx_1\ldots dx_n}  \\
  \\
&&~~~~~~~~~~~~~~~=\det   
 \left[I+K_n(x,y)\sum^m_{i=1}\lb_i\chi_{_{E_i}}(y)\right].
\eean
\end{lemma}

\proof Upon setting
$$\tilde K(x,y):=K_n(x,y)\sum^m_1\lb_i\chi_{_{E_i}}(y),
$$
and upon using the fact that $\tilde K(x,y)$ has rank
$n$ in view of its special form (\ref{7.4.2}),
the Fredholm determinant can be computed as follows:
\bean
\lefteqn{\det(I+\tilde K(x,y))}\\
&=&\sum^{\iy}_{\ell=0}\frac{1}{\ell!}\int_{\BR^{\ell}}
\det\Bigl[\tilde K(x_i,x_j)\Bigr]_{1\leq
i,j\leq\ell}~dx_1...dx_{\ell}\\
\\
&=&\sum^{n}_{\ell=0}\frac{1}{\ell!}\int_{\BR^{\ell}}\det
\Bigl[\tilde K(x_i,x_j)\Bigr]_{1\leq
i,j\leq\ell}~dx_1...dx_{\ell}
\eean\bean &=&
\sum^{n}_{\ell=0}\frac{1}{\ell!}\int_{\BR^{\ell}}
 \det\left[K_n(x_i,x_j)\sum^m_{k=1}
  \lb_k\chi_{_{E_k}}(x_j)
\right]_{1\leq
i,j\leq\ell}dx_1...dx_{\ell}\\
\\
&=& \sum^{n}_{\ell=0}\frac{1}{\ell!}
 \sum_{1\leq s_1,\ldots,s_{\ell}\leq m}
 \lb_{s_1} ...\lb_{s_{\ell}}
 \int_{\BR^{\ell}}dx_1...dx_{\ell}
  \prod^{\ell}_{r=1}\chi_{_{E_{s_r}}}(x_r)
  \det[K_n(x_i,x_j)]_{1\leq
i,j\leq\ell}
 \\
 \\
&=&\int_{\BR^n}\sum^{n}_{\ell=0}
  \sum_{1\leq s_1,\ldots,s_{\ell}\leq m}
\lb_{s_1} ...\lb_{s_{\ell}}
 \chi_{ _{E_{s_1}}}(x_1) ...\chi_{_{E_{s_{\ell}}}}(x_{\ell})\\
\\ 
& &\hspace*{1.6cm}\frac{1}{\ell!(n-\ell)!}
 \det[K_n(x_i,x_j)]_{1\leq i,j\leq n}dx_1...dx_n
 ,~~\mbox{using Lemma \ref{Lemma 7.6},}
 \\
\\
&=&\int_{\BR^n}\sum^{n}_{\ell=0}\left(\begin{array}{c}n\\\ell
\end{array}\right)
\sum_{1\leq s_1,\ldots,s_{\ell}\leq m}
 \lb_{s_1}
...\lb_{s_{\ell}}
 \chi_{_{E_{s_1}}}(x_1)...  \chi_{_{E_{s_{\ell}}}}
(x_{\ell})P_n(x)dx_1...dx_n\\
\\
&=& \int_{\BR^n}\sum^{n}_{\ell=0} \sum_{1\leq
s_1,\ldots,s_{\ell}\leq m}
 ~~\sum_{1\leq i_1< \ldots <i_{\ell}\leq n}
\lb_{s_1} ...\lb_{s_{\ell}}\chi_{_{E_{s_1}}}(x_{i_1})...
 \chi_{_{E_{s_{\ell}}}}(x_{i_{\ell}})\\
 &&\hspace{9cm}P_n(x)dx_1...dx_n
\\
&=&\int_{\BR^n}\prod^n_{k=1}\left(1+\sum^m_{i=1}
 \lb_i\chi_{_{E_i}}(x_k)\right)P_n(x)dx_1...dx_n,
\eean
 establishing Lemma \ref{Lemma 7.8}.\qed

\begin{proposition} \label{Proposition 7.9} The Fredholm determinant
is a generating function for the probabilities:
\bea  \lefteqn{
 P(\mbox{exactly $n_i$ eigenvalues
 $\in E_i, 1\leq i\leq m$})
  }\no\\
   \no\\
  &=&\prod^m_1\frac{1}{n_i!}
  \left(\frac{\pl}{\pl\lb_i}\right)^{n_i}
   \det
   \left[I+\sum^m_1 \lb_i K_n(x,y)  \chi_{_{E_i}}(y)\right]
  \Bigg\vert_{\mbox{
   all~}\lb_i=-1}.
 \label{7.4.4}\eea
 In particular
  \be   {
 P(\mbox{no eigenvalues
 $\in E_i, 1\leq i\leq m$}
  })=\det
   \left[I-    K_n(x,y)  \chi_{_{\cup_1^m E_i}}(y)\right]
  .
 \label{7.4.5}\ee
\end{proposition}

\proof The first equality below follows from Lemma
\ref{Lemma 7.8}. Concerning the second equality below,
in order to carry out the differentiation
$\prod^m_1\frac{1}{n_i!}
  \left(\frac{\pl}{\pl\lb_i}\right)^{n_i}$, one chooses
  (keeping in mind the usual product rule of differentiation)
  a first group of $n_1$ factors,
  a (distinct) second group of $n_2$ factors, ... , a $m$th group of $n_m$ factors and finally the last group of
  $n-n_1-\ldots-n_m$ remaining factors among the
  product $\prod_{k=1}^n
\left(1+\sum^m_{i=1}\lb_i \chi_{_{E_i}}(x_k)\right)$.
Then one differentiates the first $m$ groups, leaving
untouched the last group, where one sets $\lb_i=-1$.
This explains the second equality below. Let
$C^n_{n_1,...,n_m}$ be the set of distinct committees of
size $n_1,n_2,...,n_m,n_{m+1}:=n-n_1-...-n_m$ formed
with people $1,2,...,n$:

\bean \lefteqn{\prod^m_{i=1}\frac{1}{n_i!}
 \left(\frac{\pl}{\pl\lb_i}\right)^{n_i}\mbox{
det}
\left[I+\sum^m_1\lb_iK_n(x,y)\chi_{_{E_i}}(y)\right]
 \Bigg\vert_{\mbox{\rm all~}\lb_i=-1}}
 \\
 \\
&=&\prod^m_{i=1}\frac{1}{n_i!}
 \left(\frac{\pl}{\pl\lb_i}\right)^{n_i}\int_{\BR^n}\prod_{k=1}^n
\left(1+\sum^m_{i=1}\lb_i \chi_{_{E_i}}(x_k)\right)
 P_n(x)dx_1...dx_n\Bigg\vert_{\mbox{\rm all~}\lb_i=-1}
 \\
  \\
&=&\sum_{\sg\in
C^n_{n_1,...,n_m}}\int_{\BR^n}\prod_{i=1}^{n_1}
 \chi_{_{E_1}}(x_{\sg(i)})\ldots\prod^{n_1+...+n_m}_{
i=n_1+...+n_{m-1}}\chi_{_{E_m}}(x_{\sg(i)})
\\
\\
&
&\hspace*{1cm}\prod^n_{i=n_1+...+n_m+1}\left(1-\sum^m_{\ell=1}
\chi_{_{E_{\ell}}}(x_{\sg(i)})\right)P_n(x)~dx_1...dx_n\\
\\
&=&\left(
\begin{array}{c}n\\
n_1,...,n_m,n-\displaystyle{
\sum_1^m}n_i\end{array}\right)\int_{\BR^n}
 \prod^{n_1}_{i=1}\chi_{_{E_1}}(x_i)\ldots
  \prod^{n_1+...+n_m}_{i=n_1+...+n_{m-1}}\chi_{_{E_m}}(x_i)\\
\\
\\
& &\hspace*{1cm}
 \prod^n_{i=n_1+...+n_m+1}\left(1-\sum^m_{\ell=1}
 \chi_{_{E_{\ell}}} (x_i)\right)P_n(x)dx_1...dx_n\\
\\
&=&P(\mbox{exactly $n_i$ eigenvalues $\in E_i, 1\leq
i\leq m$}), \eean as follows from (\ref{7.4.3}), thus
establishing identity (\ref{7.4.4}), whereas
 (\ref{7.4.5}) follows from setting $n_1=\ldots=n_m=0$,
 completing the proof
of Proposition \ref{Proposition 7.9}.\qed

\newpage

\section{The distribution of Hermitian matrix ensembles}

\subsection{Classical Hermitian matrix ensembles}

 \subsubsection{The Gaussian Hermitian matrix ensemble (GUE)}

Let $\HR_n$ be the Hermitian ensembles
$$
\HR_n=\{n\times n\mbox{~matrices $M$ satisfying
$M^{\top}=\bar M\}$}.
$$
The real and imaginary parts of the entries $M_{ij}$ of
the $n\times n$ Hermitian matrix $(\bar M=M^{\top})$ are
all independent and Gaussian; the variables $M_{ii}$,
$1\leq i\leq n$, $\Re e ~M_{ij}$ and $\Im m ~M_{ij}$,
$1\leq i<j\leq n$, which parametrize the full matrix
have the following distribution (set ${\cal H}_n\ni
H=(u_{ij})$, with real $u_{ii}$ and off-diagonal
elements $u_{ij}:=v_{ij}+iw_{ij}$)

$$
\begin{array}{llll}
P(M_{ii}\in
du_{ii})&=\displaystyle{\frac{1}{\sqrt{\pi}}}e^{-u^2_{ii}}du_{ii}&
&1\leq
i\leq n\\
\\
P(\Re e~M_{jk}\in
dv_{jk})&=\displaystyle{\frac{2}{\sqrt{\pi}}}e^{-2v^2_{jk}}dv_{jk}&
&1\leq
j<k\leq n\\
\\
P(\Im m~M_{jk}\in
dw_{jk})&=\displaystyle{\frac{2}{\sqrt{\pi}}}e^{-2w^2_{jk}}dw_{jk}&
&1\leq j<k\leq n.
\end{array}
$$
Hence, using Haar measure (\ref{Haar measure}),
 \bea
  \lefteqn{P(M\in dH)}\no\\
  &=&\prod_1^nP(M_{ii}\in
du_{ii})\prod_{1\leq j<k\leq n}P(\Re e~M_{jk}\in
dv_{jk})P(\Im m~M_{jk}\in dw_{jk})
 \no\\
&=&c_n\prod^n_1e^{-u^2_{ii}}\prod_{1\leq j<k\leq
n}e^{-2(v^2_{jk}+w^2_{jk})}\prod^n_1du_{ii}\prod_{1\leq
j<k\leq
n}dv_{jk}dw_{jk}\no\\
\no\\
&=&c_ne^{-\sum_{1\leq i,j\leq
n}|u_{ij}|^2}\prod_1^ndu_{ii}\prod_{1\leq j<k\leq
n}d~\Re
e~u_{jk}~d~\Im m~u_{jk}\no\\
\no\\
&=&c_ne^{-\Tr H^2}dH
 =c_n\Dt_n^2(z)
 \prod_1^ne^{-z^2_i}dz_i,~~~~~\quad H\in\HR_n
 , \label{8.1.1}\eea
  using the representation of Haar measure in terms of spectral
  variables $z_i$ (see Proposition \ref{Proposition 7.2}) and where
$$
c_n=\left(\frac{2}{\pi}\right)^{n^2/2}\frac{1}{2^{n/2}}.
$$
This constant can be computed by representing the
integral of (\ref{8.1.1}) over the full range $\BR^n$ as
a determinant of a moment matrix (as in
(\ref{tau-function})) of Gaussian integrals.

\subsubsection{ Estimating covariances of complex
Gaussian populations and the Laguerre Hermitian
ensemble}\label{subsect 8.2}

Consider the complex Gaussian population
$\vec x=(x_1,\ldots,x_p)^{\top}$,
 with mean and covariance matrix given by
  $$
\vec\mu=\BE(\vec x)=(\mu_1,\ldots,\mu_p)^{\top},
 ~~~~~~~ ~\Sg=\bigl(\BE(x_i-\mu_i)(  \bar x_j-
\bar\mu_j)\bigr)_{1\leq i,j\leq p}
 $$
 and density (for the complex inner-product $\la~,~\ra$),
$$\frac{1}{(2\pi)^{p/2} (\det \Sg)^{1/2}}
e^{-\frac{1}{2} \left\la \vec x-\vec\mu ,{ \Sg^{-1}(\vec
x-\vec\mu )}\right\ra }
$$
Let $\lb_1\geq \lb_2\geq\ldots \geq \lb_p>0$ be the
eigenvalues of $\Sigma$. Taking $n$ samples of
$(x_1,\ldots,x_p)^{\top}$, consider the normalized
$p\times n$ sample matrix:
$$
x=\left(\begin{array}{cclc}
x_{11}-\frac{1}{n}\bigl(\sum_1^n x_{1i}\bigr) &
x_{12}-\frac{1}{n}\bigl(\sum_1^n x_{1i}\bigr)
 & \ldots& x_{1n}-\frac{1}{n}\bigl(\sum_1^n x_{1i}\bigr)\\
x_{21}-\frac{1}{n}\bigl(\sum_1^n x_{2i}\bigr)&x_{22}
-\frac{1}{n}\bigl(\sum_1^n x_{2i}\bigr)&\ldots&x_{2n}
-\frac{1}{n}\bigl(\sum_1^N x_{2i}\bigr)\\
\vdots&\vdots& &\vdots\\
x_{p1}-\frac{1}{n}\bigl(\sum_1^n
x_{pi}\bigr)&x_{p2}-\frac{1}{n}\bigl(\sum_1^n
x_{pi}\bigr)&\ldots&x_{pn}-\frac{1}{n}\bigl(\sum_1^n
x_{pi}\bigr)
\end{array}
\right)
$$
and the $p\times p$ sample covariance matrix,
$$ S:= \frac{1}{N-1} x\bar x^{\top} ,~~~\mbox{with
eigenvalues}
 ~~ z_1,..., z_p >0,
  $$
  which is
an unbiased estimator of $\Sigma$. It is a classical
result that when $\Sg=I$,
   the eigenvalues $z_1,\ldots,z_p> 0$ of $S$
 have the Wishart distribution, a special case of the
 Laguerre Hermitian ensemble (see Hotelling
  \cite{Hotelling} and also Muirhead \cite{Muirhead})
 \bean
 \BP_{n,p}(S\in dM)
  = c_{np}\Dt^2_p(z)\prod_1^pe^{-z_i }z_i^{n-p-1
}dz_i dU =e^{- \Tr M}(\det M)^{n-p-1 }dM .\eean

\subsubsection{Estimating the canonical correlations
between two Gaussian populations and the Jacobi
Hermitian ensemble}

In testing the statistical independence of two complex
Gaussian populations, one needs to know the distribution
of {\em canonical correlation
 coefficients}. I present here the case of real Gaussian
 populations, not knowing whether the complex case has been
 worked out, although it should proceed in the same way.
  To set up the problem, consider $p+q$ normally
distributed random variables $(X_1,...,X_p)^{\top}$ and
$(Y_1,...,Y_q)^{\top}$
 ($p\leq q$) with mean zero and
covariance matrix

\smallbreak

$\hspace*{66mm}\stackrel{p}{\longleftrightarrow}
\quad\stackrel{q}{\longleftrightarrow}$
$$\mbox{cov}\MAT{1}X\\Y\mat :=\Sigma=
 \left(\begin{array}{cc}\Sigma_{11}&\Sigma_{12}\\
\Sigma_{12}^{\top}&\Sigma_{22}\end{array}\right)
 \begin{array}{l}\updownarrow p\\ \updownarrow q\end{array}
.$$
 The method proposed by
Hotelling \cite{Hotelling} is to find linear
transformations $U=L_{1}X$ and $V=L_2Y$ of $X$ and $Y$
having the property that the correlation between the
first components $U_1$ and $V_1$ of the vectors $U$ and
$V$ is maximal subject to the condition that Var $U_1=$
Var $V_1=1$; moreover, one requires the second
components $U_2$ and $V_2$ to have maximal correlation
subjected to
$$ \left\{\begin{array}{ll} {\rm (i)}\quad \mbox{Var\,}
U_2=\mbox{Var\,}V_2=1\\ {\rm (ii)}\quad U_2
\mbox{\,\,and\,\,} V_2 \mbox{\,\,are uncorrelated with
both\,\,} U_1 \mbox{\,\,and\,\,} V_1,
\end{array}\right.
$$   etc\ldots

\noindent Then there exist orthogonal matrices $O_p\in
O(p)$, $O_q\in O(q)$ such that
$$
\Sigma_{11}^{-1/2}\Sigma_{12}\Sigma_{22}^{-1/2}=O^{\top}_{p}P
\,O_q ,$$ where $P$ has the following form:
 \newpage
 $${{q}\atop{\longleftarrow\longrightarrow}}$$ $$
\left\{\begin{array}{lll} P=\left(
\begin{array}{lllllll|l}
\rho_1& & & & & & & \\
 &\ddots& & O& & & &\\
 & &\rho_k& & & & &O\\
& & &\rho_{k+1}& & & &  \\ & O& & & &\ddots& &\\ & & & &
& & \rho_{p}&
\end{array}
\right){\Bigg\updownarrow} p,\quad
k=\mbox{\,rank\,}\Sigma_{12},\\
\hspace*{3cm}{{\longleftarrow\longrightarrow}\atop{p}}
\\
1\geq\rho_1\geq\rho_2\geq\ldots\geq\rho_k
>0,~\rho_{k+1}= ...=\rho_p=0\quad\mbox{(canonical
correlation coefficients),}\\   \\
 \rho_i \mbox{\,\,are
solutions ($\geq 0$) of\,}
\det(\Sigma_{11}^{-1}\Sigma_{12}\Sigma_{22}^{-1}\Sigma_{12}^{\top}-\rho^2I)=0.
\end{array}
\right. $$

\bigbreak

\noindent Then the covariance matrix of the vectors $$
U=L_1X:=O_p\Sigma_{11}^{-1/2}X\quad\mbox{and}\quad
V=L_2Y:=O_q\Sigma_{22}^{-1/2}Y $$ has the canonical form
($\det \Sigma_{can}=\prod_1^p (1-\rho_i^2)$) $$
\mbox{cov}\MAT{1}U\\V\mat=\Sigma_{{can}}=
\MAT{2}I_{p}&P\\P^{\top}&I_q\mat ,
$$ with
 $$\mbox{spectrum}~\Sigma_{{can}}=
 \underbrace{1,\ldots ,1}_{q-p},1-\rho_1,1+\rho_1,\ldots,
  1-\rho_p,1+\rho_p
  $$
 and inverse $$
\Sigma_{can}^{-1}=\frac{1}{\prod_1^p(1-\rho_i^2)^2}
\MAT{2}I_{p}&-P\\-P^{\top}&I_q\mat.
 $$
From here on, we may take $\Sigma=\Sigma_{can}$.
 The $n$ ($n\geq p+q$) independent samples
$(x_{11},\ldots,x_{1p},y_{11},\ldots,y_{1q})^{\top},\ldots$,
$(x_{n1},\ldots,x_{np},y_{n1},\ldots,y_{nq})^{\top}$,
arising from observing $\MAT{1}X\\Y\mat$ lead to a
matrix $\MAT{1}x\\y\mat$ of size $(p+q,n)$, having the
normal distribution \cite{Muirhead} (p. 79 and p. 539)
\bean &&(2\pi)^{-n(p+q)/2}(\det \Sigma)^{-n/2}
\exp ~{-\frac{1}{2}\Tr~ (x^{\top}~y^{\top})
\left(\begin{array}{cc}\Sigma_{11}&\Sigma_{12}\\
\Sigma_{12}^{\top}&\Sigma_{22}\end{array}\right)^{-1}
\MAT{1}x\\y\mat}\\ &&~~~= (2\pi)^{-n(p+q)/2}(\det
\Sigma)^{-n/2}
e^{ -\frac{1}{2}\Tr~\left(x^{\top}(\Sigma^{-1})_{11}x +
 y^{\top}(\Sigma^{-1})_{22}y
 +2y^{\top}(\Sigma^{-1})_{12}^{\top}x\right)}
\eean The conditional distribution of $p\times n$ matrix
$x$ given the $q\times n$ matrix $y$ is also normal: \be
(\det 2\pi\Omega)^{-n/2}
 e^{-\frac{1}{2}\Tr\Omega^{-1}
  (x-Py)(x-Py)^{\top}}\label{8.3.1}\ee with
\bean \Omega&=&\Sigma_{11}-\Sigma_{12}\Sigma_{22}^{-1}
\Sigma_{21} =\diag(1-\rho_1^2,\ldots,1-\rho_p^2)\\
P&=& \Sigma_{12}\Sigma_{22}^{-1}. \eean
 Then the
 maximum likelihood estimates $r_i$ of the $\rho_i$
  satisfy the determinantal equation
\be
  \det
(S_{11}^{-1}S_{12}S_{22}^{-1}S_{12}^{\top}-r^2I)= 0,
 \label{8.3.2}\ee
 corresponding to
  $$
  S=\MAT{2}S_{11}&S_{12}\\S_{12}^{\top}&S_{22}\mat
  :=
  \displaystyle{\left(\begin{tabular}{ll}$xx^{\top}$&$xy^{\top}$ \\
 $yx^{\top}$&$yy^{\top}$
   \end{tabular}\right)},
   $$
where $S_{ij}$ are the associated submatrices of the
{\em sample} covariance matrix $S$.

 \remark  The $r_i$
can also be viewed as $r_i=\cos \theta_i$ , where the
$\theta_1,...,\theta_{p}$
 are the {\em critical} angles
 between two planes in $\BR^n$:

 (i) a
 $p$-dimensional plane $=$ span $\{(x_{11},...,x_{n1}),...,
 (x_{1p},...,x_{np})\}$

 (ii) a $q$-dimensional plane $=$ span $\{
 (y_{11},...,y_{n1})^{\top},...,~
 (y_{1{q}},...,y_{n{q}})\}$.


\bigbreak

\noindent Since the $(q,n)$-matrix $y$ has rank$(y)=q$,
there exists a matrix $H_n\in O(n)$ such that
$yH_n=(y_1\,\,\Big |\,\, O)$; therefore acting on $x$
with $H_n$ leads to

$\hspace*{35mm}\stackrel{q}{\leftrightarrow}\,\,\,\,\,\,
\stackrel{n-q}{\leftrightarrow}\hspace{3cm}
\stackrel{q}{\leftrightarrow}\,\,\,\,\,\,\stackrel{n-q}{\leftrightarrow}$

 \vspace{-.4cm}

\be yH= (y_1\,\,\,\big |\,\,\,O)\updownarrow q, \qquad
xH_n=(u\,\,\,\big | \,\,\,v)\updownarrow p~~.
 \label{8.3.3}\ee
 With this in mind,
 \bean
  \lefteqn{S_{12}
S_{22}^{-1}S_{12}^{\top}-r^2S_{11}}\\
 &=& xy^{\top}(yy^{\top})^{-1}
 yx^{\top}-r^2xx^{\top}  \\
&=&xH(yH)^{\top}(yH(yH)^{\top})^{-1}
 yH(xH)^{\top}-r^2(xH)(xH)^{\top}\\
 &=&(u~\big |
~v)\MAT{1}y^{\top}_1\\O\mat
 \left((y_1~\big |~ O)
  \MAT{1}y^{\top}_1\\O\mat\right)^{-1}
(y_1~\big | ~O) \MAT{1}u^{\top}\\v^{\top}\mat
-r^2(u~\big | ~v)\MAT{1}u^{\top}\\v^{\top}\mat\\
%
%
&=&(u~\big |~v)\MAT{2}I_q&O\\O&0_{n-q}\mat
 \MAT{1}u^{\top}\\v^{\top}\mat
-r^2(u~\big |~v)\MAT{1}u^{\top}\\v^{\top}\mat
 \\
 &=&uu^{\top}-r^2(uu^{\top}+vv^{\top}),
\eean
 and so the equation (\ref{8.3.2}) for the $r_i$ can be
rewritten
 \be\det ( uu^{\top}-r^2(uu^{\top}+vv^{\top})=0.
 \ee
 Then setting the forms (\ref{8.3.3}) of $x$ and $y$
 in the conditional distribution
(\ref{8.3.1}) of $x$ given $y$, one computes the
following, setting $H:=H_n$,


\medbreak

$\Tr\Omega^{-1}(x-Py)(x-Py)^ {\top}$ \bean
&=&\Tr\Omega^{-1}(xH-PyH)(xH-PyH )^{\top}\\
&=&\Tr\Omega^{-1}\left((u\,\,\big |\,\,v)-P(y_1\,\,\big
|\,\,O)\right)\left((u\,\,\big
|\,\,v)-P(y_1\,\,\big |\,\,O)\right)^{\top}\\
&=&\Tr\Omega^{-1}(u-
 Py_1)(u-
 Py_1)
^{\top} +\Tr\Omega^{-1}vv^{\top};  ~~\Omega
=\diag(1-\rho_1^2,\ldots,1-\rho_p^2); \eean this
establishes the independence of
 the normal distributions $u$ and $v$, given the
 matrix $y$,
 with
$$ u\equiv N(P
y_1,\Omega),\quad v\equiv
N(O,\Omega).~~P=\diag(\rho_1,\ldots,\rho_p). $$ Hence
$uu^{\top}$ and $vv^{\top}$ are conditionally
independent and both Wishart distributed; to be precise:

\smallbreak
\begin{itemize}
  \item The $p\times p$ matrices $vv^{\top}$ are Wishart
   distributed, given $y$, with $n-q$ degrees of freedom and
   covariance $\Omega$;

  \item The $p\times p$ matrices $uu^{\top}$ are
  non-centrally Wishart distributed, given $y$, with $q$ degrees
  of freedom, with covariance $\Omega$ and with
  non-centrality matrix
$$\frac{1}{2}Py_1y_1^{\top} P^{\top}\Omega^{-1}
 .$$

  \item The marginal distribution of the $q\times q$ matrices $yy^{\top}$
   are Wishart distributed,
 with $n$ degrees of freedom and covariance $I_q$,
 because the marginal distribution of $y$ is normal
 with covariance $I_q$.

\end{itemize}



To summarize, given the matrix $y$, the sample canonical
correlation coefficients $r_1^2>\ldots
>r_p^2$ are the roots of
\bean (r_1^2>\ldots >r_p^2)&=&\mbox{~roots
of~}\det(xy^{\top}(yy^{\top})^{-1}yx^{\top}
-r^2xx^{\top})=0\\ &=&\mbox{~roots
of~}\det(uu^{\top}-r^2(uu^{\top}+vv^{\top}))=0\\
&=&\mbox{~roots
of~}\det(uu^{\top}(uu^{\top}+vv^{\top})^{-1}-r^2I)=0.
\eean

\medbreak

Then one shows that, knowing $uu^{\top}$ and $
vv^{\top}$ are Wishart and conditiionally independent,
the conditional distribution of $r_1^2>\ldots
>r^2_p$, given the
matrix $y$ is given by \bean
\pi^{p^2/2}c_{n,p,q}e^{-\frac{1}{2}\Tr
 P yy^{\top}
  P^{\top}\Omega^{-1}}
\Delta_p(r^2) \prod^p_1(r_i^2)^{\frac{1}{2}(q-p-1)}
(1-r_i^2)^{\frac{1}{2}(n-q-p-1)}.\\
%
 \sum_{\lb \in \BY}
\frac{(n/2)_{\lb}C_{\lb}(\frac{1}{2}
 P yy^{\top}
  P^{\top} \Omega^{-1})}{(q/2)_{\lb}C_{\lb}(I_p)~|\lb|!}C_{\lb}(R^2),
\eean where\footnote{Using the standard notation,
defined for a partition $\lb$,
 $$
  (a)_{\lb} :=
\prod_i(a+ (1-i))_{\lb_i}  
 ,
\mbox{with $(x)_n:=x(x+1)\ldots (x+n-1),~x_0=1$}. $$}
$$ R^2=\diag(r_1^2,\ldots,r_p^2),\quad
c_{n,p,q}= \frac{\Gamma_p(n/2)}{\Gamma_p(q/2)
\Gamma_p((n-q)/2)\Gamma_p(p/2)}, $$ and where the
$C_{\lb}$ are proportional to Jack polynomials
corresponding to the partition $\lb$; for details see
Muirhead \cite{Muirhead} and Adler-van Moerbeke
\cite{AvM3}. By taking the expectation with regard to
$y$ or, what is the same, by integrating over the matrix
$yy^{\top}$, which is Wishart distributed (see section
\ref{subsect 8.2}), one obtains:

\begin{theorem} Let $X_1,\ldots,X_p,Y_1,\ldots,Y_q$ ($p\leq q$) be normally
distributed random variables with zero means and
covariance matrix $\Sigma
=\displaystyle{\MAT{2}\Sigma_{11}&\Sigma_{12}\\
\Sigma_{21}&\Sigma_{22}\mat}$. If
$\rho^2_1,\ldots,\rho^2_p$ are the roots of
$\det(\Sigma_{11}^{-1}\Sigma_{12}\Sigma_{22}^{-1}\Sigma_{12}^{\top}-\rho^2I)=0$,
then the maximum likelihood estimates
$r_1^2,\ldots,r^2_p$ from a sample of size $n$ ($n\geq
p+q$) are given by the roots of $$
\det(xy^{\top}(yy^{\top})^{-1}yx^{\top}-r^2xx^{\top})= 0
.$$ 
Then, assuming $\rho_1^2=\ldots =\rho_p^2=0$, the joint
density of the $z_i=r^2_i$ is given by the following
density:
 \be \pi^{p^2/2}c_{n,p,q}
  \Delta_p(z)
  \prod^p_{i=1}z_i^{(q-p-1)/2}
   (1-z_i)^{(n-q-p-1)/2}
dz_i.
 \ee
 \end{theorem}

\remark Taking complex Gaussian populations should
introduce in the formula above $\Dt^2_p(z)$ instead of
$\Dt_p(z)$ and should remove the $1/2$'s in the
exponent.


\subsection{The probability for the classical Hermitian
 random ensembles and PDE's generalizing Painlev\'e}

\subsubsection{The Gaussian ensemble (GUE)
}\label{sect 8.2}

This section deals with the Gaussian Hermitian matrix
ensemble, discussed in previous section. Given the
disjoint union of intervals
 $$
E:=\bigcup_1^{r}~[c_{2i-1},c_{2i}]\subseteq \BR,
 $$
define the algebra of differential operators
\be
  {\cal B}_k=\sum_{i=1}^{2r}
c_i^{k+1} \frac{\pl}{\pl c_i}.
 \label{boundary-operators8.2}
 \ee
 The PDE (\ref{GUE-PDE}) appearing below was
 obtained by Adler-Shiota-van Moerbeke
 \cite{ASV,ASV1,AvM1}, whereas the ODE (\ref{ODE-GUE})
 was first obtained by Tracy-Widom \cite{TW-Airy}. The
 method used here is different from the one of
 Tracy-Widom, who use the method proposed by
 Jimbo-Miwa-Mori-Sato \cite{JMMS}. John Harnad
 then shows in \cite{Harnad} the relationship between the PDE's obtained
 by Tracy-Widom and by Adler-van Moerbeke.

 \begin{theorem}  \label{Theo 8.1} The log of the probability
 $$
\BP_n:=\BP_n(\mbox{all~}z_i\in
E)=\frac{\displaystyle{\int_{E^n}\Dt_n^2(z)
 \prod_1^ne^{-z^2_i}dz_i}}{
\displaystyle{\int_{\BR^n}\Dt_n^2(z)\prod_1^ne^{-z^2_i}dz_i}}
 $$ satisfies the PDE
\be ({\cal B}^4_{-1}+8n{\cal B}^2_{-1}+12{\cal
B}^2_0+24{\cal B}_0-16{\cal B}_{-1}{\cal
B}_{1})\log\BP_n+6({\cal B}_{-1}^2\log\BP_n)^2=0.
\label{GUE-PDE}\ee
 In particular
 $$f(x):=\frac{d}{d x}
  \log \BP_n(\max_i z_i \leq x)$$
  satisfies the 3rd order ODE
\be
 f^{\prime\prime\prime} +  6
~f^{\prime 2}+4  (2n- x^2) f^{\prime} + 4  x f =0
 \label{ODE-GUE}
 ,\ee
 which can be transformed into the {\bf Painlev\'e IV}
 equation.

\end{theorem}

\proof In Theorem \ref{Theo: 5.2}, it was shown that
integral (here we indicate the $t$- and $E$-dependence)
\be \tau_n(t;E)=\frac{1}{n!}\int_{ E
^n}\Dt_n^2(z)\prod_1^ne^{-z_i^2+\sum_1^{\iy} t_kz^k_i}
 dz_i
\label{integral8}\ee satisfies the KP equation,
regardless of $E$,
\be \left(\left(\frac{\pl}{\pl
t_1}\right)^4+3\left(\frac{\pl}{\pl
t_2}\right)^2-4\frac{\pl^2}{\pl t_1\pl
t_3}\right)\log\tau_n +6\left(\left(\frac{\pl}{\pl
t_1}\right)^2\log\tau_n\right)^2=0,
\label{KPequation1}\ee and in Theorem \ref{thm: Virasoro
constraints}, $\tau_n(t;E)$ and $\tau_n(t;\BR)$ were
shown to satisfy the same Virasoro constraints, and, in
particular for the Gaussian case, the three equations
(\ref{GUE1}), with the boundary term missing in the case
$\tau_n(t;\BR)$.

Let $T_i$
 denote the pure
$t$-differentiation appearing on the
 right hand side of
(\ref{GUE1}), with ``principal symbol" $\frac{\pl}{\pl
t_{i+2}}$,
\bea  T_{-1}&:=&   \frac{\pl}{\pl
t_1}-\frac{1}{2}\sum_{i\geq 2}
it_i \frac{\pl}{\pl t_{i-1}}  \no \\
 T_0&:=&
 \frac{\pl}{\pl t_2}-\frac{1}{2} \sum_{i\geq 1} it_i \frac{\pl}{\pl
t_i} \no \\
 T_1&:=& \frac{\pl}{\pl t_3}-\frac{1}{2} \sum_{i\geq
1}it_i\frac{\pl}{\pl
t_{i+1}}-\frac{n}{2}\sum_{i\geq 2}
it_i\frac{\pl}{\pl t_{i-1}}.
 \label{Ti's}\eea
Recall the differential operators ${\cal D}_i$ in terms
of the boundary operators (\ref{boundary-operators8.2}),
appearing in the Virasoro constraints (\ref{GUE1}),
 \bea
  {\cal D}_{-1}&=& -\frac{1}{2}{\cal B}_{-1} \nonumber\\
{\cal D}_{0}&=&-\frac{1}{2}{\cal B}_0  \nonumber\\
{\cal D}_{1}&=&-\frac{1}{2}\left({\cal B}_1\!+\!
 n{\cal B}_{-1}\right)
 .   \label{GUE2}
\eea
 With this notation, the Virasoro constraints (\ref{GUE1}) can be
summarized as ($F:=\log\tau_n(t;E)$)
$$
\DR_{-1}F=T_{-1}F-\frac{nt_1}{2},\quad\DR_{0}F=T_{0}F-\frac{n^2}{2},\quad\DR_{1}F=T_{
1}F-\frac{n^2t_1}{2} .
$$
Expressing the action of $T_i$ on $t_1$, one finds
  \bean
 T_{-1}t_1&=&1-t_2 ~~~~~~~~~
  T^2_{-1}t_1=T_{-1}(1-t_2)=\frac{3}{2}t_3\\
 T_1t_1&=&-nt_2~~~~~~~~~~T^3_{-1}t_1=T_{-1}T^2_{-1}t_1=T_{-1}\bigl
 (\frac{3}{2}t_3\bigr)=-3t_4
  ,\eean
one computes consecutive powers of $\DR_{i}$ and their
products, and one notices that $\DR_{i}$ involves
differentiation with regard to the boundary terms only,
implying in particular that
 {\em $\DR_{i}$ and $T_{i}$ commute}. In view of the form of
 the KP equation, containing only certain partials, and in view of the fact that
 the ``{\em principal symbol}" of $T_i$ equals $\pl /\pl t_{i+2}$, it
 suffices to compute
\bea  & &\DR_{-1}F=T_{-1}F-\frac{nt_1}{2}
\no\\
& &\DR_{-1}^2F=\DR_{-1}T_{-1}F=T_{-1}\DR_{-1}F=
 T_{-1}\left(T_{-1}F-
\frac{nt_1}{2}\right)=T^2_{-1}F-\frac{n}{2}(1-t_2)\no\\
\no\\
&
&\DR_{-1}^3F=\DR_{-1}T^2_{-1}F=T^2_{-1}\DR_{-1}F=T^2_{-1}
\left(T_{-1}F-\frac{nt_1}{2}
\right)=T^3_{-1}F-\frac{3n}{4}t_3\no\\
\no\\
&
&\DR^4_{-1}F=\DR_{-1}T^3_{-1}F=T^3_{-1}\DR_{-1}F=T^3_{-1}\left(
T_{-1}F-\frac{nt_1}{2}\right)=T^4_{-1}F+\frac{3n}{2}t_4.
\no\\
& &\DR_1F=T_1F-\frac{n^2t_1}{2}
\no\\
& &\DR_{-1}\DR_1F=\DR_{-1}T_1F=T_1\DR_{-1}F=
T_1\left(T_{-1}F-\frac{n^2t_1}{2}\right)=T_
1T_{-1}F+\frac{n^3}{2}t_2
 \no\\ & &\DR_0F=T_0F-\frac{n^2}{2}
\no\\
& &\DR^2_0F=\DR_0T_0F=T_0\DR_0F=T_0
\left(T_0F-\frac{n^2}{2}\right)=T^2_0F
 \label{8.2.7}.\eea
Since one is actually interested in the integral
(\ref{integral8}) along the locus
$\LR:=\{\mbox{all}~t_i=0\}$, and since readily from
(\ref{Ti's}) one has $T_i|_{\LR}=\pl/\pl t_{i+2}$, one
deduces from the equations above
(\ref{8.2.7})\footnote{Notice one also needs $\DR_0F$,
because $\pl F/\pl t_2$ appears in the expressions
$\DR_0^2F$ and $\DR_{-1}\DR_1F$ below.}
\bean
  \DR_{-1}^2F\Bigl\vert_{\LR}&=&
   T^2_{-1}F\Bigl\vert_{\LR}-\frac{n}{2}=\frac{\pl^2
F}{\pl
t^2_1}\Bigl\vert_{\LR}-\frac{n}{2}\\
\\
 \DR_{-1}^4F\Bigl\vert_{\LR}&=&T^4_{-1}F\Bigl\vert_{\LR}=\frac{\pl }{\pl
t_1}T^3_{-1}F\Bigl\vert_{\LR}=T^3_{-1} \frac{\pl F}{\pl
t_1}\Bigl\vert_{\LR}=\frac{\pl ^4F}{\pl
t^4_1}\Bigl\vert_{\LR}
 \\
  \DR^2_0F\Bigl\vert_{\LR}&=&T^2_0F\Bigl\vert_{\LR}=\frac{\pl }{\pl
t_2}T_0F\Bigl\vert_{\LR}=\left(\frac{\pl^2}{\pl t^2_2}-
\frac{\pl}{\pl t_2}\right)F\Bigl\vert_{\LR}
\\
 \DR_{-1}\DR_1F\Bigl\vert_{\LR}&=&T_
1T_{-1}F\Bigl\vert_{\LR}=\frac{\pl}{\pl
t_3}T_{-1}F\Bigl\vert_{\LR}=\left(\frac{\pl^2}{\pl
t_3\pl t_1}-\frac{3}{2}\frac{\pl}{\pl
t_2}\right)F\Bigl\vert_{\LR}
\\
 \DR_0F\Bigl\vert_{\LR}
  &=&T_0F\Bigl\vert_{\LR}-\frac{n^2}{2}=\frac{\pl F}{\pl t_2}\Bigl\vert_{\LR}
   -\frac{n^2}{2}
.%
 \eean
 By solving the five expressions above linearly in terms of the left hand side, one deduces
 \bean   \frac{\pl^2F}{\pl
t^2_1}\Bigl\vert_{\LR}&=&\DR_{-1}^2F\Bigl\vert_{\LR}
 +\frac{n}{2},~~~~~~~~~~~~~~~~~~~~\quad\frac{\pl^4F}{\pl
t^4_1}\Bigl\vert_{\LR}=\DR_{-1}^4F\Bigl\vert_{\LR}\\
\\
 \frac{\pl F}{\pl t_2}\Bigl\vert_{\LR}&=&\DR_0F\Bigl\vert_{\LR}+\frac{n^2}{2}
    ~~~~~~~~~~~~~~~~~~~~~~~~~\frac{\pl^2F}{\pl
t^2_2}\Bigl\vert_{\LR}
=\DR_0^2F\Bigl\vert_{\LR}+\DR_0F\Bigl\vert_{\LR}+\frac{n^2}{2}\\
\\
 \frac{\pl^2F}{\pl t_1\pl
t_3}\Bigl\vert_{\LR}&=&
\DR_{-1}\DR_1F\Bigl\vert_{\LR}
+\frac{3}{2}\DR_0F\Bigl\vert_{\LR}+\frac{3n^2}{4}
. \eean
 So, substituting into the KP equation and expressing the ${\cal D}_i$ in terms of the
 ${\cal B}_i$ as in (\ref{GUE2}), one finds
\bean 0&= &\left(\left(\frac{\pl}{\pl
t_1}\right)^4+3\left(\frac{\pl}{\pl
t_2}\right)^2-4\frac{\pl^2}{\pl t_1\pl
t_3}\right)F+6\left(\left(\frac{\pl}{\pl
t_1}\right)^2F\right)^2\Bigl\vert_{\LR}\\
\\
&=&\DR^4_{-1}+3\left(\DR^2_0F+\DR_0F+\frac{n^2}{2}\right)-4\left(\DR_{-1}\DR_1F+\frac
{3}{2}\DR_0F+ \frac{3n^2}{4}\right)
\\
& &\hspace*{3cm}+~6\left(\DR^2_{-1}F+\frac{n^2}{2}
 \right)^2\Bigl\vert_{\LR} \\
\\
&=&\DR^4_{-1}F+6n\DR^2_{-1}F+3\DR_0^2F-3\DR_0F-
4\DR_{-1}\DR_1F+6(\DR^2_{-1}F)^2\Bigl\vert_{\LR}
\eean\bean &=&\frac{1}{16}({\cal B}^4_{-1}+8n{\cal
B}_{-1}^2+12{\cal B}_0^2+24{\cal B}_0-16{\cal
B}_{-1}{\cal B}_1)F+\frac{3}{8} ({\cal
B}_{-1}^2F)^2\Bigl\vert_{\LR} ,\eean
 which establishes (\ref{GUE-PDE}) for (remember
 notation (\ref{integral8}))
 $$
 F=\log \tau_n(0,E)=\log \BP_n(\mbox{all}~z_i\in E)+\log
  \tau_n(0,\BR).
  $$
   Since the ${\cal B}_k$ are derivations with regard to the
   boundary points of the set $E$ and since
$\log \tau_n(0,\BR)$ is independent of those points,
  the equation (\ref{GUE-PDE}) is also valid for $\log
  \BP_n$; it is an equation of order 4.

  When $E$ is a
  semi-infinite interval $(-\iy,x)$, then one has ${\cal
  B}_k=x^{k+1}\pl/\pl x$ and then, of course, PDE (\ref{GUE-PDE})
  turns into an ODE (\ref{ODE-GUE}), of an order one less, involving $f(x):=\frac{d}{d x}
  \log \BP_n(\max_i z_i \leq x)$. For the connection
  with Painlev\'e IV, see section \ref{section Chazy},
  thus
   ending the proof of Theorem \ref{Theo 8.1}.
\qed



\subsubsection{The Laguerre ensemble}\label{subsubsect
8.3.1} Given $E\subset \BR^+$ and the boundary operators
 $$
   {\cal B}_{k}:=
\sum_1^{2r}c_i^{k+2}\frac{\pl}{\pl c_i}
 ,~\mbox{for}~k=-1,0,1,\ldots ,$$
 the following statement holds: (see \cite{ASV,ASV1,AvM1} for the PDE obtained below; the ODE was first obtained by
 Tracy-Widom \cite{TW-Airy})

\begin{theorem}  \label{Theo 8.2} The log of the probability
 $$
\BP_n:=\BP_n(\mbox{all~}z_i\in
E)=\frac{\displaystyle{\int_{E^n}\Dt_n^2(z)
 \prod_1^n  z_i^ae^{-z_i}dz_i}}{
\displaystyle{\int_{(\BR^+)^n}\Dt_n^2(z)\prod_1^nz_i^ae^{-z_i}dz_i}}
 $$ satisfies the PDE
$$\left(
 \begin{array}{l}
  {\cal B}_{-1}^4 - 2 {\cal B}^3_{-1} + ( 1 - a^2)
{\cal B}_{-1}^2 - 4 {\cal B}_{1} {\cal B}_{-1}
 + 3
{\cal B}_0^2 \\
 + 2(2n + a) {\cal B}_0 {\cal B}_{-1} - 2
{\cal B}_1 - (2n + a) {\cal B}_0
 \end{array}\right) \log \BP_n
$$
\be \hspace{4cm} + 6({\cal B}_{-1}^2 \log \BP_n)^2 - 4
({\cal B}_{-1}^2 \log \BP_n) ({\cal B}_{-1} \log \BP_n)
= 0. \label{PDE-Laguerre} \ee
  In particular,
 $f(x):=x\frac{d}{d x}
   \log P_n(\max_i z_i\leq x)$
satisfies
  \be
 x^2f^{\prime\prime\prime} +
xf^{\prime\prime}
 + 6
xf^{ \prime 2} - 4 ff^{\prime}-((a-x)^2 - 4nx )
f^{\prime}- (2n + a - x) f = 0 ,
\label{ODE-Laguerre}
 \ee
  which can be transformed into
the {\bf Painlev\'e V} equation.
\end{theorem}

\subsubsection{The Jacobi ensemble}
 The Jacobi weightn
 is given by
$(1-z)^{a}(1+z)^{b}$. For $E\subset [-1,+1]$, the
boundary differential
 operators ${\cal B}_k$ are now defined by
   $${\cal
  B}_k=\sum_1^{2r}c_i^{k+1}(1-c_i^2)\frac{\pl}{\pl c_i}.$$
 Introduce the following parameters
 $$
  r=4(a^2+b^2)
,~~~s=2(a^2-b^2),~~~
 q= 2(2n+a+b )^2.
$$
$$
  r=a^2+b^2
,~~~s=a^2-b^2,~~~
 q= (2n+a+b )^2.
$$
\begin{theorem}\label{Theo 8.3}
(Haine-Semengue \cite{Haine1} and Adler-van Moerbeke
\cite{AvM1})
The following probability
 \be \BP_n:=\BP_n(\mbox{all~}z_i\in
E)=\frac{{\displaystyle\int_{E^n}\Dt_n(z)^2\prod_{k=1}^n
(1-z_i)^{a}(1+z_i)^{b}dz_i}
}{{\displaystyle\int_{[-1,1]^n} \Dt_n(z)^2\prod_{k=1}^n
(1-z_i)^{a}(1+z_i)^{b}dz_i}}
 \ee
  satisfies the PDE:
 \bea \lefteqn{\hspace*{-1cm}\left(\begin{array}{l}
  {\cal
B}_{-1}^4+(q-2r+2){\cal B}^2_{-1}
+q(3  {\cal B} _{0}^2 - 2  {\cal B}_{0}+
 2   {\cal B}_2) +4{\cal B}_0 {\cal
B}_{-1}^2  \\
-2(2q-1){\cal B}_{1}{\cal B}_{-1} +(2{\cal B}_{-1}\log
\BP_n-s) ({\cal B}_1 -{\cal B}_{-1}+2{\cal B}_{0} {\cal
B}_{-1})
 \end{array}
  \right)\log \BP_n~~~~~~~~~~~~~~~~~~~~~~~~~}
\nonumber\\ \no\\&& \hspace{2cm}
  +2{\cal B}^2_{-1}\log \BP_n\left( 2{\cal B}_{0}\log \BP_n
   +3{\cal B}^2_{-1}\log \BP_n\right)
=0\label{PDE-Jacobi}\eea
 In particular, $f(x):=(1-x^2)\frac{d}{dx}\log
\BP_n(\max_{i} \lb_i \leq x)$ for $0<x<1$ satisfies:
  \bea
 &&\hspace{-1.3cm} (x^2-1)^2f^{\prime\prime\prime}
 +2(x^2-1)\left(xf^{\prime\prime}
 -3f^{\prime 2}\right)
 +\left(8
xf-q(x^{2}-1)-2sx-2r \right)f^{\prime}
 \nonumber\\&&~~~~~~~~~~~~~~~~~~~~~~~~~~~~~-f\left(2f-qx-s \right)=0,
 \label{ODE-Jacobi}\eea
which is a version of {\bf Painlev\'e VI}.

\end{theorem}

{\medskip\noindent{\it Proof of Theorems \ref{Theo 8.2}
and \ref{Theo 8.3}:\/} }
 It goes along the same lines as Theorem \ref{Theo 8.1}
for GUE, namely using the Virasoro constraints
(\ref{6.2.8}) and (\ref{6.2.9}), together with the KP
equation (\ref{KPequation1}). This then leads to the
PDE's
(\ref{PDE-Laguerre}) and (\ref{PDE-Jacobi}). The ODE's
(\ref{ODE-Laguerre}) and (\ref{ODE-Jacobi}) are found by
simple computation. For connections with the Painlev\'e
equations see section \ref{section Chazy}.\qed

\subsection{Chazy and Painlev\'e equations} \label{section
Chazy}

Each of these three equations (\ref{ODE-GUE}),
(\ref{ODE-Laguerre}), (\ref{ODE-Jacobi}) is of the Chazy
form
 \be
f^{\prime \prime\prime}+\frac{P'}{P}f^{
\prime\prime}+\frac{6}{P}f^{\prime 2}-\frac{4P'}{P^2}ff'
+\frac{P^{\prime\prime}}{P^2} f^2
+\frac{4Q}{P^2}f'-\frac{2Q'}{P^2}f+\frac{2R}{P^2}=0,
\label{Chazy} \ee
 with $P,Q,R$ having the form:
  $$
\begin{array}{llll}
Gauss & P(x)=1&~4Q(x)=-4x^2+8n&~R=0 \\ Laguerre &
P(x)=x&~4Q(x)=-(x-a)^2+4nx
&~R=0
\\ Jacobi &
P(x)=1-x^2&~4Q(x)=-(q(x^2-1)+2sx+2r)&~
 R=0 \end{array}
 $$
 The differential equation (\ref{Chazy}) belongs to the general Chazy class $$
f^{\prime\prime\prime}=F(z,f,f^{\prime},f^{\prime\prime}),~\mbox{where
$F$ is rational in $f,f^{\prime},f^{\prime\prime}$ and
locally analytic in z,}$$ subjected to the requirement
that the general solution be free of movable branch
points; the latter is a branch point whose location
depends on the integration constants. In his
classification Chazy found thirteen cases, the first of
which is given by equation (\ref{Chazy}),
 with arbitrary polynomials $P(z), Q(z), R(z)$ of
degree $3,2,1$ respectively. Cosgrove and Scoufis
\cite{CosgroveScoufis,Cosgrove}, show that this third
order equation
 has a first integral, which is second order in $f$
and quadratic in $f^{\prime\prime}$, \bea &f^{
\prime\prime 2}& +\frac{4}{P^2}
 \left( (Pf^{\prime 2}+Q f^{\prime}+R)f^{\prime}
 - (P' f^{\prime 2}+\frac{}{}Q' f^{\prime}+R')f^{}
  \right. \nonumber\\ && \hspace{2cm} \left.
   +\frac{1}{2}(P^{\prime\prime}f^{\prime
}+Q^{\prime\prime} )f^2
 -\frac{1}{6} P^{\prime\prime\prime}f^3 +c\right)=0;
\label{8.6.2}\eea $c$ is the integration constant.
Equations of the general form $$ f^{ \prime\prime
2}=G(x,f,f^{ \prime})$$ are invariant under the map $$
x\mapsto \frac{a_1z+a_2}{a_3z+a_4}~~ \mbox{and}~~
f\mapsto \frac{a_5f+a_6z+a_7}{a_3z+a_4}.$$ Using this
map, the polynomial $P(z)$ can be normalized to $$
P(z)=z(z-1),~z,~\mbox{or} ~1.$$
 Equation (\ref{8.6.2}) is a master
Painlev\'e equation, containing the 6 Painlev\'e
equations, replacing $f(z)$ by some new variable $g(z)$. 
%
e.g,

\begin{itemize}

 \item $  g^{\prime\prime 2}=-4g^{\prime 3
}-2  g^{\prime}(zg^{\prime}-g)+A_1 $ \hspace{4cm}({\bf
Painlev\'e II})

  \item $  g^{\prime\prime 2}=-4g^{\prime 3
}+4(zg^{\prime}-g)^2+A_1g^{\prime}+A_2$
\hspace{3cm}({\bf Painlev\'e IV})

  \item $ (zg^{\prime\prime })^2=(zg^{\prime}-g)\Bigl(-4g^{\prime 2
}+A_1(zg^{\prime}-g)+A_2\Bigr)+A_3g^{\prime}+A_4
 $

 \hspace{10cm}({\bf
Painlev\'e V})
  \item $ (z(z-1)g^{\prime\prime
})^2=(zg^{\prime}-g)\Bigl(4g^{\prime
2}-4g^{\prime}(zg^{\prime}-g)
 +A_2\Bigr) +A_1g^{\prime
2}+A_3g^{\prime}+A_4 $

\hspace{10cm}({\bf Painlev\'e VI})
\end{itemize}
 Now, each of these Painlev\'e II, IV, V, VI equations can be transformed
into the standard Painlev\'e equations, which are all
differential equations of the form $$ f^{\prime\prime
}=F(z,f,f^{\prime}), \mbox{rational in $f,~f^{\prime }$,
analytic in $z$,} $$ whose general solution has no
movable critical points. Painlev\'e showed that this
requirement leads to 50 types of equations, six of which
cannot be reduced to known equations.


\newcommand{\tip}{\tilde p}
\newcommand{\limn}{\lim_{n\rg\iy}}

\section{Large Hermitian matrix ensembles}

\subsection{Equilibrium measure for GUE and Wigner's semi-circle}\label{subsect equiGUE}

Remember according to
 (\ref{equil measure}),
 the average density of eigenvalues is given by
 $K_n(z,z)dz$. Pastur and Marcenko
\cite{Pastur} have proposed a method to compute the
average density of eigenvalues (equilibrium
distribution), when $n$ gets very large. For a rigorous
and very general approach, see Johansson
\cite{Johansson-1}, who also studies the fluctuations of
the linear statistics of the eigenvalues about the
equilibrium distribution.

Consider the case of a random Hermitian ensemble with
probability defined by
$$
 \frac{1}{Z_n}\int_{\HR_n(E )}dM e^{-\frac{n}{2v^2}\Tr
(M-A)^2}dM
 ,$$  for a diagonal matrix
$A=(a_1,\ldots,a_n)$.  Consider then the spectral
function of $A$, namely $d\sigma(\lb):=
 \frac{1}{n}\sum_i \dt(\lb-a_i)$. The Pastur-Marcenko
method tells us that the Stieltjes
 transform of the equilibrium measure of $d\nu(\lb)$,
  when $n\rg \iy$, namely
  $$
   f(z)=\int_{-\iy}^{\iy} \frac{d\nu(\lb)}{\lb -z},~~\Im m
   ~z\neq 0,
    $$
     satisfies the integral equation
      \be
      f(z)=\int_{-\iy}^{\iy}
      \frac{d\sigma (\lb)}{\lb-z-v^2 f(z)}
      .\label{IntEqt}\ee
    The density of the equilibrium
    distribution is then given by
     $$
      \frac{d\nu(z)}{dz}=\frac{1}{\pi} \Im m~ f(z)
     .$$
     When $A=0$, the integral equation (\ref{IntEqt})
     becomes
    $$  f(z)=
      \frac{1}{-z-v^2 f(z)}
      $$
      with solution $f(z)=\frac{1}{2v^2}
        (-z\pm \sqrt{z^2-4v^2})$, and thus one finds the
        classical semi-circle law,
      $$
         \frac{d\nu(z)}{dz}=\frac{1}{\pi} \Im m~ f(z)
        = \left\{\begin{array}{l}
        {\displaystyle\frac{1}{2\pi v^2} \sqrt{4v^2-z^2}\mbox{~~   for  }
         -2v\leq z\leq 2v}\\ \\
         0
          \mbox{~~   for  }
         |z|\geq 2v
         ,\end{array}
         \right.
         $$
         concentrated on the interval $[-2v,
         2v]$.

   As an exercise, consider now the case, where $v=1$ and where the diagonal matrix
   $A$ has two distinct eigenvalues, namely
         $$
 A=\diag(\underbrace{\al,\ldots,\al}_{pn},\underbrace{
 \beta,\ldots,\beta}_{(1-p)n}).
  $$
  See e.g., Adler-van Moerbeke \cite{AvM-asymmPearcey}. The integral equation (\ref{IntEqt}) becomes
    $$f-{\frac {1-p}{\beta-z-f}}-{\frac
{p}{\al-z-f}}=0
 ,~~\mbox{for } 0<p<1,$$
%
which, upon clearing, leads to a cubic equation for
$g:=f+z$,
 $$
 {g}^{3}- \left(  z+\al+\beta \right)
 {g}^{2}+ \left( z(\al+\beta
 )+\al\beta +1 \right) g-\al \beta z- (1-p)\al
  -p\beta=0,
$$
having, as one checks, a quartic discriminant $D_1(z)$
in $z$. Since the roots of a cubic polynomial involve,
in particular, the square root of the discriminant, the
solution $g(z)$ of the cubic will have a non-zero
imaginary part, if and only if $D_1(z)<0$. Thus one
finds the following equilibrium density,
      $$
         \frac{d\nu(z)}{dz}=\frac{1}{\pi} \Im m~ f(z)
        = \left\{\begin{array}{l}
        {\displaystyle \frac{1}{\pi}\Im m~ g(z)
        \mbox{   for $z$ such that }
        D_1(z) <0
         }\\ \\
         0 \mbox{   for $z$ such that }
        D_1(z) \geq 0
         \end{array}
         \right.
         $$
Therefore the boundary of the support of the equilibrium
measure will be given by the real roots of $D_1(z)=0$.
Depending on the values of the parameters $\al,~\beta$
and $p$, there will be four real roots or two real
roots, with a critical situation where there are three
real roots, i.e., when two of the four real ones
collide. The critical situation occurs exactly when the
discriminant $D_2(\al,\beta,p)$ (with regard to $z$) of
$D_1(z)$ vanishes, namely when
$$D_2(\al,\beta,p)=4 p(1-p)\gamma
 (\ga^3-3\ga^2+3\ga (9p^2-9p+1)-1)^3\Bigr|_{\ga=(\al-\beta)^2}
 =0.$$
 This polynomial has a positive root $\gamma$, which can
 be given explicitly,
  the others being imaginary, and one checks that, when
 one has the relationship
  $$
   \al-\beta
   = \frac{q+1}{\sqrt {{q}^{2}-q+1}}
   ,~~~\mbox{upon using the parametrization}~
   ~p=\frac{1}{q^3\!+\!1},$$
two of the four roots of $D_1(z)$ collide. This is to
say, this is the precise point at which the support of
the equilibrium measure goes from two to one interval.
This then occurs exactly at value
$$
z=\beta+{\frac {2 q-1}{\sqrt
{{q}^{2}-q+1}}}
 $$on the real line.

\subsection{Soft edge scaling
limit for GUE and the Tracy-Widom
distribution}\label{subsect 8.5}

Consider the probability measure on eigenvalues $z_i$ of
the $n\times n$ Gaussian Hermitian ensemble (GUE)
$$\BP_n(\mbox{all}~z_i \in \tilde E)=
\frac{\displaystyle{\int_{\tilde
E^n}\Dt_n^2(z)\prod^n_1e^{-z_i^2} }dz_i}{
\displaystyle{\int_{\BR^n}\Dt_n^2(z)\prod_1^n e^{-z_i^2}
}dz_i} .$$
Given the disjoint union $ E:=\bigcup_1^{
r}[x_{2i-1},x_{2i}]\subset \BR$, define the gradient and
 the Euler operator with respect to the boundary
 points of the set $  E$: \be\nabla_x=\sum_1^{2r} \frac{\pl}{\pl x_i}~~~
 \mbox{and}~~~{\cal E}_x=\sum_1^{2r} x_i\frac{\pl}{\pl x_i}.
  \label{gradients}\ee
Remember the definition of the Fredholm determinant of a
continuous kernel $K(x,y)$, the continuous analogue of
the discrete kernel (\ref{def: Fredholm det}),
\bean
 \det(I-K(x,y) ~\raisebox{.6mm}{$\chi$}{}_{E} (y) )
=1+\sum_{n=1}^{\iy}\frac{(-1)^n}{n!}
 \int_Edz_1\ldots dz_n\det(K(z_i,z_j))_{1\leq i,j\leq n}
.\eean We now state: \cite{TW-Airy,ASV,ASV1}
\begin{theorem}
The gradient
 $$
 f(x_1,\ldots,x_{2r}):=\nabla_x\log \BP(
 E^c),
  $$
   with
    $$\BP (E^c ):=\limn
\BP_n\left(\mbox{all}~
  \sqrt{2} n^{\frac{1}{6}}\Bigl( z_{i} -
 \sqrt{2n}   \Bigr)\in E^c  \right)
  , $$
 satisfies the 3rd order non-linear PDE:
 \be \left(\nabla_{\!x}^3 - 4 ({\cal E}_x -
\frac{1}{2})\right) f + 6( \nabla_{\!x}f)^2 = 0 .
    \label{Airy-PDE}\ee
In particular, for $E=(x,\iy)$,
    \bea
     \FR(x)&:=&\lim_{n\nearrow \iy}
\BP_n\left(  \sqrt{2} n^{\frac{1}{6}}\Bigl( z_{\max} -
 \sqrt{2n}   \Bigr)\leq x \right)
 \no\\%
   &=& \det
 (I-{\bf A}~\raisebox{.6mm}{$\chi$}{}_{(x,\iy)})\no\\
  &=&\exp\left({\displaystyle
  -\int^{\iy}_{x}(\al-x)g^2(\al)d\al}\right)
, \label{TW}\eea 
 is the {\bf Tracy-Widom distribution}, with\footnote{Remember the Airy function:
 $$
A(x)=\frac{1}{\pi}\int_0^{\iy}\cos\left(\frac{u^3}{3}
 +xu\right)du
 ,~\mbox{  satisfying the ODE}~
 ~  A''(x)=xA(x).$$
 }
 $${\bf A}(x,y):=\frac{A(x)A'(y)-A'(x)A(y)}{x-y}
  =
 \int^{\iy}_{0}A(u+x)A(u+y)du
 $$
 and $g(\al)$ the Hastings-McLeod (unique) solution of
 $$
\left\{\begin{array}{l} g''=\al g+2g^3\\ g(\al)\sim
 \frac{
  e^{-\frac{2}{3}  \al^{\frac{3}{2}}}}{2\sqrt \pi \al^{1/4}}
\mbox{\,\,for\,\,}\al\nearrow \infty.
\end{array}\right. (\mbox{{\bf Painlev\'e II}})
.$$
%
%

\end{theorem}



%
 %
%
\proof 
{\it Step 1}: Applying Proposition \ref{Proposition
7.9}, it follows that
\be \BP_n(\mbox{all}~z_i \in E^c)=
\det \left(I-K_n(y,z) \chi_{_{E }}(z)\right)
,\label{GUEkernel}\ee
 where the kernel $K_n(y,z)$ is given
  by Proposition \ref{Proposition 7.4},
 \be
  K_n(y,z)=\left(\frac{n}{2}\right)^{1/2}
  e^{-\frac{1}{2}(y^2+z^2)}
  ~\frac{ p_n (y)  p_{n-1}(z)-
     p_{n-1} (y)  p_{n }(z)}{y-z}
 \label{Hermite kernel}
 . \ee
 The $p_n$'s are orthonormal polynomials with respect
 to $e^{-z^2}$ and thus proportional to the classical
 Hermite polynomials:
 \be
  p_n (y):=\frac{1}{2^{n/2}\sqrt{ n! }\pi^{1/4}} 
    H_n(y)=\frac{1}{\sqrt{h_n}}y^n+\ldots
 ,\label{Hermite}\ee
   with
   $$ ~H_n(y):= e^{y^2} \left(-\frac{d}{d
y}\right)^n e^{-y^2},~~~~h_n=\frac{\sqrt{\pi}n!}{2^n}
 .$$

\medbreak

 \noindent{\it Step 2}: The Plancherel-Rotach
 asymptotic formula (see Szeg\"o \cite{Szego}) says
 that
%
%
 \bean \left.e^{-x^2/2}
\frac{n^{1/12}H_n(x)}{2^{n/2+1/4}\sqrt{n!}
 \pi^{1/4}}\right|_{\footnotesize
x=\sqrt{2n+1} +\frac{t}{\sqrt{  2}~n^{1/6}}
 }=A (t )+O(n^{-2/3})
 ,\eean
 uniformly for $t\in\mbox{compact}~K\subset \BC $ and
 thus, in view of (\ref{Hermite}),
\be
e^{-\frac{x^2}{2}}
p_n(x)\Biggl\vert_{x=\sqrt{2n+1}+\frac{t}{\sqrt{2}~n^{1/6}}}
=2^{1/4}n^{-\frac{1}{12}}\left(A(t)+{
O}(n^{-\frac{2}{3}})\right) .\label{estimate Hermite}\ee
  Since the Hermite kernel (\ref{Hermite kernel}) also involves
   $p_{n-1}(x)$, one needs an estimate like
  (\ref{estimate Hermite}), with the same scaling
  but for $p_n$ replaced by $p_{n-1}$. So, one needs the
  following:
\bean
x&=&\sqrt{2n+1}+\frac{t}{\sqrt{2}~n^{1/6}}\\
\\
&=&\sqrt{(2n-1)(1+\frac{2}{2n-1})}
 +\frac{t}{\sqrt{2}(n-1)^{1/6}
  (1+\frac{1}{n-1})^{1/6}}\\
\\
&=&\sqrt{2n-1}\left(1+\frac{1}{2n-1}+\mbox{
O}(\frac{1}{n^2})\right)+
 \frac{t}{\sqrt{2}(n-1)^{1/6}}\left(1+\mbox{ O}(\frac{1}{n})\right)\\
\\
&=&\sqrt{2n-1}+\frac{t+\frac{1}{n^{1/3}}}
 {\sqrt{2}(n-1)^{1/6}}+\mbox{O}\left(\frac{1}{n^{7/6}}\right).
\eean
Hence, from (\ref{estimate Hermite}) it follows that
\bea
  {e^{-\frac{x^2}{2}}
  p_{n-1}(x)\Biggl\vert_{x=\sqrt{2n+1}+
  \frac{t}{\sqrt{2}~n^{1/6}}}}
 &=&e^{-\frac{x^2}{2}}
 p_{n-1}(x)\Biggl\vert_{x=\sqrt{2n-1}
  +\frac{t+n^{-1/3}}{\sqrt{2}(n-1)^{1/6}}+\ldots}
   \no\\
  \no\\
 &=&
 2^{1/4}n^{-\frac{1}{12}}
  \left(A(t+n^{-1/3})+\mbox{
O}(n^{-2/3})\right)
 .\no\\ \label{estimate Hermite1}\eea
From the definition of the Fredholm determinant, one
needs to find the limit of $K_n(y;z)dz$. Therefore, in
view of (\ref{Hermite kernel}) and using the estimates
(\ref{estimate Hermite}) and (\ref{estimate Hermite1}),
\bean &
&\limn
 K_n(y;z)dz\Biggl\vert_{y=(2n+1)^{1/2}+
\frac{t}{\sqrt{2}~n^{1/6}}\atop{
z=(2n+1)^{1/2}+\frac{s}{\sqrt{2}~n^{1/6}}}}
  \\
\\
&=&-\limn  \left(\frac{n}{2}\right)^{1/2}
 e^{-\frac{1}{2}(y^2+z^2)}
  \\
   \\
   &&
    \hspace*{.8cm}\frac{ p_n(y)( p_n(z)- p_{n-1}(z)) - p_n(z)( p_n(y)-
p_{n-1}(y))}{(y-z)\sqrt{2}~n^{1/6}}
\Biggl\vert_{y=(2n+1)^{1/2}+ \frac{t}{\sqrt{2}~n^{1/6}}
\atop{z=(2n+1)^{1/2}+\frac{s}{\sqrt{2}~n^{1/6}}}}ds
\eean\bean
&=&\limn \left(\frac{n}{2}\right)^{1/2}(2^{1/4}
n^{-\frac{1}{12}})^2n^{-1/3}
 \\ \\
 &&\hspace*{1.4cm}\frac{A(t) \left(
{A(s+n ^{-1/3})-A(s)}\right) -A(s)\left( {A(t+n^{-1/3})-
 A(t)}
 \right)}{{n^{-1/3}}(t-s)}ds
 \\ \\
&=&\frac{A(t)A'(s)-A(s)A'(t)}{t-s}ds={\bf A}(t,s)
ds\eean
 and thus for $E:=\bigcup_1^{2r}[x_{2i-1},x_{2i}]$ and
  $\tilde E:=\bigcup_1^{2r}[c_{2i-1},c_{2i}],
  $
   related by
 \be\tilde E=\sqrt{2n}
 + \frac{E}{\sqrt{2}
n^{1/6}}\label{scalingGUE},\ee  one has shown (upon
taking the limit term by term in the sum defining the
Fredholm determinant)
 \bean
  \limn   \BP_n(\mbox{all}~ z_i \in \tilde E^c)
 &=&
 \limn \BP_n(\mbox{all}~z_i \in \sqrt{2n}
 + \frac{E^c}{\sqrt{2}
n^{1/6}})
  \\
  &=& \limn \BP_n(\mbox{all}~z_i \in (2n+1)^{1/2}
 + \frac{E^c}{\sqrt{2}
n^{1/6}})
 \\
  &=& \limn \det (I-K_n
\raisebox{.6mm}{$\chi$}{}_{\tilde
  E})\Bigr|_{_{\tilde E=(2n+1)^{1/2}+
\frac{E}{\sqrt{2} n^{1/6}}}}
 \\
  &=&\det(I-{\bf
A}\raisebox{.6mm}{$\chi$}{}_{
  E})
.\eean

\noindent{\it Step 3}: From Theorem \ref{Theo 8.1},
$\BP_n(\mbox{all}~ z_i \in \tilde E^c)$ satisfies, with
regard to the boundary points $c_i$ of $\tilde E$, the
PDE (\ref{GUE-PDE}); thus setting that scaling into this
PDE yields (remember $\nabla_c$ and ${\cal E}_x$ are as
in (\ref{gradients}) and $ {\cal B}_k=\sum_i c_i^{k+1}
\frac{\pl}{\pl c_i}$)
\bean \hspace*{-5cm}0&=& \bigl({\cal B} _{-1}^4+8n{\cal
B}_{-1}^2+12{\cal B}_0^2+
 24 {\cal
B}_0-16{\cal B}_{-1}{\cal B}_1\bigr)\log \BP_n
  \\
   &&~~~~~~~~~~~~~~~~~~
   ~~~~~~~~~~~~~~~~~~~~~~+6\bigl({\cal B}^2_{-1} \log \BP_n\bigr)^2\Bigr|_{c_i=\sqrt{2n}
  +\frac{x_i}{\sqrt{2}n^{1/6}}}
 \\
 &=&4n^{2/3}\left[\left(\nabla_{\!x}^3 - 4 ({\cal E}_x -
\frac{1}{2})\right) \nabla_{\!x}\log\BP
 + 6( \nabla_{\!x}^2\log\BP)^2%
\right]+o(n^{2/3})
 .\eean
 Note that in this computation, the terms of order
 $n^{4/3}$ cancel, because the leading term in
 $$
 12{\cal B}_0^2-16 {\cal B}_{-1} {\cal B}_1
 = -4\sum_i c_i^2\left(  \frac{\pl}{\pl
 c_i}\right)^2+\ldots=-16n^{4/3} \nabla_{\!x}^2+\ldots
 $$
cancels versus the leading term in $8n{\cal
B}_{-1}^2=16n^{4/3}\nabla_{\!x}^2+\ldots$; thus only the
terms of order $n^{2/3}$
 remain. Since in {\it step 2} it was shown that the
 limit exists, the term in brackets vanishes, showing
 that $\log \BP(E^c)$ satisfies the PDE (\ref{Airy-PDE}).

\medbreak

{\it Step 4}: In particular, upon picking $E=(x,\iy)$,
the PDE (\ref{Airy-PDE}) for
 $$
f(x)= \frac{\pl}{\pl x} \log \FR(x)=
 \frac{\pl}{\pl x} \log
  \lim_{n\nearrow \iy} \BP_n\left(  \sqrt{2}
n^{\frac{1}{6}}\Bigl( z_{\max} -
 \sqrt{2n}   \Bigr)\leq x \right)
 $$
becomes an ODE:
 $$
 f'''-4xf'+2f+6f'^2=0 .
 $$
 Multiplying this equation with $f''$ and
 integrating from $-\iy$ to $x$ lead to the
 differential equation (the nature of the solution shows
 the integration constant must vanish)
  \be
  f^{\prime\prime 2}+4f'(f'^2-xf'+f)=0.
  \label{McLeod1}\ee
 Then, setting
 \be
  \left\{\begin{array}{l}
  f'=-g^2\\
  f=g'^2-xg^2-g^4
   \end{array}
   \right.
   \label{McLeod2}\ee
and, since then $f''=-2gg'$, an elementary computation
shows that (\ref{McLeod2}) is an obvious solution to
equation (\ref{McLeod1}). For (\ref{McLeod2}) to be
valid, the derivative of the right hand side of the
second expression in (\ref{McLeod2}) must equal the
derivative of the right hand side of the first
expression in (\ref{McLeod2}), i.e.
we must have:
$$
 0=(f)'-f'=(g'^2-xg^2-g^4)'+g^2=2g'(g''-2g^3-xg),
 $$
 and so
  $
  g''= 2g^3+xg.
   $
 Hence
 $$
  \frac{\pl^2}{\pl x^2}\log \FR(x)
   =f'=-g^2
 $$
 Integrating once yields (assuming that $g^2$ decays
 fast enough at $\iy$, which will be apparent later)
  $$ \frac{\pl }{\pl x }\log \FR(x)
   =\int_x^{\iy}g^2(u) du
   ;$$
integrating once more and further integrating by parts
yield \bea
 \log \FR(x)
   =\int_x^{\iy}d\al \int_\al^{\iy}g^2(u) du
   &=&\int_x^{\iy}d\al \bigl(\frac{d}{d\al}\al\bigr)\int_\al^{\iy}g^2(u) du
  \no\\& =&
   \al\int_\al^{\iy}g^2(u) du\Bigr|_x^{\iy}+
    \int_x^{\iy}  \al g^2(\al)d\al \no\\
    &=&
     \int_x^{\iy}  (\al-x) g^2(\al)d\al, \label{TW1}\eea
confirming (\ref{TW}). For $x \rightarrow \iy$, one
checks that, on the one hand, from the definition of the
Fredholm determinant of ${\bf A}$, the two leading terms
are given by
 $$ \FR(x)=\det
\left(I-{\bf
A}\raisebox{.6mm}{$\chi$}{}_{_{[x,\iy)}}\right) =
 1-\int_x^{\iy}   {\bf A}(z,z)dz+\ldots$$
 and, on the other hand, from (\ref{TW1}), the two
 leading terms are
 $$
  \FR(x)=1-\int^{\iy}_{x}(\al-x)g^2(\al) d\al+\ldots
 $$
Therefore, comparing the two expressions above,
$$
\int_x^{\iy}(\al -x)g^2(\al)d\al
 =\int_x^{\iy}dz~{\bf A}(z,z)+\ldots
$$
and upon taking two derivatives in $x$,
\bean
-g^2(x)&=&\frac{\pl}{\pl x}{\bf A}(x,x)+\ldots\\
\\
&=&\frac{\pl}{\pl x}\int_0^{\iy}A(u+x)^2du +\ldots\\
\\
&=&\frac{\pl}{\pl x}\int_x^{\iy}A(u)^2du +\ldots\\
\\
&=&-A(x)^2+\ldots , \eean showing that asymptotically
$g(x)\sim A(x)$ for $x\rg \iy$. It is classically known
that asymptotically
 $$
A(x)=\frac{e^{-\frac{2}{3}x^{3/2}}}{2\sqrt{\pi}x^{1/4}}
(1+\sum_1^{\iy}\frac{\al_i}{(x^{3/2})^i}+\ldots),~~~\mbox{as
$x\rg \iy$.}$$
 The solution $g(x)$ of $ g''= 2g^3+xg$ behaving like
the Airy function at $\iy$ is unique (Hastings-McLeod
solution). It behaves like
 \bean
 g(x)&=&\mbox{Ai}(x)
 +{\bf O}\left(\frac{e^{-\frac{4}{3}x^{3/2}}}{x^{1/4}}\right)
 \mbox{~for~}x\rg\iy\\
 \\
&=&\sqrt{\frac{-x}{2}}\left(1+\frac{1}{8x^3}-\frac{73}{128
x^6}+{\bf O}(|x|^{-9})\right)\mbox{~for~}x\rg -\iy.
\eean
The Tracy-Widom distribution $\FR$ of mean and standard
deviation (see Tracy-Widom \cite{TW-ICM})
$$E(\FR)=-1.77109 ~~\mbox{and}~~ \sigma(\FR)=0.9018$$
has a density decaying for $x\rg\iy$ as (since ${\cal
F}(x)$ tends to $1$ for $x\rg \iy$)
 \be
  \FR'(x)=\FR(x)\int_x^{\iy}g^2(\al)d\al
         \sim \int_x^{\iy}A^2(\al)d\al
          \sim   \frac{1}{8\pi
          x}e^{-\frac{4}{3}x^{3/2}}
 ~~\mbox{for}~  x\rg\iy .\label{8.7.13}\ee
  The last estimate is obtained by integration by parts:
 \bean
  \int_x^{\iy} A^2(u)du&=&
 \int_x^{\iy}\frac{-1}{8\pi u}
 (1+\sum_1^{\iy}\frac{  c _i}{(u^{3/2})^i})
d(e^{-\frac{4}{3}u^{3/2}})\\
&=&
 \frac{e^{-\frac{4}{3}x^{3/2}}}{8\pi x}
  (1+\sum_1^{\iy}\frac{  c _i}{(x^{3/2})^i} )-
  \!\frac{1}{8\pi}\!
   \int_x^{\iy}\frac{e^{-\frac{4}{3}u^{3/2}}}{ u^2}
  (1+\sum_1^{\iy}\frac{  c' _i}{(u^{3/2})^i} ).
\eean
 Then, following conjectures by
 Dyson \cite{Dyson} and Widom \cite{Widom},
 Deift-Its-Krasovsky-Zhou \cite{Deift-Its-Kras}
  and
  Baik-Buckingham-DiFranco \cite{BBD}
 give a
 representation of ${\cal F}(x)$ as an integral from
 $-\iy$ to $x$ and thus this provides
  an estimate for $x\rg -\iy$,
\bean
  {\cal F}(x)&=&2^{1/24}e^{ \zeta'(-1)}
\frac{e^{-\frac{1}{12}|x|^3}}{|x|^{1/8}}
\exp\left\{\int_{-\iy}^x\left(R(y)-\frac{1}{4}
y^2+\frac{1}{8y}\right)dy\right\} \\
&=& 2^{1/24}e^{ \zeta'(-1)}
 \frac{e^{-\frac{1}{12}|x|^3}}{|x|^{1/8}}%
 \left(1+\frac{3}{2^6|x|^3}+O(|x|^{-6}\right)
,~~\mbox{for}~x\rg
-\iy,%
 \eean
   where
  $$
 R(y)=\int_y^{\iy} g(s)^2ds=g'(y)^2-yg(y)^2-g(y)^4.
  $$

\end{document}